\documentclass[12pt]{amsart}   
\usepackage{srcltx}
\usepackage{amsmath}
\usepackage{amsxtra}
\usepackage{amscd}
\usepackage{amsthm}
\usepackage{amsfonts}
\usepackage{amssymb}
\usepackage{eucal}
\usepackage{verbatim}
\usepackage[all]{xy}
\usepackage{color}
\definecolor{darkspringgreen}{rgb}{0.09, 0.45, 0.27}
\definecolor{Red}{rgb}{0.7,0,0}
\usepackage{stmaryrd}

\usepackage[pdftex,bookmarks=false,colorlinks=true,citecolor=darkspringgreen,debug=true,  naturalnames=true,pdfnewwindow=true]{hyperref}
\usepackage{geometry}
\usepackage[utf8]{inputenc}
\usepackage[T1]{fontenc}
\usepackage[english]{babel}
\usepackage{microtype} 
\usepackage{tabularx} 
\usepackage{booktabs} 
\usepackage{amsmath, amsfonts,amssymb,amsthm,mathbbol, textcomp, yfonts}
\usepackage[all]{xy}
\usepackage{tikz,nicematrix}
\usetikzlibrary{decorations.pathreplacing}
\usepackage{appendix}
\usepackage{etoolbox}%
\usepackage[T1]{fontenc}
\usepackage{lmodern}
\patchcmd{\footnotemark}{\stepcounter{footnote}}{\refstepcounter{footnote}}{}{}
\usepackage{blkarray}
\usepackage{multirow}
\usepackage{mathtools}
\usepackage{graphicx}

\newtheorem{thm}[subsubsection]{Theorem}
\newtheorem{cor}[subsubsection]{Corollary}

\newtheorem{lem}[subsubsection]{Lemma}
\newtheorem{prop}[subsubsection]{Proposition}
\newtheorem{conj}[subsubsection]{Conjecture}
\theoremstyle{definition} 
\newtheorem{defn}[subsubsection]{Definition}
\newtheorem{rem}[subsubsection]{Remark}

\newcommand{\nc}{\newcommand}
\nc{\renc}{\renewcommand} \nc{\ssec}{\subsection}
\nc{\sssec}{\subsubsection}

 \nc{\wh}{\widehat}
\nc\ol{\overline} \nc\ul{\underline} \nc\wt{\widetilde}

\nc{\BA}{{\mathbb{A}}} \nc{\BC}{{\mathbb{C}}} \nc{\BF}{{\mathbb{F}}}
\nc{\BD}{{\mathbb{D}}} \nc{\BG}{{\mathbb{G}}} \nc{\BQ}{{\mathbb{Q}}}
\nc{\BM}{{\mathbb{M}}} \nc{\BN}{{\mathbb{N}}} \nc{\BO}{{\mathbb{O}}}
\nc{\BP}{{\mathbb{P}}} \nc{\BR}{{\mathbb{R}}}
\nc{\BZ}{{\mathbb{Z}}} \nc{\BS}{{\mathbb{S}}} \nc{\BW}{{\mathbb{W}}}
\nc{\BL}{{\mathbb{L}}}

\nc{\CA}{{\mathcal{A}}} \nc{\CB}{{\mathcal{B}}} \nc{\CalD}{{\mathcal{D}}}
\nc{\CE}{{\mathcal{E}}} \nc{\CF}{{\mathcal{F}}}
\nc{\CG}{{\mathcal{G}}} \nc{\CH}{{\mathcal{H}}}
\nc{\CI}{{\mathcal{I}}} \nc{\CK}{{\mathcal{K}}} \nc{\CL}{{\mathcal{L}}}
\nc{\CM}{{\mathcal{M}}} \nc{\CN}{{\mathcal{N}}}
\nc{\CO}{{\mathcal{\bfO}}} \nc{\CP}{{\mathcal{P}}}
\nc{\CQ}{{\mathcal{Q}}} \nc{\CR}{{\mathcal{R}}}
\nc{\CS}{{\mathcal{S}}} \nc{\CT}{{\mathcal{T}}}
\nc{\CU}{{\mathcal{U}}} \nc{\CV}{{\mathcal{V}}}  \nc{\CY}{{\mathcal Y}}
\nc{\CW}{{\mathcal{W}}} \nc{\CZ}{{\mathcal{Z}}}

\nc{\cM}{{\check{\mathcal M}}{}} \nc{\csM}{{\check{\mathcal A}}{}}
\nc{\oM}{{\overset{\circ}{\mathcal M}}{}}
\nc{\obM}{{\overset{\circ}{\mathbf M}}{}}
\nc{\oCA}{{\overset{\circ}{\mathcal A}}{}}
\nc{\obA}{{\overset{\circ}{\mathbf A}}{}}
\nc{\ooM}{{\overset{\circ}{M}}{}}
\nc{\osM}{{\overset{\circ}{\mathsf M}}{}}
\nc{\vM}{{\overset{\bullet}{\mathcal M}}{}}
\nc{\nM}{{\underset{\bullet}{\mathcal M}}{}}
\nc{\oD}{{\overset{\circ}{\mathcal D}}{}}
\nc{\obD}{{\overset{\circ}{\mathbf D}}{}}
\nc{\oA}{{\overset{\circ}{\mathbb A}}{}}
\nc{\op}{{\overset{\bullet}{\mathbf p}}{}}
\nc{\cp}{{\overset{\circ}{\mathbf p}}{}}
\nc{\oU}{{\overset{\bullet}{\mathcal U}}{}}
\nc{\ofZ}{{\overset{\circ}{\mathfrak Z}}{}}

\nc{\ff}{{\mathfrak{f}}} \nc{\fv}{{\mathfrak{v}}}
\nc{\fa}{{\mathfrak{a}}} \nc{\fb}{{\mathfrak{b}}}
\nc{\fd}{{\mathfrak{d}}} \nc{\fe}{{\mathfrak{e}}}
\nc{\fg}{{\mathfrak{g}}} \nc{\fgl}{{\mathfrak{gl}}}
\nc{\fh}{{\mathfrak{h}}} \nc{\fri}{{\mathfrak{i}}}
\nc{\fj}{{\mathfrak{j}}} \nc{\fk}{{\mathfrak{k}}} \nc{\fl}{{\mathfrak{l}}}
\nc{\fm}{{\mathfrak{m}}} \nc{\fn}{{\mathfrak{n}}}
\nc{\ft}{{\mathfrak{t}}} \nc{\fu}{{\mathfrak{u}}}
\nc{\fw}{{\mathfrak{w}}} \nc{\fz}{{\mathfrak{z}}}
\nc{\fp}{{\mathfrak{p}}} \nc{\fq}{{\mathfrak{q}}} \nc{\frr}{{\mathfrak{r}}}
\nc{\fs}{{\mathfrak{s}}} \nc{\fsl}{{\mathfrak{sl}}}
\nc{\hsl}{{\widehat{\mathfrak{sl}}}}
\nc{\hgl}{{\widehat{\mathfrak{gl}}}}
\nc{\hg}{{\widehat{\mathfrak{g}}}}
\nc{\chg}{{\widehat{\mathfrak{g}}}{}^\vee}
\nc{\hn}{{\widehat{\mathfrak{n}}}}
\nc{\chn}{{\widehat{\mathfrak{n}}}{}^\vee}

\nc{\fA}{{\mathfrak{A}}} \nc{\fB}{{\mathfrak{B}}} \nc{\fC}{{\mathfrak{C}}}
\nc{\fD}{{\mathfrak{D}}} \nc{\fE}{{\mathfrak{E}}}
\nc{\fF}{{\mathfrak{F}}} \nc{\fG}{{\mathfrak{G}}} \nc{\fH}{{\mathfrak{H}}}
\nc{\fI}{{\mathfrak{I}}} \nc{\fJ}{{\mathfrak{J}}}
\nc{\fK}{{\mathfrak{K}}} \nc{\fL}{{\mathfrak{L}}}
\nc{\fM}{{\mathfrak{M}}} \nc{\fN}{{\mathfrak{N}}}
\nc{\frP}{{\mathfrak{P}}} \nc{\fQ}{{\mathfrak{Q}}}
\nc{\fS}{{\mathfrak{S}}} \nc{\fT}{{\mathfrak{T}}} \nc{\fU}{{\mathfrak{U}}}
\nc{\fV}{{\mathfrak{V}}} \nc{\fW}{{\mathfrak{W}}}
\nc{\fX}{{\mathfrak{X}}} \nc{\fY}{{\mathfrak{Y}}}
\nc{\fZ}{{\mathfrak{Z}}}

\nc{\ba}{{\mathbf{a}}}
\nc{\bb}{{\mathbf{b}}} \nc{\bc}{{\mathbf{c}}}
\nc{\be}{{\mathbf{e}}}
\nc{\bE}{{\mathbf{E}}}
\nc{\bj}{{\mathbf{j}}} \nc{\bm}{{\mathbf{m}}}
\nc{\bn}{{\mathbf{n}}} \nc{\bp}{{\mathbf{p}}}
\nc{\bq}{{\mathbf{q}}} \nc{\br}{{\mathbf{r}}} \nc{\bt}{{\mathbf{t}}}
\nc{\bfu}{{\mathbf{u}}} \nc{\bv}{{\mathbf{v}}}
\nc{\bx}{{\mathbf{x}}} \nc{\by}{{\mathbf{y}}} \nc{\bz}{{\mathbf{z}}}
\nc{\bw}{{\mathbf{w}}} \nc{\bA}{{\mathbf{A}}}
\nc{\bB}{{\mathbf{B}}} \nc{\bC}{{\mathbf{C}}}
\nc{\bD}{{\mathbf{D}}} \nc{\bF}{{\mathbf{F}}} \nc{\bG}{{\mathbf{G}}}
\nc{\bH}{{\mathbf{H}}} \nc{\bI}{{\mathbf{I}}} \nc{\bJ}{{\mathbf{J}}}
\nc{\bK}{{\mathbf{K}}} \nc{\bM}{{\mathbf{M}}} \nc{\bN}{{\mathbf{N}}}
\nc{\bO}{{\mathbf{\bfO}}} \nc{\bS}{{\mathbf{S}}} \nc{\bT}{{\mathbf{T}}}
\nc{\bU}{{\mathbf{U}}} \nc{\bV}{{\mathbf{V}}} \nc{\bW}{{\mathbf{W}}}
\nc{\bX}{{\mathbf{X}}}
\nc{\bY}{{\mathbf{Y}}} \nc{\bP}{{\mathbf{P}}}
\nc{\bZ}{{\mathbf{Z}}} \nc{\bh}{{\mathbf{h}}}

\nc{\sA}{{\mathsf{A}}} \nc{\sB}{{\mathsf{B}}}
\nc{\sC}{{\mathsf{C}}} \nc{\sD}{{\mathsf{D}}}
\nc{\sE}{{\mathsf{E}}} \nc{\sF}{{\mathsf{F}}} \nc{\sG}{{\mathsf{G}}}
\nc{\sI}{{\mathsf{I}}} \nc{\sK}{{\mathsf{K}}} \nc{\sL}{{\mathsf{L}}}
\nc{\sfm}{{\mathsf{m}}} \nc{\sM}{{\mathsf{M}}} \nc{\sO}{{\mathsf{\bfO}}}
\nc{\sQ}{{\mathsf{Q}}} \nc{\sP}{{\mathsf{P}}}
\nc{\sT}{{\mathsf{T}}} \nc{\sZ}{{\mathsf{Z}}}
\nc{\sV}{{\mathsf{V}}} \nc{\sW}{{\mathsf{W}}}
\nc{\sfp}{{\mathsf{p}}} \nc{\sq}{{\mathsf{q}}} \nc{\sr}{{\mathsf{r}}}
\nc{\st}{{\mathsf{t}}} \nc{\sfb}{{\mathsf{b}}}
\nc{\sfc}{{\mathsf{c}}} \nc{\sd}{{\mathsf{d}}}
\nc{\sz}{{\mathsf{z}}}
\nc{\si}{{\mathsf{i}}} 
\nc{\sj}{{\mathsf{j}}}

\nc{\tA}{{\widetilde{\mathbf{A}}}}
\nc{\tB}{{\widetilde{\mathcal{B}}}}
\nc{\tg}{{\widetilde{\mathfrak{g}}}} \nc{\tG}{{\widetilde{G}}}
\nc{\TM}{{\widetilde{\mathbb{M}}}{}}
\nc{\tO}{{\widetilde{\mathsf{\bfO}}}{}}
\nc{\tU}{{\widetilde{\mathfrak{U}}}{}} \nc{\TZ}{{\tilde{Z}}}
\nc{\tx}{{\tilde{x}}} \nc{\tbv}{{\tilde{\bv}}}
\nc{\tfP}{{\widetilde{\mathfrak{P}}}{}} \nc{\tz}{{\tilde{\zeta}}}
\nc{\tmu}{{\tilde{\xi}}}

\nc{\urho}{\underline{\rho}} \nc{\uB}{\underline{B}}
\nc{\uC}{{\underline{\mathbb{C}}}} \nc{\ui}{\underline{i}}
\nc{\uj}{\underline{j}} \nc{\ofP}{{\overline{\mathfrak{P}}}}
\nc{\oB}{{\overline{\mathcal{B}}}}
\nc{\og}{{\overline{\mathfrak{g}}}} \nc{\oI}{{\overline{I}}}

\nc{\eps}{\varepsilon} \nc{\hrho}{{\hat{\rho}}}
\nc{\blambda}{{\boldsymbol{\lambda}}} \nc{\bmu}{{\boldsymbol{\xi}}} \nc{\bnu}{{\boldsymbol{(\eta, \eta')}}}

\nc{\one}{{\mathbf{1}}} \nc{\two}{{\mathbf{t}}}
\nc{\Cat}{\mathop{\operatorname{\rm 1-Cat}}}
\nc{\Sym}{\mathop{\operatorname{\rm Sym}}}
\nc{\Tot}{{\mathop{\operatorname{\rm Tot}}}}
\nc{\sprd}{{\mathop{\operatorname{\rm sprd}}}}
\nc{\Spec}{\mathop{\operatorname{\rm Spec}}}
\nc{\Ker}{{\mathop{\operatorname{\rm Ker}}}}
\nc{\Isom}{{\mathop{\operatorname{\rm Isom}}}}
\nc{\Hilb}{{\mathop{\operatorname{\rm Hilb}}}}
\nc{\Conf}{{\mathop{\operatorname{\rm Conf}}}}
\nc{\deeq}{{\mathop{\operatorname{\rm deeq}}}}

\nc{\End}{{\mathop{\operatorname{\rm End}}}}
\nc{\Ran}{{\mathop{\operatorname{\rm Ran}}}}
\nc{\Sch}{{\mathop{\operatorname{\rm Sch}}}}
\nc{\Ext}{{\mathop{\operatorname{\rm Ext}}}}
\nc{\Hom}{{\mathop{\operatorname{\rm Hom}}}}
\nc{\CHom}{{\mathop{\operatorname{{\mathcal{H}}\it om}}}}
\nc{\GL}{{\mathop{\operatorname{\rm GL}}}}
\nc{\AGL}{{\mathop{\operatorname{\rm AGL}}}}
\nc{\unit}{{\mathop{\operatorname{\rm unit}}}}
\nc{\Mir}{{\mathop{\operatorname{\rm Mir}}}}
\nc{\St}{{\mathop{\operatorname{\rm St}}}}
\nc{\oblv}{{\mathop{\operatorname{\rm oblv}}}}
\nc{\gr}{{\mathop{\operatorname{\rm gr}}}}
\nc{\Id}{{\mathop{\operatorname{\rm Id}}}}
\nc{\perf}{{\mathop{\operatorname{\rm perf}}}}
\nc{\defi}{{\mathop{\operatorname{\rm def}}}}
\nc{\length}{{\mathop{\operatorname{\rm length}}}}
\nc{\supp}{{\mathop{\operatorname{\rm supp}}}}
\nc{\colim}{{\mathop{\operatorname{\rm colim}}}}
\nc{\Fun}{{\mathop{\operatorname{\rm Funct}}}}
\nc{\Hei}{{\mathop{\operatorname{\rm Heis}}}}
\nc{\Vac}{{\mathop{\operatorname{\rm Vac}}}}
\nc{\HC}{{\mathcal H}{\mathcal C}}
\nc{\ren}{{\mathsf{ren}}}
\nc{\locc}{{\mathsf{loc.c}}}
\nc{\pr}{{\operatorname{pr}}}
\nc{\Cliff}{{\mathsf{Cliff}}}
\nc{\loc}{{\operatorname{loc}}}
\nc{\Fl}{{\mathbf{Fl}}} \nc{\Ffl}{{\mathcal{F}\ell}}
\nc{\Fib}{{\mathsf{Fib}}}
\nc{\Coh}{{\mathsf{Coh}}} \nc{\FCoh}{{\mathsf{FCoh}}}
\nc{\Perf}{{\mathsf{Perf}}}
\nc{\ShvCat}{{\mathsf{ShvCat}}}
\nc{\reg}{{\text{\rm reg}}}
\nc{\gvee}{{\mathfrak g}^{\!\scriptscriptstyle\vee}}
\nc{\tvee}{{\mathfrak t}^{\!\scriptscriptstyle\vee}}
\nc{\nvee}{{\mathfrak n}^{\!\scriptscriptstyle\vee}}
\nc{\bvee}{{\mathfrak b}^{\!\scriptscriptstyle\vee}}
       \nc{\rhovee}{\rh\bfO^{\!\scriptscriptstyle\vee}}

\nc{\cplus}{{\mathbf{C}_+}} \nc{\cminus}{{\mathbf{C}_-}}
\nc{\cthree}{{\mathbf{C}_*}} \nc{\Qbar}{{\bar{Q}}}

\nc{\Gtimes}{\vphantom{j^{X^2}}\smash{\overset{G}{\vphantom{\rule{0pt}{0.30em}}\smash{\times}}}}
\nc{\sGtimes}{\vphantom{j^{X^2}}\smash{\overset{\mathsf G}{\vphantom{\rule{0pt}{0.30em}}\smash{\times}}}}

\nc{\bOmega}{{\overline{\Omega}}}

\nc{\seq}[1]{\stackrel{#1}{\sim}}
\nc{\nilp}{{\operatorname{Nilp}}}
\nc{\Bun}{{\operatorname{Bun}}}
\nc{\aff}{{\operatorname{aff}}}
\nc{\lax}{{\operatorname{lax}}}
\nc{\gen}{{\operatorname{gen}}}
\nc{\fin}{{\operatorname{fin}}}
\nc{\mir}{{\operatorname{mir}}}
\nc{\triv}{{\operatorname{triv}}}
\nc{\ext}{{\operatorname{ext}}}
\nc{\righ}{{\operatorname{right}}}
\nc{\lef}{{\operatorname{left}}}
\nc{\forg}{{\operatorname{forg}}}
\nc{\fid}{{\operatorname{fd}}}
\nc{\modu}{{\operatorname{-mod}}}
\nc{\Mor}{{\operatorname{Mor}}}
\nc{\Gr}{{\mathbf{Gr}}}
\nc{\FT}{{\operatorname{FT}}}
\nc{\Mat}{{\operatorname{Mat}}}
\nc{\MSt}{{\operatorname{MSt}}}
\nc{\sph}{{\operatorname{sph}}}
\nc{\GR}{{\mathbf{Gr}}}
\nc{\Perv}{{\operatorname{Perv}}}
\nc{\Rep}{{\operatorname{Rep}}}
\nc{\Ind}{{\operatorname{Ind}}}
\nc{\IC}{{\operatorname{IC}}}
\nc{\Proj}{{\operatorname{Proj}}}
\nc{\colimc}{{\operatorname{colim}^{co}}}
\nc{\Stab}{{\operatorname{Stab}}}
\nc{\pt}{{\operatorname{pt}}}
\nc{\bfmu}{{\boldsymbol{\xi}}}
\nc{\bfomega}{{\boldsymbol{\omega}}}
\nc{\calM}{\mathcal M}
\nc{\calA}{\mathcal A}
\nc{\calO}{\mathcal O}
\nc{\cC}{\mathcal C}
\nc{\CC}{\mathbb C}
\nc{\calN}{\mathcal N}
\nc{\grg}{\mathfrak g}
\nc{\tslash}{/\!\!/\!\!/}
\nc\grt{\mathfrak t}
\nc\bfM{\mathbf M}
\nc\bfN{\mathbf N}

\nc\ZZ{\mathbb{Z}}
\nc\calC{\mathcal C}
\nc\calF{\mathcal F}
\nc\calX{\mathcal X}
\nc\calY{\mathcal Y}
\nc\AJ{\operatorname{AJ}}
\nc\QCoh{\operatorname{QCoh}}
\nc\FM{\operatorname{FactMod}}
\nc\SFM{\operatorname{SFactMod}}
\nc\SPerv{\operatorname{SPerv}}
\nc\SD{\operatorname{SD}}
\nc\SVect{\operatorname{SVect}}
\nc\fact{\operatorname{fact}}
\nc\diag{\operatorname{diag}}
\nc\IndCoh{\operatorname{IndCoh}}
\nc\Maps{\operatorname{Maps}}
\nc\Sect{\operatorname{Sect}}
\nc\Dmod{D-\operatorname{mod}}
\newcommand\Hecke{\operatorname{Hecke}}
\nc{\calD}{\mathcal D}
\nc\bfO{\mathbf O}
\nc\bfF{\mathbf F}

\nc\GG{\mathbb G}
\nc\calK{\mathcal K}
\nc{\calG}{\mathcal G}
\nc\RHom{\operatorname{RHom}}
\nc\Res{\operatorname{Res}}
\nc\Av{\operatorname{Av}}
\nc\Const{\operatorname{Const}}
\nc\grs{\mathfrak s}
\nc{\tilX}{\widetilde X}
\nc\calB{\mathcal B}
\nc\calS{\mathcal S}
\nc\calT{\mathcal T}
\nc\calZ{\mathcal Z}
\nc\LS{\operatorname{LocSys}}
\nc\bfL{\on{\mathbf L}}

\newcommand*\circled[1]
{*{#1}}
\makeatletter
\newcommand{\raisemath}[1]{\mathpalette{\raisem@th{#1}}}
\newcommand{\raisem@th}[3]{\raisebox{#1}{$#2#3$}}
\nc{\binlim}[2][]{\def\@tempa{#1}\@ifnextchar^{\@binlim{#2}}{\@binlim{#2}^{}}}
\def\@binlim#1^#2{\mathbin{\@ifempty{#2}{\mathop{#1}}{\mathop{#1}\@xp\displaylimits\@tempa^{#2}}}}

\makeatother

\nc\cX{{\mathcal X}}

\nc\Gm{{\mathbb G_m}}
\nc{\isoto}
    {\stackrel{\displaystyle\raisemath{-.2ex}{\sim}}\to}
\nc\cD\CalD
\nc\setm\setminus
\nc\cE\CE
\renc\Hecke{\mathit{\CH\kern-.2ex ecke}}
\nc\bsl\backslash
\nc\Fq{\mathbb F_q}
\nc\x\times
\nc\bGO{{\bG_\bO}}
\nc\bla\blambda
\nc\blt\bullet
\nc\opp{{\on{op}}}
\nc\srel\stackrel
\nc\tbx{\binlim{\widetilde\boxtimes{}}}
\nc\chkbx{\binlim{\check\boxtimes}}
\nc\into\hookrightarrow
\nc\otto\leftrightarrow
\renc\setminus\smallsetminus
\nc\phitau{\varphi\tau}
\newenvironment{i-ii-iii}{%
\begin{enumerate}%
}%
{\end{enumerate}}
\nc\ceil[1]{\lceil#1\rceil}  \nc\floor[1]{\lfloor#1\rfloor}
\nc\Lie{\on{Lie}}
 \let\arXiv\arxiv
\nc\on\operatorname
\nc\kap{\kappa}
\nc\gra{\mathfrak a}
\nc\gl{\mathfrak{gl}}
\nc\sTr{\operatorname{sTr}}
\nc\hatG{\widehat{G}}
\nc\calL{\mathcal L}
\nc\Whit{\operatorname{Whit}}
\nc\KL{\operatorname{KL}}

\makeatletter

\allowdisplaybreaks
\nc\mto{\mapsto }
\nc\en{\enspace }

\numberwithin{equation}{section}
\newtheorem*{rep@theorem}{\rep@title}
\newcommand{\newreptheorem}[2]{%
\newenvironment{rep#1}[1]{%
 \def\rep@title{#2 \ref{##1}}%
 \begin{rep@theorem}}%
 {\end{rep@theorem}}}
\makeatother
\usepackage{tikz-cd}

\usepackage{wasysym, stackengine, makebox, graphicx}

\usepackage{xcolor,verbatim}
 \newcommand{\ncmd}{\newcommand*}
\newcommand{\rncmd}{\renewcommand*}
\ncmd{\DMO}{\DeclareMathOperator}
\ncmd{\ncmdd}[2]{\ncmd{#1}{{#2}}}

\makeatletter
\ncmd{\DefOps}[1]{\def\OPERATOR@NAME##1{\DeclareMathOperator{##1}{\expandafter\@gobble\string##1}}
    \def\OPERATOR@LIST##1{\ifcat\noexpand\relax\noexpand##1\OPERATOR@NAME##1\expandafter\OPERATOR@LIST\fi}
    \OPERATOR@LIST#1.}

\ncmd{\DefRm}[1]{\def\OPERATOR@NAME##1{\ncmd{##1}{\mathrm{\expandafter\@gobble\string##1}}}
    \def\OPERATOR@LIST##1{\ifcat\noexpand\relax\noexpand##1\OPERATOR@NAME##1\expandafter\OPERATOR@LIST\fi}
    \OPERATOR@LIST#1.}

\nc\negquad{\mkern-18mu}
\nc\lhs{\quad&\negquad}
\ncmd{\phtr}[2]{\lefteqn{#1{\phantom{#2}}}#2}

\def\Alphabet{ABCDEFGHIJKLMNOPQRSTUVWXYZ}
\def\newalph#1#2{\begingroup
    \def\procL@tt@r##1{%
        \@xp\gdef\csname#1\endcsname{#2}}%
    \proc@lph@bet\endgroup}
\def\proc@lph@bet{\@xp\prlist@\Alphabet\relax}
\def\prlist@#1{\ifx#1\relax\else\procL@tt@r{#1}\@xp\prlist@\fi}

\ncmd{\SmSub}[2][]{_\bgroup #2\smsub@{#1}}
\def\smsub@#1{\@ifnextchar_{#1\@smsb}{\egroup}}
\def\@smsb_#1{#1\smsub@{}}

\ncmd{\SmSup}[2][]{^\bgroup #2\smsup@{#1}}
\def\smsup@#1{\@ifnextchar^{#1\@smsp}{\@ifnextchar'{\prime\@xp\smsup@\@xp{\@xp}\@gobble}{\egroup}}}
\def\@smsp^#1{#1\smsup@{}}

\def\@binlim#1_#2{\mathbin{\@ifempty{#2}{\mathop{#1}}{\mathop{#1}\@xp\displaylimits\@tempa_{#2}}}}

\catcode`\=\active
\ncmd{\RedNote}[1]{\Text{%
  \textcolor{Red}{\sffamily#1}}}
\ncmd\Text[1]{\ifmmode\text{#1}\else#1\fi}
\ncmd\todo[1][todo]{\RedNote{#1}%
 \ifx\undefined\@todoflag
   \global\let\@todoflag\relax
  \AtEndDocument{\par\vspace{3em}\noindent \RedNote{TO DO:}}%
 \fi
 \edef\pagenmbr{\thepage}%
 \expandafter\redendnote
 \expandafter{\pagenmbr}{#1}}
\ncmd\redendnote[2]{\AtEndDocument{\\$\bullet$\ \RedNote{page #1: #2}}}
\def{\todo}
\ncmd\hidetodos{\rncmd\todo\relax}

\theoremstyle{remark}
\newalph{c#1}{\ensuremath{{\mathcal #1}}}
\newalph{s#1}{{\mathsf #1}}
\newalph{#1#1}{{\mathbb #1}}
\newalph{e#1}{{\EuScript #1}}
\newalph{g#1}{{\mathfrak #1}}
\newalph{r#1}{{\ensuremath{\mathscr #1}}}
   
\ncmd{\Fp}{{\FF_p}}

\DMO{\cHom}{\text{\textrm{\itshape{\cH}\kern-.2ex{}om}}}
\DMO{\cEnd}{\text{\textrm{\itshape{\cE}\kern-.2ex{}nd}}}    \DMO{\cExt}{\text{\textrm{\itshape{\cE}\kern-.2ex{}xt}}}

\ncmd{\sff}{\mathsf f}

\let\Im\undefined   \let\det\undefined
\ncmd{\young}[1]{\vcenter{\begin{Young}#1\crcr\end{Young}}}
\DefOps   { \REnd  \Aut
             \det \tr \rk \curv  \Coker \Im  \ad \Ad \add}
\DefRm{ \id \Vect \Pic \Loc \Hitch \Higgs \Eu \Fr}

\ncmd{\RGam}{\text{\upshape R}\Gamma}
\DMO{\chr}{char}
\ncmd{\angs}[1]{\langle#1\rangle}
\ncmd{\hmod}{\text{\upshape-mod}}
\ncmd{\cxym}[1]{\ensuremath{\vcenter{\xymatrix{#1}}}}
 \let\arXiv\arxiv
\makeatother

\title{Twisted Gaiotto equivalence for $\GL(M|N)$}
\date{}
\author[R.Travkin]{Roman Travkin}
\address{Skolkovo Institute of Science and Technology, Moscow, Russia}
\email{roman.travkin2012@gmail.com}
\author[R.Yang]{Ruotao Yang}
\address{Skolkovo Institute of Science and Technology, Moscow, Russia}
\email{yruotao@gmail.com}

\newcommand\om\omega
\nc\tboxtimes{\mathbin{\widetilde\boxtimes}}
\begin{document}

\begin{abstract}
In $\GL_N$, a series of subgroups indexed by $0\leq M\leq N-1$ were noticed by D. Gaiotto and \cite{[C], [JPS]}. It was conjectured by D. Gaiotto that the categories of twisted D-modules on the affine Grassmannian of $\GL_N$ with equivariant structure with respect to these subgroups are equivalent to the categories of representations of quantum supergroups $U_q(\gl(M|N))$. When $M=0$ and $M=N-1$, this conjecture was proved due to \cite{[G],[BFT0]}, respectively. In this paper, we prove the other cases.

We adapt the global method originating from \cite{[GN],[G]}. In order to compare the global definition of Gaiotto category with the local definition, we generalize the local-global comparison theorem of \cite{[G3]} to a general setting. 
\end{abstract}
\maketitle


\setcounter{tocdepth}{2}

\tableofcontents
\section{Introduction}
\subsection{Notations}
We denote by $\bfF=\mathbb{C}(\!(t)\!)$ the field of Laurent series and by $\bfO=\mathbb{C}[\![t]\!]$ the ring of formal power series. Given a group scheme $G$, we denote by $G(\bfF)$ the loop group of $G$ such that $\mathbb{C}$-points of $G(\bfF)$ are given by $\mathbb{C}(\!(t)\!)$-points of $G$, and denote by $G(\bfO)$ the arc group of $G$ such that $\mathbb{C}$-points of $G(\bfO)$ are given by $\mathbb{C}[\![t]\!]$-points of $G$. Set $\Gr_G:= G(\bfF)/G(\bfO)$ the affine Grassmannian of $G$.
It is known that $\Gr_G$ is a formally smooth ind-scheme.

\subsection{Reminder on FLE and Gaiotto conjecture}
For a reductive group $G$, it is known, under the name of Geometric Satake equivalence, that the category of $G(\bO)$-equivariant perverse sheaves on $\Gr_G$ is equivalent to the abelian category of finite dimensional representations of its Langlands dual group $\check{G}$. However, this equivalence is not derived and hard to generalize to the quantum case. 

According to D. Gaitsgory and J. Lurie, a suitable replacement of the above equivalence is to replace the category on the perverse sheaf side by the Whittaker category. Namely, let $U(\bF)$ be the loop group of the unipotent radical of a Borel subgroup of $G$, and let $\chi$ be a non-degenerate character of $U(\bF)$, then the category $D^{U(\bF),\chi}(\Gr_G)$, i.e., the derived category of $(U(\bF),\chi)$-equivariant D-modules on $\Gr_G$, is equivalent to the derived category of finite dimensional representations $\Rep(\check{G})$. 

Furthermore, this equivalence is also true in the quantum case. Assume $G$ to be a simple group over $\BC$. The category $D^{U(\bF),\chi}(\Gr_G)$ has a canonical deformation $D_{\kappa}^{U(\bF),\chi}(\Gr_G)$, i.e., the category of $(U(\bF),\chi)$-equivariant $\kappa$-twisted D-modules on the determinant line bundle, for $\kappa\in \BC/\BZ$. It is equivalent to the category $\Rep_q(\check{G})$ of representations of the quantum group, ref \cite{[G]} for generic case and \cite{[CDR]} for roots of unity case. Here, the quantum parameter $q= \exp(\pi \cdot i\cdot \kappa)$.

In the case of $\GL_N$, D.Gaiotto and \cite[(2.11)]{[JPS]} and \cite[Section 2, Lecture 5]{[C]} constructed a series of subgroups with an additive character $\chi$ which are related to manifestation of S-duality of supersymmetric boundary conditions. For $M=N-1$, the subgroup constructed by D.Gaiotto is just $\GL_M$, and for $M=0$, the subgroup is the unipotent group $U_{\GL_N}$ {appearing} in the fundamental local equivalence. For $0<M<N-1$, the subgroup $\GL_M\ltimes U^-_{M,N}$ is a semi-direct product of a reductive subgroup and {a} unipotent subgroup. It is conjectured that the category of twisted D-modules on $\Gr_{\GL_N}$ with equivariant structure with respect to this subgroup is closely related to the category of finite-dimensional representations of the quantum supergroup $U_q(\gl(M|N))$.

To be more precise, it is conjectured by D. Gaiotto, cf.~\cite[Section 2]{[BFGT]}, that
\begin{conj}\label{conj}
If $\kappa$ is generic, then we have a $t$-exact braided monoidal equivalence
\[\SD_{\kappa}^{\GL_M(\bO)\ltimes U^-_{M,N}(\bF),\chi}(\Gr_N)\simeq \Rep_q(\GL(M|N)),\]
where $\SD$ denotes the derived category of D-modules with coefficients in the super vector space $\SVect$, and $q= \exp(\pi\cdot i\cdot \kappa)$.
\end{conj}

When $M=0$, it is a statement about classical reductive group, and it is proved by Gaitsgory in \cite{[G]}, and when $M=N-1$, it is proved by Braverman-Finkelberg-Travkin in \cite{[BFT0]}. The main theorem (Theorem \ref{statement of main}) of this paper says that the above conjecture also holds in the other cases, i.e., $0<M<N-1$. Also, we prove that the category $\SD_{\kappa}^{\GL_M(\bO)\ltimes U^-_{M,N}(\bF),\chi}(\Gr_N)$ is the derived category of its heart category, which is different from the usual geometric Satake equivalence.

\subsection{The strategy of the proof}
In this paper, instead of working with the category of representations, we work with the category of factorization modules on the configuration space over a certain factorization algebra $\cI$. This factorization method {follows} \cite{[BFS]}, \cite{[G]}, \cite{[L]}, etc., and turns the statement of Conjecture \ref{conj} into a purely geometric statement. In general, following J. Lurie, one can associate to a Hopf algebra a factorization algebra on the configuration space $C^{\bullet}$, and there is an equivalence between the category of representations and the corresponding category of factorization modules. 

Modulo necessary preparations in Sections \ref{section 2}-\ref{SW-Zastava space}. In Section \ref{constrct functor}, by using the SW-Zastava  (Sakellaridis-Wang Zastava) space, we construct a functor from the Gaiotto category $\SD_q^{\GL_M(\bO)\ltimes U^-_{M,N}(\bF),\chi}(\Gr_N)$ to the category of factorization modules over a certain factorization algebra $\Omega$, where $\Omega$ is defined by applying our functor to the unit object in the Gaiotto category. The main technical difficulty of showing Theorem \ref{statement of main} is to calculate the image of the unit object. To be more precise, we claim (Theorem \ref{key}) that when $\kappa$ is generic, we have an isomorphism
\[\cI\simeq \Omega.\] 

In the case $M=N-1$, the map from the SW-Zastava space to the configuration space is semismall; hence, by construction, one has $\cI\simeq \Omega$ in the generic case. In the case $M=0$, the Zastava space is flat over the configuration space, and the central fiber, whose irreducible components are the Mirković--Vilonen cycles, is equidimensional. Thus one can check $\cI\simeq \Omega$ by proving the vanishing of the sub-bottom cohomology of the MV cycles. This can be further reduced to calculating the sub-bottom cohomology of the central fiber over the point of degree $-n\alpha_i$, where $\alpha_i$ is a simple root of $\GL_N$. In this case, the central fiber is $\mathbb G_m\times \mathbb G_a^{n-1}$, and its cohomology can be computed.

In the case $0<M<N-1$, the SW-Zastava space no longer has the good properties described above. Its central fibers are intersections of $\GL_M(\bF)\ltimes U_{M,N}^-(\bF)$-orbits with semi-infinite orbits in $\Gr_M\times \Gr_N$, and are more complicated.

 Using the factorization structures of $\cI$ and $\Omega$, we can reduce the proof of Theorem \ref{key} to proving the vanishing of the cohomology $H^1$ on the central fiber of the SW-Zastava space over the configuration space, i.e., Proposition \ref{Vanishment}.

In Section \ref{Section 8}, we give a rough upper bound of the dimension of the central fiber over the point with the degree $\alpha$. Here $\alpha$ is a root of the supergroup $\GL(M|N)$, see \ref{section 2.1}. Namely, we show that the dimension of the central fiber is no more than the half of the dimension of  the SW-Zastava space $\cY^\alpha$, see Corollary \ref{cor 8.1.3}.

The cohomology $H^1$ on the central fiber can possibly be non-zero only in two cases. The first case is when the dimension of the central fiber $Y$ equals $\frac{1}{2}(\dim \cY^\alpha-1)$, another case is $\dim Y= \frac{1}{2}\dim \cY^\alpha$.  In the first case, $H^1$ is the bottom cohomology on $Y$, we prove it vanishes by proving the restriction of the D-module on $Y$ is non-constant. In the latter case, $H^1$ is the sub-bottom cohomology, and the proof is much harder due to the complexity of the central fiber. 

To achieve this, we consider several cases according to the types of even simple roots {appearing} in the expression of $\alpha$. When $\alpha$ only contains odd roots, the statement is reduced to the case $M=N-1$ in \cite{[BFT0]}, and it can be solved by the semi-smallness property in \cite{[SW]}. When $\alpha$ contains at least two different even simple roots, it can be solved using the techniques in \cite{[G]}. We only need to study the case when the decomposition of $\alpha$ contains exactly one non-zero coefficient of even simple root of $\GL(M|N)$. The latter case is difficult and needs a sharper estimate of the dimension of the central fiber. We reduce the statement to the case $M=N-2$. 

In the case $M=N-2$, we give the central fiber a stratification by intersecting with another family of semi-infinite orbits, and study the dimension of each stratum. In Proposition \ref{intersection}, we give a more precise dimension estimate for the (stratum of) central fiber. In particular, in the case when the dimension of the central fiber equals the upper bound, Proposition \ref{intersection} tells us that we can choose an open subset of the central fiber, which is $\BG_a$-invariant, and the restriction of the sheaf to this open subset is equivariant against a non-trivial character. In particular, the cohomology of the restriction to this open subset vanishes in all degrees. The complement does not contribute to the $H^1$, so we get the desired vanishing property.

After establishing the isomorphism between $\cI$ and $\Omega$, we can prove that the functor $F^{loc}$ sends irreducible objects to irreducible objects. Again using the factorization property, we only need to prove the vanishing of the bottom cohomology of the restriction to the central fiber. It is not hard to achieve. 

Both sides of Conjecture \ref{conj} acquire monoidal structures given by taking the fusion products and the functor preserves the fusion products. As a result, we obtain an isomorphism of the Grothendieck rings of the twisted Gaiotto category and the category of the factorization modules, i.e., the category of finite dimensional representations of the quantum super group $U_q(\gl(M|N))$.  

To finish the proof, we now only need to repeat \cite[Section 4.3-5]{[BFT0]} verbatim. Namely, denote {by} $\cE$ the full subcategory fusion generated by the irreducible objects corresponding to the tautological representation $V^{taut}=\BC^{M|N}$, $(V^{taut})^*$, and the trivial representation, under taking fusion products, images, kernels and cokernels. It is shown in the $loc.cit$, with a proof of P. Etingof and using the facts that we have already proved, that $\cE$ is equivalent to the category of representations of the quantum supergroup. We only need to show that the full subcategory $\cE$ is actually equivalent to the whole Gaiotto category. It can be done by choosing a collection of generators such that any object is isomorphic to the image of a morphism between these generators. 

\subsection{Local-Global comparison}
For technical reasons, we need to consider two (equivalent) definitions of the twisted Gaiotto category: local Gaiotto category and global Gaiotto category. The definition of the local Gaiotto category is straightforward, it consists of those D-modules on the affine Grassmannian $\Gr_{\GL_N}$ with $(\GL_M(\bO)\ltimes U_{M,N}^-(\bF), \chi)$-equivariant structure. However, since any orbit is infinite-dimensional, we are not supposed to use the na\"{i}ve $t$-structure. Although by choosing a collection of compact objects, we can define the {`correct'} $t$-structure on the local Gaiotto category, it is not easy to do calculations. For example, it is not easy to see that the local functor $F^{loc}$ is $t$-exact. A more crucial drawback is that we do not know how to prove the local functor commutes with the Verdier duality. 

In order to address these problems, fix a global curve $C$ with a marked point $c$. We use the global model $\cM_{\infty\cdot c}$, i.e., a variant of the Drinfeld compactification (ref \cite[Section 1]{[BG]} and \cite[Section 1.2.3]{[FGV]}), to give a definition of the Gaiotto category. It consists of those {twisted} D-modules on the global model with an $H=\GL_M\ltimes U_{M,N}^-$ generic Hecke equivariant structure against $\chi$. The idea is that $\cM_{\infty\cdot c}$ is like a compactification of $H_{C-c}\backslash {(\Gr_M\times \Gr_N)}$, and $\chi|_{H_{C-c}}=0$. So, by descent, $\SD_{\kappa}^{\GL_M(\bO)\ltimes U_{M,N}^-(\bF),\chi}(\Gr_{\GL_N})\simeq \SD_{\kappa}^{H(\bF),\chi}(\Gr_M\times \Gr_N)$ can be regarded as a {certain category of D-modules} on $H_{C-c}\backslash {(\Gr_M\times \Gr_N)}$. When $M\neq N-1$, {the stack} $H_{C-c}\backslash {(\Gr_M\times \Gr_N)}$ is not algebraic, but {$\cM_{\infty\cdot c}$} is always algebraic and the equivariant structure ensures that any equivariant D-module {has zero !-stalks outside} $H_{C-c}\backslash {(\Gr_M\times \Gr_N)}$ {considered as a constructible subset in $\cM_{\infty\cdot c}$}.

We mimic the definition of \cite{[G3]} and \cite{[GN]} to give the definition of the global Gaiotto category. Namely, we choose some points $\bar{x}$ on the curve and then consider an open substack of $\cM_{\infty\cdot c}$. We can define the category of $(H,\chi)$ Hecke equivariant {twisted} D-modules on this substack, then the global Gaiotto category is defined as the limit category. The main difference is that we consider the derived category here.

In the Whittaker case (which corresponds to the Gaiotto case when $M=0$), \cite{[G2]} proves that the global Whittaker category and the local Whittaker category are equivalent. In Theorem \ref{locglob}, we generalize the local-global comparison theorem to the Gaiotto case.


Furthermore, using exactly the same argument in Section \ref{global Gaiotto} and \ref{Ran}, one can prove a more general local-global comparison theorem, {see} Theorem \ref{general}, which says that the local-global comparison theorem actually holds for any $G, H, \chi$ with some conditions (listed in Section \ref{sec 4.5.2}). In particular, it works for the case when $H$ is reductive and $\chi=0$. 

\subsection{Organization of the paper}
In Section \ref{section 2}, we briefly recall the definition of the quantum supergroup $U_q(\gl(M|N))$. Then, we sketch the definition of the factorization algebra which is a geometric replacement of $U_q(\gl(M|N))$.

In Section \ref{Gaiotto category}, we define the twisted (local) Gaiotto category and study some {of its basic properties}, e.g., we study the set of relevant orbits and irreducible objects.

In Section \ref{global Gaiotto}, we recall the definition of the global model and define the twisted global Gaiotto category, and prove a stratawise equivalence between local and global Gaiotto categories.

In Section \ref{Ran}, we use the Ran-ified (Beilinson-Drinfeld) affine Grassmannian to give a proof of the local-global comparison theorem. That is to say, the local Gaiotto category and the global Gaiotto category are not only stratawise equivalent, but equivalent for the whole category.

In Section \ref{SW-Zastava space}, we recall the SW-Zastava space and its marked point version which will be used in the construction of the functor.

In Section \ref{constrct functor}, the functors $F^{loc}$ and $F^{glob}$ are defined. As an application of the local-global comparison theorem, we see that $F^{loc}$ and $F^{glob}$ are essentially isomorphic.

In Section \ref{Section 8}, we start the proof of the main theorem of this paper. First, we calculate the intersection of $H(\bF)$-orbits and semi-infinite orbits. As a consequence, we will give a upper bound of the dimension of the intersection. 

In Section \ref{section 9}, we prove the key proposition of this paper. Namely, we will prove a vanishing of the sub-bottom cohomology with the preparation in the last section. As a result, we obtain that there is an isomorphism of Grothendieck rings of the Gaiotto category and the factorization module category.

In Section \ref{section 10}, we finish the proof by using the proof in \cite[Section 4.3-5]{[BFT0]}

\subsection{Relation with other works and generalizations}
In this paper, we work with the case $\kappa$ is generic, and one may wonder if it is possible to prove the equivalence for arbitrary $\kappa$ like the classical case, {cf.}~\cite{[FGV],[CDR]}. However, the authors do not know how to formulate the statement of Conjecture \ref{conj} in the root of unity case in general, we need to find an appropriate integral form to define the quantum supergroup. It is very likely that once one defines the quantum supergroup in the case of root of unity, then the equivalence still holds when the order of the quantum parameter avoids small torsion. But when the order is small, we need to modify the representation category to make the equivalence hold, at least in the $q=1$ case. For example, in \cite{[BFGT]} and \cite{[TY]}, the authors proved that when $q=1$, in order to have an equivalence, we need to consider the degenerate supergroup $\underline{\GL}(M|N)$ instead of $\GL(M|N)$, it is different from the classical case {of} \cite{[FGV]}.

One may also want to consider the Gaiotto conjecture for other quantum supergroups, such as the orthosympectic supergroups $\mathfrak{osp}(2k|2l)$ and $\mathfrak{osp}(2k+1|2n)$, {cf}.\ \cite{[BFT]} and \cite{[BFT2]}. Actually, our argument can be easily modified for the $\mathfrak{osp}{(2k|2l)}$ case to prove \cite[Conjecture 3.2.1]{[BFT]}. The only difference is that the fiber description \eqref{fiber} in the proof of Proposition \ref{intersection} is more complicated, and the proof of Proposition \ref{intersection} in the orthosymplectic case is more cumbersome. For the convenience of the readers who are interested in the orthosymplectic case, we sketch another proof of Proposition \ref{intersection} in Remark \ref{another proof 9.4.2}, which works verbatim for the orthosymplectic case.

The main theorem of this paper is the quantum version of a particular case of the Periods--$L$-functions duality conjecture of D. Ben-Zvi, Y. Sakellaridis and A. Venkatesh, cf.\ \cite[Conjecture 7.5.1]{[BZSV]}. The Periods--$L$-functions duality conjecture claims that there is an equivalence between the category of $G(\bO)$-equivariant D-modules on the loop space of a $G$-spherical variety $\cX$ and the category of $G^\vee$-equivariant quasi-coherent sheaves on the dual hyperspherical variety of the (twisted) cotangent bundle of $\cX$. One may wonder if one can also consider the quantum version of the general Periods--$L$-functions duality conjecture. To make sense of the quantum generalization, we need to have a reasonable deformation of $D^{G(\bO)}(\cX(\bF))$. In our case, $D^{G(\bO)}(\cX(\bF))$ is equivalent to $D^{H(\bF),\chi}(\Gr_G)$ and the deformation is given by considering twisted D-modules on $\Gr_G$. Except the Gaiotto cases mentioned in \cite[Section 3]{[BFT2]}, we do not know any other example which has a reasonable deformation.

Another generalization direction is to consider the category of twisted Gaiotto D-modules on the affine flags, like \cite{[Yan]}, and relate it with the category $O$ of the (mixed) quantum supergroup. In the upcoming paper, we will prove that there is an equivalence between a full subcategory of the twisted Gaiotto category on $\Fl_G$ and the subcategory $\mathcal{O}$ of the  mixed quantum supergroup generated by Verma modules. 
\subsubsection{Conventions}
In this paper, the schemes and stacks are defined over $\BC$ {(the field of complex numbers)}, and we will deal with D-modules. However, the same argument also works for any sheaf theory listed in \cite[0.8.8]{[GL]}, and we need to replace the twisting $\kappa$ by a character sheaf in the corresponding sheaf theory.

Let $\Vect$ be the DG category of chain complexes over $\BC$.  Similarly, we can consider  $\SVect$, the  category of chain complexes in super vector spaces. The categories considered in this paper are $\BC$-linear DG categories. For two objects $c_1, c_2\in \cC$, we will denote by $\RHom_{\cC}(c_1,c_2)\in \SVect$ the internal Hom, and denote by $\Hom_{\cC}(c_1, c_2)= H^0(\RHom_{\cC}(c_1,c_2))$.

\subsection*{Acknowledgments}
We would like to express our sincere gratitude to Michael Finkelberg and Alexander Braverman for their invaluable contributions to this paper. Their insightful discussions and generous assistance greatly influenced the development of our work. We are also thankful to Michael Finkelberg for his constant support throughout the project, and for providing us with constructive comments that significantly improved the quality of our research. Additionally, we would like to express our appreciation to Alexander Braverman for bringing to our attention the conjecture by D. Gaiotto. 

We thank Tsao-Hsien Chen for pointing out a mistake in the former version.

This paper belongs to a long series of papers about mirabolic Hecke algebras and later about their categorification, with Victor Ginzburg and Michael Finkelberg.   We are grateful to them for sharing many ideas.
Finally, we are thankful to Dennis Gaitsgory for teaching us a lot of techniques of the geometric Langlands theory.

\section{Supergroup and factorization algebra}\label{section 2}
\subsection{Supergroup}\label{section 2.1}
Let $q\in \BC^\times$ be a quantum parameter which is not a root of unity. Given such a number, according to  \cite{[CHW],[Yam]}, we can consider the quantum supergroup $U_q(\gl(M|N))$. Following the notations of \cite{[BFT0]}, we denote by $\Rep_q(\GL(M|N))$ the derived category of finite dimensional representations of $U_q(\gl(M|N))$ equipped with a grading given by the weight lattice $\Lambda$ of $\GL(M|N)$. This category acquires a naturally defined braided tensor category structure, given by the usual tensor product of modules. 

In this paper, we will assume $M<N$. Let $\delta_1,..., \delta_M, \epsilon_1,..., \epsilon_N$ be a basis of diagonal entries weights of the diagonal Cartan subgroup of $U_q(\gl(M|N))$. In this basis, the weight lattice $\Lambda$ of $U_q(\gl(M|N))$ is exactly $ \BZ^M\times\BZ^N$. 

Unlike the case of reductive groups, the Borel super subgroups of a supergroup are not necessarily conjugate to each other.
The definition of the quantum supergroup $U_q(\gl(M|N))$ involves a choice of the Borel super subgroup. Following \cite{[BFT0]}, we choose the mixed Borel super subgroup whose number of isotropic simple roots is as large as possible ($=2M$). Namely, the positive simple roots of $U_q(\gl(M|N))$ with respect to the mixed Borel super subgroup are:
\begin{gather*} 
\alpha_1=-\epsilon_1+\delta_1,\  \alpha_2=-\delta_1+\epsilon_2,\  \alpha_3=-\epsilon_2+\delta_2,\  \cdots,\  \alpha_{2M}=-\delta_M+\epsilon_{M+1},\\
\alpha_{2M+1}=-\epsilon_{M+1}+\epsilon_{M+2},\  \alpha_{2M+2}=-\epsilon_{M+2}+\epsilon_{M+3},\ \cdots,\  \alpha_{M+N-1}=-\epsilon_{N-1}+\epsilon_N.
\end{gather*}

By a dominant weight $(\lambda, (\theta, \theta'))$ of $\GL(M|N)$, we mean a dominant weight of its even part $\GL_M\times \GL_N$, i.e., a pair of a length $M$ non-decreasing integer sequence $\lambda=(\lambda_1, ...\lambda_M)$ and a length $N$ non-decreasing integer sequence $(\theta,\theta')=(\theta_1, ..., \theta_{M+1}, \theta_1',\cdots,\theta_{N-M-1}')$. According to \cite[Exercise at page 141]{[S]}, a dominant weight $(\lambda, (\theta,\theta'))$ is the highest weight of an irreducible finite dimensional representation of $U_q(\gl(M|N))$ if and only if the following condition holds:
\begin{equation}\label{condition sig}
    \begin{split}
        \text{if } \theta_i=\theta_{i+1} \text{,  then } \theta_i+\lambda_i=0, \quad \text{for } i=1,2,..., M;\\
        \text{if } \lambda_{i-1}=\lambda_i \text{, then } \theta_i+\lambda_i=0, \quad\text{for } i=2,..., M.
    \end{split}
\end{equation}

\subsection{Configuration space}
After \cite{[BFS]}, it is known that the category of representations of a group can usually be realized as the category of factorization modules over the corresponding factorization algebra. Hence, in this paper, following the method originated from \cite{[BFS]} and \cite{[G]}, instead of comparing the twisted Gaiotto category and the category of finite dimensional representations of $U_q(\gl(M|N))$, we will construct a functor from the twisted Gaiotto category to the category of twisted factorization modules on the configuration space and show that the constructed functor is an equivalence. In this section, we will recall the definition of the configuration space.

Let $C$ be a smooth connected curve over $\BC$, and let {$\Lambda^{\textnormal{neg}}=\sum_i -\BN\alpha_i\subset\Lambda$} be the {cone} of negative roots. For a negative root $\alpha$ of $\GL(M|N)$, we denote by $C^\alpha$ the scheme which classifies $\Lambda^{\textnormal{neg}}$-colored divisors on $C$ of total degree $\alpha$. (This scheme is denoted by $\Conf(C, \Lambda^{\textnormal{neg}})^{\alpha}$ in the works such as \cite{[GL], [G2], [G3]}.) It is called the configuration space of degree $\alpha$. We denote by $C^\bullet$ the disjoint union of all $C^\alpha$. If $\alpha= -\sum_{i=1}^{M+N-1} n_i \alpha_i$, then we have
\[C^\alpha\simeq \prod_{i=1}^{M+N-1} C^{(n_i)}.\]
Here, $C^{(n_i)}$ denotes the unordered configuration space of $n_i$-points in $C$. For a D-module $\cF$ on $C^\bullet$, we denote by $\cF^\alpha$ the restriction of $\cF$ on the connected component $C^\alpha$. 

The key feature of $C^\bullet$ is that it acquires a non-unital monoid structure. Namely, {addition} defines a map:
\begin{equation}\label{2.2}
    \begin{split}
       \add: C^\bullet\times C^\bullet&\longrightarrow C^\bullet
   \\ D, D'&\mapsto D+D'.
    \end{split}
\end{equation}

It is easy to see that the map $\add$ is compatible with the degree grading, i.e., \[\add: C^\alpha\times C^{\alpha'}\longrightarrow C^{\alpha+\alpha'}.\]

With respect to this monoid structure, we can consider the notion of factorization algebras. Let $(C^\bullet\times C^\bullet)_{disj}$ be the open sub-scheme of $C^\bullet\times C^\bullet$, such that the supports of $D, D'$ are disjoint. For a D-module $\cI$ on $C^\bullet$, we call it a factorization algebra if there is an isomorphism
\begin{equation}\label{fact alg}
    \add^!(\cI)|_{(C^\bullet\times C^\bullet)_{disj}}\simeq \cI\boxtimes \cI|_{(C^\bullet\times C^\bullet)_{disj}},
\end{equation}
with higher homotopy coherence.

To define the notion of a \textit{factorization module} over a given factorization algebra defined on $C^\bullet$, we need to recall the configuration space $C_{\infty \cdot c}^\bullet$ with a marked point.

Let $c\in C$ be a fixed point and $(\xi, (\eta, \eta'))\in \Lambda$ be a weight of $\GL(M|N)$. We denote by $C^{(\xi, (\eta, \eta'))}_{\infty \cdot c}$ the space which classifies $\Lambda$-colored divisors on $C$ such that the total degree is $(\xi, (\eta, \eta'))$ and the coefficient of {any} $x\neq c$ is a negative root.  Any point of $C^{(\xi, (\eta, \eta'))}_{\infty\cdot c}$ can be written as 
\[D= (\xi_c, (\eta, \eta')_c)\cdot c+ \sum \alpha_x \cdot x=-\sum_{i=1}^{M} \delta_i\cdot \Delta_i+ \sum_{i=1}^{N} \epsilon_i\cdot E_i,\]
where $x\in C-c$, $\alpha_x\in \Lambda^{\textnormal{neg}}$ {, and $\Delta_i, E_i$ are divisors in $C$ (corresponding to the coefficients of the expression of $D$ under the basis $\delta_i, \epsilon_i$)}. The space $C^{(\xi, (\eta, \eta'))}_{\infty\cdot c}$ acquires an ind-scheme structure of {ind-}finite type.
Consider the sub-scheme $C^{(\xi, (\eta, \eta'))}_{\leq (\lambda, (\theta, \theta'))\cdot c}$ of $C^{(\xi, (\eta, \eta'))}_{\infty\cdot c}$ such that $(\xi_c, (\eta, \eta')_c)\leq (\lambda, (\theta, \theta'))$. Then, for $(\lambda, (\theta, \theta'))\leq (\tilde{\lambda}, (\tilde{\theta}, \tilde{\theta}'))$, the scheme $C^{(\xi, (\eta, \eta'))}_{\leq (\lambda, (\theta, \theta'))\cdot c}$ is a closed sub-scheme of $C^{(\xi, (\eta, \eta'))}_{\leq (\tilde{\lambda}, (\tilde{\theta}, \tilde{\theta}'))\cdot c}$, and under this transition map, we have
\[C^{(\xi, (\eta, \eta'))}_{\infty\cdot c}= \colim\  C^{(\xi, (\eta, \eta'))}_{\leq (\lambda, (\theta, \theta'))\cdot c}.\]
If $(\lambda, (\theta, \theta'))-(\xi, (\eta, \eta'))= \sum n_i \alpha_i$, then there is an isomorphism 
\[C^{(\xi, (\eta, \eta'))}_{\leq (\lambda, (\theta, \theta'))\cdot c}\simeq \prod_{i=1}^{M+N-1} C^{(n_i)}.\] 
We denote by $C^\bullet_{\infty \cdot c}$ the union of $C^{(\xi, (\eta, \eta'))}_{\infty \cdot c}$ for all $(\xi, (\eta, \eta'))\in \Lambda$.

While $C^\bullet$ is a monoid by \eqref{2.2}, $C_{\infty \cdot c}^\bullet$ is a space with an action of $C^\bullet$. Namely, the {addition} defines a degree-graded map
\begin{equation}
    \begin{split}
       \add_c: C^\bullet\times C^\bullet_{\infty \cdot c}&\longrightarrow C_{\infty \cdot c}^\bullet
   \\ D, D'&\mapsto D+D'.
    \end{split}
\end{equation}
Let $(C^\bullet\times C^\bullet_{\infty \cdot c})_{disj}$ be the open subspace with the disjoint support condition. 

A factorization module over the factorization algebra $\cI$ is a constructible sheaf $\cF$ with coefficients in $\SVect$ on $C^\bullet_{\infty \cdot c}$, such that
\begin{equation}\label{fact mod}
    \add_c^!(\cF)|_{(C^\bullet\times C_{\infty \cdot c}^\bullet)_{disj}}\simeq \cI\boxtimes \cF|_{(C^\bullet\times C_{\infty \cdot c}^\bullet)_{disj}},
\end{equation}
with higher homotopy coherence.

Following \cite{[BFS]}, in order to obtain a category which is equivalent to the category of finite dimensional representations of $U_q(\gl({M|N}))$, we need to impose the following \textit{finiteness} properties:

\begin{enumerate}
    \item $\cF^{(\xi, (\eta, \eta'))}\neq 0$ only for $(\xi, (\eta, \eta'))$ {belonging} to finitely many cosets of the root lattice.
    \item For any such coset, there is a weight $(\lambda, (\theta, \theta'))$ such that the support of $\cF^{(\xi, (\eta, \eta'))}$ lies in $C^{(\xi, (\eta, \eta'))}_{\leq (\lambda, (\theta, \theta'))\cdot c}$ for a  certain $(\xi, (\eta, \eta'))$ in this coset.
    \item There are only finitely many $(\xi, (\eta, \eta'))$ such that the singular support of $\cF^{(\xi, (\eta, \eta'))}$ contains the conormal bundle to the point $(\xi, (\eta, \eta'))\cdot c\in C^{(\xi, (\eta, \eta'))}_{(\lambda, (\theta, \theta'))\cdot c}$.
\end{enumerate}

\subsection{Factorizable line bundle}\label{sec 2.3}
In this section, we will define a line bundle on the configuration space which is compatible with respect to the factorization structure.

By \cite[Section 3.10.3]{[BD]}, we can construct a canonical line bundle $\cP_{C^\bullet}$ on $C^\bullet_{\infty \cdot c}$ and $C^\bullet$ with fiber
\begin{equation}
    \begin{split}
       \cP_D=\qquad &\bigotimes_{i=1}^{M} \det R\Gamma(C, \mathcal{O}_C(-\Delta_i))\otimes \det^{-1} R\Gamma(C, \mathcal{O}_C)\\
       \otimes &\bigotimes_{i=1}^{M+1} \det^{-1} R\Gamma(C, \mathcal{O}_C(-E_i))\det R\Gamma(C, \mathcal{O}_C)\\
       \otimes &\bigotimes_{i=M+2}^N (\det^{-1} R\Gamma(C, \omega_C^{\otimes(i-M-1)}(-E_i))\otimes \det R\Gamma(C, \omega_C^{\otimes(i-M-1)}))
    \end{split}
\end{equation}
at the point \[D= -\sum_{i=1}^{M} \delta_i\cdot \Delta_i+ \sum_{i=1}^{N} \epsilon_i\cdot E_i.\]

In other words, the fiber of $\cP_{C^\bullet}$ at 
\[D= \sum_{x\in C}\sum_{i=1}^{M} \xi_{i,x}\delta_i\cdot x+ \sum_{x\in C}\sum_{i=1}^{M+1} \eta_{i,x}\epsilon_i\cdot x+\sum_{x\in C}\sum_{i=M+2}^{N} \eta'_{i,x}\epsilon_i\cdot x\]
is 
\begin{equation}
    \begin{split}
       \cP_D=\qquad &\bigotimes_{x\in C}\bigotimes_{i=1}^{M} \omega_x^{-\xi_{i,x}(\xi_{i,x}+1)/2}\\
       \otimes &\bigotimes_{x\in C}\bigotimes_{i=1}^{M+1} \omega_x^{\eta_{i,x}(\eta_{i,x}-1)/2}\\
       \otimes &\bigotimes_{x\in C}\bigotimes_{i=M+2}^N \omega_x^{\eta'_{i,x}( \eta'_{i,x}+2(i-M-1)-1)/2}.
    \end{split}
\end{equation}

The reason for the renormalization by tensoring $\det R\Gamma(C, \omega_C^{\otimes(i-M-1)})$ is that the renormalized line bundle $\cP_{C^\bullet}$ is canonically trivial on $C^{-\alpha_i}$ for any simple root $\alpha_i$, whereas the non-renormalized line bundle is not.

By definition, it is easy to see that $\cP_{C^\bullet}$ satisfies the factorization property. That is to say,
\begin{equation*}
  \add^!(\cP_{C^\bullet})|_{(C^\bullet\times C^\bullet)_{disj}}\simeq \cP_{C^\bullet}\boxtimes \cP_{C^\bullet}|_{(C^\bullet\times C^\bullet)_{disj}}.  
\end{equation*}

\subsection{The twisted factorization algebra}
We denote by ${\SD}_{\kappa}(C^\bullet)$ and ${\SD}_{\kappa}(C^\bullet_{\infty \cdot c})$ the category of $\cP_{C^\bullet}^{\kappa}$-twisted D-modules on $C^\bullet$ and $C^\bullet_{\infty\cdot c}$ with coefficients in the category of super vector spaces $\SVect$. According to the factorization property of $\cP_{C^\bullet}$, we can consider the twisted factorization algebra and factorization module, i.e., the objects in ${\SD}_{\kappa}(C^\bullet)$ and ${\SD}_{\kappa}(C^\bullet_{\infty \cdot c})$ equipped with the factorization isomorphisms in \eqref{fact alg} and \eqref{fact mod}, c.f., \cite[3.5]{[G]}, \cite[3.3]{[GL]}, \cite[2.5]{[BFT0]}, etc.. 

For a negative root $\alpha$, the configuration space $C^\alpha$ acquires a naturally defined stratification given by the pattern (i.e., partition of $\alpha$).  We denote the unique open dense sub-scheme by $\overset{\circ}{C^\alpha}$. It consists of those $\Lambda^{\textnormal{neg}}$-colored divisors whose coefficient of any point $x\in C$ is a negative simple root, i.e., it consists of $D= -\sum \alpha_{i_x}\cdot x$, where $\alpha_{i_x}$ is a simple root. If we have $\alpha= -\sum n_i \alpha_i$, then $\overset{\circ}{C^\alpha}\simeq \prod_{i=1}^{M+N-1} \overset{\circ}{C^{(n_i)}}$. Here, $\overset{\circ}{C^{(n_i)}}$ is the open sub-scheme of $C^{(n_i)}$ removing all the diagonals. 

According to the factorization property of $\cP_{C^\bullet}$, the line bundle $\cP_{C^\bullet}$ is canonically trivial on $\overset{\circ}{C^\alpha}$. Indeed, there is
\[\add^!(\cP)|_{(C^{-\alpha_1})^{n_1}\times\cdots \times (C^{-\alpha_{M+N-1}})^{n_{M+N-1}})_{disj}}\simeq \cP\boxtimes\cdots \boxtimes \cP|_{(C^{-\alpha_1})^{n_1}\times\cdots \times (C^{-\alpha_{M+N-1}})^{n_{M+N-1}})_{disj}}\]
 and $\cP$ is canonically trivial on $C^{-\alpha_i}$ for any negative simple root $-\alpha_i$.

 \begin{defn}\label{def of I}

 We define a factorization algebra $\cI$ as follows. 
 
 \begin{itemize}
     \item If $\alpha=-n\cdot \alpha_i$, for $\alpha_i$ an odd simple root, we denote by $\cI^\alpha$ the $!*$-extension of the perverse constant $\cP_{C^\bullet}^{\kappa}$-twisted local system on $\overset{\circ}{C^{\alpha}}$ to $C^\alpha$. {We place $\cI^\alpha$ in the parity $n$  $(\on{mod} 2)$ with respect to the super-grading.}
     \item If $\alpha=-n\cdot \alpha_i$, for $\alpha_i$ an even simple root, we denote by $\cI^\alpha$ the $!*$-extension of the perverse sign $\cP_{C^\bullet}^{\kappa}$-twisted local system on $\overset{\circ}{C^{\alpha}}$ to $C^\alpha$. {We place $\cI^\alpha$ in even parity with respect to the super-grading.}
     \item By factorization property, we can obtain a twisted D-module on $\overset{\circ}{C^\alpha}$ for any $\alpha$. We denote by $\cI^\alpha$ the $!*$-extension of the resulting {twisted} D-module.
 \end{itemize}
      
 \end{defn}
 \begin{rem}
 In the work \cite{[BFT0]}, since all simple roots are odd, the factorization algebra in the $loc.cit$ is $!*$-extension of the twisted constant D-module, and in \cite{[G]}, all simple roots are even, the factorization algebra in the $loc.cit$ is the $!*$-extension of the twisted sign local system.
 \end{rem}

According to the construction, the twisted  D-module $\cI$ is factorizable. So, we can consider the category of twisted factorization modules over $\cI$.

\begin{defn}
We denote by $\cI-\FM^{fin}$ the derived category of $\cP_{C^\bullet}^{\kappa}$-twisted perverse sheaves with coefficients in $\SVect$ on $C^\bullet_{\infty\cdot c}$ which are factorizable with respect to $\cI$ and satisfy the finiteness property in \cite{[BFS]}.
\end{defn}

\begin{rem}
Note that $\cI-\FM^{fin}$ is different from the category of twisted constructible sheaves which are factorizable with respect to $\cI$. Actually, $\cI-\FM^{fin}$ is not even a full subcategory of the latter. The latter is equivalent to the category $\mathcal{O}$ of the mixed quantum supergroup whose positive part is the Lusztig version quantum supergroup and the negative part is the De Concini-Kac version.   
\end{rem}

By the proof of \cite{[BFS]}, \cite[Theorem 29.2.3]{[GL]}, \cite{[L]}, and \cite[Theorem 1.2.1]{[CF]}, one can obtain the following lemma
\begin{lem}
There is a braided tensor equivalence of categories
\[\Rep_q(\GL(M|N))\simeq \cI-\FM^{fin}.\]

Furthermore, for any dominant weight $(\lambda, (\theta, \theta'))$ which satisfies the condition \eqref{condition sig}, the irreducible $U_q(\gl(M|N))$-module $V_{(\lambda, (\theta, \theta'))}$ goes to the irreducible factorization module $\IC^{\fact}_{(\lambda, (\theta, \theta'))}$ which is the $!*$-extension of the unique twisted factorization module on $C^\bullet_{=(\lambda, (\theta, \theta'))\cdot c}:= \{D\mid D= (\lambda, (\theta, \theta'))\cdot c+{\sum_{x\ne c}{}} \alpha_x\cdot x)\}$.
\end{lem}

\begin{rem}
It is shown in \cite{[L]} that for any $E_2$-category $\cC$, one can associate it with a factorizable \textit{cosheaf} of categories. 

To attach a sheaf, rather than cosheaf, \cite[Section 6]{[CF]} constructs the categorical Verdier duality. In particular, it is shown in \cite[Proposition 6.3.8]{[CF]} that for an $E_2$ DG category $\cC$ (an infinity category over $\Vect$), one can associate it with a topological factorization category $\on{Fact}(\cC)$, such that there is an equivalence between the category of (non-unital) $E_2$-algebras in $\cC$ and the category of factorization algebras in $\on{Fact}(\cC)$, and the corresponding module categories are equivalent. 

The same constructions and arguments in \cite{[CF]} apply to a more general setting, including the super-linear setting. The constructions and proofs in \cite[Section 6]{[CF]} essentially use the fact that $\Vect$ is a presentable, stable symmetric monoidal infinity category, and the tensor product is compatible with respect to colimits. In particular, the construction in \cite{[CF]} also works when we replace the base category $\Vect$ by another presentable, stable symmetric monoidal infinity category, such that the tensor product is compatible with respect to colimits, such as the category of super vector spaces. 

In particular, in our case, \cite[Proposition 6.3.8]{[CF]} says that, for an $E_2$ DG category $\cC$ over the symmetric monoidal category of super vector spaces, the categorical Verdier-duality construction of \cite{[CF]} produces a super-linear topological factorization category, with the same properties as the usual (non-super) setting.
\end{rem}

\section{Gaiotto category}\label{Gaiotto category}
In this section, we will introduce the D-module side of the main theorem. 

Given a subgroup $K$ of $G(\bF)$ and a character 
$\chi\colon K\to \mathbb{G}_a$
, we denote by $D^{K,\chi}(\Gr_G)$ (resp. ${\SD}^{K,\chi}(\Gr_G)$) the category of $(K, \chi^!(\exp))$-equivariant D-modules on $\Gr_G$ with coefficients in $\Vect$ (resp. $\SVect$). Here, $\exp$ is the exponential D-module on $\mathbb{G}_a$. 

\subsection{Definition of the Gaiotto category}
Assuming $M<N$, we denote by ${P}^-_{M+1,1,\cdots,1}$ the negative parabolic subgroup of $\GL_N$ corresponding to the partition $(M+1, 1,1,...,1)$, and denote by ${U}^-_{M,N}(\bfF)$ the unipotent radical of ${P}^-_{M+1,1,\cdots,1}(\bfF)$. Let $H$ (resp. $G$) be the semi-direct product $\GL_M\ltimes {U}^-_{M,N}$ (resp. direct product $\GL_M\times \GL_N$). The group $H$ is a subgroup of $\GL_N$ and $G$. Its loop group $H(\bF)$ naturally acts diagonally on $\Gr_G:=\Gr_M\times \Gr_N$ via 
\begin{itemize}
    \item the projection $H(\bF)\longrightarrow \GL_M(\bF) \curvearrowright \Gr_M$,
    \item and the embedding $H(\bF)\longrightarrow \GL_N(\bF)\curvearrowright \Gr_N.$
\end{itemize}

There is a group morphism from ${U}^-_{M,N}(\bfF)$ to $\bfF$ that sends $(u_{i,j})\in {U}^-_{M,N}(\bfF)$ to $\sum_{i} u_{i+1,i}$. We compose this group morphism with taking residue
\begin{equation}
    \begin{split}
        \bfF&\longrightarrow \mathbb{C}\\
        \sum a_i t^i&\mapsto a_{-1}.
    \end{split}
\end{equation}
The resulting morphism is a character of ${U}^-_{M,N}(\bfF)$ denoted by $\chi$. The conjugation action of $\GL_M(\bfF)$ on ${U}^-_{M,N}(\bfF)$ preserves ${U}^-_{M,N}(\bfF)$. Furthermore, the character $\chi$ is stable under the conjugation action. Hence, $\chi$ gives a character on $H(\bF)$.

\subsubsection{}
The Gaiotto category is defined as ${\SD}^{H(\bF),\chi}(\Gr_G)$. Equivalently, it can also be defined as  ${\SD}^{\GL_M(\bO)\ltimes {U}^-_{M,N}(\bF),\chi}(\Gr_N)$.

Here, ${\GL}_{M}({\bfO}) \ltimes U^-_{M, N}(\bfF)$ can be written of a subgroup of $\GL_N(\bF)$ as follows.

\begin{equation}\label{1.3 GL}
{\GL}_{M}({\bfO}) \ltimes U^-_{M, N}(\bfF)=
\begin{pNiceArray}{wc{2.9em}ccwc{4em}|cccc}[last-col,code-for-last-col = \quad\,]
\Block{4-4}<\huge>{\GL_M(\bO)} & &&&  &&  & & \Block[l]{4-1}{\quad\, M} \\
\\
 & &   &&  &  & & &\\
 &&& 
 &&&&\\ 
0&\Cdots          &   &0& 1&& & &\Block[l]{4-1}{\quad\, N-M} \\
* &      \Cdots    &   &*& * &\Ddots&  & &  \\
\Vdots&          &   &\Vdots& \Vdots &\Ddots&& \\
*&		\Cdots&& *& * & \Cdots  & * &1
\CodeAfter
\UnderBrace[shorten,yshift=0.5ex]{1-1}{8-4}{M}
\UnderBrace[shorten,yshift=0.5ex]{1-5}{8-8}{N-M}
\SubMatrix .{1-8}{4-8}{\}}[right-xshift=1em]
\SubMatrix .{5-8}{8-8}{\}}[right-xshift=1em]
\end{pNiceArray}  
\end{equation}
{}
\\
\subsection{Renormalization}\label{renorm}
However, in order to match the twistings and have a canonically defined character $\chi$, we need to consider renormalized objects ${U}^{-,\omega}_{M,N}(\bF)$, $H^\omega(\bF)$, $\Gr^{\omega}_N$, and $\Gr_G^\omega$.

We denote by $\rho_{M,N}$ the length $(N-M-1)$ sequence $(1,2,..., N-M-1)$. It can be regarded as a cocharacter of $T_{N-M-1}$. Fix a square root $\omega_C^{\otimes \frac{1}{2}}$ of the dualizing sheaf $\omega_C$, we denote by $\omega_C^{\rho}$ the induced $T_{N-M-1}$-bundle of $\omega_C^{\otimes \frac{1}{2}}$ associated with 
\[2\rho_{M,N}: \BG_m\longrightarrow T_{N-M-1}.\]

Fix a point $x\in C$. For any group $K$ with a group morphism $T_{N-M-1}\to K$, such {as} $\GL_N$ and $G$, we denote by $K^{\omega}(\bO)$ (resp.\ $K^{\omega}(\bF)$) the group pro-scheme (resp.\ ind-pro-scheme) of automorphisms of the principal $K$-bundle induced with respect to $T_{N-M-1}\to K$ by the $T_{N-M-1}$-bundle $\omega_C^{\rho}$ on the formal disc 
 $\cD_{x}$ (resp.\ the punctured disc 
 $\overset{\circ}{\cD_{x}}$). 

\begin{rem}
    Sometimes, in order to emphasize the chosen point $x$, we put $x$ in the subscript of the notation.
\end{rem}

\subsubsection{}
Furthermore, if $T_{N-M-1}\to K$ splits as $1\to K'\to K\to T_{N-M-1}\to 1$, such that the composition $T_{N-M-1}\rightarrow K\rightarrow T_{N-M-1}$ is identity, then we define ${K}^{\prime,\omega}(\bO)$ (resp.\ ${K}^{\prime,\omega}(\bF)$) as the kernel of ${K^\omega}(\bO)\to {T}^{\omega}_{N-M-1}(\bO)$ (resp.\ ${K^\omega}(\bF)\to {T}^{\omega}_{N-M-1}(\bF)$). In particular, if we take $K=\GL_M\ltimes B_{M,N}^{-}$ and $B_{M,N}^-$, we can define ${H}^{\omega}(\bO)$, ${H}^{\omega}(\bF)$, ${U}^{-,\omega}_{M,N}(\bO)$, and ${U}^{-,\omega}_{M,N}(\bF)$.

The advantage of $H^\omega(\bF)$ over $H(\bF)$ is that there is a canonically defined character on ${H}^{\omega}(\bF)$ which is non-degenerate on ${U}^{-,\omega}_{M,N}(\bF)$, i.e., we do not need to choose {a} uniformizer. Namely, we define the character $\chi$ of $H^{\omega}(\bF)$ by the following map
\begin{equation}\label{chi}
    \begin{split}
        \chi: H^{\omega}(\bF)\longrightarrow  H^{\omega}(\bF)/ [H^{\omega}(\bF), H^{\omega}(\bF)]\longrightarrow U_{M,N}^{-,\omega}(\bF)/[U_{M,N}^{-,\omega}(\bF), U_{M,N}^{-,\omega}(\bF)]\\ \simeq {\Gamma(\omega_{\overset{\circ}{\cD_c}})\times\cdots \times \Gamma(\omega_{\overset{\circ}{\cD_c}})}
        \overset{\text{sum}}{\longrightarrow} {\Gamma(\omega_{\overset{\circ}{\cD_c}})} \overset{\text{residue}}{\longrightarrow} {\BG_a} .
    \end{split}
\end{equation}

We denote by ${\Gr}^{\omega}_N$ the renormalized affine Grassmannian associated with ${\GL}_N$. {I}t classifies $(\cP_N, \alpha_N)$, where $\cP_N$ is an $\GL_N$-bundle on the formal disc $\cD$, and $\alpha_N$ is an identification of $\cP_N$ with $\cP_N^\omega:= \omega_{C}^{\rho}\overset{T_{N-M-1}}{\times} \GL_N$ on the punctured disc $\overset{\circ}{\cD}$. Similarly, one can define $\Gr_G^\omega$, which classifies $(\cP_G, \alpha)=(\cP_M, \cP_N, \alpha_M, \alpha_N)$, where $\cP_G$ is a $G$-bundle on $\cD$, and $\alpha$ is an identification of $\cP_N$ with $\cP_G^\omega:= \omega_{C}^{\rho}\overset{T_{N-M-1}}{\times} G=(\cP_M^{triv}, \cP_N^\omega)$ on $\overset{\circ}{\cD}$. By definition, we have $\Gr_G^\omega= \Gr_M\times \Gr_N^\omega$.

\subsubsection{}
We consider a (relative) determinant line bundle ${\cP}_{\det}$ over $\Gr_G^\omega=\Gr_M\times \Gr_N^\omega$, such that its fiber over the point $(\cP_M, \cP_N, \alpha_M, \alpha_N)$ is 
\begin{equation}\label{3.4}
\det R\Gamma(C, \cV^M_{\cP_M})\otimes \det^{-1} R\Gamma(C, \cV^N_{\cP_N})\otimes \det^{-1} R\Gamma(C, \cV^M_{\cP^{triv}_M})\otimes \det R\Gamma(C, \cV^N_{\cP_N^\omega}),
\end{equation}
 where $\cV^M_{\cP_M}$ (resp.\ $\cV^N_{\cP_N}$) denotes the {rank $M$ (resp.\ $N$)} vector bundle associated with $\cP_M$ (resp.\ $\cP_N$). 

 We denote by ${\SD}_{\kappa}(\Gr^{\omega}_G)$ the category of ${\cP}_{\det}^\kappa$-twisted D-modules on $\Gr^{\omega}_G$ with coefficients in $\SVect$. 

 {The restriction of $\cP_{\det}$ to $\Gr_N^{\omega}$ along $i: \Gr_N^{\omega}=\{t^0\}\times \Gr_N^{\omega} \rightarrow \Gr_G^\omega$ gives rise to a relative determinant line bundle on $\Gr_N^\omega$. We denote by ${\SD}_{\kappa}(\Gr^\omega_N)$ the category of $\cP_{\det}^{\kappa}$-twisted D-modules on $\Gr^\omega_N$.}
 

\begin{defn}
    We denote by $\cC_{\kappa}^{loc}(M|N)$ the category of (${H}^{\omega}(\bF),\chi$)-equivariant ${\cP}_{\det}^\kappa$-twisted D-modules on $\Gr^{\omega}_{{G}}$. Also, we denote by $\cC_{\kappa}^{loc,lc}(M|N)$ the full subcategory of $\cC_{\kappa}^{loc}(M|N)$ consisting of those objects which become compact after applying the forgetful functor $\cC_{\kappa}^{loc}(M|N)\simeq {\SD}^{\GL_M(\bO)\ltimes U_{M,N}^{-,\omega}(\bF),\chi}_{\kappa}(\Gr_N^\omega)\rightarrow {\SD}_{\kappa}^{U^{-,\omega}_{M,N}(\bF),\chi}(\Gr_N^\omega)$. 
\end{defn}

{
\begin{rem}
    The locally compactness of the local Gaiotto category can be described more explicitly as follows: an object in $\cC_{\kappa}^{loc}(M|N)$ is locally compact if and only if it is supported on finitely many orbits, and the $!$-stalk on each orbit is finite dimensional.
\end{rem}
}

\begin{rem}
The group $H^{\omega}(\bF)$ is isomorphic to $H(\bF)$. In particular, the Gaiotto category defined for $\Gr_G$ is (in fact, canonically) equivalent to the Gaiotto category defined for the renormalized affine Grassmannian. So, later, when we need to do calculations which do not depend on the renormalization with $\omega$, we will perform the calculations in the usual affine Grassmannian. 
\end{rem}

\subsection{Relevant {orbits}}
To study this category, we firstly study ${H}^{\omega}(\bF)$-orbits in $\Gr^{\omega}_{{G}}$. That is equivalent to studying $\GL_M(\bO)\ltimes {U}^-_{M,N}(\bF)$-orbits in $\Gr_N$. The following lemma is proved in \cite{[FGT]}. 

\begin{lem}
The $\GL_M(\bO)$-orbits in $\Gr_{M+1}$ are parameterized by dominant weights $(\lambda, \theta)=(\lambda_1\leq\cdots \leq \lambda_M, \theta_1\leq\cdots \leq \theta_{M+1})$.
\end{lem}
Since we will use many {matrix} calculations in the key step (Section \ref{vanishing of H1}), it is more convenient to express the above lemma in the language of matrices.

Let us denote 
\[\BL_{(\lambda,\theta)}=\begin{pNiceMatrix}
    t^{\lambda_1+\theta_1}        & &  &  \\
    &\Ddots&& \\
          & &   t^{\lambda_M+\theta_M} & \\
         t^{\theta_1}  &\Cdots&  t^{\theta_M} & t^{\theta_{M+1}} 
\end{pNiceMatrix}\in \GL_{M+1}(\bF),\]
 where $t$ is {a} uniformizer of $\bO= \BC[\![t]\!]$. With some abuse of notations, we still denote the projection of $\BL_{(\lambda,\theta)}$ in the affine Grassmannian by the same notation. The above lemma can be rewritten as
\begin{lem}\label{lem M+1}
Any $\GL_M(\bO)$-orbit of $\Gr_{M+1}$ has a representative $\BL_{(\lambda,\theta)}\in \Gr_{M+1}$, such that, $(\lambda, \theta)$ is a dominant weight of length $(M, M+1)$. 

Also, if two dominant weights $(\lambda, \theta)\neq(\tilde{\lambda}, \tilde{\theta})$, then 
\[\BL_{(\tilde{\lambda}, \tilde{\theta})}^{-1}\GL_M(\bO)\BL_{(\lambda, \theta)}\cap \GL_{M+1}(\bO)= \emptyset.\]
\end{lem}

Now, we claim,
\begin{lem}\label{H orbit}
Any $\GL_M(\bO)\ltimes U_{M,N}^-(\bF)$-orbit in $\Gr_N$ has a unique representative of the form $\BL_{(\lambda,(\theta, \theta'))}$, such that $(\lambda,\theta)$ is a dominant weight {of  $\GL_M\times\GL_{M+1}$}, and $\theta'$ is a sequence of length $N-M-1$. Here,
\[\BL_{(\lambda,(\theta, \theta'))}=\begin{pNiceMatrix}
    t^{\lambda_1+\theta_1}       &   &  &  &&&\\
    &\Ddots&& &&&\\
          & &   t^{\lambda_M+\theta_M} &&&& \\
         t^{\theta_1}  &\Cdots&  t^{\theta_M} & t^{\theta_{M+1}}&&&\\
         &&&&t^{\theta'_1}&&\\
         &&&&&\Ddots&\\
         &&&&&&t^{\theta'_{N-M-1}}
\end{pNiceMatrix}\in \GL_{N}(\bF).\]
\end{lem}

\begin{proof}
Let $\sA= (a_{i,j})_{N\times N}\in \GL_N(\bF)$. By column transformation with coefficients in $\bO$, one can find a $N\times N$-matrix $\sA'$ of the form
\[\sA'=\begin{pNiceMatrix}
*    &  \Cdots& * &&&\\
   \Vdots &&\Vdots &&&\\
         &&*&t^{\theta'_1}&&\\
         &&&\Ddots&\Ddots&\\
         *&\Cdots&&&*&t^{\theta'_{N-M-1}}
\end{pNiceMatrix}\in \GL_{N}(\bF),\]
such that $\sA\GL_N(\bO)=\sA'\GL_N(\bO)$.

Then, by row transformation, one can find a $N\times N$-matrix $\sA''$ of the form
\[\sA''=\begin{pNiceMatrix}
   *     &\Cdots   & * &&&\\
   \Vdots &&\Vdots &&&\\
         * &  \Cdots &* &&&\\
         &&&t^{\theta'_1}&&\\
         &&&&\Ddots&\\
         &&&&&t^{\theta'_{N-M-1}}
\end{pNiceMatrix}\in \GL_{N}(\bF),\]
such that $U_{M,N}^-(\bF)\sA'=U_{M,N}^-(\bF)\sA''$.

Then, by Lemma \ref{lem M+1}, one can find a $N\times N$-matrix $\sA'''$ of the form $\BL_{(\lambda, (\theta, \theta'))}$, such that $\GL_M(\bO) \sA''\GL_{M+1}(\bO)=\GL_M(\bO) \BL_{(\lambda,(\theta, \theta'))}\GL_{M+1}(\bO)$. 

In particular, we have $(\GL_M(\bO)\ltimes U_{M,N}^-(\bF)) \BL_{(\lambda,(\theta,\theta'))} \GL_N(\bO)=(\GL_M(\bO)\ltimes U_{M,N}^-(\bF)) \sA \GL_N(\bO)$. That proves the existence of a representative of the form $\BL_{(\lambda, (\theta,\theta'))}$ in each $\GL_M(\bO)\ltimes U_{M,N}^-(\bF)$-orbit of $\Gr_N$.

To prove the uniqueness, we should check: if $(\lambda,(\theta, \theta'))\neq (\tilde{\lambda},(\tilde{\theta},\tilde{\theta}'))$, then
\[\BL^{-1}_{(\tilde{\lambda},(\tilde{\theta},\tilde{\theta}'))}
\bigl( \GL_M(\bO)\ltimes U_{M,N}^-(\bF) \bigr)
\BL_{(\lambda,(\theta,\theta'))}\cap \GL_N(\bO)=\emptyset.\]

To prove this, we assume $g\in \GL_M(\bO)\ltimes U_{M,N}^-(\bF)$, and $\BL^{-1}_{(\tilde{\lambda},(\tilde{\theta},\tilde{\theta}'))} g \BL_{(\lambda,(\theta,\theta'))}\in \GL_N(\bO)$. Since $\BL^{-1}_{(\tilde{\lambda},(\tilde{\theta},\tilde{\theta}'))} g \BL_{(\lambda,(\theta,\theta'))}$ is of the form, \[\begin{pNiceMatrix}
   *      & \Cdots  & * &&&\\
   \Vdots &&\Vdots &&&\\
         &&*&t^{{\theta}_1'-\tilde{\theta}_1'}&&\\
         &&&\Ddots&\Ddots&\\
        * &&\Cdots&&*&t^{\theta_{N-M-1}'-\tilde{\theta}_{N-M-1}'}
\end{pNiceMatrix},\]
the assumption that it belongs to $\GL_N(\bO)$ implies $\theta_i'=\tilde{\theta}_i'$ for any $1\leq i\leq N-M-1$. Furthermore, according to Lemma \ref{lem M+1}, the left-top $(M+1)\times (M+1)$-submatrix of $\BL^{-1}_{(\tilde{\lambda},(\tilde{\theta},\tilde{\theta}'))} g \BL_{(\lambda,(\theta,\theta'))}$ is not in $\GL_{M+1}(\bO)$ if $(\lambda,\theta)\neq (\tilde{\lambda},\tilde{\theta})$.
\end{proof}

Not all the $\GL_M(\bO)\ltimes U_{M,N}^-(\bF)$-orbits in $\Gr_N$ carry a non-zero $(\GL_M(\bO)\ltimes U_{M,N}^-(\bF),\chi)$-equivariant ${\cP}_{\det}^\kappa$-twisted D-module supported on it. We call those orbits which admit non-zero $(\GL_M(\bO)\ltimes U_{M,N}^-(\bF),\chi)$-equivariant ${\cP}_{\det}^\kappa$-twisted D-modules the \textit{relevant} orbits. We denote by $\BO^{(\lambda,(\theta,\theta'))}$ the $\GL_M(\bO)\ltimes U_{M,N}^-(\bF)$-orbit of $\BL_{(\lambda,(\theta,\theta'))}\in \Gr_N$.

\begin{prop}
The $\GL_M(\bO)\ltimes U_{M,N}^-(\bF)$-orbit $\BO^{(\lambda,(\theta,\theta'))}$ is relevant if and only if 
\begin{equation}\label{condition 3.4}
  \begin{split}
    \textnormal{A:\quad}    \lambda_1\leq\lambda_2\leq\cdots\leq \lambda_M,\\
        \theta_{1}\leq \theta_2\leq\cdots\leq\theta_{M+1}\leq\theta'_1\leq\cdots \leq \theta'_{N-M-1},\\
        \textnormal{and B:\quad if }  \theta_i=\theta_{i+1}, \textnormal{then } \theta_i+\lambda_i=0, \textnormal{for } i=1,2,..., M;\\
        \textnormal{if } \lambda_{i-1}=\lambda_i, \textnormal{then } \theta_i+\lambda_i=0, \textnormal{for } i=2,..., M.
    \end{split}
\end{equation}
\end{prop}
\begin{proof}
We only need to check that the condition \eqref{condition 3.4} is equivalent to $\Stab_{\GL_M(\bO)\ltimes U^-_{M,N}(\bF)}(\BL_{(\lambda,(\theta,\theta'))})\subset \Ker\ \chi$ and acts trivially on the fiber of $\overset{\circ}{\cP}_{\det, N}$ over $\BL_{(\lambda,(\theta,\theta'))}\in \Gr_N$. Since $\Stab_{U^-_{M,N}(\bF)}(\BL_{(\lambda,(\theta,\theta'))})$ is unipotent and acts trivially on the fiber, we should prove that the condition \eqref{condition 3.4} is equivalent to 
\begin{equation}
    \begin{split}
        \Stab_{U_{M,N}^-(\bF)}(\BL_{(\lambda,(\theta, \theta'))})\subset \Ker\ \chi,\\
        \textnormal{and } \Stab_{\GL_M(\bO)}(\BL_{(\lambda,(\theta,\theta'))}) \textnormal{ acts trivially on the fiber of $\overset{\circ}{\cP}_{\det, N}$ over $\BL_{(\lambda,(\theta,\theta'))}$}.
    \end{split}
    \end{equation}
Consider the point $\BL_{(\lambda,\theta)}\in \Gr_{M+1}$. It is not hard to see that $\Stab_{\GL_M(\bO)}(\BL_{(\lambda,(\theta,\theta'))})=\Stab_{\GL_M(\bO)}(\BL_{(\lambda,\theta)})$ and $\Stab_{\GL_M(\bO)}(\BL_{(\lambda,(\theta,\theta'))})$ acts trivially on the fiber of $\overset{\circ}{\cP}_{\det, N}$ over $\BL_{(\lambda,(\theta,\theta'))}$ is equivalent to $\Stab_{\GL_M(\bO)}(\BL_{(\lambda,\theta)})$ acts trivially on the fiber of $\overset{\circ}{\cP}_{\det, M+1}$ over $\BL_{(\lambda,\theta)}$, where $\Stab_{\GL_M(\bO)}(\BL_{(\lambda,\theta)})$ denotes the stabilizer group of $\BL_{(\lambda,\theta)}\in \Gr_{M+1}$.

According to \cite[Proposition 4.1.1]{[BFT0]}, the stabilizer group $\Stab_{\GL_M(\bO)}(\BL_{(\lambda,\theta)})$ acts trivially on the fiber over $\BL_{(\lambda,\theta)}$ if and only if $(\lambda,\theta)$ satisfies the condition B of \eqref{condition 3.4}. 

So, it remains to show $\Stab_{U_{M,N}^-(\bF)}(\BL_{(\lambda,(\theta,\theta'))})\subset \Ker \ \chi$ if and only if $\theta_{M+1}\leq \theta'_1\leq\cdots\leq \theta'_{N-M-1}$.

Note that \[\Stab_{U_{M,N}^-(\bF)}(\BL_{(\lambda,(\theta,\theta'))})=\BL_{(\lambda,(\theta,\theta'))} \GL_N(\bO) \BL^{-1}_{(\lambda,(\theta,\theta'))} \cap U_{M,N}^-(\bF).\] 

$\Stab_{U_{M,N}^-(\bF)}(\BL_{(\lambda,(\theta,\theta'))})\subset \ker \ \chi$ if and only if for any $\sA=(a_{i,j})_{N\times N}$ belongs to the stabilizer group, the  $t^{-1}$-coefficient of $\sum_{M+1}^{N-1} a_{i+1, i}$ is zero. In fact, it is even equivalent to the condition that the $t^{-1}$-coefficient of $a_{i+1, i}$ is zero for any $i=M+1, \cdots, N-1$ and $\sA\in \Stab_{U_{M,N}^-(\bF)}(\BL_{(\lambda,(\theta,\theta'))})$. By writing $\BL_{(\lambda,(\theta,\theta'))} \GL_N(\bO) \BL^{-1}_{(\lambda,(\theta,\theta'))}$ as matrix explicitly, one can see that the latter condition is equivalent to 
\begin{equation}
    \begin{split}
        t^{-\theta_{M+1}+\theta'_1}\in \bO;\\
        t^{-\theta'_1+\theta'_2}\in \bO;\\
        \cdots\\
        t^{-\theta'_{N-M-2}+\theta'_{N-M-1}}\in \bO,
    \end{split}
\end{equation}  
which implies the condition A of \eqref{condition 3.4}.
\end{proof}

\subsubsection{}Note that there is an isomorphism
\begin{equation}\label{8.1 eq}
    \GL_{M}(\bO)\ltimes U_{M,N}^-(\bF)\backslash \Gr_N\simeq H(\bF)\backslash \Gr_G,
\end{equation}
Indeed, it follows from 
\[\xymatrix{(\GL_{M}(\bO)\ltimes U_{M,N}^-(\bF))\backslash \Gr_N\ar@{-}[r]^-{\sim}& \GL_M(\bF)\backslash (\GL_M(\bF)\times^{\GL_M(\bO)} U_{M,N}^-(\bF)\backslash \Gr_N)\ar@{-}[d]^-{\sim}\\ H(\bF)\backslash \Gr_G\ar@{-}[r]^-{\sim}&\GL_M(\bF)\backslash (\Gr_M\times U_{M,N}^-(\bF)\backslash \Gr_N).}\]
Alternatively, it comes from the fact that $\GL_M(\bF)$ acts transitively on $\Gr_{M}$, so we can remove the factor $\Gr_M$ and replace the equivariant group by the stabilizer subgroup.

Lemma \ref{H orbit} can be rewritten as
\begin{cor}\label{H orbit'}
Any $H(\bF)$-orbit of $\Gr_G$ has a unique representative
$\sL_{(\lambda,(\theta, \theta'))}= (t^{-\lambda}, \sD t^{(\theta, \theta')})\in \Gr_M\times \Gr_N=\Gr_G$, where
\[t^{-\lambda}=\begin{pNiceMatrix}
    t^{-\lambda_1}        & &  & \\
    &\Ddots&& \\
          & &   t^{-\lambda_M} &
\end{pNiceMatrix}\in \Gr_{M},\]
\[t^{(\theta, \theta')}=\begin{pNiceMatrix}
    t^{\theta_1}        & &  &  &&&\\
    &\Ddots&& &&&\\
         &  &   t^{\theta_{M+1}}&&&\\
         &&&t^{\theta'_1}&&\\
         &&&&\Ddots&\\
         &&&&&t^{\theta'_{N-M-1}}
\end{pNiceMatrix}\in \Gr_{N},\]
and \[\sD=\begin{pNiceMatrix}
    1      &  & &  &  &&&\\
          & 1 &  &  &  \\
    &&\Ddots&& &&&\\
          & &  & 1 &&&& \\
        1 &1  &\Cdots& 1 & 1&&&\\
         &&&&&1&&\\
         &&&&&&\Ddots&\\
         &&&&&&&1
\end{pNiceMatrix}\in \GL_{N},\]
such that $\lambda$ is a sequence of length $M$, $(\theta, \theta')$ is a sequence of length $N$. Here, the $(M+1)$-th row of $\sD$ is $(1,1,\cdots,1, 0,\cdots,0)$. {It is relevant if and only if the condition \ref{condition 3.4} holds.}
\end{cor}
We denote the $H(\bF)$-orbit of $\sL_{(\lambda,(\theta, \theta'))}$ by $\sO^{(\lambda,(\theta, \theta'))}$.

\subsubsection{}
We denote by $B:= B^-_M\times \sD B_N \sD^{-1}$, and by $U$ its unipotent radical. Here, $B_M^{-}$ is the lower-triangular Borel subgroup of $\GL_M$ and $B_N$ is the upper-triangular Borel subgroup of $\GL_N$. We denote by $\sS^{(\xi, (\eta,\eta'))}\subset \Gr_G$ the $U(\bF)$-orbit of $\sL_{(\xi, (\eta,\eta'))}=(t^{-\xi}, \sD t^{(\eta, \eta')})$. 
\subsection{Irreducible objects}


We denote $2\rho^{\circ}:= (\lambda^\circ, (\theta^\circ, \theta^{\circ{}'}))$, where 
\[\lambda^\circ=(-M+1, -M+3,\cdots, M-3, M-1),\]
\[\theta^\circ=(-M, -M+2, \cdots, M-2, M),\]
\[\theta^{\circ{}'}=(M+2, M+4, \cdots, 2N-M-2).\]

Any $\GL_M(\bO)$-orbit of $\Gr_N$ is finite dimensional. Actually, one can check that, for relevant $(\lambda, \theta)$, the dimension of the $\GL_M(\bO)$-orbit of $\BL_{(\lambda, (\theta,\theta'))}\in \Gr_N$ is $\langle (\lambda^\circ, \theta^\circ), (\lambda, \theta)\rangle$. 

However, any $H(\bF)$-orbit (or, $U_{M,N}^-(\bF)$-orbit) in $\Gr_G$ is of infinite dimensional, we cannot use the perverse $t$-structure of the category of D-modules on $\Gr_G$, otherwise, some distinguished object, such as $\Av_!^{H(\bF)/\GL_M(\bF), \chi}(\delta_0)$, lives in the infinite negative degree.

Mimicking the definition of the $t$-structure of the Whittaker category in \cite{[GL]}, we can consider a collection of compact objects in $\cC^{loc}_{\kappa}(M|N)$ and define the $t$-structure with respect to those compact objects. It gives rise to a well-behaved $t$-structure which is compatible with the perverse $t$-structure on the global model in Section \ref{global Gaiotto}. 

First, we need the following lemma.
\begin{lem}\label{lem 3.4.1}
    The (partially) left adjoint functor $\Av_!^{H^\omega(\bF),\chi}$ {of the forgetful functor $\oblv_{H^\omega(\bF),\chi}\colon {\SD}_{\kappa}^{H^\omega}(\bF)(\Gr_G^\omega)\rightarrow {\SD}_{\kappa}(\Gr_G^\omega)$} is well-defined on the skyscraper D-module $\delta_{\sL_{(\lambda,(\theta,\theta'))}}$, and {$\Av_!^{H^\omega(\bF),\chi}(\delta_{\sL_{(\lambda,(\theta,\theta'))}})$} corresponds to $\Av_!^{\GL_M(\bO)\ltimes U_{M,N}^{-,\omega}(\bF),\chi}(\delta_{\BL_{(\lambda,(\theta,\theta'))}})$ under the equivalence
    \begin{equation}\label{eq 3.9}
        {\SD}_{\kappa}^{H^\omega(\bF),\chi}(\Gr_G^\omega)\simeq {\SD}_{\kappa}^{\GL_M(\bO)\ltimes U_{M,N}^{-,\omega}(\bF),\chi}(\Gr_N^\omega).
    \end{equation}
\end{lem}
\begin{proof}
    First, we choose a sequence of closed group schemes of pro-finite type $H_i= \GL_{M}(\bO)\ltimes U^-_{M,N,i}$ of $\GL_M(\bO)\ltimes U_{M,N}^{-,\omega}(\bF)$, such that $H_i\subset H_{i'}$ if $i<i'$, and $\GL_M(\bO)\ltimes U_{M,N}^{-,\omega}(\bF)= \colim H_i$. For a fixed $i$, {we can choose a finite dimensional $H_i$-invariant subscheme $Y_i$ of $\Gr_N^\omega$, which contains $\BL_{(\lambda,(\theta,\theta'))}$ such that $H_i$ acts on $Y_i$ via a finite type group $\bar{H}_i$ and the kernel of $H_i\to \bar{H}_i$ is pro-unipotent.} {Since} $\delta_{\BL_{(\lambda,(\theta,\theta'))}}$ is holonomic, the {object} $\Av_!^{\bar{H}_i,\chi}(\delta_{\BL_{(\lambda,(\theta,\theta'))}})$ is well-defined. According to the assumption on $\bar{H}_i$, $\Av_!^{\bar{H}_i,\chi}(\delta_{\BL_{(\lambda,(\theta,\theta'))}})$ has a $(H_i,\chi)$-equivariant structure and it gives the object $\Av_!^{{H}_i,\chi}(\delta_{\BL_{(\lambda,(\theta,\theta'))}})$. Then, the object $\Av_!^{\GL_M(\bO)\ltimes U_{M,N}^{-,\omega}(\bF),\chi}(\delta_{\BL_{(\lambda,(\theta,\theta'))}})$ {can be obtained as} the colimit of $\Av_!^{H_i,\chi}(\delta_{\BL_{(\lambda,(\theta,\theta'))}})$.

    Then, we need to show that the object in ${\SD}_{\kappa}^{H^\omega(\bF),\chi}(\Gr_G^\omega)$ corresponding to $\Av_!^{\GL_M(\bO)\ltimes U_{M,N}^{-,\omega}(\bF),\chi}(\delta_{\BL_{(\lambda,(\theta,\theta'))}})$ is the object $\Av_!^{H^\omega(\bF),\chi}(\delta_{\sL_{(\lambda,(\theta,\theta'))}})$. The equivalence \eqref{eq 3.9} is given by taking $!$-pullback along the inclusion $i: \Gr_N^\omega= \{t^0\}\times \Gr_N^\omega \longrightarrow \Gr_G^\omega$, we denote its inverse by $(i^!)^{-1}$. 
    
    We only need to show that for any $\cF\in {\SD}_{\kappa}^{H^\omega(\bF),\chi}(\Gr_G^\omega)$, {we have}
    \begin{equation}\label{eq 3.10}
    \begin{split}
 \RHom_{{\SD}_{\kappa}(\Gr_G^\omega)}(\delta_{\sL_{(\lambda,(\theta,\theta'))}}, \oblv_{H^\omega(\bF),\chi}(\cF)
        )=\\ =\RHom_{{\SD}_{\kappa}^{H^\omega(\bF),\chi}(\Gr_G^\omega)}((i^!)^{-1}(\Av_!^{\GL_M(\bO)\ltimes U_{M,N}^{-,\omega}(\bF),\chi}(\delta_{\BL_{(\lambda,(\theta,\theta'))}})), \cF).        
    \end{split}
    \end{equation}

    Indeed, by {the} equivalence \eqref{eq 3.9}, the right-hand side of \eqref{eq 3.10} equals 
    \begin{equation}\label{eq 3.11}
        \begin{split}
            \RHom_{{\SD}_{\kappa}^{\GL_M(\bO)\ltimes U_{M,N}^{-,\omega}(\bF),\chi}(\Gr_N^\omega)}(\Av_!^{\GL_M(\bO)\ltimes U_{M,N}^{-,\omega}(\bF),\chi}(\delta_{\BL_{(\lambda,(\theta,\theta'))}}), i^!\cF)= \\ =\RHom_{{\SD}_{\kappa}(\Gr_N^\omega)}(\delta_{\BL_{(\lambda,(\theta,\theta'))}}, \oblv_{\GL_M(\bO)\ltimes U_{M,N}^{-,\omega}(\bF),\chi} i^!\cF).
        \end{split}
    \end{equation}
    The left-hand side of \eqref{eq 3.10} is the $!$-stalk of $\cF$ at $\sL_{(\lambda,(\theta,\theta'))}$. Since $\cF$ is $\GL_M(\bF)$-equivariant, the $!$-stalk at $\sL_{(\lambda,(\theta,\theta'))}$ equals the $!$-stalk at $(t^0, \BL_{(\lambda,(\theta,\theta'))})$. The latter is exactly the right-hand side of \eqref{eq 3.11}.
\end{proof}

{F}or relevant $(\lambda,(\theta,\theta'))$, we {denote by} \
\begin{equation}\label{compact t-str}
    \tilde{\Delta}_{loc}^{(\lambda,(\theta, \theta'))}:= \Av_!^{H^\omega(\bF),\chi}(\delta_{\sL_{(\lambda,(\theta,\theta'))}})[-\langle 2\rho^\circ, (\lambda, (\theta,\theta'))\rangle].
\end{equation}

\begin{defn}\label{$t$-structure}
    We define a $t$-structure of $\cC^{loc}_{\kappa}(M|N)${,} such that{,} $\cF\in \cC^{loc}_{\kappa}(M|N)$ is coconnective if and only if 
    \[\Hom_{{\SD}_{\kappa}^{H^{\omega}(\bF),\chi}(\Gr^{\omega}_G)}(\tilde{\Delta}_{loc}^{(\lambda,(\theta, \theta'))}[k], \cF)=0,\]
 for any $k>0$ and relevant $(\lambda,(\theta,\theta'))$.    
\end{defn}

{For a relevant} $(\lambda,(\theta,\theta'))$, we denote by $\omega_{\GL_M(\bF)\sL_{(\lambda, (\theta, \theta'))}}$ the $!$-extension of the twisted dualizing D-module on the $\GL_M(\bF)$-orbit of $\sL_{(\lambda, (\theta, \theta'))}$. It acquires a natural $\GL_M(\bF)$-equivariant structure, then we define \[\Delta_{loc}^{(\lambda, (\theta, \theta'))}:= \Av_!^{H^{\omega}(\bF)/\GL_M(\bF),\chi}(\omega_{\GL_M(\bF)\sL_{(\lambda, (\theta, \theta'))}})[-\langle 2\rho^\circ, (\lambda, (\theta,\theta'))\rangle].\]
{Here, $\Av_!^{H^{\omega}(\bF)/\GL_M(\bF),\chi}$ is the relative $!$-averaging functor, i.e., the (partially-defined) left adjoint functor of $\oblv_{H^{\omega}(\bF)/\GL_M(\bF),\chi}:= {\SD}_{\kappa}^{H^\omega(\bF),\chi}(\Gr_G^{\omega})\rightarrow {\SD}_{\kappa}^{\GL_M(\bF)}(\Gr_G^\omega)$.}

The above object is well-defined by the same proof {as in} Lemma \ref{lem 3.4.1}{. Actually, it is the $!$-extension of the rank $1$ local system on $\sO^{(\lambda,(\theta,\theta'))}\subset \Gr_G^\omega$ corresponding to $\chi$}. {Under the equivalence \eqref{eq 3.9}, it corresponds to $ \Av_!^{\GL_M(\bO)\ltimes U_{M,N}^{-,\omega}(\bF)/\GL_M(\bO),\chi}(\delta_{\BL_{(\lambda, (\theta, \theta'))}})[-\langle 2\rho^\circ, (\lambda, (\theta,\theta'))\rangle]$, it is the $!$-extension of the rank $1$ local system on $\BO^{(\lambda,(\theta,\theta'))}\subset 
\Gr_N^{\omega}$.}

\begin{rem}
{The} object  $\Delta_{loc}^{(\lambda, (\theta, \theta'))}$ is locally compact, i.e., it belongs to $\cC_{\kappa}^{loc, lc}(M|N)$, but it is not always compact in ${\SD}_{\kappa}^{H^{\omega}(\bF),\chi}(\Gr^\omega_G)$. {Actually, it is not hard to see that $\Delta_{loc}^{(\lambda, (\theta, \theta'))}$ generates $\cC_{\kappa}^{loc, lc}(M|N)$.}

\end{rem}

{L}et $\sj_{(\lambda, (\theta, \theta'))}$ be the locally closed embedding $\sO^{(\lambda, (\theta, \theta'))}\rightarrow \Gr^{\omega}_G$, we {define}
\[\nabla_{loc}^{(\lambda, (\theta, \theta'))}:=\sj_{(\lambda, (\theta, \theta')),*}\circ \sj_{(\lambda, (\theta, \theta'))}^*(\Delta^{(\lambda, (\theta, \theta'))}_{loc}).\]




Note that there is a canonical map
\begin{equation}\label{map of perverse}
\Delta_{loc}^{(\lambda, (\theta,\theta'))}\longrightarrow \nabla_{loc}^{(\lambda, (\theta,\theta'))}. 
\end{equation}
    We let $\IC_{loc}^{(\lambda, (\theta,\theta'))}$ be the image of $H^0$ (with respect to the $t$-structure defined in Definition \ref{$t$-structure}) of the above map.

\section{Global Gaiotto category}\label{global Gaiotto}
In order to prove that the twisted Gaiotto category is equivalent to the category of ${\cP}_{\det}^\kappa$-twisted factorization modules over $\cI$, we need to define and consider a global version of {the} Gaiotto category. In this section, we assume $C=\BP^1$.


A na\"{i}ve approach is to use the quotient prestack $H_{C-c}\backslash \Gr_G$ to define the global Gaiotto category. Namely, for a global curve $C$, since the character $\chi$ is trivial on $H_{C-c}:= \Maps(C-c, H)\subset H(\bF)$, we expect that the category of $(H(\bF),\chi)$-equivariant D-modules on $\Gr_G$ can be realized as a category of D-modules on $H_{C-c}\backslash \Gr_G$ with a certain equivariant structure. However, the prestack $H_{C-c}\backslash \Gr_G$ is not algebraic when $H$ is not reductive, for example, when $H$ is the maximal unipotent subgroup, see \cite[Remark 4.2.3]{[G3]}. In general, $H_{C-c}\backslash \Gr_G$ is just a subprestack that corresponds to a constructible set of points. The solution is to consider an algebraic stack which contains $H_{C-c}\backslash \Gr_G$, and we impose a Hecke equivariant structure such that any D-module satisfying this condition {is supported} on $H_{C-c}\backslash \Gr_G$.

\subsection{Global model}
In \cite{[GN]} and \cite{[SW]}, the authors generalized the definition of the Drinfeld compactification in \cite{[BG]} and defined global models associated with spherical varieties. Let $\cX= \overline{H\backslash G}^{\aff}=\overline{ (\GL_M\ltimes U_{M,N}^-)\backslash (\GL_M\times \GL_N)}^{\aff}$ be the affine closure of $\overset{\circ}{\cX}=H\backslash G$. It is a spherical variety. Namely, we consider the parabolic subgroup $\GL_M\times P^-_{M+1,1,\cdots,1}$ in $\GL_M\times \GL_N$, where $P^-_{M+1,1,\cdots,1}$ corresponds to the partition $(M+1,1,\cdots,1)$. The Levi subgroup is $\GL_M\times \GL_{M+1}\times \BG_m^{N-M-1}$ and $\GL_M$ is a reductive spherical subgroup of the Levi by the diagonal embedding. Now, the claim follows from \cite[Section 2.3]{[Sa]}. More directly, one can check that the Borel-orbit through the class of $(1,\ 1+\sum_{i=1}^M e_{M+1}\otimes e_i^*)$ in the open $G$-orbit $H\backslash G\subset \cX$ is open dense.

The global model associated with $\cX$ is the algebraic stack  $\Maps_{\gen}(C, \cX/G\supset \overset{\circ}{\cX}/G)$.  Here, $\Maps_{\gen}(C, \cX/G\supset \overset{\circ}{\cX}/G)$ is the open substack of the mapping stack $\Maps(C, \cX/G)$ such that the map generically lands in $\overset{\circ}{\cX}/G\simeq H\backslash G/G\simeq H\backslash \pt$.

In other words, a point of $\Maps_{\gen}(C, \cX/G\supset \overset{\circ}{\cX}/G)$ is $(\cP_G, \sigma)$, where $\cP_G$ is a principal $G$-bundle on $C\times S$ and $\sigma$ is a section $C\times S\longrightarrow \overline{H\backslash G}^{\aff}\overset{G}{\times} \cP_G$ such that $\sigma$ generically lands in ${(H\backslash G)}\overset{G}{\times} \cP_G$. Note that the generic requirement gives $\cP_G$ a generic $H$-reduction.

The following lemma is from \cite[Proposition 3.1.2]{[SW]}.

\begin{lem}
    The global model $\Maps_{\gen}(C, \cX/G\supset \overset{\circ}{\cX}/G)$ is an algebraic stack locally of finite type.
\end{lem}


\begin{defn}\label{def 4.1.2}
 We define the algebraic stack $\cM$ as the $\omega$-renormalized version of $\Maps_{\gen}(C, \cX/G\supset \overset{\circ}{\cX}/G)$, i.e., it is the fiber product \[\Maps_{\gen}(C, T_{N-M-1}\backslash \cX/G\supset T_{N-M-1}\backslash\overset{\circ}{\cX}/G)\underset{\Bun_{T_{N-M-1}}(C)}{\times} \pt.\] Here, the map $\pt\rightarrow \Bun_{T_{N-M-1}}(C)$ is given by $\omega^\rho_C$.

More precisely, an $S$-point of the algebraic stack $\cM$ is a triple $(\cP_M,\cP_N, \kappa_{V})$, where
\begin{enumerate}
    \item $\cP_M$ is a principal $\GL_M$-bundle on $C\times S$,
    \item $\cP_N$ is a principal $\GL_N$-bundle on $C\times S$,
    \item $\kappa_{\cV}$ is an injective map of coherent sheaves 
    \begin{equation}
        \kappa_{\cV}: \cV_{\tilde{\cP}_M}\longrightarrow \Ind (\cV)_{\cP_N}, \textnormal{for any } \cV\in \Rep(\GL_{M+1}\times T_{N-M-1}).
    \end{equation}
Here, $\tilde{\cP}_M:=(\cP_M\overset{\GL_M}{\times}\GL_{M+1})\times_{(C\times S)} (\omega_{C\times S}^{\rho})$ is the $(\GL_{M+1}\times T_{N-M-1})$-bundle obtained from $\cP_M\overset{\GL_M}{\times}\GL_{M+1}$ and $\omega_{C\times S}^{\rho}$, we regard $\cV$ as a $(\GL_{M+1}\ltimes B^-_{M,N})$-module by the natural projection $\GL_{M+1}\ltimes B^-_{M,N}\longrightarrow \GL_{M+1}\times T_{N-M-1}$, and $\Ind(\cV)$ denotes the induced $\GL_N$-module of the $(\GL_{M+1}\ltimes B^-_{M,N})$-module $\cV$.     
\end{enumerate}

We require that the maps $\kappa_{\cV}$ satisfy the Pl\"{u}cker relations, namely, for any $\cV_1, \cV_2\in \Rep(\GL_{M+1}\times T_{N-M-1})$, the following diagram
\begin{equation}
    \xymatrix{
\cV_{1, \tilde{\cP}_M}\otimes \cV_{2, \tilde{\cP}_M}\ar[r]^{\kappa_{\cV_1}\otimes \kappa_{\cV_2}}\ar[d]& \Ind (\cV)_{1,\cP_N}\otimes \Ind (\cV)_{2,\cP_N}\ar[d]\\
(\cV_1\otimes \cV_2)_{ \tilde{\cP}_M}\ar[r]^{\kappa_{\cV_1\otimes \cV_2}}& \Ind (\cV_1\otimes \cV_2)_{\cP_N},
}
\end{equation}
commutes.

\end{defn}

\begin{rem}
The above two descriptions are the same according to \cite[Theorem 1.1.2]{[BG]}. In the non-renormalized case, for example, note that there is an isomorphism of schemes
    \begin{equation}\label{iso of sph}
        \overline{U_{M,N}^-\backslash \GL_N}^{\aff}\simeq \overline{H\backslash G}^{\aff},
    \end{equation}
so the mapping stack \[\Maps_{\gen}(C, \cX/G\supset \overset{\circ}{\cX}/G)= \Maps_{\gen}(C, \GL_M\backslash \overline{U_{M,N}^-\backslash \GL_N}^{\aff}/\GL_N \supset \overset{\circ}{\cX}/G)\] classifies $(\cP_M, \cP_N, \sigma)$, where $\sigma$ is a section $C\rightarrow \cP_M\overset{\GL_M}{\times} \overline{U_{M,N}^-\backslash \GL_N}^{\aff} \overset{\GL_N}{\times} \cP_N$. Using \cite[Theorem 1.1.2]{[BG]}, the map $\sigma$ is equivalent to a collection of maps $\kappa_{\cV}$ for any $\cV\in \Rep(\GL_{M+1}\times T_{N-M-1})$ with Pl\"{u}cker relations.

Through this description, we see that $\Maps_{\gen}(C, \cX/G\supset \overset{\circ}{\cX}/G)$ is a fibration over $\Bun_M(C)$ and its fiber over the trivial bundle is isomorphic to a renormalized Drinfeld compactification of the moduli stack of principal $U_{M,N}^{-,\omega}$-bundles.
\end{rem}


\begin{rem}
If $\kappa_{\cV}$ is a map of vector bundles for any $\cV\in \Rep(\GL_{M+1}\times T_{N-M-1})$, the above data gives rise to a $\omega$-renormalized $H$-reduction $\cP_H$ of $\cP_N$ and an isomorphism of the induced $\GL_M$-bundle of $\cP_H$ with respect to the projection
\[H=\GL_M\ltimes U^{-,\omega}_{M,N}\longrightarrow \GL_M\]
with $\cP_M$. It is equivalent to a data of a point in the moduli stack $\Bun^\omega_H$ of $H^\omega$-bundles.
\end{rem}

\begin{rem}
    In the case $M=0$, the algebraic stack $\cM=\overline{\Bun}^{\omega}_{U^-_N}$ is exactly the Drinfeld compactification of the ($\omega$-renormalized) moduli stack of principal $U^-_N$-bundles in \cite{[BG]}. 

In the case $M=N-1$, since $H=\GL_M$ is reductive, the spherical variety $H\backslash G$ is affine. In particular, $\cX=\overset{\circ}{\cX}$ and there is no generic condition in the definition of $\Maps_{\gen}({C}, \cX/G\supset \overset{\circ}{\cX}/G)$. We have that $\cM= \Bun_{\GL_M}$ is the moduli stack of principal $\GL_M$-bundles on $C$.
\end{rem}
\subsection{Global model with a marked point}\label{sec 4.2}
We can also consider the global model with a marked point $c\in C$. In the non-renormalized setting, it classifies $(\cP_G, \sigma)$, where $\cP_G\in \Bun_G(C)$ and $\sigma$ is a section $C-c\longrightarrow \cX\overset{G}{\times}\cP_G$ generically lands in $\overset{\circ}{\cX}\overset{G}{\times}\cP_G$. 

{More precisely, to describe the $\om$-renormalized global model, we start by considering the $\om^\rho$-twist of $\cX/G$ defined by $(\cX/G)_{\om^\rho} := (\cX/G) \times^{T_{N-M-1}} \om^\rho$; it has an open substack $(\overset\circ\cX/G)_{\om^\rho}$ defined similarly.  Then our global model stack is defined as the stack of sections 
\[\cM:= \Sect_{\gen}(C,(\cX/G)_{\om^\rho} \supset (\overset\circ\cX/G)_{\om^\rho}).\]}

By the isomorphism \eqref{iso of sph}, the above stack also admits a description like Definition \ref{def 4.1.2}.

Namely, in the $\omega$-renormalized setting, we fix a point $c$ and denote by $\cM_{\infty\cdot c}$ the algebraic (ind)-stack which classifies the triples
$(\cP_M,\cP_N, \kappa_{V})$, where
\begin{enumerate}
    \item $\cP_M$ is a principal $\GL_M$-bundle on $C\times S$,
    \item $\cP_N$ is a principal $\GL_N$-bundle on $C\times S$,
    \item $\kappa_{\cV}$ is an injective map of coherent sheaves
    \begin{equation}\label{4.3}
        \kappa_{\cV}: \cV_{\tilde{\cP}_M}\longrightarrow \Ind (\cV)_{\cP_N}(\infty \cdot c), \textnormal{for any } \cV\in \Rep(\GL_{M+1}\times T_{N-M-1}),
    \end{equation}
    which satisfies the Pl\"{u}cker relation. 
\end{enumerate}    

In the case $M=0$, the algebraic ind-stack $\cM_{\infty\cdot c}=(\overline{\Bun}^{\omega}_{U_N})_{\infty\cdot c}$ is the algebraic ind-stack of the Drinfeld compactification with a possible pole at $c$. 

In the case $M=N-1$, $\cM_{\infty\cdot c}$ is the algebraic ind-stack classifies the triples $(\cP_M, \cP_{M+1}, \alpha)$, where $\cP_M$ is a principal $\GL_M$-bundle on $C$, $\cP_{M+1}$ is a principal $\GL_{M+1}$-bundle on $C$, and $\alpha: \cP_M\overset{\GL_M}{\times}\GL_{M+1}|_{C-c}\simeq \cP_{M+1}|_{C-c}$ is an isomorphism over $C-c$. 

The following lemma is proved {using the same method as in} \cite[Lemma 3.5.2]{[BFT0]}.
\begin{lem}
    The global model $\cM_{\infty\cdot c}$ is an algebraic ind-stack locally of ind-finite type.
\end{lem}


\subsubsection{}
We consider a relative determinant line bundle over $\cM_{\infty\cdot c}$ {such that} its fiber over the point $(\cP_M, \cP_N, \kappa_\cV)$ is 
\[\det R\Gamma(C, \cV^M_{\cP_M})\otimes \det^{-1} R\Gamma(C, \cV^N_{\cP_N})\otimes \det^{-1} R\Gamma(C, \cV^M_{\cP^{triv}_M})\otimes \det R\Gamma(C, \cV^N_{\cP_N^\omega}).\]

It is easy to see that it is compatible with the determinant line bundle ${\cP}_{\det}$ on $\Gr_G^{\omega}$ defined in \eqref{3.4}, so we still denote it by the same notation ${\cP}_{\det}$.
We denote by ${\SD}_{\kappa}(\cM_{\infty\cdot c})$ the category of $\cP_{\det}^{\kappa}$-twisted D-modules with coefficients in $\SVect$ on $\cM_{\infty\cdot c}$.

\subsection{Global Gaiotto category}
Now, we are going to mimic the method in \cite{[G]}, \cite{[GN]}, etc., to define the twisted global Gaiotto category. 

\subsubsection{}\label{4.3.1}
First, let us recall some notations. Given a smooth connected curve (not necessarily complete) $C$, we denote by $\Ran_C$ the \textit{Ran space} associated with it. By definition, for a test scheme $S$, the {set of} $S$-points of $\Ran_C$ classifies {finite non-empty subsets} of $\Maps(S, C)$. 

We denote by $H(\bF)_{\bar{x}}$ (resp. $H(\bO)_{\bar{x}}$) the ind-pro group scheme of sections of $H$ on the punctured formal disc (resp. formal disc) around $\bar{x}$, {see} \cite[Section 1.5.5]{[GL]}. Similar{ly} to Section \ref{renorm}, we can define the $\omega$-renormalized $H^\omega(\bF)_{\bar{x}}$ and $H^{\omega}(\bO)_{\bar{x}}$, and a canonical character \[\chi_{\bar{x}}: H^{\omega}(\bF)_{\bar{x}}\longrightarrow \BG_{a,S}\] which is trivial on $H^\omega(\bO)_{\bar{x}}$. 

\begin{rem}\label{rem 4.3.2}
    In particular, the $!$-pullback of the exponential D-module on $\BG_{a}$ to $H^{\omega}(\bF)_{\bar{x}}$ descends to a D-module on the local Hecke stack $\Hecke_{loc, \bar{x}}^H$ which classifies two $H$-bundles on the formal disc of $\bar{x}$ equipped with an isomorphism on the punctured formal disc.  
\end{rem}

If $\bar{x}$ is the disjoint union of different $S$-points $x_1, x_2,\cdots, x_n$, then, we have a canonical isomorphism
\[H^{\omega}(\bF)_{\bar{x}}\simeq \prod_i H^\omega(\bF)_{x_i},\]
and under the above isomorphism, we have $\chi_{\bar{x}}=\sum_{i} \chi_{x_i}$.

Actually, $H^\omega(\bF)_{\bar{x}}$ (resp. $H^\omega(\bO)_{\bar{x}}$, $\chi_{\bar{x}}$) assembles to a prestack $H^\omega(\bF)_{\Ran_C}$ (resp. $H^\omega(\bO)_{\Ran_C}$, $\chi_{\Ran_C}$) over $\Ran_C$, whose fiber over $\bar{x}$ is given by $H^\omega(\bF)_{\bar{x}}$ (resp. $H^\omega(\bO)_{\bar{x}}$, $\chi_{\bar{x}}$).

\subsubsection{}\label{4.3.3}
Let $\bar{x}\in \Ran_{C-c}(S)$. We denote by $\cM_{\bar{x}, \infty\cdot c}$ the open substack of $\cM_{\infty\cdot c}$ where the maps $\kappa_{\cV}$ are injective maps of vector bundles near $\bar{x}${. In other words, each $\kappa$ is injective fiberwise, i.e., the cokernel is $S$-flat near $\bar{x}$}. Since $\kappa_\cV$ is an injective vector bundle map near $\bar{x}$, the collection of maps $\kappa_\cV$ gives rise to an $H^\omega$-bundle on the formal neighborhood $\cD_{\bar{x}}$ of $\bar{x}$, which can be understood as  an $(\GL_M\ltimes B^-_{M,N})$-bundle $\cP_{\GL_M\ltimes B^-_{M,N}}$ with an identification of the induced $T_{N-M-1}$-bundle with $\omega_{C\times S}^{\rho}$.

We need to recall the {notion of} generic Hecke equivariant structure in \cite{[GN]}. Namely, let $\Hecke_{\bar{x}}$ be the stack which classifies the data $(\cP_G, \sigma, \cP_H, \alpha)$, where $(\cP_G, \sigma)$ belongs to $\cM_{\bar{x}, \infty\cdot c}$, $\cP_H$ is a principal $H^\omega$-bundle on the formal neighborhood of $\bar{x}$, and $\alpha$ is an isomorphism of $\cP_H$ and the $H^\omega$-bundle deduced from $(\cP_G, \sigma)$ on the punctured formal neighborhood of $\bar{x}$. Given such a point, we can form another point $(\cP_G',\sigma')$ in $\cM_{\bar{x}, \infty\cdot c}$ by gluing $\sigma|_{(C-c)\times S-\bar{x}}$ and $\cP_H$ along the punctured formal neighborhood of $\bar{x}$. 

{The stack $\Hecke_{\bar{x}}$ admits a map to the local Hecke stack $\Hecke_{loc, \bar{x}}^H$. We denote the $!$-pullback of the {descended} D-module in Remark \ref{rem 4.3.2} to $\Hecke_{\bar{x}}$ {by} $\chi_{\bar{x}}$.}

We consider the following diagram,
\[\xymatrix{
&\Hecke_{\bar{x}}\ar[ld]_{\overset{\gets}{h}}\ar[rd]^{\overset{\to}{h}}&\\
\cM_{\bar{x}, \infty\cdot c}&&\cM_{\bar{x}, \infty\cdot c},
}\]
where $\overset{\gets}{h}$ sends $(\cP_G, \sigma, \cP_H, \alpha)$ to $(\cP_G, \sigma)$ and ${\overset{\to}{h}}$ sends $(\cP_G, \sigma, \cP_H, \alpha)$ to $(\cP'_G, \sigma')$. 

\begin{defn}
An object of $\cC_{\kappa}^{glob}(M|N)_{\bar{x}}$ is $(\cF, \epsilon_{\bar{x}})$, where $\cF\in {\SD}_{\kappa}(\cM_{\bar{x}, \infty\cdot c})$, and $\epsilon_{\bar{x}}$ {is an isomorphism $\overset{\gets}{h}^! (\cF)\simeq \overset{\to}{h}^!(\cF)\overset{!}{\otimes} \chi$ with higher homotopy coherence.    }
\end{defn}

\subsubsection{Second description of $\cC_{\kappa}^{glob}(M|N)_{\bar{x}}$}\label{second description}
Similarly to \cite{[G3]}, one can also interpret $\cC_{\kappa}^{glob}(M|N)_{\bar{x}}$ as $(H^{\omega}(\bF)_{\bar{x}}, \chi_{\bar{x}})$-invariants of the category of {twisted} D-modules with coefficients in $\SVect$ on an $H^{\omega}(\bO)_{\bar{x}}$-torsor over $\cM_{\bar{x}, \infty\cdot c}$. 



We denote by $\tilde{\cM}_{\bar{x}, \infty\cdot c}$ the algebraic stack which classifies a data $(\cP_M,\cP_N,\kappa_{\cV})$ of $\cM_{\bar{x},\infty\cdot c}$ with an isomorphism of $\cP_{\GL_M\ltimes B^-_{M,N}}$ with $\omega_{C\times S}^{\rho}\overset{T_{N-M-1}}{\times}(\GL_M\ltimes B^-_{M,N})$ on $\cD_{\bar{x}}$, which is compatible with the existing identification of the corresponding $T_{N-M-1}$-bundles. Since any principal bundle on the formal neighborhood is \'{e}tale locally trivial, $\tilde{\cM}_{\bar{x}, \infty\cdot c}$ is an $H^\omega(\bO)_{\bar{x}}$-bundle on $\cM_{\bar{x},\infty\cdot c}$. By the standard gluing argument, such as \cite[Lemma 3.2.7]{[FGV]}, the $H^{\omega}(\bO)_{\bar{x}}$-action on $\tilde{\cM}_{\bar{x},\infty\cdot c}$ extends to an action of $H^{\omega}(\bF)_{\bar{x}}$. Then, we have $\cC_{\kappa}^{glob}(M|N)_{\bar{x}}= {\SD}_{\kappa}^{H^{\omega}(\bF)_{\bar{x}}, \chi_{\bar{x}}}(\tilde{\cM}_{\bar{x}, \infty\cdot c})$.  



\subsubsection{}
First, we note that any object in $\cC_{\kappa}^{glob}(M|N)_{\bar{x}}$ only lives on the subspace where the generalized $H^\omega$-reduction is genuine, i.e.,
\begin{lem}\label{vanish bound}
    Assume the generalized $H^\omega$-reduction of $(\cP_M,\cP_N,\kappa_{\cV})\in \cM_{\bar{x},\infty\cdot c}$ has a defect away from $c$ and $\bar{x}$, then any non-zero object $\cF\in \cC_{\kappa}^{glob}(M|N)_{\bar{x}}$ has zero stalk at $(\cP_M,\cP_N,\kappa_{\cV})$.
\end{lem}
\begin{proof}
    The proof is essentially the same as \cite[Lemma 6.2.8]{[FGV]}, we include for completeness.

    First, note that it is enough to check the claim for geometric point $\bar{x}=\{x_1,\cdots,x_m\}\in \Ran_{C-c}$. Assume that the defect locus of $(\cP_M,\cP_N,\kappa_{\cV})$ is contained in $c\cup\bar{y}=\{c,y_1,\cdots,y_n\}$ (which is away from $\bar{x}$ by definition), and the degree of the defect at these points are given by $(\lambda_c, (\theta_c, \theta_c')), (\lambda_1, (\theta_1, \theta_1')),\cdots, (\lambda_n, (\theta_n, \theta_n'))$ in \cite[Section 3.3]{[GN]}, which will also be recalled in Section \ref{4.6}. Then, the point $(\cP_M,\cP_N,\kappa_{\cV})$ is contained in the locally closed substack
    \begin{equation}
        H^\omega_{C-c-\bar{x}-\bar{y}}\backslash (\sO_{c}^{(\lambda_c,(\theta_c,\theta_c'))}\times \prod^m_{i=1} H^{\omega}(\bF)_{{x}_i}/H^{\omega}(\bO)_{{x}_i} \times \prod^n_{i=1} \sO_{y_i}^{(\lambda_i,(\theta_i,\theta_i'))}),
    \end{equation}
    and it is enough to show that there is no non-zero $(H^\omega(\bF)_{\bar{x}},\chi_{\bar{x}})$-equivariant object in
        \begin{equation}
        H^\omega_{C-c-\bar{x}-\bar{y}}\backslash (\sO_{c}^{(\lambda_c,(\theta_c,\theta_c'))}\times \prod^m_{i=1} H^{\omega}(\bF)_{{x}_i} \times \prod^n_{i=1} \sO_{y_i}^{(\lambda_i,(\theta_i,\theta_i'))}),
    \end{equation}
    i.e., there is no non-zero $(H^\omega_{C-c-\bar{x}-\bar{y}}, \chi_{\bar{x}})$--equivariant object in
        \begin{equation}\label{loc des}
        \sO_{c}^{(\lambda_c,(\theta_c,\theta_c'))}\times \prod^n_{i=1} \sO_{y_i}^{(\lambda_i,(\theta_i,\theta_i'))}.
    \end{equation}
    We denote a point in \eqref{loc des} which maps to $(\cP_M,\cP_N,\kappa_{\cV})$ by $(z_c,z_1,\cdots,z_n)$.
    
    According to the fact that the unital $H^\omega(\bF)$-orbit in $\Gr_G^\omega$ is minimal relevant with respect to the additive character $\chi$, i.e., there is no relevant orbit in the boundary of the unital orbit which supports a non-zero $(H^\omega(\bF),\chi)$-equivariant sheaf. Indeed, it follows from Corollary \ref{H orbit'} and Proposition \ref{closure}. This, together with the fact that the conjugation action of $\GL_M(\bF)$ on $U_{M,N}^{-,\omega}(\bF)$, implies that there is no $(U_{M,N}^{-,\omega}(\bF),\chi)$-relevant orbit in the boundary of the unital $H^\omega(\bF)$-orbit in $\Gr_G^\omega$. 

    Thus, one can find an element $h$ in $  (U_{M,N}^{-,\omega})(\bF)_{y_1}$, which belongs to the stabilizer of $z_1$ but does not belong to the kernel of the character $\chi_{y_1}$. 

    Fix $N>0$, we can choose an element in $(U_{M,N}^{-,\omega})_{C-c-\bar{x}-\bar{y}}$ which equals $h$ at $y_1$ and equals $1$ at $c$ and $y_i$, after modulo the \mbox{$N$-th} congruence subgroup, where $i=2,\cdots,y_n$. When $N$ is large enough, this point belongs to the stabilizer of $(z_c,z_1,\cdots,z_n)$, but by {the} residue theorem, we have $\chi_{\bar{x}}=-\chi_{y_1}\neq 0$.
\end{proof}

\subsubsection{}
If $\bar{x}\subset \bar{x}'\in \Ran_{C-c}(S)$, there is a functor 
\begin{equation}\label{4.4}
    \cC_{\kappa}^{glob}(M|N)_{\bar{x}'}\longrightarrow \cC_{\kappa}^{glob}(M|N)_{\bar{x}}.
\end{equation}

To see this, first, we note that the category of $(H^\omega(\bF)_{\bar{x}},\chi_{\bar{x}})$-equivariant {twisted} D-modules on $\tilde{\cM}_{\bar{x}, \infty\cdot c}$ is equivalent to the category of $(H^\omega(\bF)_{\bar{x}},\chi_{\bar{x}})$-equivariant {twisted} D-modules on $\tilde{\cM}_{\bar{x}, \infty\cdot c}|_{\cM_{\bar{x}', \infty\cdot c}}$, see Lemma \ref{vanish bound}. Here, we regard $\cM_{\bar{x}', \infty\cdot c}$ as an open substack of $\cM_{\bar{x}, \infty\cdot c}$ and we denote by $\tilde{\cM}_{\bar{x}, \infty\cdot c}|_{\cM_{\bar{x}', \infty\cdot c}}$ the restriction of $\tilde{\cM}_{\bar{x}, \infty\cdot c}$ to $\cM_{\bar{x}, \infty\cdot c}$. 

Then, given a point $(\cP_G, \sigma, \cP_H, \alpha)$ in $\Hecke_{\bar{x}}|_{\cM_{\bar{x}', \infty\cdot c}}$, we can extend $\cP_H$ (resp.\ $\alpha$) to a principal $H^\omega$-bundle on the formal neighborhood of $\bar{x}'$ (resp.\ an isomorphism on the punctured formal neighborhood of $\bar{x}'$) by $(\cP_G, \sigma)$. In particular, we have a map $l_{\bar{x}\to \bar{x}'}$ from $\Hecke_{\bar{x}}|_{\cM_{\bar{x}', \infty\cdot c}}$ to $\Hecke_{\bar{x}'}$ which makes the following diagram {commutative:}
\[\xymatrix{
&\Hecke_{\bar{x}'}\ar[ld]\ar[rd]&\\
\cM_{\bar{x}', \infty\cdot c}&&\cM_{\bar{x}', \infty\cdot c}\\
&\Hecke_{\bar{x}}|_{\cM_{\bar{x}', \infty\cdot c}}\ar[lu]\ar[ru]\ar[uu]_{l_{\bar{x}\to \bar{x}'}}&.
}\]

The functor which sends $(\cF, \epsilon_{\bar{x}'})$ to $(\cF, l_{\bar{x}\to \bar{x}'}^!\epsilon_{\bar{x}'})$ defines the desired functor \eqref{4.4}. 


Furthermore, if $S_1$ and $S_2$ are two schemes and $f:S_1\rightarrow S_2$, then any $S_2$-point $\bar{x}$ of $\Ran_{C-c}$ can be regarded as an $S_1$-point $\bar{x}\underset{S_2}{\times} S_1$ by composing with $f$. We have a transition functor 
\begin{equation}\label{4.5}
 \cC^{glob}_{\kappa}(M|N)_{\bar{x}} \longrightarrow   \cC_{\kappa}^{glob}(M|N)_{\bar{x}\underset{S_2}{\times} S_1}.
\end{equation}


{
Let us consider the functor from $(\Sch^\aff)^\opp$ to the category of posets, which sends $S$ to $\Ran_{C-c}(S)$ with the partial order given by containment. By Grothendieck construction, we obtain a co-Cartesian fibration over $(\Sch^\aff)^{\opp}$ whose fiber category over $S$ is $\Ran_{C-c}(S)$.
}
\begin{defn}\label{def 4.3.4}
    We define
    \[\cC_{\kappa}^{glob}(M|N):= \lim_{\bar{x}\in \Ran_{C-c}(S), S\in (\Sch^\aff)^\opp} \cC_{\kappa}^{glob}(M|N)_{\bar{x}}.\]
    Here, the index category is understood as the above co-Cartesian fibration over $(\Sch^\aff)^\opp$, and the transition functors are generated by the functors \eqref{4.4} and \eqref{4.5} defined above.
\end{defn}

{\begin{rem}
    In \cite{[G3]}, the definition of the global Whittaker category only involves the geometric points of $\Ran_{C-c}$. However, in the Gaiotto case, the category $\lim_{\bar{x}\in \Ran_{C-c}(\Spec \BC)} \cC_{\kappa}^{glob}(M|N)_{\bar{x}}$ is not the correct category to consider, its hom space is different from the expected one.
\end{rem}}




In the next section, we will give several different {frameworks and motivations  for  the above limit.}
\subsection{Continuous limit and lax sheaf of categories}\label{cont limit}
This section is a long remark to Definition \ref{def 4.3.4}, and we aim to give different descriptions of the limit. This section is irrelevant of the goal of this paper, and the readers can skip it safely.
\subsubsection{Continuous limit}
{The global category is expected to be the limit of $\cC_\kappa^{glob}(M|N)_{\bar{x}}$ for all schematic point $\bar{x}\in \Ran_{C-c}$, and the transition functors are the functors \eqref{4.4} and \eqref{4.5} defined above. To define the correct limit, in the universal property defining the limit, we should use the prestack structure of $\Ran_{C-c}$: namely, we require the family of functors to $\cC_{\kappa}^{glob}(M|N)_{\bar{x}}$ to depend ``algebraically'' on $\bar x$.}

Given a category $\cD$ internal to the category of prestacks of sets, and a sheaf of categories $\cC$ on $\cD$ {(i.e., a sheaf of categories on the object of objects of $\cD$, and a collection of functors corresponding to the object of morphisms with higher homotopy coherence)}. We can mimic the definition of limit {with a} usual category as the index category, and define the limit of $\cC$ with the internal category $\cD$ as the index category.

To be more precise, we have {the ``constant presheaf''} functor from $\Cat$ to the category of sheaves on $\cD$, which sends a category to the corresponding constant {presheaf} of categories. We define the limit $\lim_{\cD}\cC$ as the image of $\cC$ under the right adjoint functor of the constant {presheaf} functor.

In our particular case, the object of objects of $\cD$ is the prestack $\Ran_{C-c}$, and the object of morphisms of $\cD$ is $(\Ran_{C-c}\times \Ran_{C-c})^\subset$. Here, $(\Ran_{C-c}\times \Ran_{C-c})^\subset$ denotes the sub-prestack of $\Ran_{C-c}\times \Ran_{C-c}$ which consists of $(\bar{x},\bar{x}')$ such that $\bar{x}\subset \bar{x}'$. 

We let $\cC$ be the sheaf of categories which assigns $\bar{x}\in \Ran_{C-c}(S)$ the $D(S)$-module category  $\cC_{\kappa}^{glob}(M|N)_{\bar{x}}$. One can check that the limit of $\cC$ over the internal category $\cD$ is the limit in \eqref{def 4.3.4}. That is to say, we have 
\begin{equation}
    \lim_{\cD} \cC\simeq \lim_{\bar{x}\in \Ran_{C-c}(S), S\in (\Sch^\aff)^\opp} \cC_{\kappa}^{glob}(M|N)_{\bar{x}}=\colon \cC_{\kappa}^{glob}(M|N).
\end{equation}

\subsubsection{Lax sheaf of categories}
{The limit in Definition \ref{def 4.3.4} can be calculated as a limit over the category of twisted arrows of $\Sch^\aff$. }

Consider the functor 
\begin{equation}
    \ShvCat: (\Sch^\aff)^{\opp}\longrightarrow \Cat,
\end{equation}
which sends $S$ to $D(S)\modu$, the category of $D(S)$-module categories. By the Grothendieck construction, we can regard it as a Cartesian {fibration} over $\Sch^\aff$, whose fiber over $S\in \Sch^\aff$ is $D(S)\modu$. {By a \emph{lax sheaf of categories} we will mean} a section of this fibration, and the category of lax sheaves of categories {will be} denoted by $\ShvCat^{\lax}$.

\begin{rem}
    A presheaf of categories is a \emph{sheaf} precisely when it is  a co-Cartesian section.
\end{rem}

In other words, a lax sheaf of categories assigns {to} each {affine} scheme $S$ a $D(S)$-module category $\cC_S$, and  {for every morphism} $S_1\to S_2$ {in $\Sch^\aff$}, there is a morphism in $D(S_1)\modu$,
\begin{equation}
    \cC_{S_2}\otimes_{D(S_2)} D(S_1)\longrightarrow \cC_{S_1},
\end{equation}
which satisfies higher associativity.

Note that there is a naturally defined functor
\begin{equation}\label{eq 4.9}
 { \Const\colon   \Cat\longrightarrow \ShvCat^{\lax},}
\end{equation}
which sends $\cC$ to the constant sheaf of categories $\{\cC_S:= \cC\otimes D(S)\}$.

\begin{defn}
    Given a lax sheaf of categories $\{\cC_S\}\in \ShvCat^\lax$, we define its limit $\lim \{\cC_S\}$ as $\Const^R(\{\cC_S\})\in \Cat$, here $\Const^R$ denotes the right adjoint functor of $\Const$.
\end{defn}

We can calculate the limit of a lax sheaf of categories as
\begin{equation}\label{eq 4.10}
    \lim \{\cC_S\}\simeq \lim_{(S_1\to S_2)} \cC_{S_2}\otimes_{D(S_2)} D(S_1),
\end{equation}
where the right-hand side is understood as the usual limit in $\Cat$ over the category of twisted arrows {of $\Sch^\aff$}.

Indeed, it is enough to show that the data of a functor from a category $\cC'$ to $\lim_{S_1\to S_2} \cC_{S_2}\otimes_{D(S_2)} D(S_1)$ is equivalent to the data of a functor $\{\cC'\otimes D(S)\}$ to $\{\cC_S\}$. The inverse map is obvious. Given a functor $\cC'$ to $\lim_{S_1\to S_2} \cC_{S_2}\otimes_{D(S_2)} D(S_1)$, it gives rise to $\cC'\to \cC_{S_2}= \cC_{S_2}\otimes_{D(S_2)} D(S_2)$, which is equivalent to a morphism $\cC'\otimes D(S_2)\to \cC_{S_2}$ in $D(S_2)\modu$. On the other direction, the functor $\cC'\otimes D(S)\to \cC_{S}$ gives rise to $\cC'\to \cC_{S}$. Furthermore, for any $S_1\to S_2$, we have a functor $\cC'\to \cC_{S_2}\otimes_{D(S_2)} D(S_1)$ given by the composition $\cC'\to \cC_{S_2}$ and $\cC_{S_2}\to \cC_{S_2}\otimes_{D(S_2)} D(S_1)$.

According to \eqref{eq 4.10}, if we let $\cC_S:=\lim_{\bar{x}\in \Ran_{C-c}(S)} \cC_{\kappa}^{glob}(M|N)_{\bar{x}}$, then
\begin{equation}
     \lim \{\cC_S\}\simeq \lim_{(S_1\to S_2)} \cC_{S_2}\otimes_{D(S_2)} D(S_1)= \cC_{\kappa}^{glob}(M|N).
\end{equation}

\subsection{Local-Global comparison}
Similarly to \cite{[G3]}, we hope to establish an equivalence between $\cC_{\kappa}^{loc}(M|N)$ and $\cC_{\kappa}^{glob}(M|N)$. 

Let $(\cP_M,\alpha_M, \cP_N, \alpha_N)$ be a point of $\Gr_M\times \Gr^{\omega}_N$ and let $\cV$ be a representation of $\GL_{M+1}\times T_{N-M-1}$, we define the map $\kappa_\cV$ as the composition of 
\begin{equation}
    \cV_{\tilde{\cP}_M}\overset{\alpha_M}{\longrightarrow} \cV_{\tilde{\cP}^{triv}_M}\longrightarrow \Ind(\cV)_{\cP^{\omega}_N}\overset{\alpha_N}{\longrightarrow} \Ind(\cV)_{\cP_N}.
\end{equation}

One checks that the map $\kappa_{\cV}$ defined above satisfies the Pl\"{u}cker relations. So, the above assignment defines a map:
\begin{equation}\label{eq 4.14}
    \pi: \Gr^{\omega}_{{G}}=\Gr_M\times {\Gr}^{\omega}_N\longrightarrow \cM_{\infty\cdot c}
\end{equation}

By definition, the determinant line bundles on $\cM_{\infty\cdot c}$ and $\Gr^{\omega}_{{G}}$ are compatible with respect to the map $\pi$. That is to say, there is a functor
\begin{equation}
    \pi^!: {\SD}_{\kappa}(\cM_{\infty\cdot c})\longrightarrow {\SD}_{\kappa}(\Gr^{\omega}_{{G}}).
\end{equation}

We claim 

\begin{thm}\label{locglob}
A). The functor $\pi^!$ induces a functor 
\begin{equation}\label{local-global}
    \pi^!: \cC_{\kappa}^{glob}(M|N)\longrightarrow \cC_{\kappa}^{loc}(M|N).
\end{equation}

B). The functor above \eqref{local-global} is an equivalence.
\end{thm}

\subsubsection{}\label{sec 4.5.2}
For our main goal, we only perform the proof of the above theorem. However, just {repeating} the proof verbatim, one {can actually} prove the local-global comparison theorem in a more general setting. 

 Let $G$ be a reductive group, $H=H_r\ltimes H_u$ be a connected subgroup of $G$, such that $H_u$ is the unipotent radical of a parabolic subgroup $P$ and $H_r$ is reductive and belongs to the Levi subgroup of $P$. We also assume $\overset{\circ}{\cX}=H\backslash G$ is a strongly quasi-affine spherical variety, i.e., open in its canonical affine closure $\cX:= \Spec \BC[G]^H$. Let $\mathfrak{c}^-_{\cX}$ be the index set of $G(\bO)$-orbits in $\overset{\circ}{\cX}(\bF)\cap \cX(\bO)$, it is identified with a subset of the index set of $H(\bF)$-orbits in $\Gr_G$.

 We define $\cM_{\cX, \infty \cdot c}$ as the ind-stack which classifies $(\cP_G, \sigma)$, where $\cP_G$ is a $G$-bundle on $C$, and $\sigma$ is a section $C-c\longrightarrow \cX\overset{G}{\times}\cP_G$ which generically lands in $\overset{\circ}{\cX}\overset{G}{\times}\cP_G$.

{Let $\chi\colon H\to \BG_a$ be an additive character.  Note that $\chi|_{H_r}\equiv0$ necessarily.  By a slight abuse of notation, we will also denote by $\chi$ the composition of the induced map}  
   $H(\bF)\to \bF$ with the residue map {$\bF\to\BG_a$}. Then the global category corresponding to ${\SD}^{H(\bF),\chi}(\Gr_G)$ is defined as the $(H,\chi)$ generic Hecke equivariant D-modules on $\cM_{\cX, \infty\cdot c}$. 

\begin{thm}\label{general} 
    Assume that there is no non-zero $(H(\bF),\chi)$-equivariant D-module on those $H(\bF)$-orbits corresponding to $\mathfrak{c}^-_{\cX}-0$, and $H_r\backslash \overline{H_u\backslash G}^{\aff}$ is a scheme.
    
    Then the $!$-pullback from $\cM_{\cX, \infty\cdot c}$ to $\Gr_G$ induces an equivalence between the category of $(H,\chi)$-generic Hecke equivariant D-modules on $\cM_{\cX, \infty \cdot c}$ and the category of $(H(\bF),\chi)$-equivariant D-modules on $\Gr_G$. 

   { Furthermore, if the determinant line bundle on $\Gr_G$ admits an $H(\bF)$-equivariant structure, then, the $!$-pullback induces an equivalence between the corresponding twisted categories. }

\end{thm}

In the Gaiotto case, $G=\GL_M\times \GL_N$, $H=\GL_M\ltimes U_{M,N}^-$, $H_r=\GL_M$ and $H_u=U_{M,N}^-$.

\begin{rem}
The former condition in the above theorem ensures that there is no generic Hecke equivariant D-module supported on the complement of $H_{C-c}\backslash \Gr_G$, and the latter condition says that the global model is a fibration over $\Bun_{H_r}(C-c)$, whose fiber over $\cP_{H_r}\in \Bun_{H_r}(C-c)$ is the ($\cP_{H_r}$-renormalized) Drinfeld compactification for $H_u$ with a possible pole at $c$.
\end{rem}
\begin{rem}
   If $H$ is connected and reductive, then $G,H,\chi=0$ satisfies the assumption in the theorem. Indeed, we have $\cX=\overset{\circ}{\cX}=H\backslash G$ and $\mathfrak{c}_{\cX}=0$. Furthermore, the Whittaker case i.e., $H$ is unipotent, $\chi$ is non-degenerate, also satisfies the above assumption. 
\end{rem}


\subsection{Stratification of the Drinfeld compactification}\label{4.6}
First, we consider a stratification of $\cM_{\infty\cdot c}$. 

The stratification of $\cM_{\infty\cdot c}$ (also, $\cM_{\bar{x},\infty\cdot c}$) is given by the degree of the defect. Since any ${G}$-bundle on the formal disc $\cD_c$ is trivial, we can fix an isomorphism between $\cP_G$ and $\omega^{\rho}\overset{T_{N-M-1}}{\times} G$ on $\cD_c$. Given a geometric point $(\cP_G, \sigma)\in \cM_{\infty\cdot c}$, taking restriction of $\sigma$ to $\overset{\circ}{\cD}_c$, it factors through ${(\GL_M\ltimes B^-_{M,N})\backslash (\GL_M\times \GL_N})\overset{G}{\times} \cP_G$ such that the induced $T_{N-M-1}$-bundle is isomorphic to $\omega^{\rho}$. In particular, it gives rise to a point in $(H\backslash {G})^\omega(\bF)_c$, Different trivialization of $\cP_G$ on $\cD_c$ changes the resulting point by ${G}^\omega(\bO)_c$. As a consequence, for each geometric point in $\cM_{\infty\cdot c}$, there is a well-defined point in $H^{\omega}(\bF)_c\backslash \Gr^{\omega}_{{G}}$. Note that according to Corollary \ref{H orbit'}, $H^{\omega}(\bF)_c$-orbits of $\Gr^{\omega}_{{G}}$ are indexed by $(\lambda, (\theta, \theta'))$. We denote by {$\cM_{\leq (\lambda, (\theta, \theta'))\cdot c}$} the closed substack of $\cM_{\infty\cdot c}$ whose defect at $c$ is no more than $(\lambda, (\theta, \theta'))$.

There is an open substack of {$\cM_{\leq (\lambda, (\theta, \theta'))\cdot c}$} consisting of those points whose defect at $c$ is exactly $(\lambda, (\theta, \theta'))$ and do not have non-zero defect at any point $x\in C-c$, i.e., the morphism $\kappa_{\cV}$ is an injective bundle map for any $\cV\in \Rep(\GL_{M+1}\times T_{N-M-1})$. We denote it by ${\cM_{= (\lambda, (\theta, \theta'))\cdot c}}$. In the existing literature, such as \cite{[G3]}, the condition that there is no defect at any point $x\in C-c$ is called \textit{good elsewhere}.


As the first step to the proof of Theorem \ref{locglob}, we claim that pulling-back along $\pi$ defines a stratawise equivalence.
\begin{prop}\label{strata eq}
For any $(\lambda, (\theta,\theta'))$, there is an equivalence
\[\pi^!: \cC_{\kappa}^{glob}(M|N)|_{{\cM_{= (\lambda, (\theta, \theta'))\cdot c}}}\simeq {\SD}_{\kappa}^{H^{\omega}(\bF), \chi}(\sO^{(\lambda, (\theta,\theta'))})= {\SD}_{\kappa}^{\GL_M(\bO)\ltimes {U}^{-,\omega}_{M,N}(\bF), \chi}(\BO^{(\lambda, (\theta,\theta'))})\]
\end{prop}

\begin{proof}
There is an isomorphism of algebraic stacks
\[{\cM_{= (\lambda, (\theta, \theta'))\cdot c}}\simeq H^{\omega}_{C-c}\backslash \sO_c^{(\lambda,(\theta,\theta'))},\]
where the right-hand side is understood as the \'{e}tale (equivalently, fpqc) sheafification of the prestack quotient.

Furthermore, for any $\bar{x}\subset (C-c)\times S$, there is an isomorphism
\[{{\cM}_{\bar{x},= (\lambda, (\theta, \theta'))\cdot c}}\simeq H^{\omega}_{(C-c)\times S-\bar{x}}\backslash  ((H^{\omega}(\bF)_{\bar{x}}/H^{\omega}(\bO)_{\bar{x}})\times_S \sO_{c\times S}^{(\lambda,(\theta,\theta'))}).\]
Here, $H^{\omega}_{(C-c)\times S-\bar{x}}$ denotes the group ind-scheme of maps from $(C-c)\times S-\bar{x}$ to $H^{\omega}$. The action of $H^{\omega}_{(C-c)\times S-\bar{x}}$ on $H^{\omega}(\bF)_{\bar{x}}/H^{\omega}(\bO)_{\bar{x}}\times\sO_{c\times S}^{(\lambda,(\theta,\theta'))}$ is given by
\[H^{\omega}_{(C-c)\times S-\bar{x}}\longrightarrow H^{\omega}(\bF)_{\bar{x}}\times_S H^{\omega}(\bF)_{c\times S}\curvearrowright H^{\omega}(\bF)_{\bar{x}}/H^{\omega}(\bO)_{\bar{x}}\times_S\sO_{c\times S}^{(\lambda,(\theta,\theta'))}.\]

By construction and Section \ref{second description}, an object of $\cC_{\kappa}^{glob}(M|N)_{\bar{x}}|_{{\cM_{= (\lambda, (\theta, \theta'))\cdot c}}}$ is an $H^{\omega}_{(C-c)\times S-\bar{x}}$-equivariant (with respect to the diagonal action) and $(H^{\omega}(\bF)_{\bar{x}}, \chi_{\bar{x}})$-equivariant twisted D-module on $H^{\omega}(\bF)_{\bar{x}}\times_S \sO_{c\times S}^{(\lambda, (\theta, \theta'))}$, which is exactly an $H^{\omega}_{(C-c)\times S-\bar{x}}, \chi)$-equivariant (with respect to the left action) twisted D-module on $\sO_{c\times S}^{(\lambda, (\theta, \theta'))}$. Here, the character map $\chi$ of $H^{\omega}_{(C-c)\times S-\bar{x}}$ is given by restriction $H^{\omega}_{(C-c)\times S-\bar{x}}\longrightarrow H^{\omega}(\bF)_{c\times S}\overset{\chi}{\longrightarrow} \BG^1_{a,S}$.

So, we have
\begin{equation}
    \cC_{\kappa}^{glob}(M|N)_{\bar{x}}|_{{\cM_{= (\lambda, (\theta, \theta'))\cdot c}}}\simeq {\SD}_{\kappa}^{H^{\omega}_{(C-c)\times S-\bar{x}}, \chi} (\sO^{(\lambda,(\theta, \theta'))}_{c\times S}).
\end{equation}

Let $H_{gen,c}$ be the group prestack whose $S$-points classify maps from $(C-c)\times S-\bar{x}$ to $H^\omega$ for a certain $\bar{x}\in \Ran_{C-c}(S)$, and two maps $(C-c)\times S-\bar{x}\rightarrow H^\omega$ and  $(C-c)\times S-\bar{x}'\rightarrow H^\omega$ are identified if their restrictions to $(C-c)\times S-\bar{x}\cup \bar{x}'$ are the same. 

More or less tautologically from the definition, we have

\[\lim_{\bar{x}\in \Ran_{C-c}(S), S\in (\Sch^\aff)^\opp} {\SD}_{\kappa}^{H^{\omega}_{(C-c)\times S-\bar{x}}, \chi} (\sO^{(\lambda,(\theta, \theta'))}_{c\times S})\simeq {\SD}_{\kappa}^{H_{gen, c},\chi}(\sO^{(\lambda,(\theta, \theta'))}_c).\]

Let $H_{gen,c}^{sh}$ be the sheafification of $H_{gen,c}$ with respect to finite flat surjections. {The $S$-points of the stack $H_{gen,c}^{sh}$ classifies maps from $(C-c)\times S- \Gamma$ to $H^\om$ for a certain divisor $\Gamma$ in $(C-c)\times S$ which is finite over $S$, and two maps $(C-c)\times S-\Gamma\rightarrow H^\omega$ and  $(C-c)\times S-\Gamma'\rightarrow H^\omega$ are identified if their restrictions to $(C-c)\times S-(\Gamma\cup \Gamma')$ are the same. }

We only need to prove
\[{\SD}_{\kappa}^{H^\omega(\bF)_c,\chi}(\sO_c^{(\lambda, (\theta,\theta'))})\simeq {\SD}_{\kappa}^{H_{gen, c}^{sh},\chi}(\sO^{(\lambda,(\theta, \theta'))}_c).\]

It follows from the following Proposition \ref{loop H inv}


\end{proof}

\begin{prop}\label{loop H inv}
    For any dualizable category $\cC$ with an action of $H^{\omega}(\bF)$, there is an equivalence
\begin{equation}\label{h gen}
  \cC^{H^\omega(\bF)_c,\chi}\simeq \cC^{H_{gen, c}^{sh},\chi}.  
\end{equation}
\end{prop}
\begin{proof}
    For simplicity, we consider the non-renormalized groups; also, similarly to \cite[Corollary 4.3.1]{[Ber]}, we consider the case $\chi=0$. 

    By a category admitting an action of $H(\bF)$ (resp. $H_{gen, c}^{sh}$), we mean a module category over the monoidal category $(D^*(H(\bF)),\star)$ (resp. $(D^*(H_{gen, c}^{sh}),\star)$), where $D^*$ means the $*$-D-module theory in \cite{[Ras1]} and the monoidal structure is given by $*$-pushforward along the multiplication map. 

    By definition, we have 
    \begin{equation}\label{4.23 eq}
        \begin{split}
            \cC^{H_{gen, c}^{sh}}\simeq \Fun_{D^*(H_{gen, c}^{sh})}(\Vect, \cC)\\
            \simeq \Fun_{D^*(H(\bF))}(\Vect\otimes_{D^*(H_{gen, c}^{sh})}D^*(H(\bF)), \cC)\\
            \simeq \Fun_{D^*(H(\bF))}(D^*(H(\bF)/H_{gen, c}^{sh}), \cC).
        \end{split}
    \end{equation}
    Let $K_n$ be the $n$-th congruence subgroup of $H(\bF)$, we have
    \begin{equation}\label{4.24 eq}
        \begin{split}
            D^*(H(\bF)/H_{gen, c}^{sh})\simeq \on{lim}_{n,\pi_*}D^*(K_n\backslash H(\bF)/H_{gen, c}^{sh})\\
            \simeq \colim_{n,\pi^*}D^*(K_n\backslash H(\bF)/H_{gen, c}^{sh}),
        \end{split}
    \end{equation}
    where the map $\pi_{n_1, n_2}\colon K_{n_1}\backslash H(\bF)/H_{gen, c}^{sh}\rightarrow K_{n_2}\backslash H(\bF)/H_{gen, c}^{sh}$ is the projection, for $n_1\geq n_2$. The left adjoint functor $\pi^*$ exists due to the smoothness of $\pi_{n_1\to n_2}$, and the equivalence follows from the general pattern due to J.Lurie that the limit along the (possibly discontinuous) right adjoint functor is equivalent to the colimit along the left adjoint functor.

    Apply \eqref{4.24 eq} to \eqref{4.23 eq}, we obtain
    \begin{equation}
        \begin{split}
            \eqref{4.23 eq}\simeq \Fun_{D^*(H(\bF))}(D^*(H(\bF)/H_{gen, c}^{sh}), \cC)\\
            \simeq \Fun_{D^*(H(\bF))}(\Vect,\Fun(\colim_{n,\pi^*}D^*(K_n\backslash H(\bF)/H_{gen, c}^{sh}), \cC))\\
            \simeq \lim \Fun_{D^*(H(\bF))}(\Vect,\Fun(D^*(K_n\backslash H(\bF)/H_{gen, c}^{sh}), \cC)).
        \end{split}
    \end{equation}

    Since $\cC$ is dualizable, we denote by its dual by $\cC^\vee$, and for any category $\cD$, we have
    \begin{equation}
    \begin{split}
        \Fun(\cD,\cC)\simeq \Fun(\cD,\Fun(\cC^\vee,\Vect))
        \simeq \Fun(\cD\otimes \cC^\vee,\Vect)\\
        \simeq \Fun(\cC^\vee,\Fun(\cD,\Vect))
        \simeq \Fun(\cD,\Vect)\otimes \cC.
    \end{split}
    \end{equation}

    Thus, we get
    \begin{equation}
\begin{split}
     \eqref{4.23 eq}
     \simeq \lim \Fun_{D^*(H(\bF))}(\Vect,\Fun(D^*(K_n\backslash H(\bF)/H_{gen, c}^{sh}), \Vect)\otimes \cC).
\end{split}
    \end{equation}
    According to the Lemma \ref{Lemma 4.5.2}, we have
    \begin{equation}
        D^*(K_n\backslash H(\bF)/H_{gen, c}^{sh})\simeq D^*(B(K_n\cap H_{gen, c}^{sh})),
    \end{equation}
    where $B(K_n\cap H_{gen, c}^{sh})$ denotes the classifying stack of the group object $K_n\cap H_{gen, c}^{sh}$.

    By definition, we have 
    \begin{equation}
        \Fun(D^*(B(K_n\cap H_{gen, c}^{sh}), \Vect)\simeq D^!(B(K_n\cap H_{gen, c}^{sh})).
    \end{equation}

    However, by the naively $\BA^1$-contractibility of $K_n\cap H_{gen, c}^{sh}$ (cf., Lemma \ref{Lemma 4.5.2}), we have 
    \begin{equation}
        D^!(B(K_n\cap H_{gen, c}^{sh}))\simeq \Vect.
    \end{equation}

    Thus, \begin{equation}
        \begin{split}
            \eqref{4.23 eq}\simeq \lim  \Fun_{D^*(H(\bF))}(\Vect, \cC)\simeq \lim \cC^{H(\bF)}\simeq \cC^{H(\bF)}.
        \end{split}
    \end{equation}
\end{proof}

\begin{lem}\label{Lemma 4.5.2}
For any unipotent congruence subgroup $K_n$, the map \[H_{gen, c}^{sh}\longrightarrow H^{\omega}(\bF)_c/ K_n\]
is surjective and has naively $\BA^1$-contractible fibers.
\end{lem}
\begin{proof}
For simplicity, we consider the non-renormalized groups. For the surjectivity part, it just follows from the density of $\{H_{(C-c)\times S-\bar{x}}\}$ in $H(\bF)_{c\times S}$ for any $\bar{x}$ (note that $C=\BP^1$). To prove that it has naively $\BA^1$-contractible fibers, we only need to show that $H_{gen, c}^{sh}\cap K_n$ is naively $\BA^1$-contractible. We take an affine space $V$, such that the product of the open Schubert cell $U_{H_r}^- B_{H_r}$ of $H_r$ and $H_u$ is an open dense sub-scheme of $V$. Here, $H=H_r\ltimes H_u$, $H_r$ is the reductive part and $H_u$ is unipotent. Note that $H_{gen, c}^{sh}\cap K_n$ classifies rational maps from $C$ to $H$ which are regular near $c$, the value at $c$ is exactly 1, and the first $n$-order differentials at $c$ equal $0$. It equals a naively $\BA^1$-contractible space that classifies rational maps from $C$ to $V$ that are regular near $c$, the value at $c$ is exactly 1, and the first $n$ order differentials at $c$ equals $0$. 
\end{proof}

\begin{rem}
    The proof of Proposition \ref{loop H inv} relies on $\cC$ being dualizable, and we do not know the generality of the equivalence \eqref{h gen}. If $D^*(H(\bF)/H_{gen, c}^{sh})\simeq \Vect$, i.e., the $H_{gen, c}^{sh}$-coinvariants of the category of D-modules on $H(\bF)$ is equivalent to $\Vect$, then the equivalence \eqref{h gen} follows. Although we have $D^*(H(\bF)/H_{gen, c}^{sh})\simeq \colim D^*(B(K_n\cap H_{gen, c}^{sh}))$ and each $K_n\cap H_{gen, c}^{sh}$ is contractible, being different from invariants, the category of coinvariants of $\Vect$ with respect to a contractible group is not necessarily equivalent to $\Vect$.
\end{rem}

\begin{rem}
   We comment on the relation with the notion of \textit{homological contractibility} of D.~Gaitsgory. Recall that a map of prestacks $\cX_1\rightarrow \cX_2$ is called universally homologically contractible, c.f., \cite[Appendix A.1.8]{[G5]}, if for any affine test scheme $S\rightarrow \cX_2$, the $!$-pullback functor of the category of D-modules along $\cX_1\underset{\cX_2}{\times} S\rightarrow S$ is fully faithful. Consider the map of the classifying stacks
   \begin{equation}
       BH_{gen, c}^{sh}\rightarrow BH(\bF).
   \end{equation}
   Given a map from a prestack $\cX$ to $BH(\bF)$, it determines a $H(\bF)$-torsor $\tilde{\cX}$ over $\cX$, and the Proposition \ref{loop H inv} implies that if $D^!(\cX)$ is dualizable, then the $!$-pullback functor along the base change is fully faithful. In particular, it holds for any $\cX$ being a finite type affine scheme $S$. 
\end{rem}

{
\begin{rem}\label{h gen rem}
   Through the above proof, we note that the equivalence \eqref{h gen} works for an arbitrary finite type algebraic group $H= H_r\ltimes H_u$.
\end{rem}
}
\begin{rem}
Since $\chi$ is trivial on $H^{\omega}_{C-c}$, and the set of geometric points of $H^{\omega}_{C-c}\backslash \Gr^{\omega}_G$ is the union of ${\cM_{= (\lambda, (\theta, \theta'))\cdot c}}\simeq H^{\omega}_{C-c}\backslash \sO_c^{(\lambda,(\theta,\theta'))}$, we may expect that a similar argument as Proposition \ref{strata eq} works for the whole Gaiotto category. 

 However, the stack $H^{\omega}_{C-c}\backslash \Gr^{\omega}_{{G}}$ classifies a $G$-bundle on $C$ with an $H^{\omega}$-reduction over $C-c$. It is a non-locally closed constructible subspace of $\cM_{\infty\cdot c}$, and the theory of equivariant D-modules is not well-behaved. We need to apply another method to prove the equivalence for the whole Gaiotto category.
\end{rem}

\begin{rem}
    An important remark is that the Verdier duality on the global model $\cM_{\infty\cdot c}$ preserves the generic Hecke equivariance structure (up to change $\kappa$ to $-\kappa$, and $\chi$ to $-\chi$). 

Restrict the Hecke correspondence diagram in Section \ref{4.3.3} to a substack locally of finite type containing the support of $\cF$, which we can assume to be of the form $\cM_{\bar{x},\leq (\lambda,(\theta,\theta'))\cdot c}$, and consider the closed substack $\Hecke^{\leq \mu}_{\bar{x}}$ in the preimage of $\Hecke_{\bar{x}}|_{\cM_{\bar{x},\leq (\lambda,(\theta,\theta'))\cdot c}}$ given by the Schubert variety of the Beilinson-Drinfeld affine Grassmannian. We denote $i_\mu$ the closed embedding, and $\overset{\gets}{h}_\mu=\overset{\gets}{h}\circ i_\mu$, $\overset{\to}{h}_\mu=\overset{\to}{h}\circ i_\mu$. So, we should prove $\overset{\gets}{h}_\mu^! (\BD(\cF))\simeq \overset{\to}{h}_\mu^!(\BD(\cF))\overset{!}{\otimes} (-\chi)$.

Taking the Verdier dual of the assumption, we obtain $\overset{\gets}{h}^* (\BD(\cF))\simeq \overset{\to}{h}^*(\BD(\cF))\overset{*}{\otimes} \BD(\chi)\simeq \overset{\to}{h}^*(\BD(\cF))\overset{!}{\otimes} (-\chi)$. Here, we note that $\chi$ is normalized by its $!$-stalk being concentrated in degree $0$, while $\BD(\chi)$ is normalized by its $*$-stalk.

Since both $\overset{\gets}{h}_\mu$ and $\overset{\to}{h}_\mu$ are locally trivial fibrations, the $!$-pullback and $*$-pullback differ by $*$-tensoring with the relative dualizing sheaf. Hence, it remains to show $\overset{\gets}{h}_\mu^! (\bc)\simeq \overset{\to}{h}_\mu^!(\bc)$, or $\overset{\gets}{h}^! (\bc)\simeq \overset{\to}{h}^!(\bc)$. Here, $\bc$ means the constant sheaf on ${\cM_{\leq (\lambda,(\theta,\theta'))\cdot c}}$. Note that $\tilde{\cM}_{\leq (\lambda,(\theta,\theta'))\cdot c}$ is pro-smooth over ${\cM}_{\leq (\lambda,(\theta,\theta'))\cdot c}$, after choosing a dimension theory, cf., \cite{[Ber],[Ras1]}, we may regard the $!$-pullback of $\bc$ as the constant sheaf $\bc'$ on $\tilde{\cM}_{\leq (\lambda,(\theta,\theta'))\cdot c}$ up to a cohomological shift.

In particular, we only need to show the $!$-pullback of $\bc'$ along the action map $m\colon H^\omega(\bF)_{\bar{x}}{\times} \tilde{\cM}_{\bar{x},\infty\cdot c}\rightarrow \tilde{\cM}_{\bar{x},\infty\cdot c}$ is given by $\pr_2^! (\bc')$. Here, $\pr_1,\pr_2$ are projections of $H^\omega(\bF)_{\bar{x}}{\times} \tilde{\cM}_{\bar{x},\infty\cdot c}$ to its first and second factors, respectively.

Let $s$ be the isomorphism $H^\omega(\bF)_{\bar{x}}{\times} \tilde{\cM}_{\bar{x},\infty\cdot c}\rightarrow H^\omega(\bF)_{\bar{x}}{\times} \tilde{\cM}_{\bar{x},\infty\cdot c}$, which sends $(h, m)$ to $(h,hm)$. Note that $\pr_1=\pr_1\circ s$ and $m=\pr_2\circ s$. Now, the claim follows from  
\begin{equation}
\omega_{H^\omega(\bF)_{\bar{x}}}\boxtimes \bc'\simeq s^!(\omega_{H^\omega(\bF)_{\bar{x}}}\boxtimes \bc')\simeq s^!\pr_2^!(\bc')\simeq m^!(\bc').    
\end{equation}
Here, the first isomorphism follows from $\omega_{H^\omega(\bF)_{\bar{x}}}\boxtimes \bc'\simeq \pr_1^*(\omega_{H^\omega(\bF)_{\bar{x}}})$ and there is $\pr_1^*(\omega_{H^\omega(\bF)_{\bar{x}}})\simeq s^*\pr_1^*(\omega_{H^\omega(\bF)_{\bar{x}}})$.
\end{rem}

Before we give a proof of Theorem \ref{locglob}, first we need  to introduce the Ran version of twisted Gaiotto category which connects the local Gaiotto category and the global Gaiotto category.

\section{Proof of local-global comparison}\label{Ran}
\subsection{Ran space}\label{Section 5.1}
{Recall that we defined the Ran space in Section \ref{4.3.1}.} The Ran space acquires a non-unital associative commutative algebra structure with respect to taking union,
\[\Ran_C\times \Ran_C\longrightarrow \Ran_C.\]

Fix a point $c\in C$, we can also consider the Ran space $\Ran_{C, c}$ with a marked point $c$, its $S$-points classify finite subsets of $\Maps(S,C)$ containing a distinguished point $S\rightarrow \{c\}$.

Following Beilinson-Drinfeld, we consider the Ran affine Grassmannian over $\Ran_C$ and $\Ran_{C, c}$. By definition, for a group $G$, the Ran affine Grassmannian ({or, }Beilinson-Drinfeld affine Grassmannian) $\Gr^\omega_{G, \Ran}$ (resp.\ $\Gr^\omega_{G, \Ran_c}$) classifies $(\bar{x}, \cP_G, \alpha)$, where $\bar{x}$ is an $S$-point in $\Ran_C$ (resp.\ $\Ran_{C, c}$), $\cP_G$ is a $G$-bundle on the curve $C\times S$, and $\alpha$ is an isomorphism of $\cP_G$ and the $G$-bundle $\cP_G^\omega=(\cP^{triv}_M, \cP^\omega_N)$ on the open curve $C\times S-\bar{x}$.

The key property of $\Gr^\omega_{G, \Ran_C}$ (resp.\ $\Gr^{\omega}_{G, \Ran_{C, c}}$) is that it is a factorization prestack over $\Ran_C$ (resp.\ a factorization module prestack over $\Ran_{C, c}$ with respect to $\Gr^{\omega}_{G, \Ran_C}$).

We are interested in some sub-prestacks of $\Gr^{\omega}_{G, \Ran_C}$ and $\Gr^{\omega}_{G, \Ran_{C, c}}$. Let us denote by $\overline{\sO}_{\Ran_{C}}$ (resp.\ ${(\overline{\sO}_{\Ran_{C, c}})_{\infty\cdot c}}$) the closed sub-prestack of $\Gr^{\omega}_{G, \Ran_C}$ (resp.\ $\Gr^{\omega}_{G, \Ran_{C, c}}$ ) such that $(\bar{x}, \cP_G, \alpha)$ belongs to $\overline{\sO}_{\Ran_C}$ (resp.\ ${(\overline{\sO}_{\Ran_{C, c}})_{\infty\cdot c}}$) if and only if the map  
\begin{equation}\label{5.1}
    \cV_{\tilde{\cP}_M}\overset{\alpha_M}{\longrightarrow} \cV_{\tilde{\cP}^{triv}_M}\longrightarrow \Ind(\cV)_{\cP^{\omega}_N}\overset{\alpha_N}{\longrightarrow} \Ind(\cV)_{\cP_N},
\end{equation}
which is \textit{a priori} defined on $C\times S-\bar{x}$ extends to a regular map on $C\times S$ (resp.\ $(C-c)\times S$). These prestacks preserve the factorization properties and are compatible with the inclusion maps, i.e.,
\begin{equation}
\begin{split}
     \overline{\sO}_{\Ran_{C}}\underset{\Ran_C}{\times} (\Ran_C\times \Ran_C)_{disj}\simeq (\overline{\sO}_{\Ran_{C}}\times \overline{\sO}_{\Ran_{C}})\underset{\Ran_C\times\Ran_C}{\times}(\Ran_C\times\Ran_C)_{disj}\\
     {(\overline{\sO}_{\Ran_{C, c}})_{\infty\cdot c}}\underset{\Ran_{C, c}}{\times} (\Ran_C\times \Ran_{C, c})_{disj}\simeq (\overline{\sO}_{\Ran_{C}}\times {(\overline{\sO}_{\Ran_{C, c}})_{\infty\cdot c}})\underset{\Ran_C\times\Ran_{C, c}}{\times}(\Ran_C\times\Ran_{C, c})_{disj}
\end{split}
\end{equation}
\begin{rem}
    Given a point $\{x_1, x_2,\cdots, x_m\}\in \Ran_C$, the fiber of $\overline{\sO}_{\Ran_C}$ is $\prod^{m}_{i=1} \overline{\sO}^0_{x_i}$. Similarly, the fiber of ${(\overline{\sO}_{\Ran_{C, c}})_{\infty\cdot c}}$ on $\{c, x_1, x_2,\cdots, x_m\}\in \Ran_{C, c}$ is $\Gr_{G,c}^\omega\times\prod^{m}_{i=1} \overline{\sO}^0_{x_i}$.
\end{rem}

Repeating the definition of $\pi: \Gr^{\omega}_{{G}}\longrightarrow \cM_{\infty\cdot c}$, we can also construct a map from ${(\overline{\sO}_{\Ran_{C, c}})_{\infty\cdot c}}$ (resp.\ $\overline{\sO}_{\Ran_C}$) to $\cM_{\infty\cdot c}$ (resp.\ $\cM$),
\begin{equation}
    \pi_{\Ran_C}: \overline{\sO}_{\Ran_{C}}\longrightarrow \cM,
\end{equation}
and
\begin{equation}\label{eq 5.4}
    \pi_{\Ran_{C, c}}: {(\overline{\sO}_{\Ran_{C, c}})_{\infty\cdot c}}\longrightarrow \cM_{\infty\cdot c}.
\end{equation}

\subsection{Local-global comparison}
{Now, let us prove that there is an equivalence between the local Gaiotto category and the global one. Let us first explain the idea of the proof.

Since $H^{\omega}\simeq \GL_M\ltimes U^{-,\omega}_{M,N}$, we hope to write the category of $(H^\omega,\chi)$-generic Hecke equivalent {twisted} D-modules as the $\GL_M$-generic Hecke invariants of the category of $(U_{M,N}^-,\chi)$-generic Hecke equivariant {twisted} D-modules on $\cM_{\infty\cdot c}$, i.e, {``$({\SD}_{\kappa}^{U^{-,\omega}_{M,N},\Hecke,\chi}(\cM_{\infty\cdot c}))^{\GL_M,\Hecke}$''}. Then, we compare the $(U_{M,N}^-,\chi)$-generic Hecke invariant category with $(U_{M,N}^{-,\omega}(\bF),\chi)$-invariant category, and compare $\GL_M$-generic Hecke invariant category with $\GL_M(\bF)$-invariant category.

To make sense of { ``$({\SD}_{\kappa}^{U^{-,\omega}_{M,N},\Hecke,\chi}(\cM_{\infty\cdot c}))^{\GL_M,\Hecke}$''}, we should first take $(U_{M,N}^{-,\omega},\chi)$-generic Hecke invariants and then $\GL_M$-generic Hecke invariants. It indicates that we consider the following limit

\begin{equation}
  \cC_{\kappa}^{glob}(M|N)':=  \lim_{\bar{x}\in \Ran_{C-c}(S), \bar{y}\in \Ran_{(C-c)\times S,\bar{x}}(S'), (S'\to S)\in \Sch^\aff} \cC_{\kappa}^{glob}(M|N)_{\bar{x},\bar{y}}
\end{equation}
Here, $\Ran^{}_{(C-c)\times S,\bar{x}}$ is the $S$-prestack whose $S'$-points are finite subsets of $\Maps(S',C-c)$ containing $\bar{x}\times_S S'$ (which is regarded as a subset of $\Maps(S',C-c)$ by composing $\bar{x}\in \Maps(S,C-c)$ and the map $S'\to S$), the category $\cC_{\kappa}^{glob}(M|N)_{\bar{x},\bar{y}}$ is the category of ($H_{\bar{x},\bar{y}},\chi_{\bar{y}}$)-equivariant twisted D-modules on the $H^\om(\bO)_{\bar{y}}$-bundle $\tilde{\cM}_{\bar{y}, \infty\cdot c}$ in Section \ref{second description}, $H_{\bar{x}, \bar{y}}:= (\GL_M(\bO)_{\bar{y}\setminus \bar{x}\times_S S'}\times \GL_M(\bF)_{\bar{x}\times_S S'}) \ltimes U_{M,N}^{-,\om}(\bF)_{\bar{y}} \subset H^\om(\bF)_{\bar{y}}$. This group is well-defined, namely, $\GL_M(\bO)_{\bar{y}\setminus \bar{x}\times_S S'}\times \GL_M(\bF)_{\bar{x}\times_S S'}\subset \GL_M(\bF)_{\bar{y}}$ is given by the maps from the $\cD_{\bar{y}}\setminus \bar{x}\times_S S'$ to $H^\om$, and $H_{\bar{x}, \bar{y}}$ is the subgroup given by the preimage of this group under the natural projection $H^\om(\bF)_{\bar{y}}\rightarrow \GL_M(\bF)_{\bar{y}}$.

The transition functors between $\cC_{\kappa}^{glob}(M|N)_{\bar{x},\bar{y}}$ are generated by the base change functors (similar to \eqref{4.5}) for both $\bar{x}$ and $\bar{y}$, and also forgetful functors (similar to \eqref{4.4}): $\cC_{\kappa}^{glob}(M|N)_{\bar{x}', \bar{y}'}\rightarrow \cC_{\kappa}^{glob}(M|N)_{\bar{x}, \bar{y}}$, if $\bar{x}\subset \bar{x}'\in \Ran_{C-c}^{}(S)$, $\bar{y}\subset \bar{y}'\in \Ran^{}_{(C-c)\times S,\bar{x}}(S')$. 

It can be regarded as a two-step continuous limit as in Section \ref{cont limit}, namely, we can first fix $\bar{x}$ and take limit with respect to $\bar{y}$, it is a procedure of taking $(U_{M,N}^{-,\om},\chi)$-generic Hecke invariants; then we can take limit with respect to $\bar{x}$, it is a procedure of taking $\GL_M$-generic Hecke invariants.

\begin{rem}
    Note that when $\bar{y}=\bar{x}$, then $\cC_\kappa^{glob}(M|N)_{\bar{x}, \bar{y}}$ is equivalent to our {previously} defined category $\cC_{\kappa}^{glob}(M|N)_{\bar{x}}$.

    There exists a transition functor from $\cC_{\kappa}^{glob}(M|N)_{\bar{x},\bar{y}}$ to $\cC_{\kappa}^{glob}(M|N)_{\bar{x}\times_S S'}$, and a transition functor from $\cC^{glob}_{\kappa}(M|N)_{\bar{y}}$ to $\cC_{\kappa}^{glob}(M|N)_{\bar{x},\bar{y}}$.  
\end{rem}

By coinitiality of the corresponding indexing categories (in fact, they are even reflexive), there is an equivalence of categories
\begin{equation}
    \cC_{\kappa}^{glob}(M|N)\simeq \cC_{\kappa}^{glob}(M|N)'.
\end{equation}

For a fixed $\bar{x}\in \Ran_{C-c}(S)$, we denote by $\cC_{\kappa}^{glob}(M|N)_{\bar{x}, \infty}$ the limit of $\cC_{\kappa}^{glob}(M|N)_{\bar{x},\bar{y}}$ for $\bar{y}${.  We} have
\begin{equation}
    \cC_{\kappa}^{glob}(M|N)\simeq \lim_{\bar{x}\in \Ran_{C-c}(S), S\in (\Sch^\aff)^\opp} \cC_{\kappa}^{glob}(M|N)_{\bar{x}, \infty}.
\end{equation}
}




{
We claim the following two propositions.
\begin{prop}\label{prop 2}
    \begin{equation}\label{eq prop 2}
    \lim_{\bar{x}\in \Ran_{C-c}(S), S\in (\Sch^\aff)^\opp} {\SD}_{\kappa}^{\GL_{M, (C-c)\times S-\bar{x}}\ltimes U^{-,\omega}_{M,N}(\bF),\chi}(\Gr^{\omega}_{G,c\times S})\simeq {\SD}_{\kappa}^{H^{\omega}(\bF),\chi}(\Gr^{\omega}_G).
    \end{equation}
\end{prop}
\begin{proof}
\begin{equation}\label{5.19}
\begin{split}
     & \lim_{\bar{x}\in \Ran_{C-c}(S), S\in (\Sch^\aff)^\opp} {\SD}_{\kappa}^{\GL_{M, (C-c)\times S-\bar{x}}\ltimes U^{-,\omega}_{M,N}(\bF),\chi}(\Gr^{\omega}_{G,c\times S})\\
     &\simeq {\SD}_{\kappa}^{\GL_{M, gen, c}\ltimes U^{-,\omega}_{M,N}(\bF),\chi} (\Gr^{\omega}_G)\\
     &\simeq {\SD}_{\kappa}^{\GL_M(\bF)\ltimes U^{-,\omega}_{M,N}(\bF), \chi}(\Gr^{\omega}_G).
\end{split}
\end{equation}    
The third equivalence is by Proposition \ref{loop H inv}. Namely, applying the Proposition to the case $H=\GL_M$, we obtain an equivalence 
\begin{equation}\label{apply glm}
{\SD}_{\kappa}^{\GL_{M, gen, c}} (\Gr^{\omega}_G)
     \simeq {\SD}_{\kappa}^{\GL_M(\bF)}(\Gr^{\omega}_G),    
\end{equation}
 and the third equivalence above is the restriction of the equivalence \eqref{apply glm} to the full subcategories of $(U_{M,N}^{-,\omega}(\bF),\chi)$-equivariant objects.
\end{proof}

\begin{prop}\label{prop 1}
For $\bar{x}\in \Ran_{C-c}(S)$, the $!$-pullback functor along $\Gr^{\omega}_G\rightarrow \cM_{\infty\cdot c}$ induces a fully faithful functor
\begin{equation}\label{genral hu}
     \cC_{\kappa}^{glob}(M|N)_{\bar{x}, \infty}\longrightarrow {\SD}_{\kappa}^{\GL_{M, (C-c)\times S-\bar{x}}\ltimes U^{-,\omega}_{M,N}(\bF)_{c\times S},\chi}(\Gr^{\omega}_{G,c\times S}).
\end{equation}
\end{prop}

\subsubsection{}
In this section, we prove Proposition \ref{prop 1}. 


We should prove that the pullback functor induces a fully faithful functor between full subcategories
\begin{equation}\label{5.11}
   \cC_\kappa^{glob}(M|N)_{\emptyset,\infty}:= {\SD}_{\kappa}^{U^{-,\omega}_{M,N},\Hecke,\chi}( \cM_{\infty\cdot c}) \longrightarrow {\SD}_{\kappa}^{\GL_{M,C-c}\ltimes U^{-,\omega}_{M,N}(\bF),\chi}(\Gr^{\omega}_G).
\end{equation}
Then, \eqref{genral hu} follows from the facts that \eqref{genral hu} is the Hecke invariant of \eqref{5.11} with respect to $\GL_M$ at $\bar{x}$ and taking the Hecke invariants preserves the full faithfulness.

First, let us denote by $(\overline{\Bun}^{\omega}_{U_{M,N}^-})_{\infty\cdot c}$ the $\omega$-renormalized Drinfeld compactification for $U^-_{M,N}\subset G$ with a possible pole at $c$. To be more precise, it classifies $(\cP_G, \sigma)$, where $\cP_G$ is a $G$-bundle on $C$, and a section $\sigma:C-c\rightarrow \omega^\rho\overset{T_{N-M-1}}{\times}\overline{U^-_{M,N}\backslash G}^{\aff}\overset{G}{\times} \cP_G\supset \omega^\rho\overset{T_{N-M-1}}{\times}{(U^-_{M,N}\backslash G)}\overset{G}{\times} \cP_G$. We note that there is an isomorphism of algebraic stacks
\begin{equation}
    \cM_{\infty\cdot c}\simeq \GL_{M,C-c}\backslash (\overline{\Bun}^{\omega}_{U_{M,N}^-})_{\infty\cdot c}.
\end{equation}

Thus, we should prove that the pullback functor along 
\begin{equation}\label{U localglobalmap}
    \Gr^{\omega}_G\longrightarrow (\overline{\Bun}^{\omega}_{U_{M,N}^-})_{\infty\cdot c}
\end{equation}
induces a functor
\begin{equation}\label{U fun}
    {\SD}_{\kappa}^{U^{-,\omega}_{M,N},\Hecke,\chi}((\overline{\Bun}^{\omega}_{U_{M,N}^-})_{\infty\cdot c})\longrightarrow {\SD}_{\kappa}^{U^{-,\omega}_{M,N}(\bF),\chi}(\Gr^{\omega}_G),
\end{equation}
and it is fully faithful. Here,  ${\SD}_{\kappa}^{U^{-,\omega}_{M,N},\Hecke,\chi}((\overline{\Bun}^{\omega}_{U_{M,N}^-})_{\infty\cdot c})$ is the category of ($U_{M,N}^{-, \om},\chi$)-generic Hecke equivariant twisted D-modules on the Drinfeld compactification. It can be defined by the same way as for $\cC_{\kappa}^{glob}(M|N)$ in Definition \ref{def 4.3.4}, as well as Section \ref{second description}.

}

Both statements above follow from the same proof as in the main theorem of \cite{[G3]}. Indeed, since $U_{M,N}^-$ is unipotent, the category ${\SD}_{\kappa}^{U^{-,\omega}_{M,N}(\bF),\chi}(\Gr^{\omega}_G)$ is a full subcategory of ${\SD}_{\kappa}(\Gr^{\omega}_G)$. To prove that the $!$-pullback along \eqref{U localglobalmap} induces \eqref{U fun}, we only need to show that the image of the pullback functor lands in ${\SD}_{\kappa}^{U^{-,\omega}_{M,N}(\bF),\chi}(\Gr^{\omega}_G)$.  We have to check that the restriction of the image to each $U^{-,\omega}_{M,N}(\bF)$-orbit is $(U^{-,\omega}_{M,N}(\bF), \chi)$-equivariant. If we repeat the proof of Proposition \ref{strata eq}, we obtain
\begin{cor}\label{cor 521}
    For any $U^{-,\omega}_{M,N}(\bF)$-orbit $\mathfrak{O}$ of $\Gr^{\omega}_G$, the category of $(U^{-,\omega}_{M,N}(\bF),\chi)$-equivariant {twisted} D-modules on this orbit, is equivalent to, the category of $(U^{-,\omega}_{M,N},\chi)$-generic Hecke equivariant {twisted} D-modules on the corresponding substack in $(\overline{\Bun}^{\omega}_{U_{M,N}^-})_{\infty\cdot c}$ (i.e., the image of $\mathfrak{O}$ under the projection).

    In particular, the pullback functor along \eqref{U localglobalmap} gives a functor \eqref{U fun}.
\end{cor}

Now, to prove Proposition \ref{prop 1}, we should prove
\begin{lem}\label{Dennis lem}
    The functor \eqref{U fun} is fully faithful.
\end{lem}
\begin{proof}
    The proof repeats the proof of \cite[Section 6]{[G3]}, we sketch the proof here.

    Similarly to ${(\overline{\sO}_{\Ran_{C, c}})_{\infty\cdot c}}$, we can also define a closed sub-prestack ${(\overline{\sO}_{U, \Ran_{C, c}})_{\infty\cdot c}}\subset \Gr^{\omega}_{G, \Ran_{C, c}}$, such that the fiber of ${(\overline{\sO}_{U, \Ran_{C, c}})_{\infty\cdot c}}$ over a point ${\bar{x}=\{c, x_1,x_2,\cdots,x_m\}}\in \Ran_{C, c}$ is just 
 \[\prod_{i=1}^m \overline{\mathfrak{O}}_{{x_i}}^0\times \Gr^{\omega}_{G, c}.\]

Here, $\overline{\mathfrak{O}}_{{x}}^0$ is the closure of the unital $U^{-,\omega}_{M,N}(\bF)$-orbit in $\Gr^{\omega}_G$.

We consider the following diagram
\begin{equation}\label{U diagram}
    \xymatrix{
&\Ran_{C, c}\times \Gr^{\omega}_{G,c}\ar[ld]^{\pr}\ar[rd]^{\unit_U}&\\
\Gr^{\omega}_{G,c}\ar[rd]^{\pi_U}&&{(\overline{\sO}_{U, \Ran_{C, c}})_{\infty\cdot c}}\ar[ld]^{\pi_{U, \Ran_{C, c}}}\\
&(\overline{\Bun}^{\omega}_{U_{M,N}^-})_{\infty\cdot c},&}
\end{equation}
{where, the map $\unit_U$ is defined by sending a point $(\bar{x}, (\cP_G,\alpha))$ to $(\bar{x}, \cP_G,\alpha|_{C-\bar{x}})$, $\pr$ is the projection, $\pi_U$ and $\pi_{U, \Ran_{C, c}}$ are defined similarly as \eqref{eq 4.14} and \eqref{eq 5.4}.}

    To prove that the pullback along $\pi_U$ is fully faithful, we note that the Ran space $\Ran_{C, c}$ is universally homologically contractible, so we only need to check that the pullback functor 
    \[\pr^!\pi_U^!: {\SD}_{\kappa}^{U^{-,\omega}_{M,N}, \Hecke,\chi}((\overline{\Bun}^{\omega}_{U_{M,N}^-})_{\infty\cdot c})\longrightarrow {\SD}_{\kappa}^{U^{-,\omega}_{M,N}(\bF),\chi}(\Ran_{C, c}\times \Gr^{\omega}_{G,c}) \]
    is fully faithful. Here,  $U^{-,\omega}_{M,N}(\bF)$ acts on the second factor of $\Ran_{C, c}\times \Gr^{\omega}_{G,c}$.

   Let $U^{-,\omega}_{M,N}(\bF)_{\Ran_{C, c}}$ be the factorization loop group prestack associated with $U^{-,\omega}_{M,N}$, it has a character $\chi_{\Ran_{C,c}}:U^{-,\omega}_{M,N}(\bF)_{\Ran_{C, c}}\rightarrow \BG_a$ whose fiber is the character map \eqref{chi}. We can consider ${\SD}_{\kappa}^{U^{-,\omega}_{M,N}(\bF)_{\Ran_{C, c}},\chi_{\Ran_{C, c}}}({(\overline{\sO}_{U, \Ran_{C, c}})_{\infty\cdot c}})$, the category of $(U_{M,N}^-(\bF)_{\Ran_{C, c}},\chi_{\Ran_{C, c}})$-equivariant twisted D-modules on ${(\overline{\sO}_{U, \Ran_{C, c}})_{\infty\cdot c}}$. By {the same proof as in} \cite[Theorem 6.2.5]{[G3]}, we have
   \begin{equation}
       {\SD}_{\kappa}^{U^{-,\omega}_{M,N}(\bF)_{\Ran_{C, c}},\chi_{\Ran_{C, c}}}({(\overline{\sO}_{U, \Ran_{C, c}})_{\infty\cdot c}})\simeq {\SD}_{\kappa}^{U^{-,\omega}_{M,N}(\bF),\chi}(\Ran_{C, c}\times \Gr^{\omega}_{G,c}).
   \end{equation}


Indeed, the same proof as in \cite[Section 6.4]{[G3]} repeats verbatim, and we just need to replace the group $\overset{\circ}{I}{}^j$ in \cite[2.5.1]{[G3]} (also, \cite[2.2]{[Ras]}) by the preimage of $U_{M,N}^{-,\omega}(\bO)/(U_{M,N}^{-,\omega})_j$ under the natural projection from $G^\omega(\bO)$ to its quotient by $G^\omega(\bO)/G^\omega_j$, here $(U_{M,N}^{-,\omega})_j$ and $G^\omega_j$ are the corresponding $j$-th congruence subgroups. Then, following the notations in \cite{[G3]}, we denote by $I^j:=\Ad_{t^{-j\rho_{M,N}}}(\overset{\circ}{I}{}^j)$,\footnote{it is denoted by $\overset{\circ}{I}_j$ in \cite{[Ras]}.} where $\rho_{M,N}$ is defined in Section \ref{renorm}. In principle, the only thing one uses about $\rho_{M,N}$ is that the conjugation $\Ad_{t^{-\rho_{M,N}}}$ stabilizes $L_{M+1,1,\cdots,1}^-(\bO)$ and enlarges $U_{M,N}^{-,\omega}(\bO)$, where $L_{M+1,1,\cdots,1}^-,U_{M,N}^{-,\omega}$ are the Levi subgroup and unipotent radical of $P_{M+1,1,\cdots,1}^-$. With the adaptation of the group $\overset{\circ}{I}{}^j$, the only argument that does not follow immediately from the proof as in \cite[Theorem 6.2.5]{[G3]} is that one needs an analogous theorem \cite[Theorem 2.7.1 (a)]{[Ras]}, and it reduces to proving \cite[Lemma 2.12.1]{[Ras]}. However, it follows from the general Lemma \ref{sam lem} below and the fact that the unital $U_{M,N}^{-,\omega}(\bF)$-orbit in $\Gr_G^\omega$ is minimal relevant.

As in the one point case, the $!$-pullback functor along $\pi_{U,\Ran_{C, c}}$ defines a functor
   \begin{equation}\label{compo 1}
       {\SD}_{\kappa}^{U^{-,\omega}_{M,N}, \Hecke,\chi}((\overline{\Bun}^{\omega}_{U_{M,N}^-})_{\infty\cdot c})\longrightarrow  {\SD}_{\kappa}^{U^{-,\omega}_{M,N}(\bF)_{\Ran_{C, c}},\chi_{\Ran_{C, c}}}({(\overline{\sO}_{U, \Ran_{C, c}})_{\infty\cdot c}}).
   \end{equation}

   Since the diagram \eqref{U diagram} is commutative, to prove that $\pr^!\pi_U^!$ is fully faithful, we only need to prove \eqref{compo 1} is fully faithful.

{For a test scheme $S$, we call an open subset $U$ of $C\times S$ a domain, if the fiber of $U\subset C\times S$ over any point in $S$ is dense in $C$. Following \cite[Appendix A]{[G5]}, we define $\Bun_{G}^{U_{M,N}^-,gen}$ to be the prestack which classifies a $G$-bundle on $C\times S$ with an $U^{-,\omega}_{M,N}$-reduction on some domain, and two points are identified if the reductions coincide on the intersection of the domains.}
 According to \cite[Theorem A.1.10]{[G5]}, the map
   \begin{equation}\label{uhc}
       \Gr^{\omega}_{{G}, \Ran_{C, c}}\longrightarrow \Bun_G^{U_{M,N}^-, gen}
   \end{equation}
   is universally homologically contractible. 
   
  Note that the map 
   \[{(\overline{\sO}_{U, \Ran_{C, c}})_{\infty\cdot c}}\longrightarrow (\overline{\Bun}^{\omega}_{U_{M,N}^-})_{\infty\cdot c}\]
   is the base change of \eqref{uhc} along $(\overline{\Bun}^{\omega}_{U_{M,N}^-})_{\infty\cdot c}\longrightarrow \Bun_G^{U_{M,N}^-, gen}$, hence the $!$-pullback functor 
   \begin{equation}
{\SD}_{\kappa}((\overline{\Bun}^{\omega}_{U_{M,N}^-})_{\infty\cdot c})\longrightarrow  {\SD}_{\kappa}({(\overline{\sO}_{U, \Ran_{C, c}})_{\infty\cdot c}})
   \end{equation}
   is fully faithful. Also, since $U_{M,N}^-$ is unipotent, $ {\SD}_{\kappa}^{U^{-,\omega}_{M,N}, \Hecke,\chi}((\overline{\Bun}^{\omega}_{U_{M,N}^-})_{\infty\cdot c})$ and ${\SD}_{\kappa}^{U_{M,N}^-(\bF)_{\Ran_{C, c}},\chi_{\Ran_{C, c}}}({(\overline{\sO}_{U, \Ran_{C, c}})_{\infty\cdot c}})$ are full subcategories of  $ {\SD}_{\kappa}((\overline{\Bun}^{\omega}_{U_{M,N}^-})_{\infty\cdot c})$ and ${\SD}_{\kappa}({(\overline{\sO}_{U, \Ran_{C, c}})_{\infty\cdot c}})$ respectively. We obtain the proposition.
\end{proof}


The following lemma is essentially an analog of \cite[Lemma 2.12.1]{[Ras]},
\begin{lem}\label{sam lem}
Let $P$ be a parabolic subgroup of $G$, its unipotent radical is $U_P$, and its Levi subgroup is $L$. Let $\rho_P$ be a cocharacter
   $\rho_P\colon \BG_m\rightarrow G$,
which lands in the center of $L$, and its pairing with every root in $U_P$ is strictly positive. Given $j\geq 0$, denote by $G_j$ the $j$-th congruence subgroup, $G_j^j:=\Ad_{t^{-j\cdot\rho_P}}G_j$ and $G(\bO)^j:=\Ad_{t^{-j\cdot\rho_P}}G(\bO)$.

Assume $\chi$ is an (additive) character of $U_P$. Then, if the unital $U_P(\bF)$-orbit in $\Gr_{G}$ is minimal relevant,\footnote{
   In particular, if $H=H_r\ltimes U_P$, $H_r\subset L$ and its conjugation action on $H_u$ preserves the character $\chi$, then the minimal relevance condition for the unital $H(\bF)$-orbit in $\Gr_G$ implies the minimal relevance condition for the unital $U_P(\bF)$-orbit in $\Gr_G$.
} then, any ($\kappa$-twisted) $(U_P(\bF),\chi)$-equivariant D-module on $U_P(\bF)G(\bO)^j/G_j^j$ extends cleanly to its closure.

\begin{proof}
    The twisting $\kappa$ does not influence the argument, we assume $\kappa=0$. One should show that for any $U_P(\bF)$-orbit in $\overline{U_P(\bF)G(\bO)^j}/G_j^j\setminus {U_P(\bF)G(\bO)^j}/G_j^j$ is irrelevant. Here, $\overline{U_P(\bF)G(\bO)^j}$ is the pullback of the closure of the unital $U_P(\bF)$-orbit in $G(\bF)/G(\bO)^j$ along the projection $G(\bF)\rightarrow G(\bF)/G(\bO)^j$.

    We assume the closure $\overline{U_P(\bF)G(\bO)^j}$ has the decomposition 
    \begin{equation}\label{dec uo}
        \overline{U_P(\bF)G(\bO)^j}=U_P(\bF)G(\bO)^j\cup \bigcup_g U_P(\bF)\cdot g\cdot G(\bO)^j.
    \end{equation}
    Note that since the geometric points of $\Gr_G$ and $\Gr_P$ are the same, we may find a representative in $L(\bF)$ for any $U_P(\bF)$-orbit in $\Gr_G$. By conjugating with $t^{-j\rho_P}$, we may assume any representative $g$ in \eqref{dec uo} belongs to $L(\bF)$. 

    Then, as in \cite[Lemma 2.12.1]{[Ras]}, we have
    \begin{equation}
    \begin{split}
        \SD^{U_P(\bF),\chi}(U_P(\bF)\cdot g\cdot G(\bO)^j/G_j^j)\simeq \SD^{U_P(\bF),\chi}(g U_P(\bF)\overset{U_P(\bO)^j}{\times} G(\bO)^j/G_j^j)\\ \simeq \SD^{U_P(\bO)^j,\chi_g}(G(\bO)^j/G_j^j),
    \end{split}
    \end{equation}
    where, $U_P(\bO)^j:=\Ad_{t^{-j\rho_P}}U_P(\bO)$ and $\chi_g:=\chi\circ \Ad_{g}$.

    Since $G_j^j$ contains $U_P(\bO)$ and acts trivially on $G(\bO)^j/G_j^j$, we should prove that under the minimal relevance assumption of the unital orbit, and ${U_P(\bF)gG(\bO)^j}\subset \overline{U_P(\bF)G(\bO)^j}$, then $(\chi_g)|_{U_P(\bO)}\neq 0$, i.e., $U_P(\bO)\nsubseteq \Ker(\chi_g)$. The latter equals $\Ad_{g}U_P(\bO)\nsubseteq \Ker(\chi)$.

    However, by conjugating, ${U_P(\bF)gG(\bO)^j}\subset \overline{U_P(\bF)G(\bO)^j}$ implies that ${U_P(\bF)gG(\bO)}\subset \overline{U_P(\bF)G(\bO)}$. Thus, the minimal relevance condition implies that $\Stab_{U_P(\bF)}(gG(\bO)/G(\bO))=U_P(\bF)\cap \Ad_{g}G(\bO)=\Ad_{g}U_P(\bO)\nsubseteq \Ker(\chi)$. The second equality follows from the fact that $g\in L(\bF)$.
\end{proof}
\end{lem}

\subsubsection{}
Now, we can finish the proof of Theorem \ref{locglob}.
\begin{proof}
According to Proposition \ref{prop 1}, the functor $\pi^!$ induces a fully faithful functor
\[ \cC_{\kappa}^{glob}(M|N)\simeq \lim_{\bar{x}\in \Ran_{C-c}(S), S\in (\Sch^{\aff})^\opp} \cC_{\kappa}^{glob}(M|N)_{\bar{x}, \infty}\longrightarrow \]
\[\longrightarrow\lim_{\bar{x}\in \Ran_{C-c}(S), S\in (\Sch^{\aff})^\opp} {\SD}_{\kappa}^{\GL_{M, (C-c)\times S-\bar{x}}\ltimes U^{-,\omega}_{M,N}(\bF)_{c\times S},\chi}(\Gr^{\omega}_{G,c\times S}).\]

By Proposition \ref{prop 2}, the latter category is equivalent to ${\SD}_{\kappa}^{H^{\omega}(\bF),\chi}(\Gr^{\omega}_G)$. So, we proved A). of Theorem \ref{locglob}.

To prove B)., note that we have already known its fully faithfulness, we should prove the functor \eqref{local-global} is essentially surjective. Note that the Gaiotto category is generated by $*$-direct image of the $(H^\omega(\bF),\chi)$-equivariant {twisted} D-modules on relevant orbits, now the proof follows from the stratawise equivalence, i.e., Proposition \ref{strata eq}.

\end{proof}
\subsection{Objects in global Gaiotto category}

Given a relevant weight $(\lambda, (\theta, \theta'))$, let us denote by ${\IC}^{(\lambda, (\theta, \theta'))}_{glob}$ the $!*$-extension of the twisted perverse $(H^\omega,\chi)$-generic Hecke equivariant rank $1$ local system on the locally closed substack ${\cM_{= (\lambda, (\theta, \theta'))\cdot c}}$ of $\cM_{\infty\cdot c}$. Similarly, we denote by $\Delta^{(\lambda, (\theta, \theta'))}_{glob}$ (resp, $\nabla^{(\lambda, (\theta, \theta'))}_{glob}$) the $!$ (resp. $*$) extension D-module extended from ${\cM_{= (\lambda, (\theta, \theta'))\cdot c}}$. By definition, ${\IC}^{(\lambda, (\theta, \theta'))}_{glob}$ is just the image of $H^0$ of the map $\Delta^{(\lambda, (\theta, \theta'))}_{glob}\rightarrow \nabla^{(\lambda, (\theta, \theta'))}_{glob}$.

{Let $d_g:= \dim \Bun^{\omega}_H$, then the $!$-stalks of the restriction of $\IC_{glob}^{0}$ to $\Bun_H^{\om}$ lives in the degree $d_g$.

We warn that since $\Bun^{\omega}_H$ is not equi-dimensional, $d_g$ is a \emph{function} which is constant on each component. The connected components of $\Bun^{\omega}_H$ is determined by $\pi_1(H)=\BZ$ according to the degree of the determinant line bundle of the associated $\GL_M$-bundle  (i.e., $\cP_M$ in Definition \ref{def 4.1.2}), and change the degree by $1$ will change the function value by $N-M-1$. Similarly, the connected components of any relevant ${\cM_{= (\lambda, (\theta, \theta'))\cdot c}}$ are indexed by $\BZ$ according to the degree of the associated $\GL_M$-bundle. The relative dimension of the connected components of ${\cM_{= (\lambda, (\theta, \theta'))\cdot c}}$ and ${\cM_{= (\tilde{\lambda}, (\tilde{\theta}, \tilde{\theta}'))\cdot c}}$ with the same degrees is given by $\langle 2\rho^\circ, (\lambda, (\theta, \theta'))-(\tilde{\lambda}, (\tilde{\theta}, \tilde{\theta}')\rangle$.

We denote by $\tilde{\Delta}^{(\lambda, (\theta, \theta'))}_{glob}$ the object which corresponds to $\tilde{\Delta}^{(\lambda,(\theta, \theta'))}_{loc}$ under the shifted equivalence $\pi^![d_g]$ of \eqref{local-global}. We claim that 

\begin{prop}\label{check loc glob t}
The perverse $t$-structure of $\cM_{\infty\cdot c}$ is uniquely determined by the following requirement: $\cF\in \cC^{glob}_{\kappa}(M|N)$ is coconnective if and only if 
    \[\Hom_{\cC_{\kappa}^{glob}(M|N)}(\tilde{\Delta}^{(\lambda,(\theta, \theta'))}_{glob}[k], \cF)=0,\]
 for any $k>0$ and relevant $(\lambda,(\theta,\theta'))$. 
    
\end{prop}

\begin{proof}
By definition, for any $\cF\in \cC^{glob}_{\kappa}(M|N)$, we have
\begin{equation}
    \RHom_{\cC^{glob}_{\kappa}(M|N)}(\tilde{\Delta}^{(\lambda, (\theta, \theta'))}_{glob}, \cF)\simeq \RHom_{\cC^{loc}_{\kappa}(M|N)}(\tilde{\Delta}^{(\lambda, (\theta, \theta'))}_{loc}, \pi^!(\cF)[d_g]).
\end{equation}

    While the latter, by definition, is isomorphic to the shifted $!$-fiber of $\pi^!(\cF)$ at $\sL^{(\lambda, (\theta, \theta'))}\in \sO^{(\lambda, (\theta, \theta'))}$ by $d_g+ \langle 2\rho^\circ, (\lambda,(\theta, \theta'))\rangle=\dim \cM_{= (\lambda, (\theta, \theta'))\cdot c}$.
\end{proof}

\begin{rem}
    Actually, the object $\tilde{\Delta}^{(\lambda, (\theta, \theta'))}_{glob}$ can also be defined by applying the left adjoint funcor $\Av_!^{(H^\omega,\Hecke,\chi)}$ to $\pi_!\delta_{\sL_{(\lambda, (\theta, \theta'))}}[-\dim {\cM_{= (\lambda, (\theta, \theta'))\cdot c}}]$. Here, $\Av_!^{(H^\omega,\Hecke,\chi)}$ is the partially-defined left adjoint functor of the forgetful functor 
    \begin{equation}
\oblv_{H^\omega,\Hecke,\chi}: \cC^{glob}_{\kappa}(M|N)\longrightarrow {\SD}_{\kappa}(\cM_{\infty\cdot c}).
\end{equation}
    
\end{rem}

\begin{rem}
    Note that both of $\nabla^{(\lambda,(\theta, \theta'))}_{glob}$ and $\nabla^{(\lambda,(\theta, \theta'))}_{loc}$ are characterized by the properties that their $!$-stalks on the corresponding orbit is a shifted $1$-dimensional vector space and do not have $!$-stalks outside, so the shifted $!$-pullback $\pi^![d_g]$ maps $\nabla^{(\lambda,(\theta, \theta'))}_{glob}$ to $\nabla^{(\lambda,(\theta, \theta'))}_{loc}$.
\end{rem}

}
 

{ As a consequence of Proposition \ref{check loc glob t}, we have}

\begin{cor}\label{t locglob}
    $\pi^![d_g]$ is $t$-exact with respect to the $t$-structure of $\cC_{\kappa}^{loc}(M|N)$ defined in Definition \ref{$t$-structure} and the naive perverse $t$-structure of ${\SD}_{\kappa}(\cM_{\infty\cdot c})$.
\end{cor}

Furthermore, by an analog of \cite[Theorem 3.2.2]{[G2]}, we have
\begin{lem}
    \begin{equation}
    {\Delta}^{(\lambda,(\theta, \theta'))}_{loc}\simeq \pi^!({\Delta}^{(\lambda,(\theta, \theta'))}_{glob})[d_g].
\end{equation}
\end{lem}

In particular, by the definition of ${\IC}^{(\lambda,(\theta, \theta'))}_{loc}$ and ${\IC}^{(\lambda,(\theta, \theta'))}_{glob}$, there is
\begin{cor}
    \begin{equation}
    {\IC}^{(\lambda,(\theta, \theta'))}_{loc}\simeq \pi^!({\IC}^{(\lambda,(\theta, \theta'))}_{glob})[d_g].
\end{equation}
\end{cor}
    Actually, since the $\IC$ sheaves are determined by the $t$-structures, $\pi^![d_g]$ is $t$-exact implies the above corollary.

\subsubsection{Locally compact D-modules}
Consider the map
\begin{equation}
(\overline{\Bun}^{\omega}_{U_{M,N}^-})_{\infty\cdot c}\longrightarrow \cM_{\infty\cdot c}.
\end{equation}
Taking $!$-pullback along the above map induces a functor
\begin{equation}\label{5.22}
    \cC_{\kappa}^{glob}(M|N)\longrightarrow {\SD}_{\kappa}((\overline{\Bun}^{\omega}_{U_{M,N}^-})_{\infty\cdot c}).
\end{equation}

\begin{defn}
    We denote $\cC^{glob, lc}_{\kappa}(M|N)$ the full subcategory of $\cC_{\kappa}^{glob}(M|N)$ consisting of those objects which become compact after applying the forgetful functor in \eqref{5.22} to ${\SD}_{\kappa}((\overline{\Bun}^{\omega}_{U_{M,N}^-})_{\infty\cdot c})$.
\end{defn}
{
\begin{rem}
    The locally compactness can be described more explicitly as follows: an object in $\cC_{\kappa}^{glob}(M|N)$ is locally compact if and only if it is supported on finitely many strata, and the $!$-stalk on each stratum is finite dimensional.
\end{rem}
}
\begin{prop}
    The equivalence functor $\pi^![d_g]$ induces an equivalence between the category of locally compact objects. Namely, we have
    \begin{equation}
        \pi^![d_g]: \cC_{\kappa}^{glob, lc}(M|N)\simeq \cC_{\kappa}^{loc, lc}(M|N).
    \end{equation}
\end{prop}
\begin{proof}
It follows from the fact that $\cC_{\kappa}^{glob, lc}(M|N)$ is generated by ${\Delta}^{(\lambda,(\theta, \theta'))}_{glob}$ via extensions and shifts, $\cC_{\kappa}^{loc,lc}(M|N)$ is generated by ${\Delta}^{(\lambda,(\theta, \theta'))}_{loc}$ via extensions and shifts, and $\pi^![d_g]$ sends ${\Delta}^{(\lambda,(\theta, \theta'))}_{glob}$ to ${\Delta}^{(\lambda,(\theta, \theta'))}_{loc}$.

{
\begin{rem}
    The categories $\cC_{\kappa}^{glob, lc}(M|N)$ and $\cC_{\kappa}^{loc, lc}(M|N)$ are also generated by $\IC^{(\lambda,(\theta,\theta'))}_{glob}$ and  $\IC^{(\lambda,(\theta,\theta'))}_{loc}$, respectively.
\end{rem}

}
\end{proof}




\subsection{Fusion product}
In the definition of $\cC_{\kappa}^{glob}(M|N)$, we notice that we can replace the marked point $c$ by any finite set of $C$. That is to say, similar to $\cM_{\infty\cdot c}$, we define the algebraic stack $\cM_{\infty\cdot \Ran}$ over $\Ran_C$ which classifies $(\bar{c}, \cP_M, \cP_N, \kappa_\cV)$, but here we require \eqref{4.3} to be an injective map with possible poles at $\bar{c}$,
\begin{equation}\label{5.24 k}
        \kappa_{\cV}: \cV_{\tilde{\cP}_M}\longrightarrow \Ind (\cV)_{\cP_N}(\infty \cdot \bar{c}), \textnormal{for any } \cV\in \Rep(\GL_{M+1}\times T_{N-M-1}).
    \end{equation}
Its fiber over $\bar{c}\in \Ran_C$ is denoted by $\cM_{\infty\cdot \bar{c}}$.

Recall that in the definition of the global Gaiotto category, we take $C=\BP^1$. Repeating the definition of $\cC_{\kappa}^{glob}(M|N)= {{\SD}}^{H^\omega,\Hecke,\chi}_{\kappa}(\cM_{\infty\cdot c})$, we can define the category of $(H^\omega,\chi)$ generic Hecke equivariant {twisted} D-modules on $\cM_{\infty\cdot \Ran}$, which is a category over $\Ran_C$.  To be more precise, for $\bar{x}\in \Ran_C(S)$, we consider the stack $\cM_{\bar{x},\infty\cdot \Ran}$ (over $S$) whose $S'$-points classify $(\bar{c}, \cP_G, \kappa_{\cV})$, where $\bar{c}\in \Ran_C(S')$, $\cP_G$ is a $G$-bundle on $C\times S'$, and $\kappa_{\cV}$ is a collection of maps like \eqref{5.24 k} with Pl\"{u}cker relations, we require $\bar{c}\cap (\bar{x}\underset{S}{\times}S')=\emptyset$ and $\kappa_{\cV}$ is a injective bundle map near $\bar{x}\underset{S}{\times}S'$. We consider the category of Hecke equivariant {twisted} D-modules on $\cM_{\bar{x},\infty\cdot \Ran}$, and $\cC_{\kappa}^{glob}(M|N)_{\Ran}$ is defined as the limit category with transition functors defined similarly to \eqref{4.4} and \eqref{4.5}. Similarly, we can define the category of Hecke equivariant {twisted} D-modules on $\cM_{\infty\cdot \bar{c}}$, which is canonically equivalent to the fiber of $\cC_{\kappa}^{glob}(M|N)_{\Ran}$ over $\bar{c}\in \Ran_C$.

By repeating the proof of Theorem \ref{locglob}, we can prove that for any point $\bar{c}\in \Ran_C$, the category $\cC_{\kappa}^{glob}(M|N)_{\bar{c}}$ is canonically equivalent to ${\SD}_{\kappa}^{H^{{\omega}}(\bF)_{\bar{c}},\chi_{\bar{c}}}(\Gr^{\omega}_{G,\bar{c}})$. 

In fact, it is enough to show that for any finite set $I$, the $!$-pullback functor induces an equivalence
\begin{equation}
    \cC_{\kappa}^{glob}(M|N)_{\Ran}|_{C^I}\simeq {\SD}_{\kappa}^{H^{{\omega}}(\bF)_{C^I},\chi_{C^I}}(\Gr^{\omega}_{G,C^I}).
\end{equation}
Here, $H^{{\omega}}(\bF)_{C^I}$ and $\Gr^{\omega}_{G,C^I}$ denote the factorization loop group prestack associated with $H^\omega$ and the Beilinson-Drinfeld affine Grassmannian on $C^I$.

The proof of Theorem \ref{locglob} repeats verbatim, while one only needs the following modification for the full faithfulness. Namely, as in the proof of Lemma \ref{Dennis lem}, one should show that there is an analog of \cite[Theorem 6.2.5]{[G3]},
     \begin{equation}\label{dennis 6.2.5+}
       {\SD}_{\kappa}^{U^{-,\omega}_{M,N}(\bF)_{\Ran_{C, I}},\chi_{\Ran_{C, I}}}({(\overline{\sO}_{U, \Ran_{C, I}})_{\infty\cdot C^I}})\simeq {\SD}_{\kappa}^{U^{-,\omega}_{M,N}(\bF)_{C^I},\chi}(\Ran_{C, I}\times \Gr^{\omega}_{G,C^I}).
   \end{equation}
   Here, $\Ran_{C,I}$, $U^{-,\omega}_{M,N}(\bF)_{\Ran_{C, I}}$, ${(\overline{\sO}_{U, \Ran_{C, I}})_{\infty\cdot C^I}}$ are defined similarly as $\Ran_{C,c}$, $U^{-,\omega}_{M,N}(\bF)_{\Ran_{C, c}}$, ${(\overline{\sO}_{U, \Ran_{C, c}})_{\infty\cdot c}}$, while, rather than choose a marked point $c$, we define them as prestacks over $C^I$, and the corresponding point in $C^I$ (which means $I$ points in $C$) is the marked point. In particular, $\Ran_{C,c}$, $U^{-,\omega}_{M,N}(\bF)_{\Ran_{C, c}}$, ${(\overline{\sO}_{U, \Ran_{C, c}})_{\infty\cdot c}}$ are the fiber of $\Ran_{C,I}$, $U^{-,\omega}_{M,N}(\bF)_{\Ran_{C, I}}$, ${(\overline{\sO}_{U, \Ran_{C, I}})_{\infty\cdot C^I}}$ over $\{c,c,\cdots,c\}\in C^I$. For example, for a scheme $S$, the set of $S$-points of $\Ran_{C,I}$ classifies an $S$-point $\bar{c}$ in $C^I$ (equivalently, $I$ elements in $\Maps(S,C)$), and a finite subset $\bar{x}$ of $\Maps(S,C)$ which contains $\bar{c}$.

   A direct adaptation of the proof in \cite[Theorem 6.2.5]{[G3]} works here with the modification that one should just replace the group $\overset{\circ}{I}{}^j$ in \cite[2.5.1]{[G3]} by the group $(\overset{\circ}{I}{}^j)_{C^I}$ over $C^I$. Here, the group $(\overset{\circ}{I}{}^j)_{C^I}$ is a subgroup of the factorization loop group $G^\omega(\bF)_{C^I}$, and its fiber over the open locus $\overset{\circ}{C}{}^I$ is the same as the $I$-th product of $\overset{\circ}{I}{}^j$, and it extends to $C^I$ by fusion. For example, if $|I|=2$, then, the fiber over a point in the open locus $\overset{\circ}{C}{}^2=C^2\setminus \Delta$ is $\overset{\circ}{I}{}^{j}\times \overset{\circ}{I}{}^{j}$, and the fiber over a point in the diagonal is $\overset{\circ}{I}{}^{2j}$.

\subsubsection{}Assume $\bar{c}=\{c_1,c_2,\cdots,c_n\}\in \Ran_C$. Due to the factorization structure of ${\SD}_{\kappa}^{H^\omega(\bF)_{\Ran},\chi_{\Ran}}(\Gr^{{\omega}}_{G,\Ran})$, we have a canonical equivalence
 ${\SD}_{\kappa}^{H^\omega(\bF)_{\bar{c}},\chi_{\bar{c}}}(\Gr^{{\omega}}_{G,\bar{c}})\simeq \otimes_{i} {\SD}_{\kappa}^{H^{{\omega}}(\bF)_{{c_i}},\chi_{{c_i}}}(\Gr^{{\omega}}_{G,{c_i}})$. It implies that there is a canonical equivalence 
 \begin{equation}\label{5.27}
     \cC_{\kappa}^{glob}(M|N)_{\bar{c}}\simeq \cC_{\kappa}^{glob}(M|N)_{c_1}\otimes \cC_{\kappa}^{glob}(M|N)_{c_2}\otimes\cdots \otimes \cC_{\kappa}^{glob}(M|N)_{c_n}.
 \end{equation}

 Consider $\BA^1=\BP^1-\{\infty\}$ and the restriction of $\cC_{\kappa}^{glob}(M|N)_{\Ran}$ on $\BA^2$ with the diagonal $\Delta$ removed. On $\BA^2-\Delta$, we have a canonical equivalence of categories
 \begin{equation}\label{5.24}
  \cC_{\kappa}^{glob}(M|N)_{\Ran}|_{\BA^2-\Delta}\simeq (\cC_{\kappa}^{glob}(M|N)_{\Ran}|_{\BA^1}\otimes \cC_{\kappa}^{glob}(M|N)_{\Ran}|_{\BA^1})|_{\BA^2-\Delta}.
\end{equation}


\subsubsection{}
Given two objects $\cF_1, \cF_2\in \cC^{glob}_{\kappa}(M|N)$. Due to the $\BA^1$-equivariance of the category $\cC_{\kappa}^{glob}(M|N)_{\Ran}|_{\BA^1}$ and \eqref{5.24}, we can form an object in  $\cC_{\kappa}^{glob}(M|N)_{\Ran}|_{\BA^2-\Delta}$, let us denote the underlying twisted D-module on $\cM_{\infty \cdot \Ran}|_{\BA^2-\Delta}$ by $\cF_1\boxtimes \cF_2|_{\BA^2-\Delta}$. One can choose the assignment such that the functor $\cF_1, \cF_2\mapsto \cF_1\boxtimes \cF_2|_{\BA^2-\Delta}$ is $t$-exact.

We define the fusion product of  $\cF_1$ and $\cF_2$ by 
\begin{equation}
\cF_1\star \cF_2:= \pr_{\Delta,*} \psi_{x-y} (\cF_1\boxtimes \cF_2|_{\BA^2-\Delta})[-1]\in \cC_{\kappa}^{glob}(M|N)    
\end{equation}
Here, $x,y$ are two coordinates of $\BA^2$, $\pr_{\Delta}$ denotes the projection from $\cM_{\infty\cdot \Ran}|_{\Delta=\BA^1}=\BA^1\times \cM_{\infty\cdot 0}$ to $\cM_{\infty\cdot 0}$, and $\psi_{x-y}$ denotes the nearby cycles functor.

A priori, it is not clear why the fusion product provides a braided monoidal structure of the Gaiotto category (one has to check the associativity of functors defined by certain nearby cycles, which is not obvious at all!). However, as a direct corollary of the conservativity of the functor $F$, it is indeed associative. For more details, see Corollary \ref{cor fusion}.

\section{SW-Zastava space}\label{SW-Zastava space}
\subsection{SW-Zastava space}
The Zastava space was firstly introduced by Finkelberg-Mirkovi\'{c} in \cite{[FM]}. It turns out to be a very useful tool to relate the category of D-modules on the affine Grassmannian to the category of representations of groups. In \cite{[SW]}, the authors extended the definition of the Zastava space to any spherical variety.

The following definition of the Zastava space comes from \cite[Section 3.3]{[SW]}, {and we call it SW-Zastava space in the paper}.

\begin{defn}
Let $\cX$ be a spherical variety for $G$, which has an open dense $B$-orbit $\cX'$, such that $\cX'/B=\pt$. Then the Zastava space for $\cX$ is defined as:
\begin{equation}
     \Maps_{\gen} (C, \cX/B\supset \pt).
\end{equation}
The set of its $S$-points classifies the maps from $C\times S$ to $\cX/B$ which generically land in $\cX'/B=\pt$.
\end{defn}

Although the content of this section is also true in more general setting, to be more precise about the renormalization by $\om$, we focus on  the renormalized SW-Zastava space $\cY$ associated with $\cX= \overline{H\backslash G}^{\textnormal{aff}}$, where $G=\GL_M\times \GL_N$, and $H=\GL_M\ltimes U_{M,N}^{-,\om}$.

To give the precise definition of $\cY$, first, we note that for any point in $\Maps_{\gen} (C, \cX/B\supset \pt)$, the composition 
\[C\longrightarrow \cX/B\longrightarrow \pt/B\]
gives rise to a map from $\Maps_{\gen} (C, \cX/B\supset \pt)$ to $\Bun_B$. Also, by forgetting the $B$-structure, we have a map from $\Maps_{\gen} (C, \cX/B\supset \pt)$ to $\Maps_{\gen} (C, \cX/G\supset \overset{\circ}{\cX}/G)$. One can see that the composition
\[C\longrightarrow \pt/B\longrightarrow \pt/G\]
coincides with the $G$-bundle obtained from $\Maps_{\gen} (C, \cX/G\supset \overset{\circ}{\cX}/G)$. Hence, we obtain a map
\begin{equation}
\Maps_{\gen} (C, \cX/B\supset \pt)\longrightarrow \Maps_{\gen} (C, \cX/G\supset \overset{\circ}{\cX}/G)\underset{\Bun_{G}}{\times}\Bun_{B}.
\end{equation}

According to \cite[Section 3.5]{[SW]}, the above map is an open embedding whose image is the open locus such that the generalized $H$-reduction and $B$-reduction are generically transversal. 

For the $\omega$-renormalized version, we define \[\cY:= \cM\underset{\Bun_G}{\times'}\Bun_B.\]

Here, $\times'$ denotes the open subset of the relative product such that the generalized $H$-reduction and $B$-reduction are generically transversal.

{
\begin{rem}
    Similar{ly} to Section \ref{sec 4.2}, if we consider the $\om^\rho$-twist of $\cX/B$ and $\cX'/B$, we can also define $\cY$ by
    \[\cY:= \Sect_{\gen}(C,(\cX/B)_{\om^\rho} \supset (\cX'/B)_{\om^\rho}).\]

    Here, $(\cX/B)_{\om^\rho} := (\cX/B) \times^{T_{N-M-1}} \om^\rho$, and $(\cX'/B)_{\om^\rho}$ is defined similarly.
\end{rem}
}

\subsubsection{}
The transversal condition can also be interpreted in another way. First, we note that a data of a point in $\cY$ gives rise to the following maps.
\begin{itemize}
    \item For $k=1,2,\cdots , M+1$,
    \begin{equation}\label{cY 1}
       \Lambda^k \cV^N_{\cP_N}\longrightarrow \Lambda^k(\cV^M_{\cP_M}\oplus \cO),
    \end{equation}
    which can be realized as 
    \[\Lambda^k \cV^N_{\cP_N}\longrightarrow \Lambda^k \cV^M_{\cP_M},\]
    and \[\Lambda^k \cV^N_{\cP_N}\longrightarrow\Lambda^{k-1}\cV^M_{\cP_M}.\]
    \item For $k=1,2,\cdots, N-M-1$,
    \begin{equation}\label{cY 2}
       \Lambda^{M+k+1} \cV^N_{\cP_N}\longrightarrow \det \cV^M_{\cP_M} \otimes \omega_C^{\otimes 1+2+\cdots+k}.
    \end{equation}
    \end{itemize}
    \begin{itemize}
    \item For $k=1,2,\cdots, M$, 
    \begin{equation}\label{cY 3}
       \Lambda^k \cV^M_{\cP_M}\longrightarrow \cL_k ,
    \end{equation}
   which is a surjective map of vector bundles.
    \item For $k=1,2,\cdots,N$,
    \begin{equation}\label{cY 4}
       \cK_k\longrightarrow \Lambda^k \cV^N_{\cP_N},
    \end{equation}
    which is an injective map of vector bundles.
\end{itemize}

Here, the maps \eqref{cY 1} and \eqref{cY 2} arise from the generic $H$-reduction (i.e., data of $\cM$) by: for \eqref{cY 1}, we take $\cV$ to be the \textit{dual} of exterior power of the standard representation of $\GL_{M+1}$ and tensor with the trivial representation of $T_{N-M-1}$; for \eqref{cY 2}, we need to take the tensor of the dual of the determinant representation of the $\GL_{M+1}$ and the $1$-dimensional representation of $T_{N-M-1}$ corresponding to $(-1,\cdots,-1,0,\cdots,0)$, and the maps \eqref{cY 3} and \eqref{cY 4} arise from the $B$ (i.e., $B_M^-\times B_N$)-reduction of $\cP_G$.

By taking compositions of above maps, we obtain the following maps:
\begin{enumerate}
    \item  for $k=1,2,\cdots, M$,
    \begin{equation}\label{cY 5}
    \begin{split}
        \cK_k\longrightarrow \Lambda^k \cV^N_{\cP_N}\longrightarrow \Lambda^k \cV^M_{\cP_M}\longrightarrow  \cL_k,\\
        \cK_{k+1}\longrightarrow \Lambda^{k+1} \cV^N_{\cP_N}\longrightarrow \Lambda^{k}\cV^M_{\cP_M}\longrightarrow \cL_{k},\\
        \cK_{1}\longrightarrow\Lambda^{1}\cV^N_{\cP_N}\longrightarrow\cL_0:= \cO.
    \end{split}
\end{equation}
\item For $k=1,2,\cdots, N-M-1$,
\begin{equation}\label{cY 6}
   \cK_{M+k+1} \longrightarrow \Lambda^{M+k+1}\cV^N_{\cP_N}\longrightarrow \det \cV^M_{\cP_M}\otimes \omega_C^{\otimes 1+\cdots+k}\simeq \cL_M\otimes \omega_C^{\otimes 1+\cdots+k}.
\end{equation}
\end{enumerate}
The transversal condition means that the above maps are non-zero. 


\subsubsection{}
The SW-Zastava space admits a map to the Configuration space, ref \cite[Section 3.2,3.3]{[SW]}. Modulo the complexity of $\omega$-renormalization, the SW-Zastava space is $\Maps_{\gen} (C, \cX/B\supset \pt)$ and the Configuration space is $\Maps_{\gen} (C, (\cX/\!/U)/T\supset \pt)$, the map from the SW-Zastava space to the configuration space is just given by applying $\Maps(C,-/T)$ to the natural map $\cX/U\rightarrow \cX/\!/U$, and a $\omega$-renormalized version applies to our case.

Using the above description of the SW-Zastava space, we can describe this map (i.e., \cite[(3.3)]{[SW]}) more precisely.

Namely, since we require the maps \eqref{cY 5} and \eqref{cY 6} to be non-zero, we can write
\begin{equation}
    \begin{split}
        \cK_i= \cL_{i-1}(\Delta_1+\cdots+\Delta_{i-1}-E_1-\cdots-E_i)= \cL_i(\Delta_1+\cdots+\Delta_{i}-E_1-\cdots-E_i),\\
        \cL_M\otimes \omega_C^{\otimes 1+\cdots+k}(\Delta_1+\cdots+\Delta_{M}-E_1-\cdots-E_{M+k+1})=\cK_{M+k+1}.
    \end{split}
\end{equation}

Note that $\cL_0\simeq \cO$, one can determine all divisors $\Delta_i$ and $E_j$ inductively. The map sends to $D=-\sum \delta_i\cdot \Delta_i+ \sum \epsilon_i\cdot E_i$ defines a morphism
\begin{equation}
    v: \cY\longrightarrow C^\bullet.
\end{equation}

It is well-known that the map from the SW-Zastava space to the configuration space a factorizable map (e.g., \cite[Proposition 3.4.1]{[SW]}), i.e.,

\begin{equation}\label{Zastava}
      \cY\underset{C^\bullet}{\times} (C^\bullet\times C^\bullet)_{disj} \simeq (\cY\times \cY)\underset{C^\bullet\times C^\bullet}{\times} (C^\bullet\times C^\bullet)_{disj}. 
\end{equation}



\subsection{Configuration affine Grassmannian}
{To make the factorization structure more clear and define a purely local functor to the category of factorization modules}, we recall the Beilinson-Drinfeld affine Grassmannian over the Configuration spaces $C^\bullet$ and $C^\bullet_{\infty\cdot c}$.

\begin{defn}
We denote by $\Gr^{\omega}_{G, C^\bullet}$ the ind-scheme which classifies triples $(D, \cP_G, \alpha)$, where $D= -\sum \delta_i\cdot \Delta_i+ \sum \epsilon_i \cdot E_j\in C^\bullet$, $\cP_G=(\cP_M, \cP_N)$ is a $G$-bundle on $C$, and $\alpha$ is an isomorphism of $\cP_G$ with $(\cP^{triv}_M, \cP^{\omega}_N=\omega_C^{\rho}\overset{T_{N-M-1}}{\times} \GL_N)$ over $C-\supp(D)$. Here, $\supp(D)$ denotes the support of $D$.

Similarly, we denote by $\Gr^{\omega}_{G, C^\bullet_{\infty\cdot c}}$ the ind-stack which classifies $(D, \cP_G, \alpha)$, where $D\in C^{\bullet}_{\infty\cdot c}$, and $\cP_G$ and $\alpha$ are the same as above.
\end{defn}

The factorization property of $\Gr^{\omega}_{G, C^\bullet}$ says that there is an isomorphism,
\begin{equation}
    \Gr^{\omega}_{G, C^\bullet}\underset{C^\bullet}{\times} (C^\bullet\times C^\bullet)_{disj} \simeq (\Gr^{\omega}_{G, C^\bullet}\times \Gr^{\omega}_{G, C^\bullet})\underset{C^\bullet\times C^\bullet}{\times} (C^\bullet\times C^\bullet)_{disj}. 
\end{equation}

The ind-scheme $C^\bullet_{\infty\cdot c}$ is a factorization module space with respect to $\Gr^{\omega}_{G, C^\bullet}$, i.e., there is an isomorphism
\begin{equation}
    \Gr^{\omega}_{G, C^\bullet_{\infty\cdot c}}\underset{C_{\infty\cdot c}^\bullet}{\times} (C^\bullet\times C^\bullet_{\infty\cdot c})_{disj} \simeq (\Gr^{\omega}_{G, C^\bullet}\times \Gr^{\omega}_{G, C^\bullet_{\infty\cdot c}})\underset{C^\bullet\times C^\bullet_{\infty\cdot c}}{\times} (C^\bullet\times C^\bullet_{\infty\cdot c})_{disj}. 
\end{equation}

For an $S$-point $y\in \Maps_{\gen}(C, \cX/B\supset \pt)$, $y\colon C\times S\rightarrow \cX/B$, $y^{-1}(\cX'/B)= C\times S-\supp(D)$. So, $y$ defines a $B$-bundle on $C\times S$ with a section 
\[C\times S-\supp(D)\longrightarrow H\backslash \overset{\circ}{G}\overset{B}{\times} \cP_B\simeq \cP_B.\]
In particular, the above assignment defines a map 
\begin{equation}
    \Maps_{\gen}(C, \cX/B\supset \pt)\longrightarrow \Gr_{B, C^\bullet_{\infty\cdot c}}.
\end{equation}

Similarly, we have
\begin{equation}
    \cY\longrightarrow \Gr^{\omega}_{B, C^\bullet_{\infty\cdot c}}.
\end{equation}

By composing with the induction map
\[\Gr^{\omega}_{B, C^\bullet}\longrightarrow \Gr^{\omega}_{G, C^\bullet},\]
we obtain a map
\begin{equation}
    \cY\longrightarrow \Gr^{\omega}_{G, C^\bullet}.
\end{equation}

The above assignment extends to a map from the compactified SW-Zastava space $\overline{\cY}:= \cM\underset{\Bun_G}{\times'} \overline{\Bun}_B$ to $\Gr^{\omega}_{G, C^\bullet}$.

The following lemma is proved in \cite[Lemma 4.1.2]{[SW]}.
\begin{lem}
The map
\begin{equation}
    \overline{\cY}\longrightarrow \Gr^{\omega}_{G, C^\bullet}
\end{equation}
is a factorizable closed embedding.
\end{lem}

Using the above lemma, we rewrite $\overline{\cY}$ as a closed substack of $\Gr^{\omega}_{G, C^\bullet}$. Namely, $\overline{\cY}$ classifies the triples $(D, \cP_G, \alpha)$. Here,
\begin{itemize}
    \item $D=-\sum \delta_i\cdot \Delta_i+ \sum \epsilon_i\cdot E_i$,
    \item $\cP_G=(\cP_M, \cP_N)$ is a $G$-bundle on $C\times S$,
    \item $\alpha=(\alpha_M,\alpha_N)$ is an isomorphism of $G$-bundles on $C\times S-\supp(D)$,
    \begin{equation}
        \begin{split}
            \alpha_M: \cP_M|_{C\times S-\supp(D)}\simeq \cP_M^{triv}|_{C\times S-\supp(D)} \\
            \alpha_N: \cP_N|_{C\times S-\supp(D)}\simeq \cP_N^{\omega}|_{C\times S-\supp(D)}.
        \end{split}
    \end{equation}
    \item such that the composed maps
    \begin{equation}\label{6.18}
    \begin{split}
          \Lambda^k \cV^M_{\cP_M}\overset{\alpha_M}{\longrightarrow} \Lambda^k\cV^M_{\cP^{triv}_M}\longrightarrow \cL_k:= \cO(-\Delta_1-\cdots-\Delta_k), \text{for } k=1,2,\cdots, M,\\
          \cK_k:= \cO_C(-E_1-\cdots-E_k)\longrightarrow\Lambda^k \cV^N_{\cP_N^\omega}\longrightarrow\Lambda^k\cV^N_{\cP_N}, \text{for } k=1,2,\cdots, M+1,\\
          \cK_k:= \omega_C^{\otimes 1+2+\cdots+ (k-M-1)} (-E_1-\cdots-E_k)\longrightarrow\Lambda^k \cV^N_{\cP^\omega_N}\overset{\alpha_N^{-1}}{\longrightarrow}\Lambda^k\cV^N_{\cP_N}, \text{for } k=M+2,\cdots,N,
    \end{split}
    \end{equation}
    which are \textit{a priori} defined on $C\times S-\supp(D)$ extend to regular maps {of vector bundles} on $C\times S$,
    \item and $\forall\ \cV\in \Rep(\GL_{M+1}\times T_{N-M-1})$, the map
    \begin{equation}\label{6.19}
\cV_{\tilde{\cP}_M}\overset{\alpha_M}{\longrightarrow} \cV_{\tilde{\cP}^{triv}_M}\longrightarrow \Ind(\cV)_{\cP_N^\omega}\overset{\alpha^{-1}_N}{\longrightarrow} \Ind(\cV)_{\cP_N}
    \end{equation}
    which is a priori defined on $C\times S-\supp(D)$ extends to a regular map {of vector bundles} on $C\times S$.
\end{itemize}
\begin{rem}
Under the above identification, a point of $\overline{\cY}$ belongs to the algebraic stack $\cY$ if and only if the first map of \eqref{6.18} is surjective and the second, third maps of \eqref{6.18} are injective, for any $k$. 
\end{rem}


\subsubsection{}
Similar to $\cY$ and $\overline{\cY}$, one can also define the SW-Zastava space with a marked point, \[\cY_{\infty\cdot c}:= \cM_{\infty\cdot c}\underset{\Bun_G}{\times'} \Bun_B.\]
As well as the compactified SW-Zastava space with a marked point,
\[\overline{\cY}_{\infty\cdot c}:= \cM_{\infty\cdot c}\underset{\Bun_G}{\times'} \overline{\Bun}_B.\]

Both of ${\cY}_{\infty\cdot c}$ and $\overline{\cY}_{\infty\cdot c}$ are substacks of $\Gr^{\omega}_{G, C^\bullet_{\infty\cdot c}}$. Using the identifications of $\cY_{\infty\cdot c}$ and $\overline{\cY}_{\infty\cdot c}$ as substacks of $\Gr^{\omega}_{G, C^\bullet_{\infty\cdot c}}$, the algebraic ind-stack $\overline{\cY}_{\infty\cdot c}$ classifies the triples $(D, \cP_G, \alpha)$, where $D\in C^\bullet_{\infty\cdot c}$, $\cP_G$ and $\alpha$ are the same as above in $\cY$, but we require that \eqref{6.18} extends to a regular map on $(C-c)\times S$. The algebraic ind-stack $\cY_{\infty\cdot c}$ is open in $\overline{\cY}_{\infty\cdot c}$ where we require the first map of \eqref{6.18} is surjective and the second, third maps of \eqref{6.18} are injective. 

Denote by $v$ the projection map
\[v: \overline{\cY}_{\infty\cdot c}\longrightarrow C_{\infty\cdot c}^\bullet.\]
It is a factorizable map, i.e.,
\begin{equation}\label{Zastava c}
\begin{split}
     \cY_{\infty\cdot c}\underset{C_{\infty\cdot c}^\bullet}{\times} (C^\bullet\times C^\bullet_{\infty\cdot c})_{disj} \simeq (\cY\times \cY_{\infty\cdot c})\underset{C^\bullet\times C^\bullet_{\infty\cdot c}}{\times} (C^\bullet\times C_{\infty\cdot c}^\bullet)_{disj}, \\
     \overline{\cY}_{\infty\cdot c}\underset{C_{\infty\cdot c}^\bullet}{\times} ({C}^\bullet\times {C}^\bullet_{\infty\cdot c})_{disj} \simeq (\overline{\cY}\times \overline{\cY}_{\infty\cdot c})\underset{C^\bullet\times C^\bullet_{\infty\cdot c}}{\times} (C^\bullet\times C_{\infty\cdot c}^\bullet)_{disj}
\end{split}
\end{equation}

\subsection{Semi-infinite intersection}
In the last section, we realize the SW-Zastava spaces $\cY$, $\cY_{\infty\cdot c}$ and their compactifications as sub-prestacks of $\Gr^{\omega}_{G, C^\bullet}$ and $\Gr^{\omega}_{G, C^\bullet_{\infty\cdot c}}$. In this section, we will rewrite SW-Zastava spaces as intersections of two sub-prestacks of $\Gr^{\omega}_{G, C^\bullet}$ and $\Gr^{\omega}_{G, C^\bullet_{\infty\cdot c}}$, respectively.

\subsubsection{}
Let $\overline{\sO}_{C^\bullet}$ (resp.\ ${(\overline{\sO}_{C^\bullet_{\infty\cdot c}})_{\infty\cdot c}}$) be the closed sub-prestack of $\Gr^{\omega}_{G, C^\bullet}$ (resp.\ $\Gr^{\omega}_{G, C^\bullet_{\infty\cdot c}}$), such that, $(D, \cP_G, \alpha)\in \overline{\sO}_{C^\bullet}$ (resp.\ ${(\overline{\sO}_{C^\bullet_{\infty\cdot c}})_{\infty\cdot c}}$) if and only if the maps \eqref{6.19} extend to regular maps {of vector bundles} on $C\times S$ (resp.\ $(C-c)\times S$).

From the definitions, it is not hard to see that $\overline{\sO}_{C^\bullet}$ and ${(\overline{\sO}_{C^\bullet_{\infty\cdot c}})_{\infty\cdot c}}$ satisfy the factorization properties. Namely, we have isomorphisms
\begin{equation}
    \overline{\sO}_{C^\bullet}\underset{C^\bullet}{\times} (C^\bullet\times C^\bullet)_{disj} \simeq (\overline{\sO}_{C^\bullet}\times \overline{\sO}_{C^\bullet})\underset{C^\bullet\times C^\bullet}{\times} (C^\bullet\times C^\bullet)_{disj}. 
\end{equation}
\begin{equation}
    {(\overline{\sO}_{C^\bullet_{\infty\cdot c}})_{\infty\cdot c}}\underset{C_{\infty\cdot c}^\bullet}{\times} (C^\bullet\times C^\bullet_{\infty\cdot c})_{disj} \simeq (\overline{\sO}_{C^\bullet}\times {(\overline{\sO}_{C^\bullet_{\infty\cdot c}})_{\infty\cdot c}})\underset{C^\bullet\times C^\bullet_{\infty\cdot c}}{\times} (C^\bullet\times C^\bullet_{\infty\cdot c})_{disj}. 
\end{equation}

For $D=\sum_{i=1}^{l} \alpha_{x_i}\cdot x_i\in C^\bullet$, the reduced part of the fiber of $\overline{\sO}_{C^\bullet}$ over $D$ is isomorphic to the reduced part of
\[\prod_{1}^{l} \overline{\sO}_{x_i}^0\subset \prod_{1}^{l}\Gr^{\omega}_{G, x_i}.\]
Here $\overline{\sO}_{x_i}^0\subset \Gr^{\omega}_{G, x_i}$ denotes the closure of the unital $H^{\omega}(\bF)$-orbit.

For $D= (\xi, (\eta, \eta'))\cdot c+ \sum_{i=1}^{l} \alpha_{x_i}\cdot x_i\in C^\bullet_{\infty\cdot c}$, the reduced part of the fiber of ${(\overline{\sO}_{C^\bullet_{\infty\cdot c}})_{\infty\cdot c}}$ is isomorphic to the reduced part of  \[\Gr^{\omega}_{G, c}\times \prod_{i=1}^{l}\overline{\sO}_{x_i}^0\subset \Gr^{\omega}_{G,c}\times \prod_{i=1}^{l} \Gr^{\omega}_{G, x_i}.\]

\subsubsection{}
Let $\overline{\sS}_{C^\bullet}$ be the sub-prestack of $\Gr^{\omega}_{G, C^\bullet}$ such that $(D, \cP_G, \alpha)\in \overline{\sS}_{C^\bullet}$ if and only if the maps \eqref{6.18} extend to regular maps {of vector bundles} on $C\times S$. It has an open sub-prestack $\sS_{C^\bullet}$ which consists of those points such that the first map of \eqref{6.18} is surjective and the second map, third map are injective. We denote by $\overline{\sS}_{C^\bullet_{\infty\cdot c}}$ (resp.\ ${S}_{C^\bullet_{\infty\cdot c}}$) the similarly defined spaces. 

We have,
\begin{equation}
    \sS_{C^\bullet}\underset{C^\bullet}{\times} (C^\bullet\times C^\bullet)_{disj} \simeq (\sS_{C^\bullet}\times \sS_{C^\bullet})\underset{C^\bullet\times C^\bullet}{\times} (C^\bullet\times C^\bullet)_{disj},
\end{equation}
\begin{equation}
    \sS_{C^\bullet_{\infty\cdot c}}\underset{C_{\infty\cdot c}^\bullet}{\times} (C^\bullet\times C^\bullet_{\infty\cdot c})_{disj} \simeq (\sS_{C^\bullet}\times \sS_{C^\bullet_{\infty\cdot c}})\underset{C^\bullet\times C^\bullet_{\infty\cdot c}}{\times} (C^\bullet\times C^\bullet_{\infty\cdot c})_{disj}. 
\end{equation}
Similar factorization properties hold for $\overline{\sS}_{C^\bullet}$ and $\overline{\sS}_{C^\bullet_{\infty\cdot c}}$ as well.

For $D=\sum_{i=1}^{l} \alpha_{x_i}\cdot x_i\in C^\bullet$, the reduced part of the fiber of $\overline{\sS}_{C^\bullet}$ over $D$ is isomorphic to the reduced part of
\[\prod_{1}^{l} \overline{\sS}_{x_i}^{\alpha_{x_i}}\subset \prod_{1}^{l}\Gr^{\omega}_{G, x_i}.\]
Here $\overline{\sS}_{x_i}^{\alpha_{x_i}}\subset \Gr^{\omega}_{G, x_i}$ denotes the closure of the $U^{\omega}(\bF)$-orbit of $\sL^{\alpha_{x_i}}$.

For $D= (\xi, (\eta, \eta'))\cdot c+ \sum_{i=1}^{l} \alpha_{x_i}\cdot x_i\in C^\bullet_{\infty\cdot c}$, the reduced part of the fiber of $\overline{\sS}_{C^\bullet_{\infty\cdot c}}$ is isomorphic to the reduced part of  \[\overline{\sS}_{c}^{\xi, (\eta, \eta')}\times \prod_{i=1}^{l}\overline{\sS}_{x_i}^{\alpha_{x_i}}\subset \Gr^{\omega}_{G,c}\times \prod_{i=1}^{l} \Gr^{\omega}_{G, x_i}.\]

Similar fiber descriptions hold for ${S}_{C^\bullet}$ and ${S}_{C^\bullet_{\infty\cdot c}}$.

By definition, we have
\begin{equation}\label{zastava local-global}
    \begin{split}
        \cY\simeq& \overline{\sO}_{C^\bullet}\cap \sS_{C^\bullet},\\
        \overline{\cY}\simeq& \overline{\sO}_{C^\bullet}\cap \overline{\sS}_{C^\bullet},\\
        \cY_{\infty\cdot c}\simeq& {(\overline{\sO}_{C^\bullet_{\infty\cdot c}})_{\infty\cdot c}}\cap \sS_{C^\bullet_{\infty\cdot c}},\\
        \overline{\cY}_{\infty\cdot c}\simeq& {(\overline{\sO}_{C^\bullet_{\infty\cdot c}})_{\infty\cdot c}}\cap \overline{\sS}_{C^\bullet_{\infty\cdot c}},
    \end{split}
\end{equation}
and the isomorphisms are compatible with respect to the factorization structures.

\subsubsection{Factorization SW-Zastava spaces}
In the definitions of ${(\overline{\sO}_{C^\bullet_{\infty\cdot c}})_{\infty\cdot c}}$ and $\cY_{\infty\cdot c}$, we can also allow $c$ to vary in $\Ran_C$. Namely, let $C^\bullet_{\infty\cdot \Ran}$ be the scheme which classifies $(D, \bar{c}=\{c_1, c_2,\cdots, c_m\})$, where $D$ is a $\Lambda$-colored divisor on $C$, such that only the coefficients of the marked points $c_1, \cdots, c_m$ are allowed to be non-negative. The prestack ${(\overline{\sO}_{C^\bullet_{\infty\cdot \Ran}})_{\infty\cdot \Ran}}$ classifies $(\bar{c}, D, \cP_G, \alpha)$, where $(D,\bar{c})\in C^\bullet_{\infty\cdot \Ran}$, $\cP_G$ is a $G$-bundle on $C$, and $\alpha$ is an identification of $\cP_G$ and $\cP_G^\omega$ outside $\supp(D)$ such that \eqref{5.1} which is \textit{a priori} defined on $C-\supp(D)$ extends to a regular map on $C-\bar{c}$. Also, we define $\sS_{C^\bullet_{\infty\cdot \Ran}}$ as the prestack which classifies $(D,\bar{c}, \cP_G,\alpha)$, such that the maps \eqref{6.19} extend to regular maps on the whole curve.

We let ${\cY}_{C^\bullet_{\infty\cdot \Ran}}$ be the scheme over $\Ran_C$ which is defined as ${(\overline{\sO}_{C^\bullet_{\infty\cdot \Ran}})_{\infty\cdot \Ran}}\cap \sS_{C^\bullet_{\infty\cdot \Ran}}$. By the factorization property of ${(\overline{\sO}_{C^\bullet_{\infty\cdot \Ran}})_{\infty\cdot \Ran}}$ with respect to $\overline{\sO}_{C^\bullet}$ and the factorization property of $\sS_{C^\bullet_{\infty\cdot \Ran}}$ with respect to $\sS_{C^\bullet}$, we obtain that ${\cY}_{C^\bullet_{\infty\cdot \Ran}}$ is factorizable with respect to $\cY$.

\begin{rem}
    More or less by definition, we have
\begin{equation}
    {\cY}_{C^\bullet_{\infty\cdot \Ran}}\simeq \cM_{\infty\cdot\Ran}\underset{\Bun_G}{\times'} \Bun_B.
\end{equation}
\end{rem}


\section{Construction of the functor}\label{constrct functor}
In this section, we will construct a functor from the twisted Gaiotto category to the category of twisted factorization modules. 

\subsection{Global functor}
First, let us define the global functor. Consider the diagram,
\begin{center}
\[    \xymatrix {
&\cY_{\infty\cdot c}\ar[ld]_{p}\ar[rd]^{v}&\\
\cM_{\infty\cdot c}&& C^\bullet_{\infty\cdot c}.
}\]
\end{center}

Following \cite{[BFT0]}, \cite{[G]}, and \cite{[GL]}, we want to use the pullback-pushforward functor of the above correspondence diagram. But since we are in the twisted setting, in order to construct such a functor for the category of {twisted} D-modules, we need to check the compatibility of twistings.



\begin{lem}
There is a natural isomorphism of twistings
\[v^*(\cP_{C^\bullet})\simeq p^*(\cP_{\det}).\]
\end{lem}
\begin{proof}
The line bundle $p^*(\cP_{\det})$ is given by the pullback of the (renormalized relative) determinant line bundle on $\Bun_G$. Its restriction to $\Bun_B$ is isomorphic to the pullback of the (renormalized relative) determinant line bundle $\cP_{\det, T}$ on $\Bun_T$. Now, the lemma follows from the facts that $\cP_{C^\bullet}$ is isomorphic to the pullback of $\cP_{\det, T}$ along
\begin{equation*}
    \begin{split}
      \AJ:  C^\bullet_{\infty\cdot c}&\longrightarrow \Bun_T\\
        D&\mapsto (\omega_C^{\rho}\overset{T_{N-M-1}}{\times} T)(-D),
    \end{split}
\end{equation*}
and 
the composed map $\cY_{\infty\cdot c}\rightarrow \Bun_B\rightarrow \Bun_T$ equals $\AJ \circ v$.
\end{proof}

{For any $(\xi, (\eta, \eta'))$, we call the parity of $\langle2\rho^\circ, (\xi, (\eta, \eta'))\rangle$ the parity of $(\xi, (\eta, \eta'))$.}

Let us denote by $p^\bullet$ the shifted pullback functor
\begin{equation}
    {\SD}_{\kappa}(\cM_{\infty\cdot c})\longrightarrow {\SD}_{\kappa}(\cY_{\infty\cdot c}),
\end{equation}
by assigning $p^\bullet(\cF):= p^!(\cF)[-dim.rel (\Bun_B, \Bun_G)]$ and change the parity according to the parity of $(\xi, (\eta, \eta'))$. Here, the number $dim.rel(\Bun_B, \Bun_G)$ denotes the relative dimension of the connected component of $\Bun_B$ and $\Bun_G$. I.e., 
\[dim.rel(\Bun_B, \Bun_G):= \dim \Bun_B^{(\xi, (\eta, \eta'))}-\dim \Bun_G\]
on the connected component $\cY_{\infty\cdot c}^{(\xi, (\eta, \eta'))}:= \cM_{\infty\cdot c}\underset{\Bun_G}{\times'} \Bun_B^{(\xi, (\eta, \eta'))}$. Here, $\Bun_B^{(\xi, (\eta, \eta'))}$ is the connected component of $\Bun_B$ of degree $\deg(\omega^{\rho}_C)-(\xi, (\eta, \eta'))= 2(g-1)\rho_{M,N}-(\xi, (\eta, \eta'))$.

\begin{defn}\label{def glob func}
We define the functor $F^{glob}: \cC^{glob}_{\kappa}(M|N)\longrightarrow {\SD}_{\kappa}(C^\bullet_{\infty\cdot c})$ as
\[F^{glob}(\cF):= v_* p^\bullet (\cF).\]
\end{defn}

The following lemma can be proved similarly to \cite[Theorem 4.11]{[G]}.

\begin{lem}\label{clean}
For $\kappa$ generic and any $\cF\in \cC_{\kappa}^{glob}(M|N)$, the object
\[v_!p^\bullet(\cF)\in {\SD}_{\kappa}(C^\bullet_{\infty\cdot c})\] is well-defined, and the natural morphism
\begin{equation}
    v_!p^\bullet(\cF)\longrightarrow v_*p^\bullet(\cF)
\end{equation}
is an isomorphism.
\end{lem}

{\begin{rem}
    The above lemma is \textit{not} true when $\kappa$ is rational. 
\end{rem}}

\subsection{Spreading functor}
{To define a local functor, we need to first define a functor from the local twisted Gaiotto category to the category of twisted D-modules on $\Gr^\om_{G, \Ran_{C,c}}$.}

Recall the definition of ${(\overline{\sO}_{\Ran_{C, c}})_{\infty\cdot c}}$ and $\overline{\sO}_{\Ran_C}$ in Section \ref{Section 5.1}. 

Let $H'_{\Ran_{C, c}}$ be a subgroup of $ H^\omega(\bF)_{\Ran_{C, c}}$ whose fiber over $\bar{x}=\{c, x_1,\cdots, x_m\}$ is isomorphic to $H^\omega(\bF)_c \times \prod^m_i H^\omega(\bO)_{x_i}$, and consider the unit map
\begin{equation}
\begin{split}
     \unit: \Ran_{C, c}\times \Gr_{G, c}^{\omega}\longrightarrow {(\overline{\sO}_{\Ran_{C, c}})_{\infty\cdot c}}\\
     \bar{x}, (\cP_G, \alpha)\mapsto (\bar{x}, \cP_G, \alpha|_{C-\bar{x}}).
\end{split}
\end{equation}
The unit map is an $H'_{\Ran_{C, c}}$-invariant closed embedding. In particular, $\unit_!$ is well-defined and induces a functor
\begin{equation}
    \unit_!: {\SD}_{\kappa}^{H^\omega(\bF),\chi}(\Ran_{C, c}\times \Gr_{G, c}^\omega)\longrightarrow {\SD}_{\kappa}^{H'_{\Ran_{C, c}},\chi_{\Ran_{C, c}}}({(\overline{\sO}_{\Ran_{C, c}})_{\infty\cdot c}}).
\end{equation}

We denote by $\sprd_{\Ran_{C, c}}$ the following composed functor
\begin{equation}\label{eq 7.5}
\begin{split}
      \sprd_{\Ran_{C, c}}: {\SD}_{\kappa}^{H^\omega(\bF),\chi}(\Gr_{G, c}^{\omega})\overset{\pr^!}{\longrightarrow} {\SD}_{\kappa}^{H^\omega(\bF),\chi}(\Ran_{C, c}\times \Gr_{G, c}^\omega)\overset{\unit_!}{\longrightarrow} {\SD}_{\kappa}^{H'_{\Ran_{C, c}},\chi_{\Ran_{C, c}}}({(\overline{\sO}_{\Ran_{C, c}})_{\infty\cdot c}})\\ \overset{\Av_!^{H^\omega(\bF)_{\Ran_{C, c}}/ H'_{\Ran_{C, c}}, \chi_{\Ran_{C, c}}}}{\longrightarrow} {\SD}_{\kappa}^{H^\omega(\bF)_{\Ran_{C, c}},\chi_{\Ran_{C, c}}}({(\overline{\sO}_{\Ran_{C, c}})_{\infty\cdot c}}).
\end{split}
\end{equation}
Here, $\Av_!^{H^\omega(\bF)_{\Ran_{C, c}}/ H'_{\Ran_{C, c}}, \chi_{\Ran_{C, c}}}$ is the $!$-averaging with respect to $(H^\omega(\bF)_{\Ran_{C, c}},\chi_{\Ran_{C, c}})$ relative to  $H'_{\Ran_{C, c}}$.

{The relative $!$-averaging functor is well-defined by a similar argument as in \cite[Section 2.11]{[Ras]}. Namely, we need to show that $\on{act}_!(\chi_{\Ran_{C,c}}\tboxtimes \cF)$ is well-defined for $\cF\in {\SD}_{\kappa}^{H'_{\Ran_{C, c}},\chi_{\Ran_{C, c}}}({(\overline{\sO}_{\Ran_{C, c}})_{\infty\cdot c}})$, where $\chi_{\Ran_{C,c}}\tboxtimes \cF$ is the twisted D-module on $H^\omega(\bF)_{\Ran_{C, c}}\overset{H'_{\Ran_{C, c}}}{\times} {(\overline{\sO}_{\Ran_{C, c}})_{\infty\cdot c}}$ descended from $H^\omega(\bF)_{\Ran_{C, c}}{\times} {(\overline{\sO}_{\Ran_{C, c}})_{\infty\cdot c}}$, and $\on{act}$ is given by the action of $H^\omega(\bF)_{\Ran_{C, c}}$ on ${(\overline{\sO}_{\Ran_{C, c}})_{\infty\cdot c}}$. Then, similar to the definition of $H'_{\Ran_{C,c}}$, we consider a subgroup $G'_{\Ran_{C,c}}$ of $G^{\om}(\bF)_{\Ran_{C,c}}$. We note that there is an isomorphism
\begin{equation}
    H^\omega(\bF)_{\Ran_{C, c}}\overset{H'_{\Ran_{C, c}}}{\times} \Gr^\om_{G,\Ran_{C, c}}\simeq H^\omega(\bF)_{\Ran_{C, c}}G'_{\Ran_{C,c}}\overset{G'_{\Ran_{C,c}}}{\times} \Gr^\om_{G,\Ran_{C, c}}.
\end{equation}
Consider the preimage of $\overline{\sO}_{\Ran_{C}}\underset{\Ran_C}{\times} \Ran_{C,c}$ along the projection $G^{\om}(\bF)_{\Ran_{C,c}}\rightarrow \Gr^\om_{G, \Ran_{C,c}}$, and we denote it by $\overline{H^\omega(\bF)_{\Ran_{C, c}}G'_{\Ran_{C,c}}}$. The action map $\overline{\on{act}}\colon \overline{H^\omega(\bF)_{\Ran_{C, c}}G'_{\Ran_{C,c}}}\overset{G'_{\Ran_{C,c}}}{\times} \Gr^\om_{G,\Ran_{C, c}}\rightarrow \Gr^\om_{G,\Ran_{C, c}}$ is ind-proper, so $\overline{\on{act}}_!$ is well-defined and equals $\overline{\on{act}}_*$.

Due to the fact that the unital $H(\bF)$-orbit in $\Gr_G$ is minimal relevant, one can show that there is no Gaiotto equivariant twisted D-module on the complement $\overline{H^\omega(\bF)_{\Ran_{C, c}}G'_{\Ran_{C,c}}}\overset{G'_{\Ran_{C,c}}}{\times} \Gr^\om_{G,\Ran_{C, c}}- {H^\omega(\bF)_{\Ran_{C, c}}G'_{\Ran_{C,c}}}\overset{G'_{\Ran_{C,c}}}{\times} \Gr^\om_{G,\Ran_{C, c}}$. In particular, $\on{act}_!$ can be defined as the composition of $\overline{\on{act}}_!$ and a clean extension.

}

\subsubsection{Unital structure}
The image of the above functor has a canonical unital structure. We will explain it in this section.

Recall that we let $(\Ran_{C, c}\times \Ran_{C, c})^\subset$ be a sub-prestack of $\Ran_{C, c}\times \Ran_{C, c}$ which consists of $(\bar{x},\bar{x}')$ such that $\bar{x}\subset \bar{x}'$. It has two maps to $\Ran_{C, c}$: $\phi_{small}$ and $\phi_{big}$ which send $(\bar{x}, \bar{x}')$ to $\bar{x}$ and $\bar{x}'$, respectively. 

Let $\overline{\sO}_{(\Ran_{C, c}\times \Ran_{C, c})^\subset}$ be ${(\overline{\sO}_{\Ran_{C, c}})_{\infty\cdot c}}\times_{\Ran_{C, c}, \phi_{small}} (\Ran_{C, c}\times \Ran_{C, c})^\subset$, and consider the following diagram
\[\xymatrix{&\overline{\sO}_{(\Ran_{C, c}\times \Ran_{C, c})^\subset}\ar[ld]_{\phi_{\sO, small}}\ar[rd]^{\phi_{\sO, big}}&\\
{(\overline{\sO}_{\Ran_{C, c}})_{\infty\cdot c}}&&{(\overline{\sO}_{\Ran_{C, c}})_{\infty\cdot c}}},\]
where $\phi_{\sO, small}$ sends $(\bar{x}, \bar{x}', \cP_G, \alpha)$ to $(\bar{x}, \cP_G, \alpha)$ and $\phi_{\sO, big}$ sends $(\bar{x}, \bar{x}', \cP_G, \alpha)$ to $(\bar{x}', \cP_G, \alpha|_{C-\bar{x}'})$.

To introduce the definition of the unital structure in our setting, we need to first define a group prestack over $(\Ran_{C, c}\times \Ran_{C, c})^\subset$. Namely, we define $H'_{(\Ran_{C, c}\times \Ran_{C, c})^\subset}$ as the group prestack over $(\Ran_{C, c}\times \Ran_{C, c})^\subset$, whose fiber over $(\bar{x}, \bar{x}')$ is $H^\omega(\bF)_{\bar{x}}\times H^\omega(\bO)_{\bar{x}'- \bar{x}}$. To be more precise, it classifies $(\bar{x}, \bar{x}', \gamma)$, where $(\bar{x}, \bar{x}')\in (\Ran_{C, c}\times \Ran_{C, c})^\subset$ and $\gamma$ is a map from $\cD_{\bar{x}'} - \bar{x}$ to $H^\omega$. The group prestack $H'_{(\Ran_{C, c}\times \Ran_{C, c})^\subset}$ acts on $\overline{\sO}_{(\Ran_{C, c}\times \Ran_{C, c})^\subset}$ via 
\[H'_{(\Ran_{C, c}\times \Ran_{C, c})^\subset}\longrightarrow H^\omega(\bF)_{\Ran_{C, c}}\times_{\Ran_{C, c}, \phi_{small}} (\Ran_{C, c}\times \Ran_{C, c})^\subset.\]


\begin{defn}
    For an object $\cF\in {\SD}_{\kappa}^{H^\omega(\bF)_{\Ran_{C, c}},\chi_{\Ran_{C, c}}}({(\overline{\sO}_{\Ran_{C, c}})_{\infty\cdot c}})$, {a} unital structure is an isomorphism 
    \[\phi_{\sO, small}^!(\cF)\simeq \phi_{\sO, big}^!(\cF)\]
    in the category of ${\SD}_{\kappa}^{H'_{(\Ran_{C, c}\times \Ran_{C, c})^\subset},\chi_{\Ran_{C, c}}}(\overline{\sO}_{(\Ran_{C, c}\times \Ran_{C, c})^\subset})$. We denote the resulting category by ${\SD}_{\kappa}^{H^\omega(\bF)_{\Ran_{C, c}},\chi_{\Ran_{C, c}}}({(\overline{\sO}_{\Ran_{C, c}})_{\infty\cdot c}})_\unit$.
\end{defn}

The importance of the unital structure is that it ensures that any object $\cF\in {\SD}_{\kappa}^{H^\omega(\bF)_{\Ran_{C, c}},\chi_{\Ran_{C, c}}}({(\overline{\sO}_{\Ran_{C, c}})_{\infty\cdot c}})_\unit$ can be uniquely determined by its restriction $\cF_c$ {to} the fiber of ${(\overline{\sO}_{\Ran_{C, c}})_{\infty\cdot c}}$ over $\{c\}$.

Indeed, if $\bar{x}=\{c, x_1,\cdots, x_m\}\in \Ran_{C, c}$, the fiber $\overline{\sO}_{\bar{x}}$ of ${(\overline{\sO}_{\Ran_{C, c}})_{\infty\cdot c}}$ over $\bar{x}$ is isomorphic to $\Gr_{G,c}^\omega\times \overline{\sO}^0_{x_1}\times \cdots \times \overline{\sO}^0_{x_m}$. The unital structure of $\cF$ says that, the $!$-restriction of $\cF_{\bar{x}}\in {\SD}_{\kappa}^{H^(\bF)_{\bar{x}}, \chi_{\bar{x}}}(\overline{\sO}_{\bar{x}})$ to the unit section $\Gr_{G, c}^{\omega}\times \{\sL_0\}\times\cdots \times \{\sL_0\} \subset \overline{\sO}_{\bar{x}}$, is isomorphic to the $(H'_{\{c\}\subset \bar{x}}, \chi_{\bar{x}})$-equivariant D-module given by $\cF_c\in {\SD}_{\kappa}^{H^\omega(\bF)_c,\chi}(\Gr_{G,c}^\omega)$ along $\overline{\sO}_c\simeq \Gr_{G,c}^\omega\overset{\sim}{\longrightarrow} \Gr_{G, c}^{\omega}\times \{\sL_0\}\times\cdots \times \{\sL_0\} $.  This identification of the $!$-restriction of $\cF_{\bar{x}}$ with $\cF_{c}$ is in the category ${\SD}_{\kappa}^{H'_{\{c\}\subset \bar{x}}, \chi_{\bar{x}}}(\Gr_{G, c}^{\omega}\times \{\sL_0\}\times\cdots \times \{\sL_0\} )$. Note that the latter category is equivalent to ${\SD}_{\kappa}^{H^\omega(\bF)_{\bar{x}},\chi_{\bar{x}}}(\overline{\sO}_{\bar{x}})$, so the unital structure determines $\cF_{\bar{x}}$ from $\cF_c$ for any $\bar{x}\in \Ran_{C, c}$.

{It follows from the definition \eqref{eq 7.5} that} for any $\cF\in {\SD}_{\kappa}^{H^\omega(\bF),\chi}(\Gr_G^\omega)$, the image $\sprd_{\Ran_{C, c}}(\cF)$ acquires a canonically defined unital structure. Furthermore, taking $!$-restriction along $\Gr_{G,c}^\omega= \{c\}\times \Gr_{G,c}^\omega\hookrightarrow {(\overline{\sO}_{\Ran_{C, c}})_{\infty\cdot c}}$ defines the inverse functor. In other words, 
\begin{lem}
The spreading functor defines an equivalence
\begin{equation}\label{eq 7.6}
 \sprd_{\Ran_{C,c}}:   {\SD}_{\kappa}^{H^\omega(\bF),\chi}(\Gr_G^\omega) \overset{\sim}{\longrightarrow} {\SD}_{\kappa}^{H^\omega(\bF)_{\Ran_{C, c}},\chi_{\Ran_{C, c}}}({(\overline{\sO}_{\Ran_{C, c}})_{\infty\cdot c}})_\unit.
\end{equation}
\end{lem}


\begin{rem}
{To be self-complete, we mention that there is another way to assign a twisted D-module on $\Gr^\om_{G,c}$ to a twisted Gaiotto D-module on $\Gr^{\om}_{G, \Ran_{C,c}}$, although we will not use it. 
}

Let $H''_{\Ran_{C, c}}$ be the subgroup of $H^\omega(\bF)_{\Ran_{C, c}}$ whose fiber over $\bar{x}=\{c, x_1, \cdots, x_m\}\in \Ran_{C, c}$ is $H^\omega(\bF)_c\times U_{M,N}^{-,\omega}(\bF)_{x_1}\times \cdots \times U_{M,N}^{-,\omega}(\bF)_{x_m}$. Using the same proof as in \cite[Theorem 6.4.2]{[G3]} and the fact that the unital $U_{M,N}^{-,\omega}(\bF)$-orbit in $\Gr^{{\omega}}_G$ is a minimal relevant orbit, we can show with the same method as {in} \cite[Theorem 6.2.5]{[G3]} that there is an equivalence
\begin{equation}\label{eq 7.8}
    {\SD}_{\kappa}^{H^\omega(\bF)_c, \chi}(\Ran_{C, c}\times \Gr_{G,c}^\omega)\simeq {\SD}_{\kappa}^{H''_{\Ran_{C, c}},\chi_{\Ran_{C, c}}}({(\overline{\sO}_{\Ran_{C, c}})_{\infty\cdot c}}).
\end{equation}

{In particular, we obtain a functor from $\cC^{loc}_{\kappa}(M|N)$ to ${\SD}_{\kappa}^{H''_{\Ran_{C, c}},\chi_{\Ran_{C, c}}}({(\overline{\sO}_{\Ran_{C, c}})_{\infty\cdot c}})$ by composing with the $!$-pullback functor ${\SD}_{\kappa}^{H^\omega(\bF)_c, \chi}(\Gr_{G,c}^\omega)\rightarrow {\SD}_{\kappa}^{H^\omega(\bF)_c, \chi}(\Ran_{C, c}\times \Gr_{G,c}^\omega)$. Furthermore, its further composition with the forgetful functor to ${\SD}_{\kappa}({(\overline{\sO}_{\Ran_{C, c}})_{\infty\cdot c}})$ is isomorphic to the composition functor of $\sprd_{\Ran_{C,c}}$ and the forgetful functor.}
\end{rem}
{\begin{rem}
   The equivalence \eqref{eq 7.8} gives rise to an equivalence between $\cC_{\kappa}^{loc}(M|N)$ and ${\SD}_{\kappa}^{H''_{\Ran_{C, c}},\chi_{\Ran_{C, c}}}({(\overline{\sO}_{\Ran_{C, c}})_{\infty\cdot c}})_{\unit}$, i.e., the category of $(H''_{\Ran_{C, c}},\chi_{\Ran_{C, c}})$-equivariant {twisted} D-modules with {a} unital structure. To be more precise, let us consider the constant group $H^\omega(\bF)\times \Ran_{C,c}\subset H^\omega(\bF)_{\Ran_{C,c}}$ over $\Ran_{C,c}$. The unital structure for $\cF\in {\SD}_{\kappa}^{H''_{\Ran_{C, c}},\chi_{\Ran_{C, c}}}({(\overline{\sO}_{\Ran_{C, c}})_{\infty\cdot c}})$ means that there is an isomorphism
\[\phi_{\sO, small}^!(\cF)\simeq \phi_{\sO, big}^!(\cF)\]
    in the {category} ${\SD}_{\kappa}^{H^\omega(\bF)\times \Ran_{C,c},\chi}(\overline{\sO}_{(\Ran_{C, c}\times \Ran_{C, c})^\subset})$.
\end{rem}}

\subsubsection{{Factorization modules}}
The dualizing D-module $\omega_{\Ran}$ on the unit section $\Ran\subset \overline{\sO}_{\Ran_C}$ acquires a naturally defined $H^\omega(\bO)_{\Ran}$-equivariant structure. Similar{ly} to \cite{[GL]}, let us denote by $\Vac$ the unit object in ${\SD}_{\kappa}^{H^\omega(\bF)_{\Ran_C},\chi_{\Ran_C}}(\overline{\sO}_{\Ran_C})$. It is given by 
\begin{equation}\label{eq 7.3}
\Vac=\Av_!^{H^\omega(\bF)_{\Ran}/H^\omega(\bO)_{\Ran},\chi_\Ran}(\omega_{\Ran}).
\end{equation}
Here, $\Av_!^{H^\omega(\bF)_{\Ran}/H^\omega(\bO)_{\Ran},\chi_{\Ran}}$ is the $!$-averaging with respect to $(H^\omega(\bF)_{\Ran},\chi_{\Ran})$ relative to  $H^\omega(\bO)_{\Ran}$.

It is a unital factorization algebra and its restriction to the fiber of $\overline{\sO}_{\Ran_C}$ over $\bar{x}\in \Ran_C$ is just the (clean) extension of the rank $1$ local system on $\sO_{\bar{x}}^0$ corresponding to $\chi$. The functor $\sprd_{\Ran_{C, c}}$ in \eqref{eq 7.5} naturally acquires {a unital} factorization D-module structure with respect to $\Vac$. In fact, the spreading functor induces an equivalence between ${\SD}_{\kappa}^{H^\omega(\bF),\chi}(\Gr_G^\omega)$ and the category of unital $\Vac$-factorization modules (cf.~\cite[Theorem 6.13.2]{[Ras0]}).

\subsection{Local functor}
Consider the following diagram,
\begin{center}
\[    \xymatrix {
&\Ran_{C, c}\times \Gr^{\omega}_G\ar[ld]_{\pr}\ar[rd]^{\unit}&&{(\overline{\sO}_{C^\bullet_{\infty\cdot c}})_{\infty\cdot c}}\ar@{--}[ld]&\overline{\cY}_{\infty\cdot c}\ar[l]\ar[rd]^{v}&\\
\Gr^{\omega}_{G,c}&& {(\overline{\sO}_{\Ran_{C, c}})_{\infty\cdot c}}&&&C^\bullet_{\infty\cdot c}.
}\]
\end{center}

For $\cF\in {\SD}_{\kappa}^{H^{\omega}(\bF),\chi}(\Gr^\omega_G)$, we apply the spreading functor in \eqref{eq 7.6}, ${\SD}_{\kappa}^{H^\omega(\bF),\chi}(\Gr_G^\omega)\rightarrow {\SD}_{\kappa}^{H^\omega(\bF)_{\Ran_{C, c}},\chi_{\Ran_{C, c}}}({(\overline{\sO}_{\Ran_{C, c}})_{\infty\cdot c}})_\unit$. We obtain a factorization module $\sprd_{\Ran_{C, c}}(\cF)\in {\SD}_{\kappa}({(\overline{\sO}_{\Ran_{C, c}})_{\infty\cdot c}})$ over $\Vac$, which is characterized by the property
\begin{equation}
    \unit^!(\sprd_{\Ran_{C, c}}(\cF))\simeq \pr^!(\cF).
\end{equation}



\subsubsection{}\label{section 7.3}
Then, in order to obtain a D-module on $\overline{\cY}_{\infty\cdot c}$ from $\sprd_{\Ran_{C, c}}(\cF)$, we need to 'pullback' $\sprd_{\Ran_{C, c}}(\cF)$ to the configuration affine Grassmannian ${(\overline{\sO}_{C^\bullet_{\infty\cdot c}})_{\infty\cdot c}}$. It follows from the constructions in \cite[Section 18.3]{[GL]}. For the purpose of being self-complete, we briefly recall the construction in the $loc.cit$.

Consider the sub-prestack $(\Gr^{\omega}_{T,\Ran_{C, c}})_{\infty\cdot c}^{\textnormal{neg}}$ of the Beilinson-Drinfeld affine Grassmannian $\Gr^\omega_{T, \Ran_{C, c}}$ defined in \cite[Section 4.6.6]{[GL]}. Roughly speaking, it consists of those points in the Beilinson-Drinfeld affine Grassmannian with regularity and non-{redundancy} properties outside the marked point $c$ (ref. \cite[Section 4.6.2]{[GL]}). It has a map to $C^\bullet_{\infty\cdot c}$ and there is an isomorphism of prestacks
\[(\Gr^{\omega}_{T,\Ran_{C, c}})_{\infty\cdot c}^{\textnormal{neg}}\underset{\Ran_{C, c}}{\times}{{(\overline{\sO}_{\Ran_{C, c}})_{\infty\cdot c}}}\simeq (\Gr^{\omega}_{T,\Ran_{C, c}})_{\infty\cdot c}^{\textnormal{neg}}\underset{C^\bullet_{\infty\cdot c}}{\times}{{(\overline{\sO}_{C^\bullet_{\infty\cdot c}})_{\infty\cdot c}}}.\]
Furthermore, there is  an equivalence of categories
\[{\SD}_{\kappa}((\Gr^{\omega}_{T,\Ran_{C, c}})_{\infty\cdot c}^{\textnormal{neg}}\underset{C^\bullet_{\infty\cdot c}}{\times}{{(\overline{\sO}_{C^\bullet_{\infty\cdot c}})_{\infty\cdot c}}})\simeq {\SD}_{\kappa}({{(\overline{\sO}_{C^\bullet_{\infty\cdot c}})_{\infty\cdot c}}}).\]
In particular, taking $!$-pullback along  \[(\Gr^{\omega}_{T,\Ran_{C, c}})_{\infty\cdot c}^{\textnormal{neg}}\underset{\Ran_{C, c}}{\times}{{(\overline{\sO}_{\Ran_{C, c}})_{\infty\cdot c}}}\longrightarrow {{(\overline{\sO}_{\Ran_{C, c}})_{\infty\cdot c}}}\]
gives rise to a twisted D-module on ${{(\overline{\sO}_{C^\bullet_{\infty\cdot c}})_{\infty\cdot c}}}$. We denote the image of $\sprd_{\Ran_{C, c}}(\cF)\in {\SD}_{\kappa}({{(\overline{\sO}_{\Ran_{C, c}})_{\infty\cdot c}}})$ by $\sprd(\cF)\in {\SD}_{\kappa}({{(\overline{\sO}_{C^\bullet_{\infty\cdot c}})_{\infty\cdot c}}})$. 

\subsubsection{}
Recall \eqref{zastava local-global}, we have $ \overline{\cY}_{\infty\cdot c}\simeq {(\overline{\sO}_{C^\bullet_{\infty\cdot c}})_{\infty\cdot c}}\cap \overline{\sS}_{C^\bullet_{\infty\cdot c}}$. We consider the shifted dualizing D-module $\omega^{{\bullet}}_{\sS_{C^\bullet_{\infty\cdot c}}}[\deg]$ on $\sS_{C^\bullet_{\infty\cdot c}}$. Here, the cohomological shift $\deg=\langle 2\rho^\circ, (\xi,(\eta,\eta')) \rangle$ {on the corresponding connected component $\overline{\cY}_{\infty\cdot c}^{(\xi,(\eta,\eta'))}$, the supscript $\bullet$ means that its parity is given by the parity of $(\xi,(\eta,\eta'))$. In particular, the parity of the sheaf and the parity of the cohomological degree are the same}. 

\begin{rem}
    The cohomological degree $\deg=-dim.rel(\Bun^{(\xi,(\eta,\eta'))}_B, \Bun_G)-\dim \CM_{=0\cdot c}^{\xi}$, where $\CM_{=0\cdot c}^{\xi}$ denotes the connected component of $\CM_{=0\cdot c}=\Bun_H^{\om^\rho}$, such that the determinant line bundle of the associated $\GL_M$-bundle $\cP_M$ has the same degree as the degree of the determinant line bundle $\GL_M$-bundle of any point in $\Bun^{(\xi,(\eta,\eta'))}_B$.
\end{rem}

It is known (ref. \cite[Corollary 1.5.6]{[G3]}) that the $!$-extension of $\omega^{{\bullet}}_{\sS_{C^\bullet_{\infty\cdot c}}}[\deg]$ to $\overline{\sS}_{C^\bullet_{\infty\cdot c}}$ is well-defined. In the case of $\kappa$ is generic, it is isomorphic to the $*$-extension of $\omega^{{\bullet}}_{\sS_{C^\bullet_{\infty\cdot c}}}[\deg]$. In other words, the extension is clean.

Finally, we $*$ (equivalently, $!$)-pushforward the $!$-tensor product of $\sprd(\cF)$ with the clean extension of $\omega^{{\bullet}}_{\sS_{C^\bullet_{\infty\cdot c}}}[\deg]$ to $C^\bullet_{\infty\cdot c}$ along with $v$. The resulting {twisted} D-module on $C^\bullet_{\infty\cdot c}$ is denoted by $F^{loc}(\cF)$.

To summarize, the local functor $F^{loc}$ is defined as
\begin{equation}
\begin{split}
     F^{loc}: \cC^{loc}_{\kappa}(M|N)\longrightarrow {\SD}_{\kappa}(C^\bullet_{\infty\cdot c}),\\
     \cF\mapsto v_*(\sprd(\cF)\overset{!}{\otimes} \omega^{{\bullet}}_{\sS_{C^\bullet_{\infty\cdot c}}}[\deg]).
\end{split}
\end{equation}

\subsubsection{}
We claim the following comparison of local functor and global functor.
\begin{prop}\label{locglob fun}
The local functor $F^{loc}$ and the global functor $F^{glob}$ are compatible with respect to $\pi^![d_g]$, namely, there is $F^{loc}\circ \pi^![d_g]\simeq F^{glob}\colon \cC_{\kappa}^{glob}(M|N)\rightarrow {\SD}_{\kappa}(C^\bullet_{\infty\cdot c})$.
\end{prop}
\begin{proof}

Note that $\overline{\cY}_{\infty\cdot c}$ is isomorphic to ${(\overline{\sO}_{C^\bullet_{\infty\cdot c}})_{\infty\cdot c}}\cap \overline{\sS}_{C^\bullet_{\infty\cdot c}}$. By an analog of \cite[Corollary 13.4.8]{[GL]}, the $!$-pullback of the clean extension of the twisted constant D-module on $\Bun_B$ along with
\[\overline{\sS}_{C^\bullet_{\infty\cdot c}}\longrightarrow \overline{\Bun}_B\]
is isomorphic to the clean extension of the twisted shifted dualizing D-module on $\sS_{C^\bullet_{\infty\cdot c}}$. The $!$-pullback of $\cF\in \cC_{\kappa}^{glob}(M|N)$ along with $\cY_{\infty\cdot c}\rightarrow \cM_{\infty\cdot c}$ is isomorphic to the restriction of $\sprd\circ \pi^!$ on $\cY_{\infty\cdot c}$. In fact, the map $\cY_{\infty\cdot c}\rightarrow \cM_{\infty\cdot c}$ factors through ${(\overline{\sO}_{C^\bullet_{\infty\cdot c}})_{\infty\cdot c}}\rightarrow \cM_{\infty\cdot c}$, and the {twisted} D-module on ${(\overline{\sO}_{C^\bullet_{\infty\cdot c}})_{\infty\cdot c}}$ is constructed by pullback from ${(\overline{\sO}_{\Ran_{C, c}})_{\infty\cdot c}}$, so we only need to show 
\begin{equation}\label{loc-glob fun}
\sprd_{\Ran_{C, c}}\pi^!(\cF)\simeq \pi_{\Ran_{C,c}}^!(\cF).    
\end{equation}

Note that the $!$-pullback functor along $\pi_{\Ran_{C, c}}$ induces a {(in fact, equivalent)} functor
\[\cC^{glob}_{\kappa}(M|N)\longrightarrow {\SD}_{\kappa}^{H'(\bF)_{\Ran_{C, c}}, \chi}({(\overline{\sO}_{\Ran_{C, c}})_{\infty\cdot c}})_\unit,\]

Follow from the following commutative diagram,
\[\xymatrix{
&\Ran_{C, c}\times \Gr^{\omega}_{G,c}\ar[rd]^{\unit}&\\
\Gr^{\omega}_{G,c}\ar[rd]^\pi\ar[ru]^{\{c\}\times \id}&&{(\overline{\sO}_{\Ran_{C, c}})_{\infty\cdot c}}\ar[ld]^{\pi_{\Ran_{C, c}}}\\
&\cM_{\infty\cdot c},&}
\] for any $\cF\in \cC_{\kappa}^{glob}(M|N)$, we have
\[(\{c\}\times \id)^! \unit^! \pi_{\Ran_{C, c}}^!(\cF)\simeq \pi^!(\cF).\]

Now, by the equivalence between ${\SD}_{\kappa}^{H^{\omega}(\bF)_c,\chi_c}(\Gr^{\omega}_{G,c})$ and ${\SD}_{\kappa}^{H'(\bF)_{\Ran_{C, c}}, \chi}({(\overline{\sO}_{\Ran_{C, c}})_{\infty\cdot c}})_\unit$, there is
\[\sprd_{\Ran_{C,c}}(\{c\}\times \id)^! \unit^!\simeq \id\colon {\SD}_{\kappa}^{H'(\bF)_{\Ran_{C, c}}, \chi}({(\overline{\sO}_{\Ran_{C, c}})_{\infty\cdot c}})_\unit\rightarrow {\SD}_{\kappa}^{H'(\bF)_{\Ran_{C, c}}, \chi}({(\overline{\sO}_{\Ran_{C, c}})_{\infty\cdot c}})_\unit\].

{The desired isomorphism \eqref{loc-glob fun} follows. }

\end{proof}
\subsection{Factorization property of the functor}
The functor defined above sends an object in the twisted Gaiotto category to the category of twisted D-modules on the configuration space. In this section, {we will study the factorization module structure of the image of this functor.} 


Recall the factorization algebra $\Vac$ in \eqref{eq 7.3}. By Section \ref{section 7.3}, we can obtain a factorization algebra $\Vac_{C^\bullet}$ on $\overline{\sO}_{C^\bullet}$ by taking the pullback of $\Vac$. We define
\[\Omega:= v_*(\Vac_{C^\bullet}\otimes \omega_{S_{C^\bullet}})[\deg].\]

The following lemma follows from the construction.
\begin{lem}
The $\cP_{C^\bullet}^{\kappa}$-twisted D-module $\Omega$ is a factorization algebra on $C^\bullet$, and the image of $F^{loc}$ (equivalently, $F^{glob}$) has a naturally defined $\Omega$-factorization module structure, i.e.,
\begin{equation}
    \add^!(\Omega)|_{(C^\bullet\times C^\bullet)_{disj}}\simeq \Omega\boxtimes \Omega|_{(C^\bullet\times C^\bullet)_{disj}},
\end{equation}
and
\begin{equation}
    \add_c^!(F^{loc}(\cF))|_{(C^\bullet\times C_{\infty \cdot c}^\bullet)_{disj}}\simeq \Omega\boxtimes F^{loc}(\cF)|_{(C^\bullet\times C_{\infty \cdot c}^\bullet)_{disj}}.
\end{equation}
\end{lem}

That is to say, the functor $F^{loc}$ (equivalently, $F^{glob}$) factors through,
\begin{equation}\label{7.10}
    F^{loc}: \cC^{loc}_{\kappa}(M|N)\longrightarrow \Omega-\FM.
\end{equation}
Here, we denote by $\Omega-\FM$ the category of twisted factorization modules on $C^\bullet_{\infty\cdot c}$ with respect to $\Omega$. 


The main difficulty of proving the main theorem, i.e., Theorem \ref{statement of main}, of this paper is to identify $\Omega$ with the factorization algebra in Definition \ref{def of I}.
\begin{thm}\label{key}
There is a canonical isomorphism of factorization algebras on $C^\bullet$:
\begin{equation}
    \cI\simeq \Omega
\end{equation}
\end{thm}

The following sections (Section \ref{Section 8}, \ref{section 9}) are aimed to prove Theorem \ref{key}. We first need to study $H^{\omega}(\bF)$-orbits, $U^{\omega}(\bF)$-orbits of $\Gr^{\omega}_G$ and their intersections.
\section{Proof I: Orbits calculation}\label{Section 8}
Note that $\Gr_G^{\omega}$ and $\Gr_G$ are isomorphic, i.e., $\omega$-renormalization do not change the geometry. So, we only consider the usual (non-renormalized) groups and affine Grassmannian in this section.

\subsection{Intersection of $H(\bF)$-orbits and $U(\bF)$-orbits}
Recall that in Section \ref{Gaiotto 
category}, we studied $H(\bF)$-orbits of $\Gr_G$. We claim,
\begin{prop}\label{OS}
 If $\sS^{(\xi, (\eta, \eta'))}\cap \sO^{(\lambda, (\theta, \theta'))}\neq \emptyset$, then, we have
 \[(\lambda, (\theta,\theta'))\geq (\xi, (\eta,\eta')),\]
 i.e., $(\lambda, (\theta,\theta')-(\xi, (\eta,\eta'))$ belongs to $\Lambda^{pos}$ which is spanned by $\alpha_1, \alpha_2, \cdots, \alpha_{M+N-1}$.
 
\end{prop}
\begin{proof}
Let $g\in \GL_M(\bF)$, and $u\in U_{M,N}^-(\bF)$, if there is
\[g\begin{pNiceMatrix}
    t^{-\lambda_1}       &  &    \\
    &\Ddots& \\
          &   & t^{-\lambda_M} 
\end{pNiceMatrix}\in \begin{pNiceMatrix}
    t^{-\xi_1}       &  & &   \\
   \Vdots &\Ddots&\\
         * &\Cdots &  t^{-\xi_M} 
\end{pNiceMatrix}\GL_M(\bO),\]
then
\[g\in \begin{pNiceMatrix}
    t^{-\xi_1}       &  &   \\
   \Vdots &\Ddots& \\
       *   &\Cdots &  t^{-\xi_M} 
\end{pNiceMatrix}\GL_M(\bO) \begin{pNiceMatrix}
    t^{\lambda_1}       &  &    \\
    &\Ddots& \\
          &   & t^{\lambda_M} 
\end{pNiceMatrix}.\]
So, we can rewrite our assumption $\sS^{(\xi, (\eta, \eta'))}\cap \sO^{(\lambda, (\theta, \theta'))}\neq \emptyset$ as
\begin{equation}\label{intersect 8.1}
    \sD^{-1} U^-_{M}(\bF)t^{-\xi} \GL_M(\bO) U_{M,N}^-(\bF) \BL_{(\lambda, (\theta, \theta'))}\cap U_{N}(\bF) t^{(\eta, \eta')}\GL_N(\bO)\neq \emptyset.
\end{equation}

{First, we note that the above intersection is non-empty only if
\begin{equation}\label{eq triv}
    \sum_{i=1}^{M} \lambda_i+\sum_{i=1}^{M+1}\theta_{i}+\sum_{i=1}^{N-M-1}\theta_{i}'=\sum_{i=1}^{M} \xi_i+ \sum_{i=1}^{M+1} \eta_i+\sum_{i=1}^{N-M-1}\eta_{i}'.  
\end{equation}
}

{Within this proof, we denote by $\fA$ the former subset of \eqref{intersect 8.1} and by $\fB$ the latter subset. }


For $i\leq N-M-1$, the minimal degree of the determinant of $i\times i$-minors in lines $N-i+1, N-i+2, \cdots, N$ of any $\sB\in \fB$ is 
\[\eta'_{N-M-i}+\eta'_{N-M-i+1}+\cdots+\eta'_{N-M-1}.\]

{On the other hand, by considering the $i\times i$-minor $\sA_{\{N-i+1,\cdots, N\},\{N-i+1,\cdots, N\}}$ of any element $\sA\in\fA$,} we obtain \[\det\begin{pNiceMatrix}
    t^{\theta'_{N-M-i}}       &  & &  & \\
   \Vdots &\Ddots&& \\
   *        &\Cdots  & t^{\theta'_{N-M-1}}&
\end{pNiceMatrix}\in t^{\eta'_{N-M-i}+\eta'_{N-M-i+1}+\cdots+\eta'_{N-M-1}}\bO. \]
It implies that for $i\leq N-M-1$, we have
\begin{equation}\label{eq 8.2}
    \eta'_{N-M-i}+\eta'_{N-M-i+1}+\cdots+\eta'_{N-M-1}\leq \theta'_{N-M-i}+\theta'_{N-M-i+1}+\cdots+\theta'_{N-M-1}.
\end{equation}

For $1 \leq i\leq M$, the minimal degree of the determinant of $i \times i $-minors in the lines $1,\cdots,i$ of $\sA$, { for any $\sA\in \fA$} is \[-\xi_1-\xi_2-\cdots-\xi_i+\lambda_1+\theta_1+\lambda_2+\theta_2+\cdots+\lambda_i+\theta_i.\]

In particular, the degree of any $i\times i$-minor in the lines $1,\cdots, i$ of any element in \eqref{intersect 8.1} belongs to $t^{-\xi_1-\xi_2-\cdots-\xi_i+\lambda_1+\theta_1+\lambda_2+\theta_2+\cdots+\lambda_i+\theta_i}\bO$.

By multiplying by an element in $\GL_N(\bO)$, we obtain an element in $U_{N}(\bF)t^{\{\eta,\eta'\}}$, whose left-top $i\times i$-minor's determinant is $t^{\eta_1+\eta_2+\cdots+\eta_i}$. So, we have
\begin{equation}
\begin{split}
\eta_1+\eta_2+\cdots+\eta_i \geq -\xi_1-\xi_2-\cdots-\xi_i+\lambda_1+\theta_1+\lambda_2+\theta_2+\cdots+\lambda_i+\theta_i,
\end{split}
\end{equation}
i.e.,
\begin{equation}\label{eq 8.4} 
\begin{split}
\xi_1+\eta_1+\xi_2+\eta_2+\cdots+\xi_i+\eta_i \geq \lambda_1+\theta_1+\lambda_2+\theta_2+\cdots+\lambda_i+\theta_i.
\end{split}
\end{equation}

Now, we rewrite \eqref{intersect 8.1} as
\[ U^-_{M}(\bF)t^{-\xi} \GL_M(\bO) U_{M,N}^-(\bF) \BL_{(\lambda, \theta, \theta')}\GL_N(\bO)\cap \sD U_{N}(\bF) t^{(\eta, \eta')}\neq \emptyset, \]
and {we denote the first subset by $\sA'$, the second subset by $\sB'$}.

For $0\leq i\leq M$, we consider the rows $1,2,\cdots,i, M+1$. {For any $\sA\in \fA'$}, the minimal degree of the determinant of $(i+1)\times (i+1)$-minors of those rows is $-\xi_1-\xi_2-\cdots-\xi_i+\lambda_1+\theta_1+\cdots+\lambda_{i}+\theta_{i}+\theta_{i+1}$. Note that here we use the assumption $\lambda_1\leq \lambda_2\leq\cdots\leq \lambda_{i}$ and $\theta_1\leq\theta_2\leq\cdots\leq\theta_{M+1}$.

For any $\sB\in \fB'$, we consider the $(i+1)\times (i+1)$-minor 
\[\sB_{\{1,2,\cdots,i, M+1\}, \{1,2,\cdots,i+1\}}=\begin{pNiceMatrix}
    t^{\eta_{1}}       &*  & \Cdots&  &* \\
       & t^{\eta_{2}} &&&\Vdots\\
   &&\Ddots&\Ddots& \\
       & & & t^{\eta_{i}}&*\\
       v_1&v_2&\Cdots& v_i& v_{i+1}
\end{pNiceMatrix}.\]

Since $(v_1, v_2,\cdots, v_{i+1})$ equals the sum of row vectors of the lines $1, 2, \cdots,i$ with $(0, 0, \cdots, 0, t^{\eta_{i+1}})$, there is
\begin{equation}
\begin{split}
    &\det \sB_{\{1,2,\cdots,i, M+1\}, \{1,2,\cdots,i+1\}}\\= &\det \begin{pNiceMatrix}
    t^{\eta_{1}}       & * &\Cdots & & *\\
          & t^{\eta_{2}} &\Ddots&&  \\
   &&\Ddots&&\Vdots \\
       & &  & t^{\eta_{i}}&*\\
       &&& & t^{{\eta}_{i+1}}\\
\end{pNiceMatrix}\\=& t^{\eta_1+\eta_2+\cdots+\eta_{i}+\eta_{i+1}.}
\end{split}
\end{equation}
So, \begin{equation}
\begin{split}
\eta_1+\eta_2+\cdots+\eta_{i}+\eta_{i+1}\\ \geq -\xi_1-\xi_2-\cdots-\xi_i+\lambda_1+\theta_1+\cdots+\lambda_{i}+\theta_{i}+\theta_{i+1},
\end{split}
\end{equation}
i.e.,
\begin{equation}\label{eq 8.7}
\begin{split}
    \xi_1 +\eta_1+\xi_2+\eta_2+\cdots+\xi_i+\eta_i+\eta_{i+1}\\ \geq \lambda_1+\theta_1+\lambda_2+\theta_2+\cdots+\lambda_i+\theta_i+\theta_{i+1}.
\end{split}
\end{equation}

{The formulas \eqref{eq triv}, \eqref{eq 8.2}, \eqref{eq 8.4}, \eqref{eq 8.7} are equivalent to $(\lambda, (\theta,\theta'))\geq (\xi, (\eta,\eta'))$.}
\end{proof}

\begin{prop}\label{one point}
If $(\lambda,(\theta,\theta'))=(\xi,(\eta,\eta'))$, then the intersection $\sO^{(\lambda,(\theta,\theta'))}\cap \sS^{(\xi,(\eta,\eta'))}$ consists of exactly one point.
\end{prop}
\begin{proof}
It only contains the point $(t^{-\xi}, \sD t^{(\eta,\eta')})\in \Gr_M\times \Gr_N$.
\end{proof}
\subsection{Closure relation of  relevant orbits}
In this section, we will study the closure relation of $H^{\omega}(\bF)$-orbits in $\Gr^{\omega}_G$. We only need to consider the closure relation of $\GL_M(\bO)\ltimes U_{M,N}^-(\bF)$-orbits in $\Gr_N$. 

We claim
\begin{prop}\label{closure}
For two $\GL_M(\bO)\ltimes U_{M,N}^-(\bF)$-orbits $\BO^{(\lambda, (\theta, \theta'))}$ and $\BO^{(\tilde{\lambda}, \tilde{\theta}, \tilde{\theta'})}$ of $\Gr_N$,
\begin{equation}
    \begin{split}
        &\BO^{(\lambda, (\theta, \theta'))}\subset \overline{\BO}^{(\tilde{\lambda}, (\tilde{\theta}, \tilde{\theta}'))}\\
        & \Longleftrightarrow\\
        & (\lambda, (\theta, \theta'))\leq (\tilde{\lambda}, (\tilde{\theta}, \tilde{\theta}')).
   \end{split}
\end{equation}
\end{prop}
\begin{proof}
By \cite[Proof of Lemma 3.6.2 (b)]{[BFT0]} and the usual semi-infinite orbit closure relation, we obtain the \textit{if} direction. 

Now, we prove the \textit{only if} direction.

For $1\leq i\leq N-M-1$, the minimal degree of the determinant of $(M+1+i)\times (M+1+i)$-minors in lines $1,2,\cdots, M+1+i$ of any (representative of) element in $\BO^{(\tilde{\lambda}, (\tilde{\theta}, \tilde{\theta}'))}$ is $\tilde{\lambda}_1+\tilde{\theta}_1+\cdots+\tilde{\lambda}_M+\tilde{\theta}_M+\tilde{\theta}_{M+1}+\tilde{\theta_1}'+\tilde{\theta}_2'+\cdots+ \tilde{\theta}_i'$, the minimal degree of the determinant of $(M+1+i)\times (M+1+i)$-minors in lines $1,2,\cdots, M+1+i$ of any element in $\BO^{({\lambda}, {\theta}, {\theta}')}$ is ${\lambda}_1+{\theta}_1+\cdots+{\lambda}_M+{\theta}_M+\theta_{M+1}+{\theta_1}'+{\theta}_2'+\cdots+ {\theta}_i'$. So, we have
\begin{equation}
   \tilde{\lambda}_1+\tilde{\theta}_1+\cdots+\tilde{\lambda}_M+\tilde{\theta}_M+\tilde{\theta}_{M+1}+\tilde{\theta_1}'+\tilde{\theta}_2'+\cdots+ \tilde{\theta}_i'\leq {\lambda}_1+{\theta}_1+\cdots+{\lambda}_M+{\theta}_M+\theta_{M+1}+{\theta_1}'+{\theta}_2'+\cdots+ {\theta}_i'.
\end{equation}

For $1\leq i\leq M$, the minimal possible minimal degree of the determinant of $i\times i$-minors in lines $1,2,\cdots, i$ of all (representatives of) elements of $\BO^{\tilde{\lambda}, (\tilde{\theta}, \tilde{\theta}')}$ is $\tilde{\lambda}_1+\tilde{\theta}_1+\cdots+\tilde{\lambda}_i+\tilde{\theta}_i$, and  the minimal possible minimal degree of the determinant of $i\times i$-minors in lines $1,2,\cdots, i$ of all elements of $\BO^{{\lambda}, {\theta}, {\theta}'}$ is ${\lambda}_1+{\theta}_1+\cdots+{\lambda}_i+{\theta}_i$. So, we have
\begin{equation}
    \tilde{\lambda}_1+\tilde{\theta}_1+\cdots+\tilde{\lambda}_i+\tilde{\theta}_i\leq {\lambda}_1+{\theta}_1+\cdots+{\lambda}_i+{\theta}_i.
\end{equation}

For $0\leq i\leq M$, the minimal possible minimal degree of the determinant of $(i+1)\times (i+1)$-minors in lines $1,2,\cdots, i, M+1$ of all elements of $\BO^{\tilde{\lambda}, (\tilde{\theta}, \tilde{\theta}')}$ is $\tilde{\lambda}_1+\tilde{\theta}_1+\cdots+\tilde{\lambda}_i+\tilde{\theta}_i+\tilde{\theta}_{i+1}$, and  the minimal possible minimal degree of the determinant of $(i+1)\times (i+1)$-minors in lines $1,2,\cdots, i, M+1$ of all elements of $\BO^{{\lambda}, {\theta}, {\theta}'}$ is ${\lambda}_1+{\theta}_1+\cdots+{\lambda}_i+{\theta}_i+\theta_{i+1}$. So, we have
\begin{equation}
   \tilde{\lambda}_1+\tilde{\theta}_1+\cdots+\tilde{\lambda}_i+\tilde{\theta}_i+\tilde{\theta}_{i+1}\leq {\lambda}_1+{\theta}_1+\cdots+{\lambda}_i+{\theta}_i+\theta_{i+1}.
\end{equation}
It finishes the proof of the only if direction.
\end{proof}

\subsection{Dimension estimate, I}
For the convenience of future use, we introduce the following notations:
\begin{itemize}
    \item For $1 \leq i\leq M-1$, we denote $\alpha_{i, \GL_M}:= \alpha_{2i}+\alpha_{2i+1}=-\delta_{i}+\delta_{i+1}$.
    \item For $ 1\leq i\leq M$, we denote $\alpha_{i, \GL_N}:= \alpha_{2i-1}+\alpha_{2i}= -\epsilon_{i}+\epsilon_{i+1}$.
    \item For $M+1\leq i\leq N-1$, we denote $\alpha_{i, \GL_N}:= \alpha_{i+M}$.
\end{itemize}

Applying the argument of \cite[Proposition 6.1.1]{[SW]} to the above propositions, and using the fact that 
\begin{equation}
    \overline{\sS}^{(\lambda,(\theta, \theta'))}= \bigsqcup_{(\lambda,(\theta, \theta'))\underset{G}{\leq} (\tilde{\lambda}, (\tilde{\theta}, \tilde{\theta}'))}{\sS}^{(\tilde{\lambda}, (\tilde{\theta}, \tilde{\theta}'))},
\end{equation}
we obtain the following corollary
\begin{cor}\label{cor 8.3.1}
\begin{equation}
    \dim \overline{\sO}^{(\lambda, (\theta, \theta'))} \cap \overline{\sS}^{(\xi, (\eta, \eta'))} \leq \max\{\sum_{i=1}^{M-1} a_i+ \sum_{j=1}^{M} b_j+ \sum_{k=M+1}^{N-1} c_k\},
\end{equation}
where the maximum is taken over the set
\[\{(\tilde{\xi}, (\tilde{\eta}, \tilde{\eta}')): (\xi, (\eta, \eta'))\underset{G}{\leq} (\tilde{\xi}, (\tilde{\eta}, \tilde{\eta}'))\leq (\lambda, (\theta, \theta'))\},\]
and 
\begin{equation}
    \begin{split}
        (\tilde{\xi}, (\tilde{\eta}, \tilde{\eta}'))-({\xi}, {\eta}, {\eta}')\\
        = \sum_{i=1}^{M-1} a_i \alpha_{i, \GL_M}+ \sum_{j=1}^{M} b_j \alpha_{j, \GL_N}+ \sum_{k=M+1}^{N-1} c_k \alpha_{k, \GL_N}.
    \end{split}
\end{equation}

\end{cor}

We denote by ${\cY}_{= (\lambda, (\theta, \theta'))\cdot c}^{(\xi, (\eta, \eta'))}$ the preimage of $\CM_{= (\lambda, (\theta, \theta'))\cdot c}$ in $\CY_{\infty\cdot c}^{(\xi, (\eta, \eta'))}$. It is a smooth locally closed subscheme of $\CY_{\infty\cdot c}^{(\xi, (\eta, \eta'))}$, ref. \cite[Proposition 3.6.1]{[SW]}. Let $\overset{\circ}{\cY}{}_{= (\lambda, (\theta, \theta'))\cdot c}^{(\xi, (\eta, \eta'))}$ be the open subset of ${\cY}_{= (\lambda, (\theta, \theta'))\cdot c}^{(\xi, (\eta, \eta'))}$, which is the preimage of the open locus $\overset{\circ}{C}{}_{= (\lambda, (\theta, \theta'))\cdot c}^{(\xi, (\eta, \eta'))}$, i.e.,
\[\xymatrix{
{\cY}_{= (\lambda, (\theta, \theta'))\cdot c}^{(\xi, (\eta, \eta'))}\ar[d]& \overset{\circ}{\cY}{}_{= (\lambda, (\theta, \theta'))\cdot c}^{(\xi, (\eta, \eta'))}\ar[l]\ar[d]\\
C_{\leq (\lambda, (\theta, \theta'))\cdot c}^{(\xi, (\eta, \eta'))}& \overset{\circ}{C}{}_{= (\lambda, (\theta, \theta'))\cdot c}^{(\xi, (\eta, \eta'))}\ar[l]
}\]

    Although the connected components of $\CM_{= (\lambda, (\theta, \theta'))\cdot c}$ are indexed by $\BZ$, it is not hard to see that there only exists at most one connected component $\CM^{\xi}_{= (\lambda, (\theta, \theta'))\cdot c}$ such that the relative fiber product with $\Bun_B^{(\xi, (\eta,\eta'))}$ over $\Bun_G$ is non-empty. Here, $\CM^{\xi}_{= (\lambda, (\theta, \theta'))\cdot c}$ denotes the connected component of $\CM_{= (\lambda, (\theta, \theta'))\cdot c}$ such that the degree of the determinant line bundle of the associated $\GL_M$-bundle has the same degree as the degree of the determinant line bundle of the associated $\GL_M$-bundle of any point in $\Bun_B^{(\xi, (\eta,\eta'))}$. 

If we have 
\begin{equation}
    (\lambda, (\theta, \theta'))-(\xi, (\eta, \eta'))= \sum_{i=1}^{M+N-1} n_i \alpha_i,
\end{equation}
then the fiber of $\overset{\circ}{\cY}{}_{= (\lambda, (\theta, \theta'))\cdot c}^{(\xi, (\eta, \eta'))}$ over the point 
\[D= (\lambda, (\theta, \theta'))\cdot c+\sum_{x\in \supp(D)} \alpha_{i_x}\cdot x\]
is isomorphic to
\[(\sO^{(\lambda,(\theta,\theta'))}\cap \sS^{(\lambda,(\theta,\theta'))})\times \prod_{x} (\sO^{0}\cap \sS^{\alpha_{i_x}}_x)\simeq \prod_{x} (\sO^{0}\cap \sS^{\alpha_{i_x}}_x). \]

If $1\leq i\leq 2M$, then $\sO^{0}\cap \sS^{\alpha_{i_x}}_x=\pt$, if $2M+1 \leq i\leq M+N-1$, then $\sO^{0}\cap \sS^{\alpha_{i_x}}_x$ is 1-dimensional. So, the dimension of fiber over $D$ is $\sum_{i=2M+1}^{N+M-1} n_i$, i.e., the sum of coefficients of even simple roots of the supergroup.

So, we have
\begin{equation}\label{eq 8.17}
    \dim \overset{\circ}{\cY}{}_{= (\lambda, (\theta, \theta'))\cdot c}^{(\xi, (\eta, \eta'))}= \sum_{i=1}^{2M} n_i+ \sum_{i=2M+1}^{M+N-1} 2n_i, 
\end{equation}
i.e., the sum of coefficients of odd simple roots plus two times sum of coefficients of even simple roots.

\begin{rem}\label{8.3.2}
    Actually, through the same analysis as above, one can show that the dimension of the complement of $\overset{\circ}{\cY}{}_{= (\lambda, (\theta, \theta'))\cdot c}^{(\xi, (\eta, \eta'))}$ in ${\cY}{}_{= (\lambda, (\theta, \theta'))\cdot c}^{(\xi, (\eta, \eta'))}$ is strictly smaller than $\dim \overset{\circ}{\cY}{}_{= (\lambda, (\theta, \theta'))\cdot c}^{(\xi, (\eta, \eta'))}$. Namely, given a decomposition $(\xi,(\eta,\eta'))=(\xi_c,(\eta_c,\eta_c'))+\sum_{k} \alpha_k$ and consider the corresponding stratum consisting of the points of the form $(\xi_c,(\eta_c,\eta_c'))\cdot c+\sum \alpha_{k}\cdot x$ in ${C}_{\leq (\lambda, (\theta, \theta'))\cdot c}^{(\xi, (\eta, \eta'))}$. The dimension of its preimage in $\cY_{= (\lambda, (\theta, \theta'))\cdot c}^{(\xi, (\eta, \eta'))}$ equals $\dim \overset{\circ}{\cY}{}_{= (\lambda, (\theta, \theta'))\cdot c}^{(\xi, (\eta, \eta'))}$ only if the following two conditions hold: (1). $\dim (\sO^0\cap \sS^{\alpha_{k}})+1=\dim (\sO^0\cap \sS^{\alpha_{k,1}})+\cdots+\dim (\sO^0\cap \sS^{\alpha_{k,m}})+m$, for any $k$. Here, $\alpha_{k}=\sum \alpha_{k,i}$ and each $\alpha_{k,i}$ is a negative simple root of the supergroup; (2). $(\xi_c,(\eta_c,\eta_c'))=(\lambda, (\theta, \theta'))$. For the right-hand side of the first condition, it equals the number of the odd simple roots plus twice of the number of the even simple roots appearing in $\alpha_k$. On the other hand, $\dim (\sO^0\cap \sS^{\alpha_{k}})$ is no more than half of the number of the odd simple roots plus the number of the even simple roots appearing in $\alpha_k$. So, the equality implies that the number of odd roots plus twice of the number of even roots is no more than $2$, i.e., if $\alpha_k$ is not a negative simple root, then it must be a sum of two odd negative simple roots $\alpha_{k,1},\alpha_{k,2}$. However, in this case, $\sO^0\cap \sS^{\alpha_1+\alpha_2}$ is a (possibly empty) set of points (the case when $\alpha_{k,1}+\alpha_{k,2}$ is a simple root of $\GL_M\times \GL_N$ is proved in \cite[Lemma 5.4.2]{[SW]}, the other cases are direct consequences of Corollary \ref{cor 8.3.1}). Thus, the condition (1) is not satisfied. 
\end{rem}

From Corollary \ref{cor 8.3.1} and \eqref{eq 8.17}, we obtain the following corollary.
\begin{cor}\label{cor 8.1.3}
For any irreducible component $Y$ of $\overline{\sO}^{(\lambda,(\theta,\theta'))}\cap \overline{\sS}^{(\xi, (\eta, \eta'))}$, we have
\begin{equation}\label{eq 8.18}
    \dim \overset{\circ}{\cY}{}_{= (\lambda, (\theta, \theta'))\cdot c}^{(\xi, (\eta, \eta'))}\geq 2 \dim Y.
\end{equation}

Furthermore, if the equality holds, then there are non-negative integers $a_i, b_j, c_k$, such that
\begin{equation}
    (\lambda, (\theta, \theta'))-(\xi, (\eta, \eta'))= \sum_{i=1}^{M-1} a_i \alpha_{i, \GL_M}+\sum_{j=1}^M b_j \alpha_{j, \GL_N}+ \sum_{k=M+1}^{N-1} c_k \alpha_{k, \GL_N},
\end{equation}
and if there is $Y$, such that \begin{equation}
    \dim \overset{\circ}{\cY}{}_{= (\lambda, (\theta, \theta'))\cdot c}^{(\xi, (\eta, \eta'))}= 2 \dim Y+1,
\end{equation}
then there are non-negative integers $a_i, b_j, c_k$, such that
\begin{equation}
    (\lambda, (\theta, \theta'))-(\xi, (\eta, \eta'))= \alpha+ \sum_{i=1}^{M-1} a_i \alpha_{i, \GL_M}+\sum_{j=1}^M b_j \alpha_{j, \GL_N}+ \sum_{k=M+1}^{N-1} c_k \alpha_{k, \GL_N},
\end{equation}
here $\alpha$ is an odd simple root of the supergroup.
\end{cor}
\begin{proof}
Let $(\tilde{\xi}, \tilde{\eta},\tilde{\eta}')$ be an element from the set
\[\{(\tilde{\xi}, (\tilde{\eta}, \tilde{\eta}')): (\xi, (\eta, \eta'))\underset{G}{\leq} (\tilde{\xi}, (\tilde{\eta}, \tilde{\eta}'))\leq (\lambda, (\theta, \theta'))\},\]
such that
\begin{equation}
    \begin{split}
        (\tilde{\xi}, (\tilde{\eta}, \tilde{\eta}'))-({\xi}, ({\eta}, {\eta}'))\\
        = \sum_{i=1}^{M-1} a_i \alpha_{i, \GL_M}+ \sum_{j=1}^{M} b_j \alpha_{j, \GL_N}+ \sum_{k=M+1}^{N-1} c_k \alpha_{k, \GL_N}
    \end{split}
\end{equation}
takes the maximum. Since $(\tilde{\xi}, (\tilde{\eta}, \tilde{\eta}'))\leq (\lambda, (\theta, \theta'))$, we have
\begin{equation}
    \begin{split}
        \sum_{i=1}^{2M} n_i\geq 2(\sum_{i=1}^{M-1} a_i+\sum_{j=1}^{M} b_j )\\
        \sum_{i=2M+1}^{N+M-1} n_i\geq \sum_{k=M+1}^{N-1} c_k,
    \end{split}
\end{equation}
By summing up, we get \eqref{eq 8.18}. The latter two statements follow immediately when considering the cases where the equality holds.
\end{proof}

\begin{rem}
Note that 
\begin{equation}
\begin{split}
     \dim \overset{\circ}{\cY}{}_{= (\lambda, (\theta, \theta'))\cdot c}^{(\xi, (\eta, \eta'))}= \dim \cM^{\xi}_{=(\lambda, (\theta, \theta'))\cdot c}+ dim.rel (\Bun_B^{(\xi, (\eta, \eta'))}, \Bun_G),
\end{split}
\end{equation}
so the above corollary can be rewritten as
\begin{equation}
    \dim \cM^\xi_{=(\lambda, (\theta, \theta'))\cdot c}+ dim.rel (\Bun_B^{(\xi, (\eta, \eta'))}, \Bun_G)\geq 2\dim Y,
\end{equation}
for any irreducible component of $\overline{\sO}^{(\lambda,(\theta,\theta'))}\cap \overline{\sS}^{(\xi, (\eta, \eta'))}$.
\end{rem}

\section{Proof II: Image of irreducible objects}\label{section 9}
With the preparation in the last section, now, we are going to show that the irreducible objects in the twisted Gaiotto category go to the irreducible objects (i.e., twisted IC D-modules on the configuration space) in the category of factorization modules. 

\subsection{Basic properties of the functor}
First, we show that the functor $F^{glob}$ is $t$-exact with respect to the perverse $t$-structures on both sides.
\begin{prop}\label{p tex}
The shifted pullback functor $p^\bullet: \cC_{\kappa}^{glob}(M|N)\longrightarrow {\SD}_{\kappa}(\cY_{\infty\cdot c})$ is $t$-exact.
\end{prop}

\begin{proof}
First, we note that for $(\xi, (\eta, \eta'))$ large enough, i.e., $(\xi, (\eta, \eta'))$ is deep enough in the dominant chamber, $\Bun_B^{(\xi, (\eta, \eta'))}\rightarrow \Bun_G$ is smooth  (note that according to our notation, it means the degree is very anti-dominant). So, as an open subspace of $\cM_{\infty\cdot c}\underset{\Bun_G}{\times} \Bun_B^{(\xi, (\eta,\eta'))}$, the SW-Zastava space $\cY^{(\xi,(\eta,\eta'))}_{\infty\cdot c}$ is smooth over $\cM_{\infty\cdot c}$. So, the restriction of $p^\bullet$ to those components is $t$-exact. According to Proposition \ref{locglob fun} and the factorization property of $\sprd(\pi^!(\cF))$ for $\cF\in \cC_{\kappa}^{glob}(M|N)$, we obtain that $p^\bullet(\cF)\in {\SD}_{\kappa}(\cY_{\infty\cdot c})$ is factorizable with respect to the factorizable structure of the SW-Zastava space (see \eqref{Zastava c}). 

Now, we assume  $(\tilde{\xi}+\xi, (\tilde{\eta}+\eta, \tilde{\eta}'+\eta'))$ and $({\xi}, ({\eta}, {\eta}'))$  to be large enough.

According to \cite[Theorem 16.2.1]{[GN]} and \cite[Lemma 3.5.4]{[SW]}, there is a correspondence 
\begin{equation}
    \xymatrix{
    \cY_{\infty\cdot c}^{(\tilde{\xi}, (\tilde{\eta}, \tilde{\eta}'))}&  (\cY_{\infty\cdot c}^{(\tilde{\xi}, (\tilde{\eta}, \tilde{\eta}'))}\times \cY^{(\xi, (\eta,\eta'))}_{\text{genu}})_{disj}\ar[r]\ar[l]& \cY_{\infty\cdot c}^{(\tilde{\xi}+\xi, (\tilde{\eta}+\eta, \tilde{\eta}'+\eta'))}\ar[d]\\&& \cM_{\infty\cdot c}.
    }
\end{equation}
Here, $\cY^{(\xi, (\eta,\eta'))}_{\text{genu}}$ is the  preimage of $\Bun^{\omega}_H\subset \cM$ in $\cY^{(\xi, (\eta,\eta'))}$, {it is isomorphic to ${\cY^{(\xi, (\eta,\eta'))}_{=0\cdot c}}$}.\footnote{Here, the subscript 'genu' means that the generic $H$-structure is genuine.} The left arrow is smooth surjective, and the right composed arrow is smooth.

For $\cF\in \cC_{\kappa}^{glob}(M|N)$, $p^{\bullet,(\tilde{\xi}+\xi, (\tilde{\eta}+\eta, \tilde{\eta}'+\eta')) }(\cF)$ is perverse, and by the factorization property, its further pullback along the \'{e}tale map $(\cY_{\infty\cdot c}^{(\tilde{\xi}, (\tilde{\eta}, \tilde{\eta}'))}\times \cY^{(\xi, (\eta,\eta'))}_{\text{genu}})_{disj}\longrightarrow \cY_{\infty\cdot c}^{(\tilde{\xi}+\xi, (\tilde{\eta}+\eta, \tilde{\eta}'+\eta'))}$ is perverse and isomorphic to the restriction of $p^{\bullet, (\tilde{\xi},(\tilde{\eta},\tilde{\eta}'))}(\cF)\boxtimes p^{\bullet, (\xi, (\eta,\eta'))}(\IC_{glob}^{0})$. 

Since we require $(\xi,(\eta,\eta'))$ to be large enough, $p^{\bullet, (\xi, (\eta,\eta'))}(\IC_{glob}^{0})$ is a perverse local system on a smooth variety. It implies that the restriction of  $p^{\bullet, (\tilde{\xi},(\tilde{\eta},\tilde{\eta}'))}(\cF)\boxtimes \IC_{\cY^{(\xi,(\eta,\eta')),\circ}}$ to  $(\cY_{\infty\cdot c}^{(\tilde{\xi}, (\tilde{\eta}, \tilde{\eta}'))}\times \cY^{(\xi, (\eta,\eta'))}_{\text{genu}})_{disj}$ is perverse. 

Now, since the left arrow is smooth surjective, and a D-module is concentrated in degree 0 if and only if it is concentrated in degree 0 in a smooth atlas (for example, see \cite[Lemma 4.1]{[LO]}), the D-module $p^{\bullet, (\tilde{\xi},(\tilde{\eta},\tilde{\eta}'))}(\cF)$ is perverse.
\end{proof}

Applying the same method as above, we can obtain the following corollaries.


\begin{cor}\label{cor 9.1.2}
The shifted pullback $p^\bullet$ along $\cY_{\infty\cdot c}\rightarrow \cM_{\infty\cdot c}$ sends the irreducible object $\IC_{glob}^{(\lambda,(\theta,\theta'))}\in \cC_{\kappa}^{glob}(M|N)$ to the $!*$-extension of a twisted perverse D-module (i.e., the restriction of $p^\bullet(\IC_{glob}^{(\lambda,(\theta,\theta'))})$) on $\cY_{=(\lambda,(\theta,\theta'))\cdot c}$. 
\end{cor}

\begin{rem}
    The above corollary implies that the dimension of each connected component in the smooth stack $\cY_{=(\lambda,(\theta,\theta'))\cdot c}^{(\xi,(\eta,\eta'))}$ has the same dimension. In particular, according to Remark \ref{8.3.2}, we obtain that $\cY_{=(\lambda,(\theta,\theta'))\cdot c}^{(\xi,(\eta,\eta'))}$ is connected.
\end{rem}

\begin{cor}\label{p ver}
The shifted pullback functor $p^\bullet: \cC_{\kappa}^{glob}(M|N)\longrightarrow {\SD}_{\kappa}(\cY_{\infty\cdot c})$ commutes with the Verdier duality functor.
\end{cor}

As a combination of Lemma \ref{clean}, Proposition \ref{p tex}, and Corollary \ref{p ver}, we have
\begin{cor}\label{t exact}
For $\kappa$ generic, the functor $F^{glob}: \cC_{\kappa}^{glob}(M|N)\longrightarrow {\SD}_{\kappa}(C^\bullet_{\infty\cdot c})$ is $t$-exact, and commutes with the Verdier duality functor. 
\begin{proof}
The map $v: \cY_{\infty\cdot c}\rightarrow C^\bullet_{\infty\cdot c}$ is affine. So, $v_*p^\bullet:\cC_{\kappa}^{glob}(M|N)\longrightarrow {\SD}_{\kappa}(C^\bullet_{\infty\cdot c}) $ is right $t$-exact, i.e., for $\cF\in \cC_{\kappa}^{glob}(M|N)$, $v_*p^\bullet(\cF)$ is concentrated in non-positive degrees. Also, $v_!p^\bullet:\cC_{\kappa}^{glob}(M|N)\longrightarrow {\SD}_{\kappa}(C^\bullet_{\infty\cdot c}) $ is left $t$-exact, i.e., for $\cF\in \cC_{\kappa}^{glob}(M|N)$, we have $v_!p^\bullet(\cF)$ is concentrated in non-negative degrees. According to Lemma \ref{clean}, $v_*\simeq v_!$, so we obtain that $F^{glob}$ is $t$-exact.

The commutativity with the Verdier duality functor follows from Corollary \ref{p ver} and $\BD v_*\simeq v_!\BD$.
\end{proof}
\end{cor}

{\begin{rem}
    Actually, for any $G, H, \chi$ mentioned in the Theorem \ref{general}, if the determinant line bundle on $\Gr_G$ admits an $H(\bF)$-equivariant structure, then, we can still define the functor from the category of $\cP_{\det}^\kappa$-twisted $(H(\bF),\chi)$-equivariant D-modules on $\Gr_G$ to the category of $\cP_{C^\bullet}$-twisted D-modules on the configuration space, and Proposition \ref{p tex}, Corollary \ref{cor 9.1.2}, \ref{p ver}, \ref{t exact} still hold.
\end{rem}}

\subsubsection{}
Consider the intersection $\sO^{(\lambda,(\theta,\theta'))}\cap \sS^{(\xi, (\eta,\eta'))}\subset \Gr_G$. For relevant $(\lambda, (\theta,\theta'))$, there is a unique $H(\bF)$-equivariant trivialization of the restriction of $\cP_{\det}$ to $\sO^{(\lambda,(\theta,\theta'))}$, and there is a unique $U(\bF)$-equivariant trivialization of $\cP_{\det}$ on $\sS^{(\xi, (\eta,\eta'))}$. We denote by $\psi_{(\lambda,(\theta,\theta')),(\xi, (\eta,\eta'))}$ the ratio of the above two trivializations. 

There is an action of 
\[T_{N-M-1}(\bO):=\{\begin{pNiceMatrix}
    1       &   &  & &&\\
    &\Ddots&& &&\\
          &   & 1 &&&\\
          &&&a_1&&\\
          &&&&\Ddots&\\
          &&&&&&a_{N-M-1}
\end{pNiceMatrix}\in \GL_N(\bO), a_i\in \bO^\times\}\]
on $\sO^{(\lambda,(\theta,\theta'))}\cap \sS^{(\xi, (\eta,\eta'))}\subset \Gr_G$. We have the following lemma which is an analog of \cite[Lemma 5.5]{[G]}.
\begin{lem}
The function $\psi_{(\lambda,(\theta,\theta')),(\xi, (\eta,\eta'))}$ intertwines the natural $T_{N-M-1}(\bO)$-action on $\sO^{(\lambda,(\theta,\theta'))}\cap \sS^{(\xi, (\eta,\eta'))}$ and the action on $\BG_m$ given by the character
\[T_{N-M-1}(\bO)\longrightarrow T_{N-M-1}\overset{(\theta'-\eta',-)}{\longrightarrow} \BG_m.\]
\end{lem}

Let $\Psi_{(\lambda,(\theta,\theta')), (\xi,(\eta,\eta'))}$ be the $!$-pullback of Kummer local system on $\BG_m$ corresponding to $\kappa$ along the map $\psi_{(\lambda,(\theta,\theta')),(\xi, (\eta,\eta'))}$, whose $!$-stalks are concentrated in degree $0$. We have the following lemma which is an analog of \cite[Lemma 3.7.2]{[BFT0]} and \cite[Lemma 5.5]{[G]}. Let us sketch the proof.
\begin{lem}\label{non trivial}
The local system $\Psi_{(\lambda,(\theta,\theta')), (\xi,(\eta,\eta'))}$ is not trivial on any irreducible component $Y\subset \sO^{(\lambda,(\theta,\theta'))}\cap \sS^{(\xi, (\eta,\eta'))}$ if $\dim Y=\frac{1}{2} \dim \overset{\circ}{\cY}{}_{= (\lambda, (\theta, \theta'))\cdot c}^{(\xi, (\eta, \eta'))}$ or $\dim Y=\frac{1}{2} (\dim \overset{\circ}{\cY}{}_{= (\lambda, (\theta, \theta'))\cdot c}^{(\xi, (\eta, \eta'))}-1)$ except for $(\lambda,(\theta,\theta'))=(\xi,(\eta,\eta'))$ or $(\lambda,(\theta,\theta'))=(\xi,(\eta,\eta'))+\alpha$, where $\alpha$ is an odd simple root of the supergroup.
\end{lem}
\begin{proof}
By the $T_{N-M-1}(\bO)$-equivariance against a non-trivial character, the local system $\psi_{(\lambda,(\theta,\theta')),(\xi, (\eta,\eta'))}$ is trivial only if $\eta'=\theta'$. In this case, we can reduce the question to the case $M=N-1$.

If $\dim Y=\frac{1}{2} \dim \overset{\circ}{\cY}{}_{= (\lambda, (\theta, \theta'))\cdot c}^{(\xi, (\eta, \eta'))}$, it is proved in \cite[Lemma 3.7.2]{[BFT0]}. Actually, the proof of \cite[Lemma 3.7.2]{[BFT0]} is also available for $\dim Y=\frac{1}{2} (\dim \overset{\circ}{\cY}{}_{= (\lambda, (\theta, \theta'))\cdot c}^{(\xi, (\eta, \eta'))}-1)$.

Indeed, similar{ly} to \cite{[BFT0]}, we only need to check that there is no boundary component $\partial_i \bar{Y}$ of the closure $ \bar{Y}$ which is an irreducible component of $\overline{\sO}^{(\lambda,(\theta,\theta'))}\cap \overline{\sS}^{(\xi+\alpha_{\GL_M}, (\eta,\eta')+\alpha_{\GL_N})}$.  Here $\alpha_{\GL_M}$ (resp.\ $\alpha_{\GL_N}$) is a certain simple root of $\GL_M$ (resp.\ $\GL_N$), which in turns to show the dimension of $\overline{\sO}^{(\lambda,(\theta,\theta'))}\cap \overline{\sS}^{(\xi+\alpha_{\GL_M}, (\eta,\eta')+\alpha_{\GL_N})}$ is strictly smaller than $\dim Y-1$.

In this case, according to Corollary \ref{cor 8.1.3}, we have
\begin{equation}\label{eq 9.2}
    (\lambda, (\theta, \theta'))-(\xi, (\eta, \eta'))= \alpha+ \sum_{i=1}^{M-1} a_i \alpha_{i, \GL_M}+\sum_{j=1}^M b_j \alpha_{j, \GL_N}+ \sum_{k=M+1}^{N-1} c_k \alpha_{k, \GL_N},
\end{equation}
where $\alpha$ is an odd simple root of the supergroup. 


The \eqref{eq 9.2} and Corollary \ref{cor 8.3.1} implies
\[\dim \overline{\sO}^{(\lambda,(\theta,\theta'))}\cap \overline{\sS}^{(\xi+\alpha_{\GL_M}, (\eta,\eta')+\alpha_{\GL_N})}\leq \sum a_i+ \sum b_j +\sum c_k-2.\]
Furthermore, the assumption $\dim Y=\frac{1}{2} (\dim \overset{\circ}{\cY}{}_{= (\lambda, (\theta, \theta'))\cdot c}^{(\xi, (\eta, \eta'))}-1)$ implies $\sum a_i+ \sum b_j +\sum c_k= \dim Y$.




\end{proof}

\subsection{Intersection with the unital relevant orbit}
Now, we focus on the case $(\lambda,(\theta,\theta'))=0$ and study the intersection of $\sO^0$ with $U(\bF)$-orbits.  

It is more convenient to do calculations in $\Gr_N$ than in $\Gr_M\times \Gr_N$, so we first realize the intersection $\sO^{0}\cap \sS^{(\xi, (\eta,\eta'))}$ as a subspace of $\Gr_N$.

Note that there is a projection 
\[\pr_2: \Gr_M\times \Gr_N\longrightarrow \Gr_N.\]
We claim that it induces a map from $\sO^0\cap \sS^{(\xi,(\eta,\eta'))}$ to 
\begin{equation}\label{intersection 9.4}
    \begin{split}
        W^{(\xi,(\eta,\eta'))}&= \begin{pNiceMatrix}
   t^{-\xi_1}   &  & &  &  &&&\\
       *   & &  &  &  \\
   \Vdots&\Ddots&\Ddots&& &&&\\
    *     & \Cdots& * &t^{-\xi_M} &&&& \\
       &  &&  & 1&&&\\
     *    &&\Cdots&&*&1&&\\
      \Vdots   &&&&&\Ddots&\Ddots&\\
       *  &\Cdots&&&&&*&1
\end{pNiceMatrix}\GL_N(\bO)/\GL_N(\bO)\cap\\
&\cap \sD\begin{pNiceMatrix}
    t^{\eta_{1}}       &*  &\Cdots &&&*\\
          & \Ddots&\ddots&&& \Vdots \\
       & & t^{{\eta}_{M+1}}&\ddots&&\\
       &&&t^{\eta_1'}&\ddots&\\
       &&&&\Ddots&*\\
       &&&&&t^{\eta_{N-M-1}'}
\end{pNiceMatrix}\GL_N(\bO)/\GL_N(\bO)\\
&= U_M^-(\bF) U_{M,N}^-(\bF) \diag(t^{-\xi_1},\cdots, t^{-\xi_M},1,\cdots,1)\GL_N(\bO)/\GL_N(\bO)\cap\\
&\cap \sD U_N(\bF)\diag(t^{\eta_1},\cdots, t^{\eta_{M+1}},t^{\eta_1'},\cdots, t^{\eta_{N-M-1}'})\GL_N(\bO)/\GL_N(\bO)\\
&\subset \Gr_N.
    \end{split}
\end{equation}
Recall that
\[\sO^0= \{(g \GL_M(\bO)/\GL_M(\bO) , g u\GL_N(\bO)/\GL_N(\bO))| g\in \GL_M(\bF), u\in U_{M,N}^-(\bF)\},\]
and
\begin{equation*}
    \begin{split}
\sS^{(\xi, (\eta,\eta'))}=& \sS^{-\xi}\times \sD\sS^{(\eta, \eta')}\\ =&  U_M^-(\bF) t^{-\xi}\GL_M(\bO)/\GL_M(\bO)\times \sD U_N(\bF) t^{(\eta,\eta')}\GL_N(\bO)/\GL_N(\bO)\\ \subset& \Gr_M\times \Gr_N.       
    \end{split}
\end{equation*}


If $g \GL_M(\bO)/\GL_M(\bO)\in \sS^{-\xi}$, we have $g\in  U_M^-(\bF) t^{-\xi} \GL_M(\bO)$. So, the image of $\sO^0\cap \sS^{(\xi,(\eta,\eta'))}$ under $\pr_2$ lands in the first space of \eqref{intersection 9.4}. 
By definition, it also lands in $\sD \sS^{(\eta, \eta')}$. In particular, the restriction of $\pr_2$ to $\sO^0\cap \sS^{(\xi,(\eta,\eta'))}$ induces a map
\begin{equation}\label{pr2}
    \pr_2: \sO^0\cap \sS^{(\xi,(\eta,\eta'))}\longrightarrow W^{(\xi,(\eta,\eta'))}.
\end{equation}

In fact, we have
\begin{lem}
The above map \eqref{pr2} is an isomorphism.
\end{lem}
\begin{proof}
Let us construct the inverse map of $\pr_2$. Take a point $\sB$ in the first space of \eqref{intersection 9.4}. We regard it as a point in the affine Grassmannian associated with the parabolic subgroup $P_{M,1,\cdots,1}$ corresponding to the partition $(M,1,1,\cdots,1)$. $\GL_M$ is the Levi subgroup of $P_{M,1,\cdots,1}$, we denote the image of $\sB$ in $\Gr_{M}$ by $\sC$. One can check easily that
\begin{equation}
    \begin{split}
        \iota_{\pr}: W^{(\xi,(\eta,\eta'))}&\longrightarrow \sO^0\cap \sS^{(\xi,(\eta,\eta'))}\\
        \sB&\mapsto (\sC, \sB),
    \end{split}
\end{equation}
is the inverse map of \eqref{pr2}.
\end{proof}
 
 Then, by taking transpose inverse, we see that the above intersection $W^{(\xi,(\eta,\eta'))}$ is isomorphic to 
 \begin{equation} \label{9.13}
 \begin{split}
     \begin{pNiceMatrix}
t^{\xi_1 }&*&\Cdots&*&1 &*&\Cdots&*\\
&     &\Ddots & \Vdots &  &&&\\
    &   & \Ddots   & * &\Vdots &&& \\
   &&&t^{\xi_M}&1 &\Vdots&&\Vdots\\
      &   & &&1 &*& \\
         &&&&&1&\Ddots&\\
         &&&&&&\Ddots&*\\
         &&&&&&&1
\end{pNiceMatrix}\GL_N(\bO)/\GL_N(\bO)\cap\\ 
\cap\begin{pNiceMatrix}
    t^{-\eta_{1}  }       &  & & &&&&\\
         * &  &&&& &&\\
         &\ddots&\Ddots&&&&&\\
 \Vdots & &\ddots&t^{-\eta_{M+1}  }&&&&\\
     &  &&&t^{-\eta_1'}&&&\\
     &&&&{\ddots}&&&\\
     & & &&&\ddots&\ddots&\\
      * &&\Cdots&&&&*&t^{-\eta_{N-M-1}'}.
\end{pNiceMatrix}\GL_N(\bO)/\GL_N(\bO)\\
\subset \Gr_N.
 \end{split}
\end{equation}

\begin{rem}\label{notation remark}
    From now on, we will denote by $\mathfrak{A}$ the first subspace of $\Gr_N$, and by $\mathfrak{B}$ the second subspace.
\end{rem}
\subsection{Vanishing of $H^1$}\label{vanishing of H1}
In this section, we will prove the most technical result of this paper, Proposition \ref{Vanishment}. 

To prove Theorem \ref{key}, we need to show the following two facts:
\begin{enumerate}
    \item There is an isomorphism 
    \begin{equation}
        \cI|_{\overset{\circ}{C^\bullet}}\simeq \Omega|_{\overset{\circ}{C^\bullet}}.
    \end{equation}
    \item The factorization algebra $\Omega$ is the $!*$-extension of $\Omega|_{\overset{\circ}{C^\bullet}}$.
\end{enumerate}

For the first claim, we note that both $\cI$ and $\Omega$ are factorizable, so we only need to check that there is a canonical isomorphism
\begin{equation}\label{over simple root}
    \cI|_{-\alpha_i\cdot x}\simeq \Omega|_{-\alpha_i\cdot x},
\end{equation}
for any simple root $\alpha_i$ of $\GL(M|N)$ and point $x\in C$.

If $1\leq i\leq 2M$, then the restriction of $\cY\to C^\bullet$ over $-\alpha_i\cdot C$ is an isomorphism, and the restriction of $p^\bullet(\IC_{glob}^0)$ to $\cY^{-\alpha_i}$ is the twisted constant perverse sheaf.

If $2M+1\leq i\leq M+N-1$, the isomorphism \eqref{over simple root} follows from \cite[Section 18.4.5-18.5]{[GL]} or \cite[Section 5.1]{[G]}.

To check the fact that $\Omega$ is the $!*$-extension of $\Omega|_{\overset{\circ}{C^\bullet}}$, we will use the method in \cite[Section 6]{[G]}.

Namely, by the factorization structure and induction, we may assume that the intermediate-extension property is already known away from the main diagonal
\begin{equation}
    \begin{split}
        C\hookrightarrow C^{(\xi,(\eta,\eta'))}\\
        x\mapsto (\xi,(\eta,\eta'))\cdot x,
    \end{split}
\end{equation}
and it remains to check that there is no subobject or quotient supported on the main diagonal. Since \(\Omega\) is perverse and Verdier self-dual, it is enough to verify that the \(!\)-restriction of $\Omega$ to the main diagonal lives in the perverse degree no less than $1$. As the main diagonal is one-dimensional, this is
equivalent to saying that, for a point
\(x\in C\), we need to check that the $!$-stalk of $\Omega=F^{glob}(\IC_{glob}^0)\in {\SD}_{\kappa}(C^\bullet)$ at the point $(\xi,(\eta,\eta'))\cdot x\in C^{(\xi,(\eta,\eta'))}$ is concentrated in degrees $\geq 2$. 

According to \eqref{zastava local-global}, the central fiber of $\cY$ over $(\xi,(\eta,\eta'))\cdot x\in C^{(\xi,(\eta,\eta'))}$ is the intersection $\overline{\sO}^0\cap \sS^{(\xi,(\eta,\eta'))}$. {Furthermore, due to the fact that the unital $H(\bF)$-orbit is minimal relevant, there is no Gaiotto D-module supported on the boundary $\overline{\sO}^0-\sO^0$, and we are only interested in the intersection $\sO^0\cap \sS^{(\xi, (\eta,\eta'))}$. The latter intersection is the fiber of $\cY_{\text{genu}}:= \Bun_H^\om\underset{\Bun_G}{\times'}\Bun_B \subset \cY$.}

Let us denote by $i_x$ the closed embedding of the central fiber over $(\xi,(\eta,\eta'))\cdot x$ to $\cY_{\text{genu}}$. By the base change theorem, 
\begin{equation}
    \begin{split}
        i_x^!\Omega\overset{\text{glob def}}{=}&H^\bullet(\sO^0\cap \sS^{(\xi,(\eta,\eta'))}, i_x^! p^\bullet(\IC^0_{glob}))\\
        \overset{\text{loc def}}{=}&H^\bullet(\sO^0\cap \sS^{(\xi,(\eta,\eta'))}, i_x^!(\sprd(\IC^0_{loc})\overset{!}{\otimes}\omega^{{\bullet}}_{\sS_{C^\bullet}}[\deg]).
    \end{split}
\end{equation}
 Recall the notation $\deg=\langle 2\rho^\circ, (\xi,(\eta,\eta')) \rangle$, and equals $\dim {\overset{\circ}{\cY}{}^{(\xi, (\eta,\eta'))}_{=0\cdot c}}=\dim \cM^{\xi}_{=0\cdot c} +\dim \Bun_B^{(\xi, (\eta,\eta'))}-\dim \Bun_G$, if $(\xi,(\eta,\eta'))$ is a span of negative simple roots of the supergroup.

Note that under the trivialization of twisting on the semi-infinite orbit $\sS_{C^\bullet}$, the restriction of twisted dualizing D-module $\omega^{{\bullet}}_{\sS_{C^\bullet}}[\deg]$ to the fiber $\sS^{(\xi,(\eta,\eta'))}$ is isomorphic to $\omega_{\sS^{(\xi,(\eta,\eta'))}}[-\dim {\overset{\circ}{\cY}{}^{(\xi, (\eta,\eta'))}_{=0\cdot c}}]$, which is placed in same parity as $(\xi, (\eta,\eta'))$. Also, under the $H(\bF)$-equivariant trivialization of $\cP_{\det}$ on the relevant unital orbit, the $!$-restriction of $\sprd(\IC^0_{loc})$ to $\sO^0$ is isomorphic to the rank 1 Gaiotto local system $\IC^0_{loc}$ on $\sO^0$. So, the $!$-fiber of $p^\bullet(\IC^0_{glob})$ at any point $y$ in the central fiber of $\cY_{\text{genu}}$ over $(\xi,(\eta,\eta'))\cdot x$ is $\BC[-\dim {\overset{\circ}{\cY}{}^{(\xi, (\eta,\eta'))}_{=0\cdot c}}]$. As a result, the restriction of the D-module $\sprd(\IC^0_{loc})\overset{!}{\otimes}\omega^{{\bullet}}_{\sS_{C^\bullet}}[\deg]$ to any irreducible component $Y$ of $\sO^0\cap \sS^{(\xi,(\eta,\eta'))}$ is concentrated in degrees $\geq \dim {\overset{\circ}{\cY}{}^{(\xi, (\eta,\eta'))}_{=0\cdot c}}-\dim Y$.  In particular, $H^i(Y, \IC_{loc}^0\overset{!}{\otimes} \omega_{\sS^{(\xi,(\eta,\eta'))}}[-\dim {\overset{\circ}{\cY}{}^{(\xi, (\eta,\eta'))}_{=0\cdot c}}])$ is concentrated in degrees $\geq  \dim {\overset{\circ}{\cY}{}^{(\xi, (\eta,\eta'))}_{=0\cdot c}}-2\dim Y$. Since the dimension of any irreducible component $Y$ of $\sO^0\cap \sS^{(\xi,(\eta,\eta'))}$ is no more than the half of $\dim {\overset{\circ}{\cY}{}^{(\xi, (\eta,\eta'))}_{=0\cdot c}}$. So, $H^i(\sO^0\cap \sS^{(\xi,(\eta,\eta'))}, \IC_{loc}^0\overset{!}{\otimes} \omega_{\sS^{(\xi,(\eta,\eta'))}}[-\dim {\overset{\circ}{\cY}{}^{(\xi, (\eta,\eta'))}_{=0\cdot c}}])=0$ for any $i< 0$. We have to prove the $H^1$ (resp.\ $H^0$) of the $!$-restriction of $p^\bullet(\IC^0_{glob})$ vanishes, if $(\xi,(\eta,\eta'))$ is not a negative simple root (resp.\ if $(\xi,(\eta,\eta'))\neq 0$).

That is to say, to prove Theorem \ref{key}, it suffices to show
\begin{prop}\label{Vanishment}
\begin{equation}\label{H0}
    H^0(\sO^0\cap \sS^{(\xi,(\eta,\eta'))}, \IC_{loc}^0\overset{!}{\otimes} \omega_{\sS^{(\xi,(\eta,\eta'))}}[-\dim {\overset{\circ}{\cY}{}^{(\xi, (\eta,\eta'))}_{=0\cdot c}}])=0,
\end{equation}
if $(\xi,(\eta,\eta'))\neq 0$, and 
\begin{equation}\label{H1}
    H^1(\sO^0\cap \sS^{(\xi,(\eta,\eta'))}, \IC_{loc}^0\overset{!}{\otimes} \omega_{\sS^{(\xi,(\eta,\eta'))}}[-\dim {\overset{\circ}{\cY}{}^{(\xi, (\eta,\eta'))}_{=0\cdot c}}])=0,
\end{equation}
if $(\xi,(\eta,\eta'))$ is not a negative simple root.

\end{prop}

\begin{proof}
Let $Y$ be an irreducible component of $\sO^0\cap \sS^{(\xi,(\eta,\eta'))}$, in order to prove \eqref{H1}, we only need to check that for any irreducible component (or even, its open dense subset) $Y$ of $\sO^0\cap \sS^{(\xi,(\eta,\eta'))}$ such that $\dim Y= \frac{1}{2} \dim {\overset{\circ}{\cY}{}^{(\xi, (\eta,\eta'))}_{=0\cdot c}}$ or $\frac{1}{2}(\dim {\overset{\circ}{\cY}{}^{(\xi, (\eta,\eta'))}_{=0\cdot c}}-1)$, we have
\begin{equation}\label{9.17}
     H^1(Y, \IC_{loc}^0\overset{!}{\otimes} \omega_{\sS^{(\xi,(\eta,\eta'))}}[-\dim {\overset{\circ}{\cY}{}^{(\xi, (\eta,\eta'))}_{=0\cdot c}}])=0.
\end{equation}

If $\dim {\overset{\circ}{\cY}{}^{(\xi, (\eta,\eta'))}_{=0\cdot c}}= 2\dim Y+1$, then $H^1$ is the bottom cohomology of $H^i(Y, \IC_{loc}^0\overset{!}{\otimes} \omega_{\sS^{(\xi,(\eta,\eta'))}}[-\dim {\overset{\circ}{\cY}{}^{(\xi, (\eta,\eta'))}_{=0\cdot c}}])$, so we only need to prove that the restriction of $\IC_{loc}^0\overset{!}{\otimes} \omega_{\sS^{(\xi,(\eta,\eta'))}}[-\dim {\overset{\circ}{\cY}{}^{(\xi, (\eta,\eta'))}_{=0\cdot c}}]$ to an open subset of $Y$ is not constant. Note that under the $H(\bF)$-equivariant trivialization on $\sO^0$, $\IC_{loc}^0$ is just the (shifted) local system corresponding to $\chi$, and under the $U(\bF)$-equivariant trivialization on $\sS^{(\xi,(\eta,\eta'))}$, $\omega_{\sS^{(\xi,(\eta,\eta'))}}$ is the dualizing D-module. So, the restriction of $\IC_{loc}^0\overset{!}{\otimes} \omega_{\sS^{(\xi,(\eta,\eta'))}}[-\dim {\overset{\circ}{\cY}{}^{(\xi, (\eta,\eta'))}_{=0\cdot c}}]$ to $Y$ can be rewritten as the tensor product of the local system corresponding to $\chi$, and the local system $\Psi_{0,(\xi,(\eta,\eta'))}$ which measures the ratio of two trivializations. Since the first local system is either non-tame or constant, but $\Psi_{0,(\xi,(\eta,\eta'))}$ is tame, the desired non-constant property follows from the fact that $\Psi_{0,(\xi,(\eta,\eta'))}$ is not constant unless $(\xi,(\eta,\eta'))$ is an odd negative simple root of $\GL(M|N)$ (ref. Lemma \ref{non trivial}). 

Similarly, we can show the claim \eqref{H0} about $H^0$.

Now, we focus on \eqref{H1} and the case $\dim {\overset{\circ}{\cY}{}^{(\xi, (\eta,\eta'))}_{=0\cdot c}}= 2\dim Y$. In order to simplify the notations, we denote by $d$ the dimension of $Y$. Note that according to Corollary \ref{cor 8.1.3}, we have 
\begin{equation}\label{expression}
    -(\xi, (\eta, \eta'))= \sum_{i=1}^{M-1} a_i \alpha_{i, \GL_M}+\sum_{j=1}^M b_j \alpha_{j, \GL_N}+ \sum_{k=M+1}^{N-1} c_k \alpha_{k, \GL_N}.
\end{equation}

Case I. In the expression \eqref{expression}, if there are at least two different $k_1$ and $k_2$ belong to $M+1, M+2,\cdots, N-1$, such that the coefficients of $\alpha_{k_1,\GL_N}$ and $\alpha_{k_2,\GL_N}$ are non-zero. Then, we can use the same method as \cite[6.4, 6.5 Case (2)]{[G]} to prove \eqref{9.17}.

Note that $\IC^0_{loc}$ and $\Psi_{0,(\xi,(\eta,\eta'))}$ are lisse, the $!$-restriction of $p^\bullet(\IC_{glob}^0)$ to $Y$ is a D-module concentrated in the degree $d$ and we can write it as $\chi\overset{!}{\otimes} \Psi_{0,(\xi,(\eta,\eta'))}[-2d]$, where $\chi$ is the local system on $Y$ corresponding to the character $\chi$ whose $!$-stalks are concentrated in degree $0$.

Consider the map 
\[\chi_{univ}: Y\longrightarrow \sO^0\longrightarrow \BG_a^{N-M-1}\]
induced by \eqref{chi}. For any non-degenerate linear map $l\colon \BG_a^{N-M-1}\rightarrow \BG_a^1$, the character of $H(\bF)$ in the definition of Gaiotto D-module can be taken as the composition of $\chi_{univ}$ and $l$.

Since $T_{N-M-1}(\bO)$ acts on $\sO^0\cap \sS^{(\xi,(\eta,\eta'))}$, it acts on $Y$ as well. Note that $\chi_{univ}$ is $T_{N-M-1}(\bO)$-equivariant, so the D-module $M_Y:= \chi_{univ,*}(\Psi_{0,(\xi,(\eta,\eta'))})$ is $T_{N-M-1}(\bO)$-equivariant against a non-trivial character.
\begin{equation}\label{MY}
    \begin{split}
        H^1(Y, \IC_{loc}^0\overset{!}{\otimes} \omega_{\sS^{(\xi,(\eta,\eta'))}}[-2d])\ =\ &H^1(Y, \chi\overset{!}{\otimes} \Psi_{0,(\xi,(\eta,\eta'))}[-2d])\\
        =\ & H^1(Y, (l\circ \chi_{univ})^!(\exp)\overset{!}{\otimes} \Psi_{0,(\xi,(\eta,\eta'))}[-2d])\\
        =\  &H^1(\BG_a^{N-M-1}, l^!(\exp)\overset{!}{\otimes} M_Y[-2d]).
    \end{split}
\end{equation}

Note that according to \eqref{MY}, the cohomology \eqref{9.17} shifted by $N-M-1$ is the fiber of the Fourier transform $FT(M_Y[-2d])$ of $M_Y[-2d]$ at $l\in (\BG_a^{N-M-1})^*$. By the $T_{N-M-1}(\bO)$-equivariance property of $M_Y$, $FT(M_Y[-2d])$ is lisse away from diagonals. In particular, it is lisse near $l$. So, we only need to prove that $M_Y$ is concentrated in degrees $\geq -2d+2$. Since $\Psi_{0,(\xi,(\eta,\eta'))}$ is in degree $-d$, we should show that $\chi_{univ,*}[-d+2]$ is left $t$-exact.

One can prove that $M_Y$ actually lives on $\BG_a^r$, where $r$ is the number of different $\alpha_{i,\GL_N}$, such that $M+1\leq i\leq N-1$, appear in the expression \eqref{expression} of $(\xi,(\eta,\eta'))$. 

Indeed, consider a point in the intersection \eqref{9.13},  i.e., $\fA\cap \fB$ (see, Remark \ref{notation remark}). If we take an upper matrix representative $\sA\in \fA$,
then the map $\chi_{univ}$ sends  to the $t^{-1}$-coefficients of $(\sA_{M+1,M+2}, \sA_{M+2, M+3}, \cdots, \sA_{N-1,N})$. We claim that for $M+1\leq i\leq N-1$, if the coefficients of $\alpha_{i,\GL_N}$ in \ref{expression} is zero, then $\sA_{i, i+1}\in \bO$. 

Namely, for such an $i$, since the coefficient of $\alpha_{i,\GL_N}$ is zero, we have 
\[\eta_{i-M}'+\eta_{i+1-M}'+\cdots+\eta_{N-M-1}'=0.\]


Let us take another lower matrix representative $\sB\in \fB$ of the chosen point in the intersection, and let us consider the $(N-i+1)\times(N-i+1)$-minor $\sB_{(i,i+1,\cdots, N), (i', i+1,\cdots, N)}$ for any $i'<i$. Since there exists a $(N-i+1)\times(N-i+1)$-matrix $C_{i'}$ with coefficients in $\bO$, such that \[A_{(i,i+1,\cdots, N), (i, i+1,\cdots, N)} C_{i'}=\sB_{(i,i+1,\cdots, N), (i', i+1,\cdots, N)},\]
we have 
\[\det \sA_{(i,i+1,\cdots, N), (i, i+1,\cdots, N)} \det C_{i'}=\det \sB_{(i,i+1,\cdots, N), (i', i+1,\cdots, N)}.\]

In particular, $\sB_{i,i'}\in \bO$ for all $i'<i$. 

Also, by Proposition \ref{OS}, we have 
\[-\eta_{i-M-1}'-\eta_{i-M}'-\cdots-\eta_{N-M-1}'\geq 0, \textnormal{if $M+2 \leq i$},\]
and
\[-\eta_{M+1}-\eta_1'-\eta_2'+\cdots-\eta_{N-M-1}'\geq 0.\]
So, $\sB_{i,i}\in \bO$ according to our assumption. In addition, $\sB_{i, i'}=0$ for $i'>i$. It implies that $\sA_{i, i+1}\in \bO$.

By the equivariance property of $M_Y$ against a non-trivial character, it is the clean extension of its restriction {to} the complement of the coordinate hyperplanes in $\BG_a^{r}$. Note that the dimension of the fiber of $\chi_{univ}: Y\longrightarrow \BG_a^{r}$ over this open subset is no more than $d-r$, which is no more than $d-2$ according to the assumption. Now, the desired property follows from the fact that $\chi_{univ, *}[-d+2]$ is left $t$-exact.

Case II. If there is no $M+1\leq i\leq N-1$ such that the coefficient of $\alpha_{i,\GL_N}$ in \eqref{expression} is non-zero. Then, the intersection \eqref{9.13} can be regarded as the intersection of the unital $\GL_{M}(\bF)$-orbit and the semi-infinite orbit $\sS^{(\xi, \eta)}$ in $\Gr_M\times \Gr_{M+1}$, and it is proved in \cite[Corollary 4.2.3]{[BFT0]} with \cite[Proposition 6.1.1]{[SW]}. Namely, it is proved in \cite{[SW]} that the dimension of any irreducible component $Y$ of $\sO^0\cap \sS^{(\xi,\eta)}$ is strictly smaller than the half of the dimension of the SW-Zastava space.

From now on, we assume that there is exactly one $i$ which belongs to $M+1,\cdots, N-1$, such that the coefficient of $\alpha_{i,\GL_N}$ is non-zero.

Case III. We assume that $i\geq M+2$. 

In this case, if $\eta_{N-M-1}'=0$, then one can show that the intersection \eqref{9.13} is isomorphic to the intersection of the unital relevant orbit of $\Gr_M\times \Gr_{N-1}$ and the semi-infinite orbit $\sS^{(\xi,(\eta, \tilde{\eta}'))}\subset \Gr_M\times \Gr_{N-1}$. Here $\tilde{\eta}'= (\eta_1', \eta_2',\cdots, \eta_{N-M-2}')$. So, we reduce the question from  $M,N$ case to $M, N-1$ case. 

So, without loss of generality, we can assume that $i= N-1$. In this case, the $\eta'$ which appears in the intersection \eqref{9.13} is $\eta'=(0, \cdots,0, n,-n)$. Here $n=-\eta'_{N-M-1}$ is the coefficient $c_{N-1}$ of $\alpha_{N-1,\GL_N}$ in \eqref{expression}. 

For any element in the intersection, let us take an upper triangle matrix representative $\sA\in \fA$. Since $\det \sA_{\{N-1,N\},\{N-1, N\}}=1$, we obtain that the minimal degree of the determinant of $2\times 2$-minors in the lines $N-1,N$ of any lower triangle matrix representative $\sB\in \fB$ is $0$. Since $\sB_{N,N}=t^n$ and $\sB_{N-1, N}=0$, we obtain that $\sB_{N-1, i}\in t^{-n}\bO$ for $1\leq i\leq N-1$. In particular, one can eliminate $\sB_{N-1, i}=0$ for $1\leq i\leq N-1$.

Similarly, for this new representative $\sB$, we have $\sB_{N,i}\in t^n \bO$ for $1\leq i\leq N-2$. So, one can further eliminate all $\sB_{N,i}$, for $1\leq i\leq N-2$.

Then, one can check that for $M+1< i<N-1$, there is $\sB_{i,j}\in \bO$. So, one can further eliminate $\sB_{i,j}$ for $M+1<i<N-1$ and $1\leq j<i$.

So, the intersection \eqref{9.13} can be rewritten as
\begin{equation}\label{9.21}
\begin{split}
  \fA \cap \begin{pNiceMatrix}
    t^{-\eta_{1}  }       & & & & &&&\\
     *       &&  & & &&&\\
        \Vdots&\Ddots  & \Ddots &&&&& \\
       *  &\Cdots&*&t^{-\eta_{M+1}  }&&&&\\
  && &&1&&&\\
  &     &&&&\Ddots&&\\
  &     &&&&&t^{-n}&\\
  &     &&&&&*&t^{n}
\end{pNiceMatrix}\GL_N(\bO)/\GL_N(\bO)\in \Gr_N,
\end{split}
\end{equation}
which is further isomorphic to 
\begin{equation}\label{9.22}
\begin{split}
    \begin{pNiceMatrix}
t^{\xi_1 }&*&\Cdots&*&1 &&&\\
&    &\Ddots &\Vdots  &  &&&\\
    &   & \Ddots   & * &\Vdots &&& \\
   &&&t^{\xi_M }&1 &&&\\
      &   & &&1 && \\
         &&&&&\Ddots&&\\
         &&&&&&1&*\\
         &&&&&&&1
\end{pNiceMatrix}\GL_N(\bO)/\GL_N(\bO)\\ \cap \begin{pNiceMatrix}
    t^{-\eta_{1}  }       & & & & &&&\\
     *       &&  & & &&&\\
        \Vdots&\Ddots  & \Ddots &&&&& \\
       *  &\Cdots&*&t^{-\eta_{M+1}  }&&&&\\
  && &&1&&&\\
  &     &&&&\Ddots&&\\
  &     &&&&&t^{-n}&\\
  &     &&&&&*&t^{n}
\end{pNiceMatrix}\GL_N(\bO)/\GL_N(\bO)\in \Gr_N.
\end{split}
\end{equation}

The latter intersection in $\Gr_{N}$ is isomorphic to the product of
\begin{equation}\label{9.23}
\begin{split}
      \begin{pNiceMatrix}
t^{\xi_1 }&*&\Cdots&*&1 \\
&    &\Ddots &\Vdots  &  \\
    &   & \Ddots   & * &\Vdots  \\
   &&&t^{\xi_M }&\\
   &&&&1
\end{pNiceMatrix}\GL_{M+1}(\bO)/\GL_{M+1}(\bO)\\ \cap \begin{pNiceMatrix}
    t^{-\eta_{1}  }    &   &  & \\
   *    &  &  & \\
        \Vdots &\Ddots & \Ddots & \\
       *  &\Cdots&*&t^{-\eta_{M+1}  }
\end{pNiceMatrix}\GL_{M+1}(\bO)/\GL_{M+1}(\bO) 
\end{split}    
\end{equation}
and
\begin{equation}
\begin{split}
\begin{pNiceMatrix}
1&*\\
&1
\end{pNiceMatrix}\GL_2(\bO)/\GL_2(\bO) \cap \begin{pNiceMatrix}
t^{-n}&\\
*&t^n
\end{pNiceMatrix}\GL_2(\bO)/\GL_2(\bO).
\end{split}
\end{equation}

According to our assumption \eqref{expression}, we have  $-(\xi, \eta)= \sum_{i=1}^{M-1} a_i \alpha_{i, \GL_M}+\sum_{j=1}^M b_j \alpha_{j, \GL_N}$. Then, according to \cite[Proposition 6.1.1]{[SW]}, the dimension of the former factor of the above product is no more than $\sum_{i=1}^{M-1}a_i+ \sum_{j=1}^M b_j$, and the equality holds only if $(\xi,\eta)=0$. Also, the dimension of the latter intersection is $n$. So, only if $(\xi,\eta)=0$, the dimension of any irreducible component of \eqref{9.23} has the chance to equal to $\sum_{i=1}^{M-1} a_i+\sum_{j=1}^M b_j+\sum_{k=M+1}^{N-1}c_k$. In this case, the claim \eqref{H1} can be reduced to the case $M=0, N=2$, which is proved in \cite[Section 6.5 Case (1)]{[G]}. Namely, if $n>1$, then the intersection $\{\begin{pNiceMatrix}
t^{-n}&\\
*&t^n
\end{pNiceMatrix}\}\in \Gr_2$ is isomorphic to $\BG_m\times \BG_a^{n-1}$. Since there is no non-constant map $\BG_a^{n-1}\longrightarrow \BG_m$, the difference of trivializations factors through $\BG_m$, and the exponential D-module comes from $\BG_a^{n-1}$, the cohomology of $\IC_{loc}^0\overset{!}{\otimes} \omega_{\sS^{(\xi,(\eta,\eta'))}}[-\dim {\overset{\circ}{\cY}{}^{(\xi, (\eta,\eta'))}_{=0\cdot c}}]$ on $Y$ is $0$ in all degrees.

In order to finish the proof (i.e., the case $i=M+1$), we need a better estimate of the dimension of the intersection \eqref{9.13}.

\subsection{Dimension estimate, II}
In the case $M=0$, it is well-known (ref. \cite{[MV]}) that the dimension of any irreducible component $Y$ of \eqref{9.13} equals the upper bound of the estimate in \eqref{eq 8.18}, i.e., the half of the dimension of the Zastava space. If $M=N-1$, then it is well-known (cf.~\cite[Theorem 1.5.1, or Proposition 6.1.1]{[SW]}) that any irreducible component of the intersection is strictly less than the upper bound of the estimate except the trivial case: $(\xi,(\eta,\eta'))=0$.

The proof of $H^1=0$ can be reduced to the case $M=N-2$. In this section, we study the intersection in the case $M=N-2$.

For a positive integer $-\eta_1'$, two sequences $\sj:=(1\leq j_k<j_{k-1}<\cdots< j_1\leq M+1)$ and $\si:= (i_k<i_{k-1}<\cdots<i_1<-\eta_1')$, we denote by $\fB_{\si, \sj}$ the subspace
\begin{equation}
    \begin{pNiceMatrix}
    t^{-\eta_{1}  }       &  & & &&&&\\
         * & \Ddots && &&&&\\
         \Vdots&\ddots&t^{\mathrlap{-\eta_{j_k}}}&&&&&\\
         &&\ddots&\Ddots&&&&\\
         &&&\ddots&t^{\mathrlap{-\eta_{j_1}}}&&\\
  \Vdots &&&&\ddots&\Ddots&&\\
       *&\Cdots&&\Cdots&\Cdots&*&&\\
       &&t^{i_k}\bO^\times&&t^{i_1}\bO^\times&&&t^{-\eta_{1}  '}
\end{pNiceMatrix}\GL_N(\bO)/\GL_N(\bO)
\end{equation}
of $\fB$. In other words, it consists of those points of $\fB$ which admit a representative $\sB$ such that
\begin{equation}
    \begin{split}
        \sB_{N,j_k}\in t^{i_k}\bO^\times, \sB_{N, j_{k-1}}\in t^{i_{k-1}}\bO^\times,\cdots, \sB_{N,j_1}\in t^{i_1}\bO^\times,
    \end{split}
\end{equation}
and $\sB_{N,i}=0$ for all other $1\leq i\leq M+1$.

It is not hard to see that for any geometric point $\sB$ of $\fB$, it belongs to exactly one $\fB_{\si,\sj}$. Furthermore, the subset $\fA\cap \fB_{\si,\sj}\subset \fA\cap \fB$ is non-empty only in the case $i_k=0$. So, we need to study, under the assumption $i_k=0$, when the dimension of the intersection $\fA\cap \fB_{\si,\sj}$ is equal to the upper bound of the estimate.

\begin{prop}\label{intersection}
If the dimension of $\fA\cap \fB_{\si,\sj}$ equals the upper bound of the estimate in \eqref{eq 8.18}, then
\begin{equation}\label{9.27}
    \begin{split}
        -\eta_{j_k}=i_k-i_{k-1},\\
        -\eta_{j_{k-1}}=i_{k-1}-i_{k-2},\\
        \cdots,\\
        -\eta_{j_1}=i_1+\eta_1',\\
        \eta_{i}=0, \Text{for all other\ } i\in 1,2,\cdots, M+1\\ 
        \xi_j=0, \Text{for all } j.
    \end{split}
\end{equation}

In particular, we have $-\eta_{M+1}  \leq 0$.
\end{prop}
\begin{proof}
One can check that for any element of the subspace $\fB_{\si,\sj}$, we can choose a representative of the form
\begin{equation}
\begin{split}
      &\begin{pNiceMatrix}
    t^{-\eta_{1}  }       &  & & &&&&*\\
        *  & \Ddots && &&&&\\
         \Vdots&\ddots&t^{\mathrlap{-\eta_{j_k}+i_{k-1}-i_k}}&&&&&\Vdots\\
         &&\ddots&\Ddots&&&&\\
         &&&\ddots&t^{\mathrlap{-\eta_{j_1}-\eta_1'-i_1}}&&&\Vdots\\
  \Vdots &&&&\ddots&\Ddots&&*\\
       *&\Cdots&&\Cdots&\Cdots&*&&*\\
       &&&&&&&1
\end{pNiceMatrix}\\
= U_{M+1}^-(\bF) \bF^{M+1} \diag(t^{-\eta_1},&\cdots,t^{-\eta_{j_k-1}}, t^{-\eta_{j_k}+i_{k-1}-i_k}, \cdots,t^{-\eta_{j_1}-\eta_1'-i_1},\cdots,t^{\eta_{M+1}}).
\end{split}
\end{equation}

Ignore the last row, we denote by $\bv_i$ the $i$-th column (i.e., a ($M+1$)-vector) of the above matrix. We can check that 
\begin{equation}\label{fiber}
\begin{split}
     \bv_{M+2}\in \bO\cdot \bv_{1}+\bO\cdot \bv_{2}+\cdots+\bO\cdot \bv_{j_{k}-1}\\
     +t^{i_k-i_{k-1}}\bO^\times\cdot \bv_{j_k}+t^{i_k-i_{k-1}}\bO\cdot \bv_{j_k+1}+\cdots+t^{i_k-i_{k-1}}\bO\cdot \bv_{j_{k-1}-1}\\
     +t^{i_k-i_{k-2}}\bO^\times\cdot \bv_{j_{k-1}}+t^{i_k-i_{k-2}}\bO\cdot \bv_{j_{k-1}+1}+\cdots+t^{i_k-i_{k-2}}\bO\cdot \bv_{j_{k-2}-1}\\
     +\cdots\\
     +t^{i_k-i_1}\bO^\times\cdot\bv_{j_2}+t^{i_k-i_1}\bO\cdot \bv_{j_2+1}+\cdots+t^{i_k-i_1}\bO\cdot\bv_{j_1-1}\\
     +t^{i_k+\eta_1'}\bO^\times\cdot \bv_{j_1}+t^{i_k+\eta_1'}\bO\cdot\bv_{j_1+1}+\cdots+t^{i_k+\eta_1'}\bO\cdot \bv_{M+1}.
\end{split}
\end{equation}

So, the intersection $\fA\cap \fB_{\si,\sj}$ is a fibration over 
\begin{equation}\label{9.30}
\begin{split}
      \begin{pNiceMatrix}
    t^{-\eta_{1}  }       &  & & &&&\\
     *     & \Ddots && &&&\\
         \Vdots&\ddots&t^{\mathrlap{-\eta_{j_k}+i_{k-1}-i_k}}&&&&\\
         &&\ddots&\Ddots&&&\\
         &&&\ddots&t^{\mathrlap{{-\eta_{j_1}-\eta_1'-i_1}}}&&\\
  \Vdots &&&&\ddots&\Ddots&\\
       *&\Cdots&&\Cdots&\Cdots&*&t^{-\eta_{M+1}  }
\end{pNiceMatrix}\GL_{M+1}(\bO)/\GL_{M+1}(\bO)\cap\\ \cap\begin{pNiceMatrix}
t^{\xi_1 }&*&\Cdots&*&1 \\
&     &\Ddots &\Vdots  &  \\
    &   & \Ddots   & * &\Vdots  \\
   &&&t^{\xi_M }&\\
   &&&&1
\end{pNiceMatrix}\GL_{M+1}(\bO)/\GL_{M+1}(\bO)\in \Gr_{M+1},
\end{split}
\end{equation}
with the fiber controlled by \eqref{fiber}.

Over a fixed point in \eqref{9.30}, the dimension of the fiber is no more than
\begin{equation}
    \begin{split}
        (j_{k-1}-j_k) i_{k-1}\\
        + (j_{k-2}-j_{k-1}) i_{k-2}\\
        +\cdots\\
        +(j_1-j_2) i_1\\
        -(M+2-j_1)\eta_1'.
    \end{split}
\end{equation}

It equals
\begin{equation*}
    \begin{split}
        \dim (U_{M+2}^-(\bF) \diag(t^{-\eta_{1}  },t^{-\eta_{2}},\cdots, t^{-\eta_{M+1}  }, t^{-\eta_1'})\GL_{M+2}(\bO)/\GL_{M+2}(\bO)\cap\\ \cap U_{M+2}(\bF)\diag(t^{\xi_1 },t^{\xi_2} ,\cdots, t^{\xi_M} ,1,1)\GL_{M+2}(\bO)/\GL_{M+2}(\bO))\\-\\ \dim (U_{M+1}^-(\bF)\diag(t^{-\eta_{1}  },\cdots, t^{-\eta_{j_k}+i_{k-1}-i_k},\cdots t^{-\eta_{j_1}-\eta_1'-i_1},\cdots, t^{-\eta_{M+1}}  )\GL_{M+1}(\bO)/\GL_{M+1}(\bO)\cap\\ \cap U_{M+1}(\bF)\diag(t^{\xi_1 },t^{\xi_2 },\cdots, t^{\xi_M },1)\GL_{M+1}(\bO)/\GL_{M+1}(\bO)).
    \end{split}
\end{equation*}
Here, the first intersection is taken in $\Gr_{M+2}$, and the second is in $\Gr_{M+1}$.


So, if we have 
\begin{equation*}
    \begin{split}
        \dim (\fA\cap \fB_{\si,\sj})=\dim (U_{M+2}^-(\bF) \diag(t^{-\eta_{1}  },t^{-\eta_{2}},\cdots, t^{-\eta_{M+1}  }, t^{-\eta_1'})\GL_{M+2}(\bO)/\GL_{M+2}(\bO)\cap\\ \cap U_{M+2}(\bF)\diag(t^{\xi_1 },t^{\xi_2} ,\cdots, t^{\xi_M} ,1,1)\GL_{M+2}(\bO)/\GL_{M+2}(\bO)),
    \end{split}
\end{equation*}
then, the dimension of the intersection \eqref{9.30} should be equal to 
\begin{equation}\label{9.35}
    \begin{split}
        \dim (U_{M+1}^-(\bF)\diag(t^{-\eta_{1}  },\cdots, t^{-\eta_{j_k}+i_{k-1}-i_k},\cdots t^{-\eta_{j_1}-\eta_1'-i_1},\cdots, t^{-\eta_{M+1}}  )\GL_{M+1}(\bO)/\GL_{M+1}(\bO)\cap\\ \cap U_{M+1}(\bF)\diag(t^{\xi_1 },t^{\xi_2 },\cdots, t^{\xi_M },1)\GL_{M+1}(\bO)/\GL_{M+1}(\bO)).
    \end{split}
\end{equation}
It is an intersection in $\Gr_{M+1}$, and according to \cite[Proposition 6.1.1]{[SW]}, it happens only if the condition \eqref{9.27} holds. Indeed, according to our assumption \eqref{expression}, $(-\xi_1,-\xi_2,\cdots,-\xi_M) ,\ (-\eta_{1},\cdots, -\eta_{j_k}+i_{k-1}-i_k,\cdots,-\eta_{j_1}-\eta_1'-i_1,\cdots, -\eta_{M+1})$ are spanned by $\alpha_{i,\GL_M}$ and  $\alpha_{j,\GL_{M+1}}$, respectively.  In non-trivial case, the dimension of \eqref{9.30} is strictly less than the half of the dimension of the corresponding SW-Zastava component, which is equal to the number with multiplicities of simple roots of $\GL_M\times \GL_{M+1}$ appears in $(-\xi_1,-\xi_2,\cdots,-\xi_M) ,\ (-\eta_{1},\cdots, -\eta_{j_k}+i_{k-1}-i_k,\cdots,-\eta_{j_1}-\eta_1'-i_1,\cdots, -\eta_{M+1})$. It equals $\sum a_i+\sum b_j+\sum c_k-     (j_{k-1}-j_k) i_{k-1}
        - (j_{k-2}-j_{k-1}) i_{k-2}
        -\cdots
        -(j_1-j_2) i_1
        +(M+2-j_1)\eta_1'$. It is the dimension of \eqref{9.35}.

\end{proof}

\begin{rem}\label{another proof 9.4.2}
Let us briefly sketch another proof of the above proposition, which can be generalized to the orthosymplectic group case easily.

Define $U^\prime$ as $U_M^-\times \sD U_{M+1} U_{M,M+2}^- \sD^{-1}$. In other words, $U'$ is the semi-direct product of the normalizer subgroup of $U_{M,M+2}^-$ in $U$ and the group $U_{M,M+2}^-$. Then, the intersection $\fA\cap \fB_{\si, \sj}\subset \sO^0\cap \sS^{\xi, (\eta, \eta')}$ equals the intersection of $\sO^0\cap \sS^{\xi, (\eta, \eta')}$ with a single $U'(\bF)$-orbit $\sS^{\prime, (\xi, \tilde{\eta}, \tilde{\eta}')}$. Note that both $U'$ and $H$ belong to the parabolic subgroup of normalizer subgroup of $U_{M,M+2}^-$ in $\GL_M\times \GL_{M+2}$. This parabolic subgroup has a homomorphism to $\GL_M\times
 \GL_{M+1}$. In particular, there is a map 
\begin{equation}\label{reinter prop inter}
    \fA\cap \fB_{\si, \sj}= \sO^0\cap \sS^{\xi, (\eta, \eta')}\cap \sS^{\prime, (\xi, \tilde{\eta}, \tilde{\eta}')}\hookrightarrow \sO^0\cap \sS^{\prime, (\xi, \tilde{\eta}, \tilde{\eta}')}\longrightarrow \Gr_{M}\times \Gr_{M+1}.
\end{equation}

The Proposition \ref{intersection} is equivalent to the following claim: if the dimension of $\fA\cap \fB_{\si,\sj}$ equals the upper bound of the estimate, then its image under the map \eqref{reinter prop inter} is just a point.

Using this reformulation, Proposition \ref{intersection} can be roughly reproved as follows. One can check that the image of $\sO^0\cap \sS^{\xi, (\eta, \eta')}\cap \sS^{\prime, (\xi, \tilde{\eta}, \tilde{\eta}')}$ in $\Gr_M\times \Gr_{M+1}$ is the intersection of the unital $\GL_M(\bF)$-orbit and a semi-infinite orbit $\sS^{\prime, (\xi, \tilde{\eta}, \tilde{\eta}')-\alpha}$ in $\Gr_M\times \Gr_{M+1}$. Here, $\alpha$ is a span of simple roots of $\GL_{M+2}$ and $(\xi, \tilde{\eta}, \tilde{\eta}')-\alpha$ belongs to the coweight lattice of $\GL_M\times \GL_{M+1}\subset \GL_M\times \GL_{M+2}$. The dimension of the fiber of \eqref{reinter prop inter} is no more than the number (with multiplicities) of simple roots of $\GL_M\times \GL_{M+2}$ in the decomposition of  $\alpha$, so $\dim \fA\cap \fB_{\si, \sj}$ equals the upper bound only if the image of \eqref{reinter prop inter}, as an intersection of the unital orbit and a semi-infinite orbit in $\Gr_M\times \Gr_{M+1}$, has the dimension given by the upper bound. However, according to \cite[Theorem 1.5.1, or Proposition 6.1.1]{[SW]}, it only happens when the intersection is a point.

\end{rem}
\subsubsection{}
With the preparation in the last sections, now we are going to finish the proof of \eqref{H1}. 

Case IV. Now, we assume that $\alpha_{M+1,\GL_N}$ is the only simple root among $\alpha_{M+1,\GL_N},\cdots, \alpha_{N-1, \GL_N}$ such that its coefficient of \eqref{expression} is non-zero. By applying the same analysis as Case III, we can assume that $N=M+2$. 

The intersection is 
\begin{equation}\label{9.24}
\begin{split}
       \begin{pNiceMatrix}
t^{\xi_1 }&*&\Cdots&*&1 &*\\
&     &\Ddots &\Vdots  &  &\\
    &   & \Ddots   & * &\Vdots & \\
   &&&t^{\xi_M }&1 &\Vdots&\\
      &   & &&1 &*\\
      &&&&&1
\end{pNiceMatrix}\GL_{M+2}(\bO)/\GL_{M+2}(\bO)\cap\\
\cap\begin{pNiceMatrix}
    t^{-\eta_{1}  }      & &  & & \\
   *      & &  & & \\
          &\Ddots& \Ddots && \\
  \Vdots &&&t^{-\eta_{M+1}  }&\\
       *&\Cdots&&*&t^{-\eta_1'}
\end{pNiceMatrix}\GL_{M+2}(\bO)/\GL_{M+2}(\bO)\subset \Gr_{M+2}.
\end{split}
\end{equation}

According to Proposition \ref{intersection}, we only need to assume $-\eta_{M+1}  \leq 0$.

Case IV A). Assume $-\eta_{M+1}  \leq -2$. Since the minimal degree of entries of the $N$-th row of any representative matrix $\sB\in \fB$ in the intersection is $0$, we have $\sB_{N,N-1}\in \bO$.

Consider the left action of
\begin{equation}\label{9.36}
    \BG_a=\begin{pNiceMatrix}
1&&&\\
&\Ddots&&\\
&&1& t^{-1}\BC\\
&&&1
\end{pNiceMatrix}
\end{equation}
on the intersection.

We claim that the intersection is $\BG_a$-invariant. Note that $\fA$ is obviously $\BG_a$-invariant, we only need to show that for any representative $\sB\in \fB$ of any element in the intersection, we have $c\cdot \sB\in \fB$ for any $c\in \BG_a$. 

Note that since $\sB$ lands in the intersection, $\sB_{N, i}\in \bO$. In particular, $\sB_{N-1,N-1}+t^{-1}\bO$ has the same degree as $\sB_{N-1,N-1}$. Also, left multiplication with an element in $\BG_a$ does not change the determinant of the $2\times 2$-minor $\sB_{\{N-1,N\},\{N-1,N\}}$. So, the intersection is $\BG_a$-invariant.

Furthermore, we can write $Y$ as $Y'\times \BG_a^1$. Since there is no non-constant map $\BG_a^1\rightarrow \BG_m$, the difference of trivializations factors through $Y'$. Also, the exponential D-module comes from pullback of $\exp$ on $\BG_a^1$. In particular, $\IC_{loc}^0\overset{!}{\otimes} \omega_{\sS^{(\xi,(\eta,\eta'))}}[-\dim {\overset{\circ}{\cY}{}^{(\xi, (\eta,\eta'))}_{=0\cdot c}}]$ is $(\BG_a, \exp)$-equivariant, its cohomology is $0$ in all degrees.



Case IV B). Assume $-\eta_{M+1}  =0$ and $-\eta_1'>0$ or $-\eta_{M+1}  =-1$ and $-\eta_1'>1$. We divide the intersection into two subsets $P_1$ and its complement $P_2$. Here, $P_1$ consists of those elements which  admit a representative $\sB\in \fB$, such that $\sB_{N, N-1}\in t^{-\eta_1'-\eta_{M+1}}\bO$. 

According to Proposition \ref{intersection}, the dimension of $P_2$ is strictly smaller than the critical dimension. To prove \eqref{H1}, we only need to show that the $H^1$ cohomology on  $P_1$ vanishes. For any element in $P_1$, we can choose a representative $\sB$, such that $\sB_{N,N-1}=0$. So, $Y$ is $\BG_a$-invariant with respect to the action \eqref{9.36}. Then, using similar arguments as the proof of Case IV A), we can show that the cohomology on $P_1$ vanishes in all degrees.

Case IV C). Assume $-\eta_{M+1}  =-1$ and $-\eta_1'=1$. Then, the intersection is of the maximal possible dimension only if $\xi_i=0$ and $\eta_i=0$ for any $1\leq i\leq M$. In other words, $(\xi,(\eta, \eta'))$ is a negative simple root.
\end{proof}

\begin{rem}
    When $M=N-2$, through the same analysis as above, we can show that $\Omega\simeq \cI$ for any non-degenerate $\kappa$. Here, we define $\Omega$ to be $v_*(\sprd(\IC^{0}_{loc})\overset{!}{\otimes} \IC^{\frac{\infty}{2}}_{q^{-1},C^\bullet_{\infty\cdot c}})[\deg]$, and $\IC^{\frac{\infty}{2}}_{q^{-1},C^\bullet_{\infty\cdot c}}$ denotes the twisted semi-infinite IC D-module on $\overline{\sS}_{C^\bullet_{\infty\cdot c}}$, ref \cite[Section 13]{[GL]}, \cite{[G5]}.
\end{rem}

\subsection{Grothendieck groups}
As a corollary of Theorem \ref{key}, we claim
\begin{cor}\label{cor 9.6.1}
$F^{loc}(\IC_{loc}^{(\lambda,(\theta,\theta'))})$ is isomorphic to $\IC^{\fact}_{(\lambda,(\theta,\theta'))}$.
\end{cor}

\begin{proof}
We mimic the proof of \cite[Section 6.7]{[G]}. Namely, we use inductive method to show these two twisted D-modules are isomorphic on $C^{(\xi,(\eta,\eta'))}_{\infty\cdot c}$. 

Note that according to Proposition \ref{OS}, the image of $\IC_{loc}^{(\lambda,(\theta,\theta'))}$ is supported on $C^\bullet_{\leq (\lambda,(\theta,\theta'))\cdot c}$. Because of Proposition \ref{one point}, the $!$-stalk of the image of $\IC_{loc}^{(\lambda,(\theta,\theta'))}$ at $(\lambda,(\theta,\theta'))\cdot c$ is $\BC$ up to a shift. By the factorization property and $\Omega\simeq \cI$, i.e., Theorem \ref{key}, there is an isomorphism
\begin{equation}
    F^{loc}(\IC_{loc}^{(\lambda,(\theta,\theta'))})|_{C^\bullet_{=(\lambda,(\theta,\theta'))\cdot c}}\simeq \IC^{\fact}_{\lambda,(\theta,\theta')}|_{C^\bullet_{=(\lambda,(\theta,\theta'))\cdot c}}.
\end{equation}

To show the above isomorphism extends to $C^{(\xi,(\eta,\eta'))}_{\leq (\lambda,(\theta,\theta'))\cdot c}$ for all $(\xi,(\eta,\eta'))$, we use induction on $(\lambda,(\theta,\theta'))-(\xi,(\eta,\eta'))$. 

If $(\lambda,(\theta,\theta'))=(\xi,(\eta,\eta'))$, the claim is obvious. Now, assume the claim is true for all $(\tilde{\lambda},(\tilde{\theta},\tilde{\theta}'))$ and $(\tilde{\xi},(\tilde{\eta},\tilde{\eta}'))$ such that $(\tilde{\lambda},(\tilde{\theta},\tilde{\theta}'))-(\tilde{\xi},(\tilde{\eta},\tilde{\eta}'))<(\lambda,(\theta,\theta'))-(\xi,(\eta,\eta'))$. Then, by factorization property, the isomorphism on $C^{(\xi,(\eta,\eta'))}_{=(\lambda,(\theta,\theta'))\cdot c}$ extends to the complement of $\{(\xi,(\eta,\eta'))\cdot c\}\hookrightarrow C^{(\xi,(\eta,\eta'))}_{\leq (\lambda,(\theta,\theta'))\cdot c}$.

The $!$-stalk of $F^{loc}(\IC_{loc}^{(\lambda,(\theta,\theta'))})$ at $\{(\xi,(\eta,\eta'))\cdot c\}$ is
\begin{equation}\label{9.38}
    H^\bullet(\Gr_G, \IC_{loc}^{(\lambda,(\theta,\theta'))}\overset{!}{\otimes} \omega_{\sS^{(\xi,(\eta,\eta'))}}[\deg]).
\end{equation}

Here $\IC_{loc}^{(\lambda,(\theta,\theta'))}\overset{!}{\otimes} \omega_{\sS^{(\xi,(\eta,\eta'))}}[\deg])$ is the $!$-restriction of $p^\bullet(\IC_{glob}^{(\lambda,(\theta,\theta'))})$ to the central fiber over $(\xi,(\eta,\eta'))\cdot c$. We need to show this cohomology is concentrated in degrees $\geq 1$. Note that according to Corollary \ref{t exact}, $F^{loc}(\IC_{loc}^{(\lambda,(\theta,\theta'))})$ is perverse. So, the cohomology \eqref{9.38} concentrates in non-negative degrees. 

It suffices to show that the restriction of $p^\bullet(\IC^{(\lambda,(\theta,\theta'))}_{glob})$ to the central fiber is non-constant on any irreducible component (or, its open dense subset) of the critical dimension (i.e., half of $\dim \overset{\circ}{\cY}{}_{= (\lambda, (\theta, \theta'))\cdot c}^{(\xi, (\eta, \eta'))}$). It follows from the fact that (up to a cohomology shift) $\IC_{loc}^{(\lambda,(\theta,\theta'))}\overset{!}{\otimes} \omega_{\sS^{(\xi,(\eta,\eta'))}})$ is the tensor of a non-tame ramified local system and a non-constant tame ramified local system $\Psi_{(\lambda,(\theta,\theta')),(\xi,(\eta,\eta'))}$, where the non-tame ramified local system comes from the pullback of the exponential D-module and the tamely ramified local system comes from the ratio of two trivializations on the central fiber (See Lemma \ref{non trivial}).
\end{proof}

\begin{cor}
The functor $F^{glob}$ (equally, $F^{loc}$) in Definition \ref{def glob func} induces a $t$-exact, conservative functors
\begin{equation}\label{functor}
\begin{split}
    F^{glob}: \cC_{\kappa}^{glob,lc}(M|N)\longrightarrow \cI-\FM^{fin}.\\
    F^{loc}: \cC_{\kappa}^{loc,lc}(M|N)\longrightarrow \cI-\FM^{fin},
\end{split}
\end{equation}
and the restrictions to the heart are faithful.
\end{cor}
\begin{proof}
The finiteness property of the image of $F^{glob}$ follows from that any irreducible factorization module $\IC_{(\lambda,(\theta,\theta'))}^{\fact}$ (for $(\lambda,(\theta,\theta'))$ relevant) satisfies the finiteness property in \cite{[BFS]}. Conservative and faithful properties follow from Corollary \ref{cor 9.6.1}.
\end{proof}
Based on the above propositions and the corollaries, we have

\begin{cor}\label{cor fusion}
The fusion product defines a braided monoidal structure of the twisted Gaiotto category, and the functor $F^{glob}$ (equally, $F^{loc}$) is braided monoidal.
\end{cor}
\begin{proof}
Note that both monoidal category structures of $\cC_{\kappa}^{glob}(M|N)$ and $\cI-\FM^{fin}$ are given by taking the fusion product, i.e., nearby cycles. For $(\xi, (\eta, \eta'))$ large enough, the map ${\cY}^{(\xi, (\eta, \eta'))}_{\infty\cdot c}\rightarrow \cM_{\infty\cdot c}$ is smooth, so the fusion products are compatible with the pullback on ${\cY}^{(\xi, (\eta, \eta'))}_{\infty\cdot c}$. Using the same analysis of Proposition \ref{p tex}, we obtain that the fusion products are compatible with the pullback on ${\cY}_{\infty\cdot c}$ by the factorizable property. Now the functor $F^{glob}$ preserves the fusion product follows from the fact that the map $\overline{\cY}_{\infty\cdot c}\rightarrow C^\bullet_{\infty\cdot c}$ is proper, and the proper map preserves the nearby cycles. Now, the statement follows from the conservativity of the functor $F^{glob}$.
\end{proof}
\begin{cor}
The functor \eqref{functor} induces a monoidal isomorphism of the Grothendieck rings of $\cC_{\kappa}^{loc,lc}(M|N)=\cC_{\kappa}^{glob,lc}(M|N)$ and $\cI-\FM^{fin}$.
\end{cor}

\section{Proof of equivalence}\label{section 10}
In the last section, we have already shown that the functor $F^{loc}$ sends irreducible objects in the twisted Gaiotto category to the irreducible objects in $\cI-\FM^{fin}$. So, we need to check that the functor is fully faithful. 

\subsection{End of the proof}
The full faithfulness of $F^{loc}$ follows from repeating the same proof as in \cite[Section 4.3-5]{[BFT0]} word-by-word, and the key properties are that the Gaiotto category is rigid (Lemma \ref{almost eq}) and has a collection of projective-injective generators (Corollary \ref{cor ab}). In order to be self-complete, we briefly sketch the proof.

\begin{defn}
We call a relevant orbit $\sO^{(\lambda,(\theta,\theta'))}$ \textit{typical}, if we have $\lambda_i+\theta_j\neq 0$ for $i\in 1,2,\cdots, M$ and $j\in 1,2,\cdots, M+1$.
\end{defn}
\begin{rem}
In this case,  the stabilizer in $H^{\omega}(\bF)$ of  $\sL_{(\lambda,(\theta,\theta'))}$ is pro-unipotent.

\end{rem}

By the same proof as in \cite[Proposition 4.4.1]{[BFT0]}, we have
\begin{lem}\label{typical}
If $\sO^{(\lambda,(\theta,\theta'))}$ is typical, then for any relevant $(\tilde{\lambda},(\tilde{\theta},\tilde{\theta}'))$ not equal to $(\lambda,(\theta,\theta'))$, the $!$-stalk of $\IC_{loc}^{(\tilde{\lambda},(\tilde{\theta},\tilde{\theta}'))}$ at $\sO^{(\lambda,(\theta,\theta'))}$ is zero. 
\end{lem}

Let $\lambda^\dag$ be $(-M, -M+1,\cdots, -1)$, $\theta^\dag$ be $(-M, -M+1,-M+2,\cdots,-1, 0)$, and $\theta^{\dag '}$ be $(0,0,\cdots,0)$.
\begin{cor}\label{pro-inj}
The twisted D-module $\IC_{loc}^{(\lambda^\dag,(\theta^\dag, \theta^{\dag '}))}$ is projective and injective in $\cC_{\kappa}^{loc}(M|N)^\heartsuit$.
\end{cor}
\begin{proof}
According to Proposition \ref{closure}, one can check that there is no other relevant $H(\bF)$-orbit in $\Gr_G$ in the closure of $\sO^{(\lambda^\dag,(\theta^\dag, \theta^{\dag '}))}$. In particular, $\IC^{(\lambda^\dag,(\theta^\dag, \theta^{\dag '}))}_{loc}$ is the clean extension from its restriction {to} $\sO^{(\lambda^\dag,(\theta^\dag, \theta^{\dag '}))}$.

Now, $\Ext_{D^b(\cC^{loc}_{\kappa}(M|N)^\heartsuit)}^1(\IC^{(\lambda^\dag,(\theta^\dag, \theta^{\dag '}))}_{loc}, \IC^{(\lambda,(\theta,\theta'))}_{loc})=\Ext_{\cC^{loc}_{\kappa}(M|N)}^1(\IC^{(\lambda^\dag,(\theta^\dag, \theta^{\dag '}))}_{loc}, \IC^{(\lambda,(\theta,\theta'))}_{loc})$ equals the $!$-stalk of $\IC^{(\lambda,(\theta,\theta'))}_{loc}$ at $\sO^{(\lambda^\dag,(\theta^\dag, \theta^{\dag '}))}$, which is zero by Lemma \ref{typical}. Here, $D^b(\cC^{loc}_{\kappa}(M|N)^\heartsuit)$ denotes the derived category of $\cC^{loc}_{\kappa}(M|N)^\heartsuit$. As a result, $\IC_{loc}^{(\lambda^\dag,(\theta^\dag, \theta^{\dag '}))}$ is projective. Applying the Verdier duality functor, we obtain that $\IC_{loc}^{(\lambda^\dag,(\theta^\dag, \theta^{\dag '}))}$ is injective.
\end{proof}

Denote $\IC_{taut}:= \IC_{loc}^{(0,0,\cdots,0), ((-1,0,\cdots,0), (0,0,\cdots, 0))}$ and $\IC_{taut}^*:= \IC_{loc}^{(0,0,\cdots,0), ((0,0,\cdots,0), (0,0,\cdots, 1))}$. Let $\cE$ be the full fusion subcategory of $\cC^{loc}_{\kappa}(M|N)^\heartsuit$ generated by $\IC_{loc}^0$, $\IC_{taut}$, and $\IC_{taut}^*$ via taking fusion products, images, kernels and cokernels.

The following lemma is proved in \cite[Lemma 4.3.2, Corollary 4.3.3]{[BFT0]} with a method proposed by P. Etingof.
\begin{lem}\label{almost eq}
$F^{loc}|_{\cE}: \cE\rightarrow \cI-\FM^{fin,\heartsuit}$ is an equivalence of rigid braided abelian categories. In particular, since any irreducible object is contained in $\cE$, we obtain that any irreducible object $\IC_{loc}^{(\lambda,(\theta,\theta'))}$ is dualizable.
\end{lem}

\begin{cor}\label{cor ab}
a). For any relevant $(\lambda,(\theta,\theta'))$, the fusion product $ \IC_{loc}^{(\lambda,(\theta,\theta'))}\star \IC_{loc}^{(\lambda^\dag,(\theta^\dag, \theta^{\dag '}))}$ is projective and injective in $\cC_{\kappa}^{loc}(M|N)^\heartsuit$.

b). Any irreducible object $\IC_{loc}^{(\lambda,(\theta,\theta'))}$ is a submodule (also, a quotient) of some fusion products in a).
\end{cor}
\begin{proof}
The claim a). follows from the fact that the fusion product is bi-exact, Lemma \ref{almost eq}, and Corollary \ref{pro-inj}.

b). Since \begin{equation*}
    \begin{split}
        \Hom_{\cC_{\kappa}^{loc}(M|N)^\heartsuit}(\IC_{loc}^{(\lambda,(\theta,\theta'))}\star \IC_{loc}^{(\lambda^\dag,(\theta^\dag, \theta^{\dag '})),*}, \IC_{loc}^{(\xi,(\eta,\eta'))})\\= \Hom_{\cC_{\kappa}^{loc}(M|N)^\heartsuit}(\IC_{loc}^{(\lambda,(\theta,\theta'))}, \IC_{loc}^{(\xi,(\eta,\eta'))}\star \IC_{loc}^{(\lambda^\dag,(\theta^\dag, \theta^{\dag '}))}),
    \end{split}
\end{equation*}
we only need to take a $(\xi,(\eta,\eta'))$ such that the above Hom is not zero. Then $\IC_{loc}^{(\lambda,(\theta,\theta'))}$ is a submodule of $\IC_{loc}^{(\xi,(\eta,\eta'))}\star \IC_{loc}^{(\lambda^\dag,(\theta^\dag, \theta^{\dag '}))}$. Similar for the claim about quotient.
\end{proof}
We have the following theorem which is an analog of \cite[Theorem 4.5.1]{[BFT0]}.
\begin{thm}\label{theorem 10.1.7}
The natural functor 
\begin{equation}
    {{D}}^b(\cC_{\kappa}^{loc,lc}(M|N)^\heartsuit)\longrightarrow \cC_{\kappa}^{loc,lc}(M|N)
\end{equation}
is an equivalence.
\end{thm}
\begin{proof}
According to the above corollary, we only need to show that
\begin{equation}\label{naive compare}
    \begin{split}
        \RHom_{{{D}}^b(\cC_{\kappa}^{loc,lc}(M|N)^\heartsuit)}(\IC_{loc}^{(\lambda,(\theta,\theta'))}\star \IC_{loc}^{(\lambda^\dag,(\theta^\dag, \theta^{\dag '}))}, \IC_{loc}^{(\xi,(\eta,\eta'))}\star \IC_{loc}^{(\lambda^\dag,(\theta^\dag, \theta^{\dag '}))})\longrightarrow\\ \RHom_{\cC_{\kappa}^{loc,lc}(M|N)}(\IC_{loc}^{(\lambda,(\theta,\theta'))}\star\IC_{loc}^{(\lambda^\dag,(\theta^\dag, \theta^{\dag '}))}, \IC_{loc}^{(\xi,(\eta,\eta'))}\star \IC_{loc}^{(\lambda^\dag,(\theta^\dag, \theta^{\dag '}))}),
    \end{split}
\end{equation}
is an isomorphism. 

For the Hom in ${{D}}^b(\cC_{\kappa}^{loc,lc}(M|N)^\heartsuit)$, since $\IC_{loc}^{(\lambda,(\theta,\theta'))}\star \IC_{loc}^{(\lambda^\dag,(\theta^\dag, \theta^{\dag '}))}$ is projective, it is concentrated in degree $0$. Note that $\IC_{loc}^{(\lambda,(\theta,\theta'))}$ is dualizable, the left-hand side of \eqref{naive compare} equals the multiplicity of $\IC_{loc}^{(\lambda^\dag,(\theta^\dag, \theta^{\dag '}))}$ in $\IC_{loc}^{(\lambda,(\theta,\theta')),*}*\IC_{loc}^{(\xi,(\eta,\eta'))}*\IC_{loc}^{(\lambda^\dag,(\theta^\dag, \theta^{\dag '}))}$.

Since $\IC_{loc}^{(\lambda^\dag,(\theta^\dag, \theta^{\dag '}))}$ is the clean extension from its restriction {to} $\sO^{(\lambda^\dag,(\theta^\dag, \theta^{\dag '}))}$, the right-hand side of \eqref{naive compare} equals the $!$-stalk of $\IC_{loc}^{(\lambda,(\theta,\theta')),*}*\IC_{loc}^{(\xi,(\eta,\eta'))}*\IC_{loc}^{(\lambda^\dag,(\theta^\dag, \theta^{\dag '}))}$ at $\sO^{(\lambda^\dag,(\theta^\dag, \theta^{\dag '}))}$. By Lemma \ref{typical}, it is concentrated in degree $0$ and equals the multiplicity  of $\IC_{loc}^{(\lambda^\dag,(\theta^\dag, \theta^{\dag '}))}$ in $\IC_{loc}^{(\lambda,(\theta,\theta')),*}*\IC_{loc}^{(\xi,(\eta,\eta'))}*\IC_{loc}^{(\lambda^\dag,(\theta^\dag, \theta^{\dag '}))}$.
\end{proof}

Now, we can finish the proof of the main theorem.
\begin{thm}\label{statement of main}
The functor \eqref{functor} is an equivalence of braided categories.
\end{thm}
\begin{proof}
We should prove that $\cE\rightarrow \cC_{\kappa}^{loc,lc}(M|N)^{\heartsuit}$ is essentially surjective. According to Lemma \ref{almost eq} and Corollary \ref{cor ab}, all {twisted} D-modules of the form $\IC_{loc}^{(\lambda,(\theta,\theta'))}\star \IC_{loc}^{(\lambda^\dag,(\theta^\dag, \theta^{\dag '}))}$ belong to $\cE$. Note that any object of $\cC_{\kappa}^{loc,lc}(M|N)^{\heartsuit}$ is finite length, so it is the image of a map from a direct sum of $\IC_{loc}^{(\lambda,(\theta,\theta'))}\star \IC_{loc}^{(\lambda^\dag,(\theta^\dag, \theta^{\dag '}))}$ to a direct sum of $\IC_{loc}^{(\lambda,(\theta,\theta'))}\star \IC_{loc}^{(\lambda^\dag,(\theta^\dag, \theta^{\dag '}))}$, so it is also contained in $\cE$.

\end{proof}

\end{document}